\documentclass[11pt, a4paper]{amsart} %,leqno 
\usepackage{amsmath,amssymb,amsthm}
\usepackage{mathrsfs}

\usepackage{xcolor}
\usepackage{hyperref}
\usepackage{fullpage}
\usepackage{tikz-cd}

\numberwithin{equation}{subsection}

\newtheorem{theorem}[equation]{Theorem}
\newtheorem*{theorem*}{Theorem}
\newtheorem{lemma}[equation]{Lemma}
\newtheorem*{lemma*}{Lemma}
\newtheorem{proposition}[equation]{Proposition}
\newtheorem*{proposition*}{Proposition}
\newtheorem{corollary}[equation]{Corollary}
\newtheorem*{corollary*}{Corollary}
\newtheorem{conjecture}[equation]{Conjecture}
\newtheorem*{conjecture*}{Conjecture}

\theoremstyle{definition}

\newtheorem*{definition*}{Definition}

\newtheorem*{notation*}{Notation}

\newtheorem*{choice*}{Choice}
\newtheorem{remark}[equation]{Remark}
\newtheorem*{remark*}{Remark}
\newtheorem{assumption}[equation]{Assumption}
\newtheorem*{assumption*}{Assumption}

\newtheorem*{axiom*}{Axiom}
\newtheorem{condition}[equation]{Condition}

\allowdisplaybreaks[1]

\newcommand{\abs}[1]{\left\lvert#1\right\rvert}

\DeclareMathOperator{\Irr}{Irr}
\DeclareMathOperator{\End}{End}

\DeclareMathOperator{\Frob}{Frob}

\DeclareMathOperator{\Ad}{Ad}

\DeclareMathOperator{\ind}{ind}
\DeclareMathOperator{\Ind}{Ind}

\DeclareMathOperator{\Hom}{Hom}

\DeclareMathOperator{\red}{red}
\DeclareMathOperator{\lMod}{-Mod}
\DeclareMathOperator{\rMod}{Mod-}
\DeclareMathOperator{\Res}{Res}

\DeclareMathOperator{\id}{id}

\DeclareMathOperator{\ord}{ord}
\DeclareMathOperator{\supp}{supp}

\DeclareMathOperator{\Gal}{Gal}

\DeclareMathOperator{\der}{der}
\DeclareMathOperator{\GL}{GL}
\DeclareMathOperator{\SL}{SL}

\DeclareMathOperator{\spl}{spl}
\newcommand*{\autm}{\mathrm{autm}}
\newcommand*{\qs}{\mathrm{qs}}

\newcommand*{\spec}{\mathrm{spec}}

\DeclareMathOperator{\Aut}{Aut}
\DeclareMathOperator{\PGL}{PGL}
\DeclareMathOperator{\Pic}{Pic}
\newcommand*{\op}{\mathrm{op}}
\DeclareMathOperator{\Br}{Br}

\DeclareMathOperator{\Alg}{-Alg}
\DeclareMathOperator{\Rep}{Rep}
\DeclareMathOperator{\LLC}{LLC}
\DeclareMathOperator{\Cent}{Cent}
\DeclareMathOperator{\diag}{diag}

\DeclareMathOperator{\Red}{Red}
\newcommand*{\base}{\mathrm{base}}

\newcommand{\unr}{{\textup{unr}}}

\newcommand{\isoarrow}{\stackrel{\sim}{\longrightarrow}}

\newcommand{\Fgrp}{X} % notation for finite abelian group 
\newcommand{\resi}[1][F]{k_{#1}} % residue field
\newcommand{\cts}{\mathrm{cts}}

\newcommand{\bC}{\mathbb{C}}
\newcommand*{\C}{\mathbb{C}}

\newcommand{\Z}{\mathbb{Z}}

\title{Supercuspidal blocks and the local Langlands correspondence}
\author{Noriyuki Abe}
\address[N. Abe]{Graduate School of Mathematical Sciences, the University of Tokyo, 3-8-1 Komaba, Meguro-ku, Tokyo 153-8914, Japan.}
\thanks{The first-named author was supported by JSPS KAKENHI Grant Number 23H01065.}
\email{abenori@ms.u-tokyo.ac.jp}
\author{Kazuma Ohara}
\address[K. Ohara]{Universität Bonn\\
	Mathematisches Institut\\
	Endenicher Allee 60\\
	53115 Bonn\\
	Germany}
\email{kohara@math.uni-bonn.de}

\begin{document}
\begin{abstract}
Let $G$ be a connected reductive group defined over a non-archimedean local field $F$, and let $\mathrm{Rep}^{\mathfrak{s}}(G(F))$ be a supercuspidal Bernstein block.
In this paper, we define a cohomology class $[\mu_{\mathfrak{s}}]$ represented by a 2-cocycle that encodes the equivalence class of $\mathrm{Rep}^{\mathfrak{s}}(G(F))$.
We then formulate a conjecture describing $[\mu_{\mathfrak{s}}]$ in terms of the enhanced $L$-parameter $(\varphi, \rho)$ attached to an irreducible representation in $\mathrm{Rep}^{\mathfrak{s}}(G(F))$ via the conjectural local Langlands correspondence.
Finally, we verify this conjecture for the local Langlands correspondence for non-singular supercuspidal representations constructed by Kaletha.
\end{abstract}
\maketitle

\tableofcontents

\section{Introduction}
\counterwithin{equation}{section}
Let $G$ be a connected reductive group defined over a non-archimedean local field $F$ and $G(F)$ the group of $F$-valued points.
The category of smooth complex representations of $G(F)$ decomposes as a product of indecomposable full subcategories, called Bernstein blocks.
In many cases, each Bernstein block $\Rep^{\mathfrak{s}}(G(F))$ is equivalent to the category of right unital $\mathcal{H}^{\mathfrak{s}}$-modules, where $\mathcal{H}^{\mathfrak{s}}$ is a semi-direct product of an affine Hecke algebra and a twisted group algebra.
See, for instance, a recent work \cite{2024arXiv240807801A, 2024arXiv240807805A}, where such an equivalence is established under the assumptions that $G$ splits over a tamely ramified extension of $F$, and the residue characteristic $p$ of $F$ does not divide the order of the absolute Weyl group of $G$.
Roughly speaking, the affine Hecke algebra part of $\mathcal{H}^{\mathfrak{s}}$ arises from parabolic induction, whereas the twisted group algebra part comes from the supercuspidal representations.
In particular, if the block $\Rep^{\mathfrak{s}}(G(F))$ is supercuspidal, i.e.,\  if $\Rep^{\mathfrak{s}}(G(F))$ consists of supercuspidal representations, the algebra $\mathcal{H}^{\mathfrak{s}}$ reduces to a twisted group algebra.
In this paper, we provide a conjectural description of the Morita equivalence class of $\mathcal{H}^{\mathfrak{s}}$ for a supercuspidal block $\Rep^{\mathfrak{s}}(G(F))$ in terms of the corresponding $L$-parameters via the local Langlands correspondence.
We also prove that this conjecture holds for Kaletha's local Langlands correspondence for non-singular supercuspidal representations~\cite{2019arXiv191203274K}.

We explain our results in more detail.
Let $\Rep^{\mathfrak{s}}(G(F))$ be a supercuspidal block.
We fix an irreducible representation $\pi$ in $\Rep^{\mathfrak{s}}(G(F))$.
We denote by $X_{\unr}(G)$ the group of unramified characters of $G(F)$ and define a subgroup $G(F)^{1}$ of $G(F)$ as
\[
G(F)^1 = \bigcap_{\chi \in X_{\unr}(G)} \ker(\chi).
\]
We also define $N_{\mathfrak{s}} := N_{G(F)}(\pi^{1})/G(F)^{1}$, where $\pi^{1}$ is an irreducible subrepresentation of the restriction of $\pi$ to $G(F)^{1}$, and $N_{G(F)}(\pi^{1})$ denotes the normalizer of $\pi^1$ in $G(F)$.
Then there exists a 2-cocycle $\mu \colon N_{\mathfrak{s}} \times N_{\mathfrak{s}} \to \mathbb{C}^{\times}$ such that 
\[
\Rep^{\mathfrak{s}}(G(F)) \simeq \rMod \mathcal{H}^{\mathfrak{s}},
\]
where $\mathcal{H}^{\mathfrak{s}} = \mathbb{C}[N_{\mathfrak{s}}, \mu]$ is the twisted group algebra 
\[
\mathbb{C}[N_{\mathfrak{s}}, \mu] = \bigoplus_{n \in N_{\mathfrak{s}}} \mathbb{C} \cdot b_{n},
\]
whose multiplication is defined by $b_{n_1} b_{n_2} = \mu(n_1, n_2) \cdot b_{n_1 n_2}$ for $n_1, n_2 \in N_{\mathfrak{s}}$, and $\rMod \mathcal{H}^{\mathfrak{s}}$ denotes the category of right unital $\mathcal{H}^{\mathfrak{s}}$-modules.
Therefore, understanding the structure of the block $\Rep^{\mathfrak{s}}(G(F))$ reduces to determining the Morita equivalence class of $\mathcal{H}^{\mathfrak{s}}$.
The structure of the center $Z^{\mathfrak{s}}$ of $\mathcal{H}^{\mathfrak{s}}$ is better understood and admits the following description.
We define the subgroup $X_{\unr}(G)_{\pi}$ of $X_{\unr}(G)$ as 
\[
X_{\unr}(G)_{\pi} = \left\{
\chi \in X_{\unr}(G) \mid \pi \otimes\chi \simeq \pi
\right\}.
\]
Then the center $Z^{\mathfrak{s}}$ of $\mathcal{H}^{\mathfrak{s}}$ is the group algebra $\mathbb{C}[T_{\mathfrak{s}}]$, where $T_{\mathfrak{s}}$ is the subgroup of $N_{\mathfrak{s}}$ defined as 
\[
T_{\mathfrak{s}} = \left\{
g \in G(F)/G(F)^{1} \mid \text{$\chi(g) = 1$ for all $\chi \in X_{\unr}(G)_{\pi}$}
\right\}.
\]
As a consequence, $Z^{\mathfrak{s}}$ can be identified with the ring of regular functions on the complex algebraic torus $X_{\unr}(G)/X_{\unr}(G)_{\pi}$, which parametrizes the irreducible representations in $\Rep^{\mathfrak{s}}(G(F))$.
Although the algebra $\mathcal{H}^{\mathfrak{s}}$ is in general not Morita equivalent to its center $Z^{\mathfrak{s}}$, the base change $\mathcal{H}^{\mathfrak{s}} \otimes_{Z^{\mathfrak{s}}} \mathbb{C}[N_{\mathfrak{s}}]$ to $\mathbb{C}[N_{\mathfrak{s}}]$ is Morita equivalent to $\mathbb{C}[N_{\mathfrak{s}}]$.
Algebras over $Z^{\mathfrak{s}}$ with this property are called Azumaya algebras, and their Morita equivalence classes are parametrized by $H^2\left(
(N_{\mathfrak{s}}/T_{\mathfrak{s}})^{*}, \mathbb{C}[N_{\mathfrak{s}}]^{\times}\right)$, where $(N_{\mathfrak{s}}/T_{\mathfrak{s}})^{*}$ denotes the character group of $N_{\mathfrak{s}}/T_{\mathfrak{s}}$.
We denote by $c_{\mathcal{H}^{\mathfrak{s}}} \in H^2\left(
(N_{\mathfrak{s}}/T_{\mathfrak{s}})^{*}, \mathbb{C}[N_{\mathfrak{s}}]^{\times}\right)$ the cohomology class corresponding to $\mathcal{H}^{\mathfrak{s}}$.
We also define a cohomology class $[\mu_{\mathfrak{s}}] \in H^2\left(
X_{\unr}(G)_{\pi}, \mathbb{C}^{\times}
\right)$ as follows.
For each $\chi \in X_{\unr}(G)_{\pi}$, we fix an isomorphism $I_{\chi} \colon \pi \isoarrow \pi \otimes \chi$.
For $\chi_1, \chi_2 \in X_{\unr}(G)_{\pi}$, we define $\mu_{\mathfrak{s}}(\chi_1, \chi_2) \in \bC^{\times}$ by
\[
I_{\chi_1} \circ I_{\chi_2} = \mu_{\mathfrak{s}}(\chi_1, \chi_2) \cdot I_{\chi_1 \chi_2}.
\]
Then $\mu_{\mathfrak{s}}$ gives a 2-cocycle on $X_{\unr}(G)_{\pi}$, and its cohomology class $[\mu_{\mathfrak{s}}] \in H^2\left(
X_{\unr}(G)_{\pi}, \mathbb{C}^{\times}
\right)$ is independent of $\pi$ and $I_{\chi}$.
The first main result of this paper is the comparison of $c_{\mathcal{H}^{\mathfrak{s}}}$ and $[\mu_{\mathfrak{s}}]$.
We regard $H^2\left(
(N_{\mathfrak{s}}/T_{\mathfrak{s}})^{*}, \mathbb{C}^{\times}
\right)$ as a subset of $H^2\left(
X_{\unr}(G)_{\pi}, \mathbb{C}^{\times}
\right)$ via the inflation map
\[
H^2\left(
(N_{\mathfrak{s}}/T_{\mathfrak{s}})^{*}, \mathbb{C}^{\times}
\right)
\hookrightarrow H^2\left(
X_{\unr}(G)_{\pi}, \mathbb{C}^{\times}
\right)
\]
induced by the restriction map $X_{\unr}(G)_{\pi} \twoheadrightarrow (N_{\mathfrak{s}}/T_{\mathfrak{s}})^{*}$.
\begin{theorem}[{Corollaries~\ref{corollary:twistedgroupalgebraasanazumayaalgebra} and \ref{corollarymus=muvee}}]
\label{theorem1inIntro}
We have $[\mu_{\mathfrak{s}}] \in H^2\left(
(N_{\mathfrak{s}}/T_{\mathfrak{s}})^{*}, \mathbb{C}^{\times}
\right)$.
Moreover, the cohomology class $c_{\mathcal{H}^{\mathfrak{s}}}$ is the image of $[\mu_{\mathfrak{s}}]^{-1}$ via the natural map
\[
H^2\left(
(N_{\mathfrak{s}}/T_{\mathfrak{s}})^{*}, \mathbb{C}^{\times}
\right) \to H^2\left(
(N_{\mathfrak{s}}/T_{\mathfrak{s}})^{*}, \mathbb{C}[N_{\mathfrak{s}}]^{\times}\right).
\]
\end{theorem}

Therefore, the cohomology class $[\mu_{\mathfrak{s}}]$ encodes the information of the $Z^{\mathfrak{s}}$-category $\Rep^{\mathfrak{s}}(G(F))$.
The rest of the paper is devoted to describing $[\mu_{\mathfrak{s}}]$ in terms of spectral-side data via the (refined) local Langlands correspondence.
In this paper, we adopt the refinement of the local Langlands correspondence based on rigid inner forms as developed by Kaletha and Dillery (\cite{MR3548533, MR4617942}).
We denote by $W_{F}$ the Weil group of $F$, by $\widehat{G}$ the dual group of $G$ over $\mathbb{C}$, and by ${^LG}$ the $L$-group $\widehat{G} \rtimes W_{F}$ of $G$.
In the conjectural local Langlands correspondence, it is expected that we can attach to an irreducible smooth representation $\pi$ of $G(F)$ a $\widehat{G}$-conjugacy class of enhanced $L$-parameters $(\varphi, \rho)$, where $\varphi \colon W_{F} \times \SL_{2}(\mathbb{C}) \to {^LG}$ is a group homomorphism called an $L$-parameter, and $\rho$ is an irreducible representation of the $S$-group $S_{\varphi}^{+}$, which is a modification of the centralizer $\Cent_{\widehat{G}}(\varphi)$ in $\widehat{G}$ of the image of $\varphi$.
For simplicity, in this introduction we restrict ourselves to the case where $S_{\varphi}^{+} = \Cent_{\widehat{G}}(\varphi)$.
We attach to an enhanced $L$-parameter $(\varphi, \rho)$ a cohomology class $[\mu_{(\varphi, \rho)}] \in H^2\left(
X_{\unr}(G)_{[(\varphi,\rho)]}, \mathbb{C}^{\times}
\right)$ as follows.
Let $X_{\unr}(G_{F^{\unr}})$ denote the group of unramified characters of $G(F^{\unr})$, where $F^{\unr}$ denotes the maximal unramified extension of $F$.
We can define an action of $X_{\unr}(G_{F^{\unr}})$ on the set of enhanced $L$-parameters $(\varphi, \rho)$ using the Kottwitz homomorphism.
For an enhanced $L$-parameter $(\varphi, \rho)$, we define a subgroup $N_{\varphi, \rho}^{+}$ of $\widehat{G}$ as
\[
N_{\varphi, \rho}^{+} = \left\{
g \in \widehat{G} \mid {^g(\varphi, \rho)} = \chi \cdot (\varphi, \rho) \quad \text{for some $\chi \in X_{\unr}(G_{F^{\unr}})$}
\right\}.
\]
Then we have an exact sequence 
\[
1 \to S_{\varphi}^{+} \to N_{\varphi, \rho}^{+} \to X_{\unr}(G_{F^{\unr}})_{[(\varphi, \rho)]} \to 1,
\]
where $X_{\unr}(G_{F^{\unr}})_{[(\varphi, \rho)]}$ denotes the stabilizer of the $\widehat{G}$-conjugacy class of $(\varphi, \rho)$.
The representation $\rho$ of $S_{\varphi}^{+}$ extends to a projective representation of $N_{\varphi, \rho}^{+}$, which defines a 2-cocycle $\mu_{(\varphi, \rho)} \in Z^2\left(
X_{\unr}(G_{F^{\unr}})_{[(\varphi, \rho)]}, \mathbb{C}^{\times}
\right)$.
Moreover, we prove that the cohomology class $[\mu_{(\varphi, \rho)}]$ of $\mu_{(\varphi, \rho)}$ is contained in the image of the inflation map (Proposition~\ref{prop:muforspecsideFunrvsF})
\[
H^2\left(
X_{\unr}(G)_{[(\varphi, \rho)]}, \mathbb{C}^{\times}
\right)
\hookrightarrow H^2\left(
X_{\unr}(G_{F^{\unr}})_{[(\varphi, \rho)]}, \mathbb{C}^{\times}
\right),
\]
where $X_{\unr}(G)_{[(\varphi, \rho)]}$ is the image of $X_{\unr}(G_{F^{\unr}})_{[(\varphi, \rho)]}$ via the restriction map $X_{\unr}(G_{F^{\unr}}) \twoheadrightarrow  X_{\unr}(G)$.
We regard $[\mu_{(\varphi, \rho)}] \in H^2\left(
X_{\unr}(G)_{[(\varphi, \rho)]}, \mathbb{C}^{\times}
\right)$.
The main conjecture of this paper relates the cohomology class $[\mu_{\mathfrak{s}}] \in H^2\left(
X_{\unr}(G)_{\pi}, \mathbb{C}^{\times}
\right)$ to $[\mu_{(\varphi, \rho)}]$:
\begin{conjecture}[{Conjecture~\ref{conjecture:comparisonoftwococycles}}]
\label{conjectureinIntro}
Assume that there exists a local Langlands correspondence for irreducible representations in $\Rep^{\mathfrak{s}}(G(F))$ that is compatible with unramified twists.
Fix an irreducible representation $\pi\in \Rep^{\mathfrak{s}}(G(F))$ and let $(\varphi,\rho)$ be an enhanced $L$-parameter whose $\widehat{G}$-conjugacy class corresponds to $\pi$.
The compatibility with unramified twists implies $X_{\unr}(G)_{\pi} = X_{\unr}(G)_{[(\varphi, \rho)]}$.
Then we have $[\mu_{\mathfrak{s}}] = [\mu_{(\varphi, \rho)}]$ as elements of $H^2\left(
X_{\unr}(G)_{\pi}, \mathbb{C}^{\times}
\right) = H^2\left(
X_{\unr}(G)_{[(\varphi, \rho)]}, \mathbb{C}^{\times}
\right)$.
\end{conjecture}

The main result of this paper is a proof of Conjecture~\ref{conjectureinIntro} for Kaletha's local Langlands correspondence for non-singular supercuspidal representations.
In \cite{2019arXiv191203274K}, Kaletha defined and constructed a class of supercuspidal representations called non-singular supercuspidal representations.
He also gave a candidate for the local Langlands correspondence for non-singular supercuspidal representations.
If the residue characteristic $p$ of $F$ does not divide the order of the absolute Weyl group of $G$, this correspondence involves all supercuspidal $L$-parameters, which are expected to be precisely the parameters corresponding to those discrete series $L$-packets of $G$ consisting entirely of supercuspidal representations.
Kaletha's construction of the local Langlands correspondence depends on several choices.
In particular, it depends on a choice of isomorphism $\iota \colon \Box_{\autm} \isoarrow \Box_{\spec}$, where $\Box_{\autm}$ and $\Box_{\spec}$ are groups defined on the automorphic and spectral sides, respectively.
For more details, see Section~\ref{subsec:Computation of the cohomology class attached to non-singular Lparameter}.
In this paper, we define the ``inertial versions'' $\widetilde{\Box}_{\autm, F^{\unr}}$ and $\widetilde{\Box}_{\spec, F^{\unr}}$ that contain $\Box_{\autm}$ and $\Box_{\spec}$, respectively.
The main result of this paper is the following:
\begin{theorem}[{Theorem~\ref{mainthm:comparisonoftildeautmandtildespec}, Corollary~\ref{corollary:Kaletha'sLLCsatisfiesmainconjecture}}]
\label{mainresultinIntro}
There exists an isomorphism $\widetilde{\iota} \colon \widetilde{\Box}_{\autm, F^{\unr}} \isoarrow \widetilde{\Box}_{\spec, F^{\unr}}$ whose restriction $\iota \colon \Box_{\autm} \isoarrow \Box_{\spec}$ to $\Box_{\autm}$ has the following property: the corresponding local Langlands correspondence in \cite{2019arXiv191203274K} satisfies Conjecture~\ref{conjectureinIntro}.
\end{theorem}

\subsection*{Relation to prior work}

The idea of describing a Bernstein block $\Rep^{\mathfrak{s}}(G(F))$ and its corresponding algebra $\mathcal{H}^{\mathfrak{s}}$ in terms of the spectral data was developed in a series of works by Aubert--Moussaoui--Solleveld~\cite{2017arXiv170103593A,MR3845761,MR3969871}.
Using the generalized Springer correspondence, they established a block decomposition for enhanced $L$-parameters and for each block $\mathfrak{s}^{\vee}$, they defined an associated $\mathbb{C}$-algebra $\mathcal{H}^{\mathfrak{s}^{\vee}}$.
They conjectured that each block $\mathfrak{s}^{\vee}$ corresponds to a Bernstein block $\Rep^{\mathfrak{s}}(G(F))$ and that the algebra $\mathcal{H}^{\mathfrak{s}^{\vee}}$ is closely related to the algebra $\mathcal{H}^{\mathfrak{s}}$ attached to $\Rep^{\mathfrak{s}}(G(F))$ by the work of Solleveld~\cite{sol2022}.
Several instances of the local Langlands correspondence have been established by verifying this conjecture in special cases, including the correspondence for unipotent representations~\cite{Lusztig95unipotent, Lusztig02unipotentII, MR4167790, MR4515286, Solleveld23unipotentLLC, MR4620884}, for principal series representations of quasi-split $p$-adic groups~\cite{MR4862261}, and for non-singular depth-zero representations~\cite{MR5056155}.
For the twisted group algebra part, namely for supercuspidal blocks, they only consider the center $Z^{\mathfrak{s}}$ of $\mathcal{H}^{\mathfrak{s}}$.
Although the algebra $\mathcal{H}^{\mathfrak{s}}$ is closely related to its center, the two do not coincide in general.
An example illustrating this phenomenon can be found in \cite[\S 4.4]{MR3331726}.
Since supercuspidal representations are both projective and injective in the category of representations with fixed central character, the category of finite-length representations in a fixed supercuspidal block is equivalent to the category of finite-dimensional modules over the center $Z^{\mathfrak{s}}$ of $\mathcal{H}^{\mathfrak{s}}$. However, this equivalence no longer holds once the finite-length condition is removed.
In this paper, we study the Morita equivalence class of $\mathcal{H}^{\mathfrak{s}}$ as an Azumaya algebra over $Z^{\mathfrak{s}}$.
In this sense, our work may be viewed as complementary to the works by Aubert--Moussaoui--Solleveld.

Furthermore, a categorical version of the local Langlands correspondence has recently been formulated by Fargues and Scholze (\cite{2021arXiv210213459F}).
This categorical correspondence aims to describe not only the set of irreducible smooth representations of $G(F)$, but also the entire category $\Rep(G(F))$ of smooth representations.
On the spectral side, it involves the derived category of coherent sheaves on the stack of $L$-parameters $\mathrm{Par}_{G}$.
The authors expect that the cohomology class $[\mu_{(\varphi, \rho)}]$ constructed in this paper can be naturally interpreted within the framework of $\mathrm{Par}_{G}$.

\subsection*{Structure of the paper}

We give a brief overview of the following sections of this paper.

We collect necessary notation in Section~\ref{sec:Notation and assumptions}.

In Section~\ref{subsec:Basic facts about Azumaya algebras}, we record basic facts about Azumaya algebras, including their classification in terms of cohomology classes of 2-cocycles.
In Section~\ref{subsec:Twisted group algebra as an Azumaya algebra}, we apply these results to a special class of Azumaya algebras, namely twisted group algebras.
Finally in Section~\ref{subsection:Azumaya algebra attached to a supercuspidal block}, we prove Theorem~\ref{theorem1inIntro} using the results of previous subsections.

In Section~\ref{sec:aconjecture}, we attach to an enhanced $L$-parameter $(\varphi, \rho)$ a cohomology class $[\mu_{(\varphi, \rho)}]$ represented by a 2-cocycle and formulate Conjecture~\ref{conjectureinIntro}.

In Section~\ref{sec:The case of non-singular supercuspidal representations}, we prove Theorem~\ref{mainresultinIntro}.
In Section~\ref{subsec:Computation of the cohomology class attached to a non-singular supercuspidal representation}, we recall Kaletha's construction of non-singular supercuspidal representations and describe $[\mu_{\mathfrak{s}}]$ as the cohomology class attached to a projective representation of the group $\widetilde{\Box}_{\autm, F^{\unr}}$.
In Section~\ref{subsec:Computation of the cohomology class attached to non-singular Lparameter}, we recall Kaletha's local Langlands correspondence for non-singular supercuspidal representations and describe $[\mu_{(\varphi, \rho)}]$ as the cohomology class attached to a projective representation of the group $\widetilde{\Box}_{\spec, F^{\unr}}$.
In Sections~\ref{subsec:deformationogtildeBoxautm} and \ref{subsec:deformationofwidetildespec}, we describe $\widetilde{\Box}_{\autm, F^{\unr}}$ and $\widetilde{\Box}_{\spec, F^{\unr}}$ as deformations of simpler groups $\widetilde{\Box}_{\autm, F^{\unr}, \base}$ and $\widetilde{\Box}_{\spec, F^{\unr}, \base}$, respectively.
Using these deformations, in Section~\ref{subsec:comparisonofextensions}, we prove the main result Theorem~\ref{mainresultinIntro}.

In Appendix~\ref{sec:A cohomology class attached to a projective representation}, we give a construction of a certain cohomology class used in the main part of this paper.
In Appendix~\ref{sec:An example}, we directly compute the cohomology classes defined in this paper for the example considered in~\cite{MR3331726}.

\subsection*{Acknowledgments}

The authors would like to express their gratitude to the Max Planck Institute for Mathematics for its hospitality and financial support.
They would also like to thank Peter Dillery for explaining rigid inner forms over local function fields to them. The authors gratefully acknowledge OpenAI for providing free access to ChatGPT through the ChatGPT for Academic Researchers program.

\subsection*{Use of AI}
During the final proofreading stage, the authors used GPT-5.6 Sol and GPT-6 Astra (OpenAI) to assist in identifying typographical errors and potential mistakes.
We have not used any AI-generated text without careful verification.
The text of this paper has been written under the full responsibility of the authors.

\section{Notation and assumptions}\label{sec:Notation and assumptions}
\counterwithin{equation}{subsection}
Let $F$ be a non-archimedean local field with residue characteristic $p$.
We write $\resi$ for the residue field of $F$.
We fix a separable closure $\overline{F}$ of $F$.
When we refer to a separable field extension $E/F$, we always assume that $E \subset \overline{F}$.
For an algebraic extension $E$ of $F$, let $\mathcal{O}_{E}$ denote the ring of integers of $E$.
Let $E/F$ be a finite extension.
We denote by $N_{E/F}$ the norm map, by $E^{\unr}$ the maximal unramified extension of $E$, and by $\ord_{E} \colon E^{\times} \to \mathbb{Z}$ the discrete valuation normalized by $\ord_{E}(E^{\times}) = \mathbb{Z}$.
We also denote by $\ord_{E^{\unr}}$ the unique extension of $\ord_{E}$ to $E^{\unr, \times}$.
For a positive integer $n$, we denote by $E_{n}$ the unramified extension of $E$ of degree $n$.

We denote by $\Gamma_F$ the absolute Galois group $\Gal(\overline{F}/F)$ of $F$ and by $W_F$ the Weil group of $F$.
We also denote by $I_F$ the inertia subgroup of $\Gamma_F$ and by $P_F$ the wild inertia subgroup of $\Gamma_F$.
We fix a lift $\Frob_{F} \in \Gamma_{F}$ of the geometric Frobenius element in $\Gamma_F/I_F$.
Let $\ord_{W_{F}} \colon W_{F} \to \mathbb{Z}$ be the group homomorphism defined by $\ord_{W_{F}}(\Frob_{F}^{i} \cdot x) = i$ for $i \in \mathbb{Z}$ and $x \in I_{F}$.
For a finite Galois extension $E/F$, we denote by $\Gamma_{E/F}$ the corresponding Galois group and by $\Frob_{E/F}$ the image of $\Frob_{F}$ in $\Gamma_{E/F}$.

For a linear algebraic group $G$ over a field $k$ and an algebraic field extension $k'$ of $k$, we write $G_{k'}$ for the base change of $G$ to $k'$.
For a connected reductive group $G$, we denote by $G_{\der}$ the derived subgroup of $G$, by $G_{\mathrm{sc}}$ the simply connected cover of $G_{\der}$, and by $Z(G)$ the center of $G$.
For a closed subgroup $H$ of $G$, we denote by $N_{G}(H)$ (resp.\ $\Cent_{G}(H)$) the normalizer (resp.\ centralizer) of $H$ in $G$.
If $H = T$ is a maximal torus of $G$, we also denote by $W_{G}(T) := N_{G}(T)/T$ the Weyl group of $G$ with respect to $T$.

Suppose that $G$ is a connected reductive group defined over $F$.
We denote by $\widehat{G}$ the dual group of $G$ over $\mathbb{C}$ and by ${^LG}$ the $L$-group $\widehat{G} \rtimes W_{F}$ of $G$.
For a maximal torus $T$ of $G$, we denote by $X^{*}(T)$ and $X_{*}(T)$ the (algebraic) character group and cocharacter group of $T_{\overline F}$ equipped with the actions
of the absolute Galois group $\Gamma_{F}$, respectively.
Let $\Phi(G, T) \subset X^{*}(T)$ denote the root system of $G_{\overline{F}}$ with respect to $T_{\overline F}$.
For a root $\alpha \in \Phi(G, T)$, we denote by $\alpha^{\vee} \in X_{*}(T)$ the corresponding coroot.

Let $E = F$ or $F^{\unr}$.
We denote by $X^{*}_{E}(G)$ the group of algebraic characters that are defined over $E$.
We define a subgroup $G(E)^{1}$ of $G(E)$ as
\[
G(E)^{1} = \left\{
g \in G(E) \mid \text{$\ord_{E}(\lambda(g)) = 0$ for all $\lambda \in X^{*}_{E}(G)$}
\right\}.
\]
Let $X_{\unr}(G_{E})$ denote the group of characters of $G(E)$ that are trivial on $G(E)^{1}$.
We call an element $\chi$ of $X_{\unr}(G_{E})$ an \emph{unramified character} of $G(E)$.

We denote by $\mathcal{B}(G, E)$ the enlarged Bruhat--Tits building of $G_{E}$, and for a maximal $E$-split $F$-torus $T$ of $G$, we denote by $\mathcal{A}(G, T, E)$ the apartment of $T_{E}$ in $\mathcal{B}(G, E)$.
We also write $\mathcal{B}^{\red}(G, E)$ for the reduced building of $G_{E}$ and $\mathcal{A}^{\red}(G, T, E)$ for the apartment of $T_{E}$ in $\mathcal{B}^{\red}(G, E)$.
Let $x \in \mathcal{B}(G, F)$ and we write $[x]$ for the image of $x$ in $\mathcal{B}^{\red}(G, F)$.
We denote by $G(E)_{[x]}$ the stabilizer of $[x]$ in $G(E)$, by $G(E)_{x,0}$ the group of $\mathcal{O}_{E}$-points of the connected parahoric group scheme of $G$ associated to the point $x$, and by $G(E)_{x,0+}$ the pro-$p$ radical of $G(E)_{x,0}$.
We also denote by $\mathsf{G}_{[x]}$ the reductive quotient of the special fiber of the smooth integral model of $G$ whose group of $\mathcal{O}_{F^{\unr}}$-points is $G(F^{\unr})_{[x]}$.
Let $\mathsf{G}_{[x]}^{\circ}$ be its neutral connected component, which is the reductive quotient of the special fiber of the parahoric group scheme of $G$ associated to $[x]$.
Thus, we have $\mathsf{G}_{[x]}(\resi) = G(F)_{[x]}/G(F)_{x, 0+}$ and $\mathsf{G}^{\circ}_{[x]}(\resi) = G(F)_{x,0} / G(F)_{x, 0+}$.
In the case that $G$ is a torus, we omit the notation $x$ and write $G(E)_{x, 0} = G(E)_{0}$, $G(E)_{x, 0+} = G(E)_{0+}$, $\mathsf{G}_{[x]} = \mathsf{G}$, and $\mathsf{G}^{\circ}_{[x]} = \mathsf{G}^{\circ}$, which do not depend on $x$.

We denote by $\Rep(G(F))$ the category of smooth complex representations of $G(F)$.
Let $\mathfrak{I}(G)$ denote the set of inertial equivalence classes, i.e.,\ equivalence classes $[L,\sigma]_G$ of cuspidal pairs $(L,\sigma)$ in $G$, where $L$ is a Levi subgroup of $G$, $\sigma$ is an irreducible supercuspidal representation of $L(F)$, and where the equivalence is given by conjugation by $G(F)$ and twisting by unramified characters of $L(F)$.
For each $\mathfrak{s} \in \mathfrak{I}(G)$, let $\Rep^{\mathfrak{s}}(G(F))$ denote the corresponding Bernstein block of $\Rep(G(F))$, namely the full subcategory of $\Rep(G(F))$ consisting of smooth representations $(\pi, V)$ such that any irreducible subquotient $\pi'$ of $\pi$ occurs as a subquotient of a parabolically induced representation from an irreducible supercuspidal representation $\sigma$ of a Levi subgroup $L(F)$ of $G(F)$, where $(L, \sigma) \in \mathfrak{s}$.
If $\mathfrak{s} = [G, \pi]_{G}$ for an irreducible supercuspidal representation $\pi$ of $G(F)$, we call the corresponding Bernstein block $\Rep^{\mathfrak{s}}(G(F))$ a \emph{supercuspidal block}.

For a group $G$, we denote by $G^{*}$ the group of characters $\chi \colon G \to \mathbb{C}^{\times}$.
For $g \in G$, we denote by $\Ad(g)$ the adjoint action $x \mapsto gxg^{-1}$ of $g$ on $G$.
We denote by $\mathbb{C}[G]$ the group algebra of $G$ over $\mathbb{C}$.
For $g \in G$, we denote by $b_{g} \in \mathbb{C}[G]$ the corresponding basis element.
Suppose that $G$ acts on a set $X$.
We denote by $X^{G}$ the subset of $G$-fixed points:
\[
X^G = \left\{
x \in X \mid \text{$g \cdot x = x$ for all $g \in G$}
\right\}.
\]
If $X$ carries an algebraic structure (e.g.,\ a ring, module, algebra, etc.) and the $G$-action preserves this structure, then $X^{G}$ inherits the same type of structure.
For a $G$-module $X$, we also denote by $X_{G}$ the group of $G$-coinvariants:
\[
X_{G} = X/\left \langle 
g \cdot x - x \mid \text{$g \in G$ and $x \in X$}
\right \rangle.
\]
Suppose that $H$ is a subgroup of a group $G$ and $g \in G$. We denote $gHg^{-1}$ by $^gH$. If $\rho$ is a representation of $H$, we write $^g\rho$ for the representation $x\mapsto \rho(g^{-1}xg)$ of $^gH$. 
We define
\[
N_{G}(H) = \left\{
g \in G \mid {^gH} = H
\right\}
\quad
\text{ and }
\quad
N_{G}(\rho) = \left\{
g \in N_{G}(H) \mid {^g\rho} \simeq \rho
\right\}.
\]
We denote by $\Irr(G, \rho)$ the set of irreducible representations of $G$ whose restriction to $H$ contains $\rho$.
We note that if $H$ is a normal subgroup of $G$ of finite index and $\rho$ is normalized by $G$, the set $\Irr(G, \rho)$ is the set of irreducible representations of $G$ whose restriction to $H$ is $\rho$-isotypic.

Suppose that $K$ is an open subgroup of a locally profinite group $G$.
For a smooth representation $(\rho, V_{\rho})$ of $K$, we denote by 
\[
\left(
\ind_{K}^{G} (\rho), \ind_{K}^{G} (V_{\rho})
\right)
\]
the compactly induced representation of $(\rho, V_{\rho})$.

For a representation $(\rho, V_{\rho})$ of a group $G$, we identify $\rho$ with its representation space $V_{\rho}$ by abuse of notation. 
For any vector space or module $V$, we write $\id_V$ for the identity map on $V$.

For a unital ring $A$, we denote by $A \lMod$ the category of left unital $A$-modules and by $\rMod A$ the category of right unital $A$-modules.

For an algebraic group $G$, we denote by $G^{\circ}$ its neutral connected component and by $\pi_0(G)$ the group of components.
For a finitely generated abelian group $X$, we denote by $X_{\mathrm{tor}}$ its torsion subgroup.

Throughout the paper, we fix a prime number $\ell \neq p$ and an isomorphism $\mathbb{C} \simeq \overline{\mathbb{Q}_{\ell}}$.
All representations in this paper are on vector spaces over $\mathbb{C}$ (or over $\overline{\mathbb{Q}_{\ell}}$, which we identify with $\mathbb{C}$ via the fixed isomorphism above).

\section{A cohomology class attached to a supercuspidal block}
\label{sec:A cohomology class attached to a supercuspidal block}

Let $G$ be a connected reductive group over $F$ and let $\Rep^{\mathfrak{s}}(G(F))$ be a supercuspidal block.
Then we have a $\mathbb{C}$-algebra $\mathcal{H}^{\mathfrak{s}}$ and an explicit equivalence of categories 
\[
\Rep^{\mathfrak{s}}(G(F)) \simeq \rMod \mathcal{H}^{\mathfrak{s}}.
\]
In this section, we realize $\mathcal{H}^{\mathfrak{s}}$ as an Azumaya algebra over its center $Z^{\mathfrak{s}}$.
Morita equivalence classes of Azumaya algebras over $Z^{\mathfrak{s}}$ are parametrized by certain cohomology classes of 2-cocycles.
We provide an explicit method for computing the cohomology class corresponding to the Morita equivalence class of $\mathcal{H}^{\mathfrak{s}}$.

\subsection{Basic facts about Azumaya algebras}
\label{subsec:Basic facts about Azumaya algebras}

In this subsection, we recall basic facts about Azumaya algebras.
In particular, we explain the classification of Morita equivalence classes of Azumaya algebras in terms of cohomology classes of 2-cocycles.

Let $R_{0}$ be a noetherian commutative ring and $A$ be an $R_{0}$-algebra that is a finite $R_{0}$-module.
We call $A$ an \emph{Azumaya algebra} over $R_{0}$ if there exists a faithfully flat commutative $R_{0}$-algebra $R$ such that $R \otimes_{R_0} A$ is isomorphic, as $R$-algebras, to a matrix algebra $M_{m}(R)$ for some $m \ge 1$.
In this case, we say that the algebra $R$ \emph{splits} the Azumaya algebra $A$.
For an Azumaya algebra $A$ over $R_0$, we denote by $A^{\op}$ the opposite algebra to $A$.

Let $R$ be a noetherian commutative ring whose Picard group $\Pic(R)$ is trivial.
Let $\Fgrp$ be a finite abelian group acting on $R$.
We write $R_0 := R^{\Fgrp}$ and assume that the inclusion map $R_0 \to R$ is faithfully flat.
We also assume the following descent conditions.
\begin{assumption}
\label{assumptionfordescent}
\begin{enumerate}
\item
\label{assumptionfordescent1}
For any $R_0$-module $M_0$, the map 
\[
M_0 \to R \otimes_{R_0} M_0, \qquad m_0 \mapsto 1 \otimes m_{0}
\]
induces an isomorphism of $R_0$-modules $M_0 \isoarrow (R \otimes_{R_0} M_0)^{\Fgrp}$, where 
\[
(R \otimes_{R_0} M_0)^{\Fgrp}
:=
\left\{
m \in R \otimes_{R_0} M_0 \mid \text{$(\chi \otimes \id_{M_0})(m) = m$ for all $\chi \in \Fgrp$}
\right\}.
\]
\item
\label{assumptionfordescent2}
For any $R$-module $M$ with a semilinear action of $\Fgrp$, the map
\[
R \otimes_{R_0} M^{\Fgrp} \to M, \qquad r \otimes m \mapsto rm
\]
is an isomorphism of $R$-modules.
\end{enumerate}
\end{assumption}
In the following subsections, we will apply the result of this subsection to the following setting. 
\begin{lemma}
\label{lemmahowtoapplysubsec1.1in1.2}
Let $N$ be a finitely generated abelian group and $T$ be a subgroup of $N$ of finite index.
We set $R = \mathbb{C}[N]$.
We define an action of the finite group $\Fgrp := (N/T)^{*}$ on $R$ by
\[
 \chi(b_{n}) = \chi(n) \cdot b_{n}, \qquad \chi \in (N/T)^{*}, n \in N.
 \]
 In this setting, the map $R_0 \to R$ is faithfully flat, $\Pic(R)$ is trivial, and Assumption~\ref{assumptionfordescent} is satisfied.
\end{lemma}
\begin{proof}
By the definition of the action of $\Fgrp$ on $R$, we have $R_{0} = \mathbb{C}[T]$.
Hence, as an $R_0$-module, we have 
\[
R = \bigoplus_{w \in N/T} R_0 \cdot b_{\widetilde{w}},
\]
where $\widetilde{w} \in N$ is a lift of $w$.
In particular, $R$ is a free $R_0$-module, which implies that $R_0 \to R$ is faithfully flat.

There exists an isomorphism $N\simeq \Z^{r}\times N_{0}$ where $N_{0}$ is a finite abelian group.
Hence $\C[N]\simeq \C[\Z^{r}]\otimes_{\C}\C[N_{0}]\simeq \C[\Z^{r}]\otimes_{\C} \prod_{\chi\in N_{0}^{*}}\C\simeq\prod_{\chi\in N_{0}^{*}} \C[\Z^{r}]$.
Since $\C[\Z^{r}]$ is isomorphic to a Laurent polynomial ring, we have $\Pic(\C[\Z^{r}])=0$.
Hence $\Pic(\C[N]) = 0$.

We prove that Assumption~\ref{assumptionfordescent} is satisfied in this setting.
Let $M_0$ be an $R_0$-module.
Then we have an isomorphism of $R_0$-modules
\[
I \colon R \otimes_{R_0} M_0 \simeq \bigoplus_{w \in N/T} M_0, \qquad \sum_{w \in N/T} b_{\widetilde{w}} \otimes m_{w} \mapsto (m_{w})_{w \in N/T}.
\]
Under this isomorphism, the action $\chi \mapsto \chi \otimes \id_{M_0}$ of $\Fgrp$ on $R \otimes_{R_0} M_0$ corresponds to the action of $\Fgrp$ on $\bigoplus_{w \in N/T} M_0$ defined as 
\[
\chi \colon (m_{w})_{w \in N/T} \mapsto (\chi(w) \cdot m_{w})_{w \in N/T}.
\]
Hence, we conclude that 
\begin{align*}
(R \otimes_{R_0} M_0)^{\Fgrp} &= I^{-1}\left(
\left\{
(m_{w})_{w \in N/T} \mid \text{$m_{w} = 0$ for all $w \neq 0$}
\right\}
\right) \\
&= \left\{
1 \otimes m_0 \mid \text{$m_0 \in M_0$}
\right\}.
\end{align*}
Thus, the first condition of Assumption~\ref{assumptionfordescent} is satisfied.

Let $M$ be an $R$-module with a semilinear action of $\Fgrp$.
Then we have
\[
R \otimes_{R_0} M^{\Fgrp} = \bigoplus_{w \in N/T} R_{0} \cdot b_{\widetilde{w}} \otimes_{R_0} M^{\Fgrp}
\simeq \bigoplus_{w \in N/T} b_{\widetilde{w}} \cdot M^{\Fgrp}.
\]
Since the action of $\Fgrp$ on $M$ is semilinear, we have
\[
b_{\widetilde{w}} \cdot M^{\Fgrp} = \left\{
m \in M \mid \text{$\chi \cdot m = \chi(w) \cdot m$ for all $\chi \in \Fgrp$}
\right\}.
\]
Since $\Fgrp = (N/T)^{*}$ is a finite abelian group, we have the eigenspace decomposition
\[
M = \bigoplus_{w \in \left(
(N/T)^{*}
\right)^{*} = N/T} \left\{
m \in M \mid \text{$\chi \cdot m = \chi(w) \cdot m$ for all $\chi \in \Fgrp$}
\right\} = \bigoplus_{w \in N/T} b_{\widetilde{w}} \cdot M^{\Fgrp}.
\]
Thus, the second condition of Assumption~\ref{assumptionfordescent} is satisfied.
\end{proof}

Let $A$ be an Azumaya algebra over $R_0$ that is split by $R$.
We attach an element $c_{A} \in H^2(\Fgrp, R^{\times})$ to $A$ as follows.
We fix an isomorphism $i_{A} \colon R \otimes_{R_0} A \simeq M_{m}(R)$.
For $\chi \in \Fgrp$, we define a ring automorphism
$
M_{m}(\chi) \colon M_{m}(R) \to M_{m}(R)
$
by applying $\chi$ to each entry of $M_{m}(R)$.
We also define a group automorphism $\PGL_{m}(\chi) \colon \PGL_{m}(R) \to \PGL_{m}(R)$ by applying $\chi$ to each entry.
Since $\Pic(R) = \{0\}$, according to \cite[Chapter~3, Corollary~4.8]{MR279159}, we have $\Aut_{R \Alg}\left(
M_{m}(R)
\right) = \PGL_{m}(R)$, where $\Aut_{R \Alg}\left(
M_{m}(R)
\right)$ denotes the group of $R$-algebra automorphisms of $M_{m}(R)$.
Under this identification, we have 
\[
M_{m}(\chi) \circ g \circ (M_{m}(\chi))^{-1} = \PGL_{m}(\chi)(g) 
\]
for $g \in \PGL_{m}(R)$ and $\chi \in \Fgrp$.
We regard $M_{m}(R)$ and $\PGL_{m}(R)$ as $\Fgrp$-modules via $M_m(\chi)$ and $\PGL_{m}(\chi)$, respectively.
We define another action $a_{i_{A}}$ of $\Fgrp$ on $M_{m}(R)$ as 
\[
a_{i_{A}}(\chi) =  i_{A} \circ (\chi \otimes \id_{A}) \circ i_{A}^{-1}
\]
for $\chi \in \Fgrp$.
Using this, we define an element $z_{i_{A}} \in Z^1\left(\Fgrp, \PGL_{m}(R)\right)$ by
\[
a_{i_{A}}(\chi) = z_{i_{A}}(\chi) \circ M_{m}(\chi)
\]
for $\chi \in \Fgrp$.
We claim that the cohomology class $[z_{i_{A}}] \in H^1\left(\Fgrp, \PGL_{m}(R)\right)$ of $z_{i_{A}}$ does not depend on the choice of an isomorphism $i_{A}$.
Let $i'_{A} \colon R \otimes_{R_0} A \simeq M_{m}(R)$ be another isomorphism. 
Then there exists $g \in \PGL_{m}(R)$ such that $i'_{A} = g \circ i_{A}$.
By the definitions of $a_{i_{A}}$ and $z_{i_{A}}$, we obtain that 
\[
a_{i'_{A}}(\chi) = g \circ a_{i_{A}}(\chi) \circ g^{-1}
\qquad
\text{and}
\qquad
z_{i'_{A}}(\chi) = g \circ z_{i_{A}}(\chi) \circ \PGL_{m}(\chi)(g)^{-1}
\]
for $\chi \in \Fgrp$.
Thus, we have $[z_{i_{A}}] = [z_{i'_{A}}]$, as desired.
Now, we define the element $c_{A} \in H^2(\Fgrp, R^{\times})$ as the image of $[z_{i_{A}}]$ via the connecting homomorphism $H^1\left(\Fgrp, \PGL_{m}(R)\right) \to H^2(\Fgrp, R^{\times})$ obtained from the short exact sequence 
\[
1 \to R^{\times} \to \GL_{m}(R) \to \PGL_{m}(R) \to 1
\]
of $X$-modules.
\begin{lemma}
\label{lemma:cdefinesgrouphom}
Let $A_1$ and $A_2$ be Azumaya algebras over $R_0$ that are split by $R$.
Then we have $c_{A_1^{\op}} = c_{A_1}^{-1}$ and $c_{A_1 \otimes_{R_0} A_2} = c_{A_1} \cdot c_{A_2}$.
\end{lemma}
\begin{proof}
The claims follow from the definition of $c_{A}$.
\end{proof}

We say that two Azumaya algebras $A_1$ and $A_2$ over $R_0$ are \emph{Morita equivalent} if the $R_0$-categories $\rMod A_1$ and $\rMod A_2$ are equivalent.
Following \cite[Chapter~3, Section~4]{MR279159}, let $\Br(R_0)$ denote the set of Morita equivalence classes of Azumaya algebras over $R_0$.
According to \cite[Chapter~3, Proposition and Definition~4.10]{MR279159}, the tensor product $\otimes_{R_0}$ induces a structure of an abelian group on $\Br(R_0)$.
We call this abelian group the \emph{Brauer group} of $R_0$.
For an Azumaya algebra $A$ over $R_0$, we denote by $[A] \in \Br(R_0)$ its Morita equivalence class.

\begin{proposition}
\label{prop:MoritaequivclassofAzumayaalgebras}
Let $A_1$ and $A_2$ be Azumaya algebras over $R_0$ that are split by $R$.
Then we have $c_{A_1} = c_{A_2}$ if and only if $A_1$ and $A_2$ are Morita equivalent.
\end{proposition}

Although this proposition follows from \cite[Corollary~1]{MR586870} after a suitable translation of the terminology into our setting, we include a brief proof for completeness.
We first consider the following special case.

\begin{lemma}
\label{lemmaforprop:MoritaequivclassofAzumayaalgebras}
Let $A$ be an Azumaya algebra over $R_0$ that is split by $R$.
Then we have $c_{A} = 1$ if and only if $A$ is Morita equivalent to $R_0$.
\end{lemma}
\begin{proof}
Suppose that $A$ is Morita equivalent to $R_0$.
According to \cite[Chapter~3, Proposition and Definition~4.10]{MR279159}, there exists a faithfully projective $R_{0}$-module $P_{0}$ such that $A \simeq \End_{R_0}(P_0)$.
Here, a faithfully projective $R_0$-module means a finitely generated projective module with trivial annihilator.
We prove that $R \otimes_{R_0} P_0 \simeq R^{\oplus m}$ for some $m$.
Since $P_0$ is a finitely generated projective $R_0$-module and $R$ splits $A$, we have 
\[
\End_{R}\left(
R \otimes_{R_0} P_0
\right) \simeq R \otimes_{R_0} \End_{R_0}(P_0) \simeq M_{m}(R)
\]
for some $m \ge 1$.
For $1 \le i \le m$, let $e_{i}$ denote the endomorphism of $R \otimes_{R_0} P_0$ corresponding to $E_{ii} \in M_{m}(R)$, where $E_{ii}$ denotes the matrix whose $(i, i)$-entry is $1$ and all other entries are $0$.
We define a submodule $P_{i}$ of $R \otimes_{R_0} P_0$ as $P_i = e_{i}(R \otimes_{R_0} P_0)$.
Then we have 
\[
R \otimes_{R_0} P_0 = \bigoplus_{1 \le i \le m} P_{i}.
\]
Moreover, for each $1 \le i \le m$, we have 
\[
\End_{R}(P_i) = e_i \End_{R}(R\otimes_{R_{0}}P_{0}) e_i \simeq E_{ii} M_{m}(R) E_{ii} \simeq R.
\]
Hence, $P_{i}$ is a projective $R$-module of rank one.
Since $\Pic(R) = \{0\}$, we have $P_{i} \simeq R$ for each $1 \le i \le m$.
Thus, we obtain that $R \otimes_{R_0} P_0 \simeq R^{\oplus m}$.
We fix an isomorphism $i_{P_0} \colon R \otimes_{R_0} P_0  \isoarrow R^{\oplus m}$.
The isomorphism $i_{P_0}$ induces an isomorphism 
\[
i_{A} \colon R \otimes_{R_0} \End_{R_0}(P_0) \simeq \End_{R}\left(
R \otimes_{R_0} P_0
\right) \isoarrow \End_{R}\left(
R^{\oplus m}
\right) = M_{m}(R).
\]
 We define an element $z_{P_0} \in Z^1(\Fgrp, \GL_{m}(R))$ by
\[
i_{P_0} \circ (\chi \otimes \id_{P_0}) \circ i_{P_0}^{-1} = z_{P_0}(\chi) \circ \chi^{\oplus m}
\]
for $\chi \in \Fgrp$.
Then the image of $z_{P_0}$ in $Z^1(\Fgrp, \PGL_{m}(R))$ agrees with the element $z_{i_{A}}$ attached to $i_{A}$.
Since the sequence
\[
H^1(\Fgrp, \GL_m(R)) \to H^1(\Fgrp, \PGL_m(R)) \to H^2(\Fgrp, R^{\times})
\]
is exact, we obtain that $c_{A} = 1$ in this case.

Conversely, suppose that $c_{A} = 1$.
We fix an isomorphism $i_{A} \colon R \otimes_{R_0} A \simeq M_{m}(R)$.
By Assumption~\ref{assumptionfordescent}~\eqref{assumptionfordescent1}, $i_{A}$ restricts to an isomorphism of $R_{0}$-algebras $A \simeq M_{m}(R)^{(\Fgrp, z_{i_{A}})}$, where 
\[
M_{m}(R)^{(\Fgrp, z_{i_{A}})} := \left\{
x \in M_{m}(R) \mid \text{$(z_{i_{A}}(\chi) \circ M_m(\chi))(x) = x$ for all $\chi \in \Fgrp$}
\right\}.
\]
Since $c_{A} = 1$, there exists a lift $\widetilde{z}_{i_{A}} \in Z^1\left(
\Fgrp, \GL_{m}(R)
\right)$ of $z_{i_{A}} \in Z^1\left(
\Fgrp, \PGL_{m}(R)
\right)$.
We define an $R_{0}$-module $P_0$ as 
\[
P_0 = \left\{
x \in R^{\oplus m} \mid \text{$(\widetilde{z}_{i_{A}}(\chi) \circ \chi^{\oplus m})(x) = x$ for all $\chi \in \Fgrp$}
\right\}.
\]
By Assumption~\ref{assumptionfordescent}~\eqref{assumptionfordescent2}, the map 
\[
R \otimes_{R_0} P_0 \to R^{\oplus m}, \quad r \otimes x \mapsto rx
\]
is an isomorphism of $R$-modules, which induces an isomorphism 
\[
R \otimes_{R_0} \End_{R_0}(P_0) \simeq \End_{R}\left(
R \otimes_{R_0} P_0
\right) \isoarrow \End_{R}\left(
R^{\oplus m}
\right) = M_{m}(R).
\]
By the definition of $P_0$, the action of $\Fgrp$ on $R \otimes_{R_0} \End_{R_0}(P_0)$ induced by its action on the first factor corresponds to the action of $\Fgrp$ on $ M_{m}(R)$ given by $\chi \mapsto z_{i_{A}}(\chi) \circ M_m(\chi)$.
Using Assumption~\ref{assumptionfordescent}~\eqref{assumptionfordescent1} again, we obtain an isomorphism of $R_0$-algebras 
\[
\End_{R_0}(P_0) \simeq \left(
R \otimes_{R_0} \End_{R_0}(P_0)
\right)^{\Fgrp} \isoarrow M_{m}(R)^{(\Fgrp, z_{i_{A}})}  \simeq A.
\]
Since $R \otimes_{R_0} P_0 \simeq R^{\oplus m}$, we obtain from the faithfully flat descent that $P_0$ is a finitely generated projective $R_0$-module.
Moreover, since the annihilator of $P_0$ is contained in the annihilator of $R \otimes_{R_0} P_0$, we also obtain that $P_0$ has zero annihilator.
Thus, the $R_0$-module $P_0$ is faithfully projective, and the $R_0$-algebra $A \simeq \End_{R_0}(P_0)$ is Morita equivalent to $R_0$ by \cite[Chapter~3, Proposition and Definition~4.10]{MR279159}.
\end{proof}

\begin{proof}[Proof of Proposition~\ref{prop:MoritaequivclassofAzumayaalgebras}]

Assume that $A_1$ and $A_2$ are Morita equivalent.
Then we have 
\[
[A^{\op}_{1} \otimes_{R_0} A_{2}] = [A^{\op}_{1} \otimes_{R_0} A_{1}].
\]
According to Lemma~\ref{lemma:cdefinesgrouphom}, we have $c_{A^{\op}_{1} \otimes_{R_0} A_{1}} = 1$.
Hence, we obtain from Lemma~\ref{lemmaforprop:MoritaequivclassofAzumayaalgebras} that 
\[
[A^{\op}_{1} \otimes_{R_0} A_{2}] = [A^{\op}_{1} \otimes_{R_{0}} A_{1}] = [R_0].
\]
Thus, applying Lemmas~\ref{lemma:cdefinesgrouphom} and \ref{lemmaforprop:MoritaequivclassofAzumayaalgebras} again, we conclude that 
\[
c_{A_1}^{-1} \cdot c_{A_2} = c_{A^{\op}_{1} \otimes_{R_0} A_{2}} = 1, 
\]
equivalently, $c_{A_1} = c_{A_2}$.

Conversely, assume that $c_{A_1} = c_{A_2}$.
According to Lemmas~\ref{lemma:cdefinesgrouphom} and \ref{lemmaforprop:MoritaequivclassofAzumayaalgebras}, we have $[A_1 \otimes_{R_0} A_1^{\op}] = [R_0]$.
Thus, we obtain that 
\[
[A_1 \otimes_{R_0} A_1^{\op} \otimes_{R_0} A_2] = [R_0 \otimes_{R_0} A_2] = [A_2].
\]
Since $c_{A_1} = c_{A_2}$, according to Lemma~\ref{lemma:cdefinesgrouphom}, we have $c_{A_1^{\op} \otimes_{R_0} A_2} = 1$.
Applying Lemma~\ref{lemmaforprop:MoritaequivclassofAzumayaalgebras} again, we have 
\[
[A_1 \otimes_{R_0} A_1^{\op} \otimes_{R_0} A_2] = [A_1 \otimes_{R_0} R_0] = [A_1].
\]
Thus, we conclude that $[A_1] = [A_2]$.
\end{proof}

\subsection{Twisted group algebra as an Azumaya algebra}
\label{subsec:Twisted group algebra as an Azumaya algebra}

In this subsection, we realize a twisted group algebra as an Azumaya algebra over its center and compute the cohomology class of a 2-cocycle corresponding to its Morita equivalence class. 
The results of this subsection will be used in Section~\ref{subsection:Azumaya algebra attached to a supercuspidal block} to compute the cohomology class attached to a supercuspidal block.

Let $N$ be a finitely generated abelian group and $\mu \colon N \times N \to \mathbb{C}^{\times}$ be a 2-cocycle.
We define the twisted group algebra $\mathbb{C}[N, \mu]$ as the $\mathbb{C}$-vector space 
\[
\mathbb{C}[N, \mu] = \bigoplus_{x \in N} \mathbb{C} \cdot b_{x}
\]
equipped with the $\mathbb{C}$-algebra structure defined by $b_{x} b_{y} = \mu(x, y) \cdot b_{xy}$ for $x, y \in N$.
We note that the isomorphism class of $\mathbb{C}[N, \mu]$ only depends on the cohomology class $[\mu]$ of $\mu$.
For $x, y \in N$, we define $c(x, y) = c_{\mu}(x, y) := \mu(y, x) \mu(x, y)^{-1}$, which only depends on the cohomology class $[\mu]$.
\begin{lemma}
\label{lemma:alternatingformdeterminestwococycle}
The map $c_{\mu}$ is an alternating bilinear form on $N \times N$.
Moreover, the map $\mu \mapsto c_{\mu}$ induces a bijection between $H^2(N, \mathbb{C}^{\times})$ and the set of alternating bilinear forms $N \times N \to \mathbb{C}^{\times}$.
\end{lemma}
\begin{proof}
The first claim is \cite[Theorem~2.8 (i)]{MR2914869} and the second claim follows from \cite[Theorem~2.8 (ii)--(iv)]{MR2914869}.
\end{proof}

We denote by $T$ the radical of $c$.
Then the subalgebra $\mathbb{C}[T, \mu] := \bigoplus_{x \in T} \mathbb{C} \cdot b_{x}$ is the center of $\mathbb{C}[N, \mu]$.
From now on, we assume the following condition.
\begin{assumption}
We have $\#(N/T) < \infty$.
\end{assumption}
The alternating form $c$ induces a symplectic form $(N/T) \times (N/T) \to \mathbb{C}^{\times}$, which we also denote by $c$.
Let $N/T = X^{+} \oplus X^{-}$ be a decomposition into Lagrangian subgroups; that is, $X^{+}$ and $X^{-}$ are subgroups such that $c$ is trivial on both $X^{+}$ and $X^{-}$.
The existence of such a decomposition follows from \cite[Theorem~4.1]{MR2598895}.
Hence, we have $c(x_1^{+}, x_2^{+}) = c(x_1^{-}, x_2^{-}) = 1$ for all $x_1^{+}, x_{2}^{+} \in X^{+}$ and $x_{1}^{-}, x_2^{-} \in X^{-}$, and $c$ induces a perfect pairing on $X^{+} \times X^{-}$.
We denote by $N^{+}$ (resp.\ $N^{-}$) the preimage of $X^{+}$ (resp.\ $X^{-}$) in $N$.
We define a $\mathbb{C}$-algebra $\mathbb{C}[N^{+}] \otimes^{c}_{\mathbb{C}[T]} \mathbb{C}[N^{-}]$ as follows:
\begin{enumerate}
\item
As a $\mathbb{C}$-vector space, we have $\mathbb{C}[N^{+}] \otimes^{c}_{\mathbb{C}[T]} \mathbb{C}[N^{-}] = \mathbb{C}[N^{+}] \otimes_{\mathbb{C}[T]} \mathbb{C}[N^{-}]$.
\item
The $\mathbb{C}$-algebras $\mathbb{C}[N^{+}]$ and $\mathbb{C}[N^{-}]$ are subalgebras.
\item
For $n^{+} \in N^{+}$ and $n^{-} \in N^{-}$, we have 
\[
(b_{n^+} \otimes 1) (1 \otimes b_{n^-}) = (b_{n^+} \otimes b_{n^-})
\qquad 
\text{and}
\qquad
(1 \otimes b_{n^-}) (b_{n^+} \otimes 1) = c(n^+, n^-) \cdot (b_{n^+} \otimes b_{n^-}).
\]
\end{enumerate}
In other words, $\mathbb{C}[N^{+}] \otimes^{c}_{\mathbb{C}[T]} \mathbb{C}[N^{-}]$ is the twisted group algebra $\mathbb{C}[N, \mu']$ attached to the 2-cocycle $\mu'$ on $N \times N$ defined as follows:
For $n^{+}_{1}, n^{+}_{2} \in N^{+}$ and $n^{-}_{1}, n^{-}_{2} \in N^{-}$, we have
\[
\mu'(n^{+}_{1}n^{-}_{1}, n^{+}_{2}n^{-}_{2}) = c(n^{+}_{2}, n^{-}_{1}).
\]

\begin{lemma}
\label{lemma:replacingmu}
We have $[\mu] = [\mu']$.
Hence, we have an isomorphism of $\mathbb{C}$-algebras
\[
\mathbb{C}[N^{+}] \otimes^{c}_{\mathbb{C}[T]} \mathbb{C}[N^{-}] \simeq \mathbb{C}[N, \mu].
\]
\end{lemma}
\begin{proof}
For $n^{+}_{1}, n^{+}_{2} \in N^{+}$ and $n^{-}_{1}, n^{-}_{2} \in N^{-}$, we have
\begin{align*}
c_{\mu'}(n^{+}_{1}n^{-}_{1}, n^{+}_{2}n^{-}_{2}) &= \mu'(n^{+}_{2}n^{-}_{2}, n^{+}_{1}n^{-}_{1}) \mu'(n^{+}_{1}n^{-}_{1}, n^{+}_{2}n^{-}_{2})^{-1} \\
&= c(n^{+}_{1}, n^{-}_{2}) c(n^{+}_{2}, n^{-}_{1})^{-1} \\
&= c(n^{+}_{1}, n^{-}_{2}) c(n^{-}_{1}, n^{+}_{2}) \\
&= c(n^{+}_{1}n^{-}_{1}, n^{+}_{2}n^{-}_{2}) = c_{\mu}(n^{+}_{1}n^{-}_{1}, n^{+}_{2}n^{-}_{2}).
\end{align*}
Thus, the claim follows from Lemma~\ref{lemma:alternatingformdeterminestwococycle}.
\end{proof}

We prove that $\mathbb{C}[N^{+}] \otimes^{c}_{\mathbb{C}[T]} \mathbb{C}[N^{-}]$ is an Azumaya algebra over its center $\mathbb{C}[T]$ that is split by $\mathbb{C}[N]$.
Let $V^{+}$ denote the $\mathbb{C}$-vector space whose basis $\{v_{x^{+}}\}$ is indexed by $X^{+}$.
We set $V^{+}_{\mathbb{C}[N]} := \mathbb{C}[N] \otimes_{\mathbb{C}} V^{+}$ and define a morphism of $\mathbb{C}[N]$-modules
\[
i_{c} \colon \mathbb{C}[N] \otimes_{\mathbb{C}[T]} \left(
\mathbb{C}[N^{+}] \otimes^{c}_{\mathbb{C}[T]} \mathbb{C}[N^{-}]
\right)
 \to \End_{\mathbb{C}[N]}\left(
 V^{+}_{\mathbb{C}[N]}
 \right)
\]
by
\[
i_{c}(b_{n} \otimes (b_{n^+} \otimes b_{n^-}))(v_{x^{+}}) = c(x^{+}, \overline{n}^{-}) b_{n n^{+} n^{-}} v_{\overline{n}^{+} x^{+}}
\]
for $n \in N$, $n^{+} \in N^{+}$, $n^{-} \in N^{-}$, $x^{+} \in X^{+}$, where $\overline{n}^{+}$ and $\overline{n}^{-}$ denote the images of $n^{+}$ and $n^{-}$ in $N/T$, respectively.
By the definition of the algebra structure on $\mathbb{C}[N^{+}] \otimes^{c}_{\mathbb{C}[T]} \mathbb{C}[N^{-}]$, $i_{c}$ is a morphism of $\mathbb{C}[N]$-algebras.
We index the elements of $X^{+}$ and $X^{-}$ as $X^{+} = \left\{x^{+}_1, x^{+}_2, \ldots , x^{+}_m\right\}$ and $X^{-} = \left\{x^{-}_1, x^{-}_2, \ldots , x^{-}_m\right\}$, where $m = \# X^+ = \# X^-$.
Using this, we identify $V^{+} \simeq \mathbb{C}^{\oplus m}$, $V^{+}_{\mathbb{C}[N]} \simeq \mathbb{C}[N]^{\oplus m}$, $\End_{\mathbb{C}}(V^{+}) \simeq M_{m}(\mathbb{C})$, and $\End_{\mathbb{C}[N]}\left(V^{+}_{\mathbb{C}[N]}\right) \simeq M_{m}\left(\mathbb{C}[N]\right)$.
For each $1 \le i \le m$, we fix lifts $n^{+}_{i} \in N^{+}$ of $x^{+}_{i}$ and $n^{-}_{i} \in N^{-}$ of $x^{-}_{i}$.
Then the set $\left\{
b^{-1}_{n^{+}_i n^{-}_{j}} \otimes b_{n^{+}_{i}} \otimes b_{n^{-}_{j}} \mid 1 \le i, j \le m
\right\}$ is a $\mathbb{C}[N]$-basis of $\mathbb{C}[N] \otimes_{\mathbb{C}[T]} \left(
\mathbb{C}[N^{+}] \otimes^{c}_{\mathbb{C}[T]} \mathbb{C}[N^{-}]
\right)$.
We denote by $V_{c}$ the $\mathbb{C}$-subspace of $\mathbb{C}[N] \otimes_{\mathbb{C}[T]} \left(
\mathbb{C}[N^{+}] \otimes^{c}_{\mathbb{C}[T]} \mathbb{C}[N^{-}]
\right)$ spanned by $\left\{
b^{-1}_{n^{+}_i n^{-}_{j}} \otimes b_{n^{+}_{i}} \otimes b_{n^{-}_{j}} \mid 1 \le i, j \le m
\right\}$.
\begin{proposition}
\label{prop:CNmuisazumaya}
For $1 \le i, j \le m$ and $x^{+} \in X^{+}$, we have 
\[
i_{c}(b^{-1}_{n^{+}_i n^{-}_{j}} \otimes b_{n^{+}_{i}} \otimes b_{n^{-}_{j}})(v_{x^{+}}) = c(x^{+}, x^{-}_{j}) \cdot v_{x^{+}_{i} x^{+}}.
\]
Moreover, the induced map
\[
i_{c}\restriction_{V_{c}} \colon V_{c} \to \End_{\mathbb{C}}(V^{+}) 
\]
is an isomorphism of $\mathbb{C}$-vector spaces.
Thus, the $\mathbb{C}[N]$-algebra map $i_{c}$ is also an isomorphism.
\end{proposition}
\begin{proof}
\begin{subequations}
The first claim follows from the definition of $i_{c}$.
We claim that 
\begin{equation}
\label{eq:linearlyindepset}
\left\{
i_{c}(b^{-1}_{n^{+}_i n^{-}_{j}} \otimes b_{n^{+}_{i}} \otimes b_{n^{-}_{j}}) \mid 1 \le i, j \le m
\right\} \subset \End_{\mathbb{C}}(V^+)
\end{equation}
is linearly independent.
For each $1 \le i \le m$, we denote by $V^{+}_{i}$ the $\mathbb{C}$-subspace of $\End_{\mathbb{C}}(V^+)$ spanned by $i_{c}(b^{-1}_{n^{+}_i n^{-}_{j}} \otimes b_{n^{+}_{i}} \otimes b_{n^{-}_{j}}) $ for all $1 \le j \le m$.
We obtain from the first claim that the sum $\sum_{1 \le i \le m} V^{+}_i$ is direct.
Moreover, since $c$ induces a perfect pairing on $X^{+} \times X^{-}$, the set 
\[
\left\{
\left(
c(x^{+}_{k}, x^{-}_{j}) 
\right)_{1 \le k \le m}
\mid 1 \le j \le m
\right\}
\subset \mathbb{C}^{m}
\]
is linearly independent.
Hence, for each $1 \le i \le m$, the set 
\[
\left\{
i_{c}(b^{-1}_{n^{+}_i n^{-}_{j}} \otimes b_{n^{+}_{i}} \otimes b_{n^{-}_{j}}) \mid 1 \le j \le m
\right\} \subset \End_{\mathbb{C}}(V^+)
\]
is linearly independent.
Thus, we conclude that the set in \eqref{eq:linearlyindepset} is linearly independent.
Since $\dim_{\mathbb{C}}\left(
\End_{\mathbb{C}}(V^{+})
\right) = m^2$, it is a basis of $\End_{\mathbb{C}}(V^{+})$.
Therefore, we obtain that the morphism of $\mathbb{C}$-vector spaces 
\[
i_{c}\restriction_{V_{c}} \colon V_{c} \to \End_{\mathbb{C}}(V^{+}) 
\]
is an isomorphism.
The last claim follows from the fact that $i_{c}$ is the base change of $i_{c}\restriction_{V_{c}}$ to $\mathbb{C}[N]$.
\end{subequations}
\end{proof}
According to Proposition~\ref{prop:CNmuisazumaya}, $\mathbb{C}[N^{+}] \otimes^{c}_{\mathbb{C}[T]} \mathbb{C}[N^{-}]$ is an Azumaya algebra over $\mathbb{C}[T]$, and $\mathbb{C}[N]$ splits $\mathbb{C}[N^{+}] \otimes^{c}_{\mathbb{C}[T]} \mathbb{C}[N^{-}]$.
Applying the construction of $c_{A}$ in Section~\ref{subsec:Basic facts about Azumaya algebras} to $R = \mathbb{C}[N]$, $\Fgrp = (N/T)^{*}$, which acts on $\mathbb{C}[N]$ by 
\[
\chi(b_{n}) = \chi(n) \cdot b_{n}, \qquad \chi \in (N/T)^{*}, n \in N,
\]
we attach an element 
\[
c_{\mathbb{C}[N^{+}] \otimes^{c}_{\mathbb{C}[T]} \mathbb{C}[N^{-}]} \in H^2\left(
(N/T)^{*}, \mathbb{C}[N]^{\times}
\right)
\]
to the Azumaya algebra $\mathbb{C}[N^{+}] \otimes^{c}_{\mathbb{C}[T]} \mathbb{C}[N^{-}]$ over $\mathbb{C}[T] = \mathbb{C}[N]^{\Fgrp}$ (see Lemma~\ref{lemmahowtoapplysubsec1.1in1.2}).
In order to compute $c_{\mathbb{C}[N^{+}] \otimes^{c}_{\mathbb{C}[T]} \mathbb{C}[N^{-}]}$, we define $[\mu]^{\vee} \in H^2\left(
(N/T)^{*}, \mathbb{C}^{\times}
\right)$ as follows.
Since $c$ is a perfect pairing on $(N/T) \times (N/T)$, we have an isomorphism of groups 
\[
\phi_{c} \colon N/T \isoarrow (N/T)^{*}, \qquad x \mapsto [y \mapsto c(x, y)].
\]
We define $[\mu]^{\vee}$ as the image of $[\mu] = [\mu'] \in H^2(N/T, \mathbb{C}^{\times})$ via the isomorphism $H^2(
N/T, \mathbb{C}^{\times}) \isoarrow H^2\left(
(N/T)^{*}, \mathbb{C}^{\times}
\right)$ induced by $\phi_{c}$.

\begin{proposition}
\label{prop:cohomologyclassofCN+TN-}
The cohomology class $c_{\mathbb{C}[N^{+}] \otimes^{c}_{\mathbb{C}[T]} \mathbb{C}[N^{-}]}$ agrees with the image of $[\mu]^{\vee}$ under the natural map $H^2\left(
(N/T)^{*}, \mathbb{C}^{\times}
\right) \to H^2\left(
(N/T)^{*}, \mathbb{C}[N]^{\times}
\right)$.
\end{proposition}
\begin{proof}
Let  $\chi \in (N/T)^{*}$.
We write $\phi^{-1}_{c}(\chi) = \chi^{+} \cdot \chi^{-}$, where $\chi^{+} \in X^{+}$ and $\chi^{-} \in X^{-}$.
We define $g_{\chi} \in \Aut_{\mathbb{C}[N]}\left(
V^{+}_{\mathbb{C}[N]}
\right)$ as 
\[
g_{\chi}(v_{x^+}) = c(x^{+}, \chi^{-}) \cdot v_{\chi^{+} x^{+}}
\]
for $x^{+} \in X^{+}$.
Then for $1 \le i, j \le m$ and $x^{+} \in X^{+}$, we have 
\begin{align*}
&\quad  \left(
g_{\chi} \circ i_{c}(b^{-1}_{n^{+}_i n^{-}_{j}} \otimes b_{n^{+}_{i}} \otimes b_{n^{-}_{j}}) \circ g_{\chi}^{-1}
\right)(v_{x^{+}}) \\
&= 
c((\chi^{+})^{-1} x^{+}, \chi^{-})^{-1} \cdot 
\left(
g_{\chi} \circ i_{c}(b^{-1}_{n^{+}_i n^{-}_{j}} \otimes b_{n^{+}_{i}} \otimes b_{n^{-}_{j}})
\right)(v_{(\chi^{+})^{-1} x^{+}}) \\
&= c((\chi^{+})^{-1} x^{+}, \chi^{-})^{-1} c((\chi^{+})^{-1} x^{+}, x^{-}_{j}) \cdot g_{\chi}(v_{x^{+}_{i} (\chi^{+})^{-1} x^{+}}) \\
&= c((\chi^{+})^{-1} x^{+}, \chi^{-})^{-1} c((\chi^{+})^{-1} x^{+}, x^{-}_{j}) c(x^{+}_{i} (\chi^{+})^{-1} x^{+}, \chi^{-}) \cdot v_{x^{+}_{i} x^{+}} \\
&= c((\chi^+)^{-1}, x^{-}_{j}) c(x^{+}, x_j^{-}) c(x_i^{+}, \chi^{-}) \cdot v_{x^{+}_{i} x^{+}} \\
&= c(\chi^+, x^{-}_{j})^{-1} c(\chi^{-}, x_i^{+})^{-1} c(x^{+}, x_j^{-}) \cdot v_{x^{+}_{i} x^{+}} \\
&= c(\chi^{+} \chi^{-}, x_i^{+} x_j^{-})^{-1} \cdot i_{c}(b^{-1}_{n^{+}_i n^{-}_{j}} \otimes b_{n^{+}_{i}} \otimes b_{n^{-}_{j}})(v_{x^+}) \\
&= \chi(x_i^{+} x_{j}^{-})^{-1} \cdot i_{c}(b^{-1}_{n^{+}_i n^{-}_{j}} \otimes b_{n^{+}_{i}} \otimes b_{n^{-}_{j}})(v_{x^+}) \\
&= 
i_{c}\left(
(\chi \otimes \id_{\mathbb{C}[N^{+}] \otimes^{c}_{\mathbb{C}[T]} \mathbb{C}[N^{-}]})
\left(
b^{-1}_{n^{+}_i n^{-}_{j}} \otimes b_{n^{+}_{i}} \otimes b_{n^{-}_{j}}
\right)
\right)(v_{x^{+}}).
\end{align*}
According to Proposition~\ref{prop:CNmuisazumaya} and the definition of $z_{i_{c}}$ in Section~\ref{subsec:Basic facts about Azumaya algebras}, we have $z_{i_{c}}(\chi) = \overline{g}_{\chi}$, where $\overline{g}_{\chi}$ denotes the image of $g_{\chi} \in \Aut_{\mathbb{C}[N]}\left(
V^{+}_{\mathbb{C}[N]}
\right) \simeq \GL_{m}\left(
\mathbb{C}[N]
\right)$ in $\PGL_{m}\left(
\mathbb{C}[N]
\right)$.
Hence, $c_{\mathbb{C}[N^{+}] \otimes^{c}_{\mathbb{C}[T]} \mathbb{C}[N^{-}]}$ is the cohomology class of the following 2-cocycle:
\[
\left(
(N/T)^{*}
\right)^2 \ni
(\chi_1, \chi_2) \mapsto g_{\chi_1} g_{\chi_2} g_{\chi_1 \chi_2}^{-1} = c(\chi_2^{+}, \chi_1^{-}) = \mu'(\chi^{+}_1 \chi^-_1, \chi^+_2 \chi^-_2).
\]
Thus, we obtain that $c_{\mathbb{C}[N^{+}] \otimes^{c}_{\mathbb{C}[T]} \mathbb{C}[N^{-}]} = [\mu']^{\vee} = [\mu]^{\vee}$, as desired.
\end{proof}

\begin{corollary}
\label{corollary:twistedgroupalgebraasanazumayaalgebra}
The algebra $\mathbb{C}[N, \mu]$ is an Azumaya algebra over its center $\mathbb{C}[T, \mu]$ that is split by $\mathbb{C}[N]$.
Moreover, the cohomology class $c_{\mathbb{C}[N, \mu]} \in H^2\left(
(N/T)^{*}, \mathbb{C}[N]^{\times}
\right)$ attached to it agrees with the image of $[\mu]^{\vee} \in H^2\left(
(N/T)^{*}, \mathbb{C}^{\times}
\right)$ under the natural map $H^2\left(
(N/T)^{*}, \mathbb{C}^{\times}
\right) \to H^2\left(
(N/T)^{*}, \mathbb{C}[N]^{\times}
\right)$.
\end{corollary}
\begin{proof}
The corollary follows from Lemma~\ref{lemma:replacingmu} combined with Propositions~\ref{prop:CNmuisazumaya} and \ref{prop:cohomologyclassofCN+TN-}.
\end{proof}

\subsection{Azumaya algebra attached to a supercuspidal block}
\label{subsection:Azumaya algebra attached to a supercuspidal block}

Let $G$ be a connected reductive group over $F$ and let $\Rep^{\mathfrak{s}}(G(F))$ be a supercuspidal Bernstein block.
In this subsection, we attach to $\Rep^{\mathfrak{s}}(G(F))$ a cohomology class $[\mu_{\mathfrak{s}}]$ of a 2-cocycle, which encodes the information of $\Rep^{\mathfrak{s}}(G(F))$ up to equivalence.

Let $\pi^{1}$ be an irreducible subrepresentation of $\pi \restriction_{G(F)^{1}}$, where $\pi$ is any irreducible supercuspidal representation in $\Rep^{\mathfrak{s}}(G(F))$.
We define a subgroup $N$ of $G(F)$ as $N = N_{G(F)}(\pi^{1})$.
For $n \in N$, we denote by $\overline{n}$ the image of $n$ in $N/G(F)^{1}$.
By Schur's lemma, for each $n \in N$, the space $\Hom_{G(F)^{1}}\left(
{^n\pi^{1}}, \pi^{1}
\right)$ is one-dimensional.
We fix a family of non-zero elements 
\[
\mathcal{T} = 
\left\{
T_{n} \in \Hom_{G(F)^{1}}\left(
{^n\pi^{1}}, \pi^{1}
\right)
\right\}_{n \in N}
\]
that satisfies the following condition.
\begin{condition}
\label{conditionforcT}
For all $n \in N$ and $g \in G(F)^{1}$, we have $T_{ng} = T_{n} \circ \pi^{1}(g)$.
\end{condition}
By the construction in Appendix~\ref{sectionA cohomology class attached to a projective representation}, the family $\mathcal{T}$ determines a 2-cocycle $\mu_{\mathcal{T}} \in Z^2(N/G(F)^{1}, \mathbb{C}^{\times})$, whose cohomology class $[\mu_{\mathcal{T}}] \in H^2(N/G(F)^{1}, \mathbb{C}^{\times})$ does not depend on the choice of a family $\mathcal{T}$.
Let $\mathcal{H}(G(F), \pi^{1})$ be the $\C$-algebra of functions $\varphi \colon G(F) \to \End_{\mathbb{C}}(\pi^{1})$ that are compactly supported modulo $G(F)^{1}$ and satisfy $\varphi(g_{1} g g_{2}) = \pi^{1}(g_{1}) \circ \varphi(g) \circ \pi^{1}(g_{2})$ for all $g\in G(F)$ and $g_{1}, g_{2} \in G(F)^{1}$.
The multiplication is given by the convolution product
\[
(\varphi_{1} * \varphi_{2})(g) = \sum_{h\in G(F)/G(F)^{1}} \varphi_{1}(h) \circ \varphi_{2}(h^{-1} g).
\]
\begin{proposition}[\cite{MR2508719}]
\label{propfromRochestructureofscblock}
We have an equivalence of categories 
\[
\Rep^{\mathfrak{s}}(G(F)) \simeq \rMod \mathcal{H}(G(F), \pi^{1}).
\]
Moreover, we have an isomorphism of $\mathbb{C}$-algebras
\[
\mathcal{H}(G(F), \pi^{1}) \simeq \mathbb{C}[N/G(F)^{1}, \mu_{\mathcal{T}}].
\]
\end{proposition}
\begin{proof}
The first claim follows from \cite[\S 1.5.4]{MR2508719}.
For $n \in N$, let $\varphi_{n} \in \mathcal{H}(G(F), \pi^{1})$ be the element such that $\supp(\varphi_{n}) = nG(F)^{1}$, and $\varphi_{n}(n) = T_{n}$.
By Condition~\ref{conditionforcT}, the element $\varphi_{n}$ only depends on the image of $n$ in $N/G(F)^{1}$.
For $w \in N/G(F)^{1}$, we write $\varphi_{w} := \varphi_{n}$, where $n$ is any element of $N$ with $\overline{n} = w$.
According to \cite[\S 1.6.1]{MR2508719}, we have 
\[
\mathcal{H}(G(F), \pi^{1}) = \bigoplus_{w\in N/G(F)^{1}} \mathbb{C} \cdot \varphi_{w},
\]
and for $w_1, w_2 \in N/G(F)^{1}$, we have $\varphi_{w_1} * \varphi_{w_2} = c_{w_1, w_2} \cdot \varphi_{w_1 w_2}$ for some $c_{w_1, w_2} \in \mathbb{C}^{\times}$.
We prove that $c_{w_1, w_2} = \mu_{\mathcal{T}}(w_1, w_2)$.
We fix lifts $n_1$ and $n_2$ of $w_1$ and $w_2$ in $N$, respectively.
Then we have 
\begin{align*}
c_{w_1, w_2} \cdot \varphi_{w_1 w_2}(n_1 n_2) &=
(\varphi_{w_1} * \varphi_{w_2})(n_1 n_2) \\
&= \sum_{g \in G(F)/G(F)^{1}} \varphi_{w_1}(g) \circ \varphi_{w_2}(g^{-1} n_1 n_2) \\
&= \varphi_{w_1}(n_1) \circ \varphi_{w_2}(n_1^{-1} n_1 n_2) \\
&= \varphi_{n_1}(n_1) \circ \varphi_{n_2}(n_2) \\
&= T_{n_1} \circ T_{n_2} \\
&= \mu_{\mathcal{T}}(\overline{n}_1, \overline{n}_2) \cdot T_{n_1 n_2} \\
&= \mu_{\mathcal{T}}(w_1, w_2) \cdot \varphi_{w_1 w_2}(n_1 n_2).
\end{align*}
Since $\varphi_{w_1 w_2}(n_1 n_2) = T_{n_1 n_2} \neq 0$, we obtain that $c_{w_1, w_2}  = \mu_{\mathcal{T}}(w_1, w_2)$.
\end{proof}

By Proposition~\ref{propfromRochestructureofscblock}, in order to understand the supercuspidal Bernstein block $\Rep^{\mathfrak{s}}(G(F))$, it suffices to determine the Morita equivalence class of $\mathbb{C}[N/G(F)^{1}, \mu_{\mathcal{T}}]$.
According to Corollary~\ref{corollary:twistedgroupalgebraasanazumayaalgebra}, the algebra $\mathbb{C}[N/G(F)^{1}, \mu_{\mathcal{T}}]$ is an Azumaya algebra over its center, and the corresponding cohomology class can be computed in terms of $\mu_{\mathcal{T}}$.
In the remainder of this subsection, we provide a description of this cohomology class in terms of irreducible representations in $\Rep^{\mathfrak{s}}(G(F))$.
We remark that the fact that $\mathbb{C}[N/G(F)^{1}, \mu_{\mathcal{T}}]$ is an Azumaya algebra over its center was already proved in \cite[Proposition~8.1]{MR1981032} and \cite[Proposition~9.4]{MR4538041}.

We define a 2-cocycle $\mu_{\mathfrak{s}} \in Z^2\left(
X_{\unr}(G)_{\pi}, \mathbb{C}^{\times}
\right)$, where $\pi$ is an irreducible representation in $\Rep^{\mathfrak{s}}(G(F))$, and 
\[
X_{\unr}(G)_{\pi} = \left\{
\chi \in X_{\unr}(G) \mid \pi \otimes\chi \simeq \pi
\right\}.
\]
For each $\chi \in X_{\unr}(G)_{\pi}$, we fix an isomorphism 
\[
I_{\chi} \colon \pi \isoarrow\pi  \otimes \chi.
\]
Let $\chi_1, \chi_2 \in X_{\unr}(G)_{\pi}$.
Then as subspaces of $\End_{\mathbb{C}}(\pi)$, we have 
\[
\Hom_{G(F)}\left(
\pi, \pi \otimes \chi_1
\right) = \Hom_{G(F)}\left(
\pi \otimes \chi_2, \pi \otimes \chi_1 \chi_2
\right).
\]
Hence, we can compose the elements $I_{\chi_1}$ and $I_{\chi_2}$, and the composition $I_{\chi_1} \circ I_{\chi_2}$ is a non-zero element of $\Hom_{G(F)}\left(
\pi, \pi \otimes \chi_1 \chi_2
\right)$.
By Schur's lemma, there exists $\mu_{\mathfrak{s}}(\chi_1, \chi_2) \in \bC^{\times}$ such that
\[
I_{\chi_1} \circ I_{\chi_2} = \mu_{\mathfrak{s}}(\chi_1, \chi_2) \cdot I_{\chi_1 \chi_2}.
\]
It is easy to check that $\mu_{\mathfrak{s}}$ defines a 2-cocycle $X_{\unr}(G)_{\pi} \times X_{\unr}(G)_{\pi} \to \mathbb{C}^{\times}$.
While the 2-cocycle $\mu_{\mathfrak{s}}$ depends on the choice of $I_{\chi}$, its cohomology class $[\mu_{\mathfrak{s}}]$ does not.
Moreover, the cohomology class $[\mu_{\mathfrak{s}}]$ is independent of the choice of $\pi$:
\begin{lemma}
The group $X_{\unr}(G)_{\pi}$ and the cohomology class $[\mu_{\mathfrak{s}}] \in H^2\left(
X_{\unr}(G)_{\pi}, \mathbb{C}^{\times}
\right)$ do not depend on the choice of an irreducible representation $\pi$ in $\Rep^{\mathfrak{s}}(G(F))$.
\end{lemma}
\begin{proof}
Let $\pi'$ be another irreducible representation in $\Rep^{\mathfrak{s}}(G(F))$.
Then there exists $\chi' \in X_{\unr}(G)$ such that $\pi' \simeq \pi \otimes \chi'$.
As subspaces of $\End_{\mathbb{C}}(\pi)$, we have 
\[
\Hom_{G(F)}\left(
\pi, \pi \otimes \chi
\right) = \Hom_{G(F)}\left(
\pi \otimes \chi', \pi \otimes \chi \chi'
\right) \simeq \Hom_{G(F)}\left(
\pi', \pi' \otimes \chi
\right).
\]
Thus, the conditions $\pi' \otimes \chi \simeq \pi'$ and $\pi \otimes \chi \simeq \pi$ are equivalent, and the cohomology class $[\mu_{\mathfrak{s}}]$ is independent of the choice of $\pi$.
\end{proof}

We call $[\mu_{\mathfrak{s}}]$ the \emph{cohomology class attached to the supercuspidal block $\Rep^{\mathfrak{s}}(G(F))$}.
We relate the cohomology class $[\mu_{\mathfrak{s}}]$ to the cohomology class $[\mu_{\mathcal{T}}]^{\vee}$ defined in Section~\ref{subsec:Twisted group algebra as an Azumaya algebra}.

As in Section~\ref{subsec:Twisted group algebra as an Azumaya algebra}, we define an alternating bilinear form $c = c_{\mu_{\mathcal{T}}}$ on $(N/G(F)^{1}) \times (N/G(F)^{1})$ by $c_{\mu_{\mathcal{T}}}(x, y) := \mu_{\mathcal{T}}(y, x) \mu_{\mathcal{T}}(x, y)^{-1}$ for $x, y \in N/G(F)^{1}$.
The relation between $c$ and the family $\mathcal{T}$ is described as follows.

\begin{lemma}
\label{lemmacomputationofc}
For $n_{1}, n_{2} \in N$, we have 
\[
T_{n_2} \circ T_{(n_2)^{-1} n_1 n_2} = c(\overline{n}_1, \overline{n}_2) \cdot T_{n_1} \circ T_{n_2}.
\]
\end{lemma}
\begin{proof}
We have 
\[
T_{n_1} \circ T_{n_2} = \mu_{\mathcal{T}}(\overline{n}_1, \overline{n}_2) \cdot T_{n_1 n_2}
\]
and
\[
T_{n_2} \circ T_{(n_2)^{-1} n_1 n_2} = \mu_{\mathcal{T}}(\overline{n}_2, \overline{n}_2^{-1} \overline{n}_1 \overline{n}_2) \cdot T_{n_1 n_2} = \mu_{\mathcal{T}}(\overline{n}_2, \overline{n}_1) \cdot T_{n_1 n_2}.
\]
Thus, we conclude that 
\begin{align*}
T_{n_2} \circ T_{(n_2)^{-1} n_1 n_2} &= \mu_{\mathcal{T}}(\overline{n}_2, \overline{n}_1) \cdot T_{n_1 n_2} \\
&= \mu_{\mathcal{T}}(\overline{n}_2, \overline{n}_1) \mu_{\mathcal{T}}(\overline{n}_1, \overline{n}_2)^{-1} \cdot T_{n_1} \circ T_{n_2} \\
&= c(\overline{n}_1, \overline{n}_2) \cdot T_{n_1} \circ T_{n_2}.
\qedhere
\end{align*}
\end{proof}

The radical of $c$ is described as follows.
Following \cite[\S 1.6.3]{MR2508719}, we define a subgroup $T$ of $G(F)$ as 
\[
T = \bigcap_{\chi \in X_{\unr}(G)_{\pi}} \ker(\chi).
\]
As explained in \cite[\S 1.6.3]{MR2508719}, we have $X_{\unr}(G)_{\pi} = (G(F)/T)^{*}$, and it is a finite group.
Moreover, according to \cite[Proposition~1.6.3.2]{MR2508719}, the subalgebra $\mathcal{H}(T, \pi^{1})$ is the center of the $\mathbb{C}$-algebra $\mathcal{H}(G(F), \pi^{1}) \simeq \mathbb{C}[N/G(F)^{1}, \mu_{\mathcal{T}}]$.
Hence, the group $T/G(F)^{1}$ is the radical of the alternating form $c$ on $(N/G(F)^{1}) \times (N/G(F)^{1})$.
We also use the same notation $c$ for the symplectic form $(N/T) \times (N/T) \to \mathbb{C}^{\times}$ induced by $c$.
We fix a decomposition $N/T = N^{+}/T \oplus N^{-}/T$ into Lagrangian subgroups.
According to Lemma~\ref{lemma:replacingmu}, we may and do replace $\mathcal{T}$ so that the following holds.
\begin{condition}
For $x_1^{+}, x_2^{+} \in N^{+}/G(F)^{1}$ and $x_1^{-}, x_2^{-} \in N^{-}/G(F)^{1}$, we have 
\[
\mu_{\mathcal{T}}(x^{+}_{1}x^{-}_{1}, x^{+}_{2}x^{-}_{2}) = c(x^{+}_{2}, x^{-}_{1}).
\]
\end{condition}
Since $\mu_{\mathcal{T}}$ is trivial on $(N^{+}/G(F)^{1}) \times (N^{+}/G(F)^{1})$, the family $\left\{T_{n}\right\}_{n \in N^+}$ defines an extension of $\pi^{1}$ to $N^{+}$, which we denote by $\pi^{+}$.
More precisely, for $n \in N^{+}$, we define $\pi^{+}(n) := T_{n}$.
We also write $\pi_{T} := \pi^{+} \restriction_{T}$.
In order to relate $[\mu_{\mathfrak{s}}]$ and $[\mu_{\mathcal{T}}]^{\vee}$, we prepare several lemmas about $\pi^{+}$.
Recall that $\phi_{c}$ denotes the isomorphism $N/T \isoarrow (N/T)^{*}$ given by $\phi_{c}(x)(y) = c(x, y)$ for $x, y \in N/T$.
For $x, y \in N$, we write $\phi_{c}(x) = \phi_{c}(xT)$ and $\phi_{c}(x)(y) =  \phi_{c}(xT)(yT)$ by abuse of notation.

\begin{lemma}
\label{lemma:relationofnconjandphicntwist}
For $n \in N$, we have
\[
T_{n} \in \Hom_{N^{+}}\left(
{^n\pi^{+}} \otimes \phi_{c}(n)\restriction_{N^{+}}, \pi^{+}
\right).
\]
In particular, we have ${^n\pi^{+}} \otimes \phi_{c}(n)\restriction_{N^{+}} \simeq \pi^{+}$.
\end{lemma}
\begin{proof}
According to Lemma~\ref{lemmacomputationofc}, for $n^{+} \in N^{+}$, we have 
\[
T_{n} \circ {^n\pi^{+}}(n^{+}) = T_{n} \circ T_{n^{-1} n^{+} n} = c(\overline{n}^{+}, \overline{n}) \cdot T_{n^{+}} \circ T_{n} = \phi_{c}(n)(n^{+})^{-1} \cdot \pi^{+}(n^{+}) \circ T_{n}.
\]
Thus, $T_{n}$ intertwines $^n\pi^{+}$ and $\pi^{+} \otimes \phi_{c}(n)^{-1} \restriction_{N^+}$. 
Equivalently, $T_{n}$ intertwines ${^n\pi^{+}} \otimes \phi_{c}(nT)\restriction_{N^{+}}$ and $\pi^{+}$.
Since ${^n\pi^{+}} \otimes \phi_{c}(nT)\restriction_{N^{+}}$ and $\pi^{+}$ are irreducible, we obtain the last claim. 
\end{proof}

\begin{lemma}
The compactly induced representation $\ind_{N^+}^{G(F)}(\pi^{+})$ is an irreducible representation in $\Rep^{\mathfrak{s}}(G(F))$.
\end{lemma}
\begin{proof}
First, we prove that $\ind_{N^+}^{G(F)}(\pi^{+})$ is irreducible.
It suffices to show that $N_{G(F)}(\pi^{+}) = N^{+}$.
Suppose that $n \in G(F)$ normalizes $\pi^{+}$.
Since $\pi^{+}$ is an extension of the representation $\pi^{1}$, we have $n \in N$.
Since $\Hom_{G(F)^{1}}\left(
{^n\pi^{1}}, \pi^{1}
\right)$ is one-dimensional, we have 
\[
\{0\} \neq \Hom_{N^{+}}\left(
{^n\pi^{+}}, \pi^{+}
\right) = \Hom_{G(F)^{1}}\left(
{^n\pi^{1}}, \pi^{1}
\right).
\]
In particular, we have $T_{n} \in \Hom_{N^{+}}\left(
{^n\pi^{+}}, \pi^{+}
\right) $.
Hence, for $n^{+} \in N^{+}$, we have 
\[
T_{n} \circ T_{n^{-1} n^+ n} = T_{n} \circ {^n\pi^{+}}(n^+) = \pi^{+}(n^{+}) \circ T_n = T_{n^{+}} \circ T_{n}.
\]
Hence, we obtain from Lemma~\ref{lemmacomputationofc} that $c(\overline{n}^{+}, \overline{n}) = 1$ for all $n^{+} \in N^{+}$.
Thus, we obtain that $n \in N^{+}$, as desired.

We prove that $\ind_{N^{+}}^{G(F)}(\pi^{+}) \in \Rep^{\mathfrak{s}}(G(F))$.
Let $\pi$ be an irreducible representation in $\Rep^{\mathfrak{s}}(G(F))$.
By the Mackey theory, we have
\[
\ind_{N^{+}}^{G(F)}(\pi^{+}) \restriction_{G(F)^{1}} = \bigoplus_{g \in G(F)/N^{+}} {^g\pi^{1}}.
\]
Hence, we obtain that 
\[
\Hom_{G(F)^{1}}\left(
\ind_{N^{+}}^{G(F)}(\pi^{+}) \restriction_{G(F)^{1}}, \pi\restriction_{G(F)^{1}}
\right) \supseteq \Hom_{G(F)^{1}}\left(
\pi^{1}, \pi\restriction_{G(F)^{1}}
\right) \neq \{0\}.
\] 
Thus, we obtain from \cite[Proposition~1.4.1.1]{MR2508719} that $\ind_{N^{+}}^{G(F)}(\pi^{+}) \simeq \pi \otimes \eta \in \Rep^{\mathfrak{s}}(G(F))$ for some $\eta \in X_{\unr}(G)$.
\end{proof}

We now compute the 2-cocycle $\mu_{\mathfrak{s}}$ using the irreducible representation $\pi = \ind_{N^{+}}^{G(F)}(\pi^{+})$.
For $\chi \in X_{\unr}(G)_{\pi}$, we choose an isomorphism $I_{\chi} \colon \pi \isoarrow\pi  \otimes \chi$ as follows.
We fix a lift $n \in N$ of $\phi_{c}^{-1}(\chi \restriction_{N}) \in N/T$.
Then we have an isomorphism 
\[
\ind^{G(F)}_{N^{+}}(\pi^{+}) \isoarrow \ind^{G(F)}_{N^{+}}\left(
^{n^{-1}}\pi^{+}
\right), \qquad f \mapsto [g \mapsto f(n g)].
\]
According to Lemma~\ref{lemma:relationofnconjandphicntwist}, we have 
\[
T_{n} \in \Hom_{N^{+}}\left(
{^n\pi^{+}} \otimes \phi_{c}(n)\restriction_{N^{+}}, \pi^{+}
\right) = \Hom_{N^{+}}\left(
\pi^{+} \otimes \chi\restriction_{N^{+}}, {^{n^{-1}}\pi^{+}}
\right).
\]
Hence, we have an isomorphism 
\[
\ind^{G(F)}_{N^{+}}\left(
^{n^{-1}}\pi^{+}
\right) \isoarrow \ind^{G(F)}_{N^{+}}\left(
\pi^{+} \otimes \chi\restriction_{N^+}
\right), \qquad f \mapsto [g \mapsto T_{n}^{-1}(f(g))].
\]
Finally, we have an isomorphism 
\[
\ind^{G(F)}_{N^{+}}\left(
\pi^{+} \otimes \chi\restriction_{N^+}
\right) \isoarrow \ind^{G(F)}_{N^{+}}(\pi^{+}) \otimes \chi, \qquad f \mapsto [g \mapsto \chi(g)^{-1} \cdot f(g)].
\]
Now, we define the isomorphism $I_{\chi}$ as the composition of the isomorphisms above.
Explicitly, for $f \in \ind^{G(F)}_{N^{+}}(\pi^{+})$ and $g \in G(F)$, we have
\[
I_{\chi}(f)(g) = \chi(g)^{-1} \cdot T_{n}^{-1}(f(n g)).
\]
We note that the isomorphism $I_{\chi}$ does not depend on the choice of a lift $n$ of $\phi_{c}^{-1}(\chi \restriction_{N})$.

\begin{proposition}
\label{prop:forrelationofmuandmuautm}
Let $\chi_1, \chi_2 \in X_{\unr}(G)_{\pi}$ and let $n_{1} \in N$ (resp.\ $n_{2} \in N$) be a lift of $\phi_{c}^{-1}\left(
\chi_1 \restriction_{N}
\right)$ (resp.\ $\phi_{c}^{-1}\left(
\chi_2 \restriction_{N}
\right)$).
Then we have 
\[
I_{\chi_1} \circ I_{\chi_2} = \mu_{\mathcal{T}}(\overline{n}_1, \overline{n}_2)^{-1} \cdot I_{\chi_1 \chi_2}.
\]
\end{proposition}
\begin{proof}
Note that $n_2 n_1$ is a lift of $\phi_{c}^{-1}(\chi_1 \chi_2\restriction_{N}) = \phi_c^{-1}(\chi_2 \chi_1\restriction_{N})$.
Let $f \in \ind^{G(F)}_{N^{+}}(\pi^{+})$.
Then for $g \in G(F)$, we have 
\[
I_{\chi_1 \chi_2}(f)(g) = \chi_1 \chi_2(g)^{-1} \cdot T^{-1}_{n_2 n_1}(f(n_2 n_1 g)).
\]
On the other hand, we have 
\begin{align*}
\left(
I_{\chi_1} \circ I_{\chi_2}
\right)(f)(g) &= \chi_1(g)^{-1} \cdot T_{n_1}^{-1}\left(
\left(
I_{\chi_2}(f)
\right)
(n_1 g)
\right) \\
&= \chi_1(g)^{-1} \chi_2(n_1 g)^{-1} \cdot (T_{n_2} \circ T_{n_1})^{-1}(f(n_2 n_1 g)).
\end{align*}
Thus, we obtain that
\begin{align*}
\left(
I_{\chi_1} \circ I_{\chi_2}
\right)(f)(g) &= \chi_{2}(n_1)^{-1} \mu_{\mathcal{T}}(\overline{n}_2, \overline{n}_1)^{-1} \cdot I_{\chi_1 \chi_2}(f)(g) \\
&= c(\overline{n}_2, \overline{n}_1)^{-1} \mu_{\mathcal{T}}(\overline{n}_2, \overline{n}_1)^{-1} \cdot I_{\chi_1 \chi_2}(f)(g) \\
&= \mu_{\mathcal{T}}(\overline{n}_1, \overline{n}_2)^{-1} \cdot I_{\chi_1 \chi_2}(f)(g).
\qedhere
\end{align*}
\end{proof}

We regard $H^2\left(
(N/T)^{*}, \mathbb{C}^{\times}
\right)$ as a subset of $H^2\left(
X_{\unr}(G)_{\pi}, \mathbb{C}^{\times}
\right)$ via the inflation map
\[
H^2\left(
(N/T)^{*}, \mathbb{C}^{\times}
\right)
\hookrightarrow H^2\left(
(G(F)/T)^{*}, \mathbb{C}^{\times}
\right) = H^2\left(
X_{\unr}(G)_{\pi}, \mathbb{C}^{\times}
\right),
\]
which is injective by the inflation-restriction exact sequence and the injectivity of $\mathbb{C}^{\times}$.

\begin{corollary}
\label{corollarymus=muvee}
We have $[\mu_{\mathfrak{s}}] \in H^2\left(
(N/T)^{*}, \mathbb{C}^{\times}
\right)$.
Moreover, we have $[\mu_{\mathcal{T}}]^{\vee} = [\mu_{\mathfrak{s}}]^{-1}$.
\end{corollary}
\begin{proof}
The claim follows from Proposition~\ref{prop:forrelationofmuandmuautm} and the definitions of $[\mu_{\mathcal{T}}]^{\vee}$ and $[\mu_{\mathfrak{s}}]$.
\end{proof}

According to Proposition~\ref{prop:MoritaequivclassofAzumayaalgebras}, Corollary~\ref{corollary:twistedgroupalgebraasanazumayaalgebra}, and Corollary~\ref{corollarymus=muvee}, the Morita equivalence class of the Azumaya algebra $\mathbb{C}[N/G(F)^{1}, \mu_{\mathcal{T}}]$ over its center $\mathbb{C}[T/G(F)^{1}]$ is determined by the cohomology class $[\mu_{\mathcal{T}}]^{\vee} =  [\mu_{\mathfrak{s}}]^{-1}$.
Thus, the Bernstein block $\Rep^{\mathfrak{s}}(G(F))$ is determined, up to equivalence, by its center $\mathbb{C}[T/G(F)^{1}]$ and the cohomology class $[\mu_{\mathfrak{s}}]$.
The center $\mathbb{C}[T/G(F)^{1}]$ of the Bernstein block $\Rep^{\mathfrak{s}}(G(F))$ admits an explicit description as the ring of regular functions on the complex algebraic torus parametrizing the irreducible representations in $\Rep^{\mathfrak{s}}(G(F))$ (\cite[Proposition~1.6.4.1]{MR2508719}).
In the following section, we formulate a conjecture describing the cohomology class $[\mu_{\mathfrak{s}}]$ in terms of ``spectral-side data'' via the local Langlands correspondence.

\section{A conjecture}
\label{sec:aconjecture}

In this section, we attach to an $L$-parameter $\varphi \colon W_{F}\times \SL_{2}(\C) \to {^LG}$ and its ``enhancement'' $\rho$ a cohomology class $[\mu_{(\varphi, \rho)}]$ represented by a 2-cocycle.
We then formulate a conjecture relating $[\mu_{(\varphi, \rho)}]$ to the cohomology class $[\mu_{\mathfrak{s}}]$ attached to a supercuspidal block $\Rep^{\mathfrak{s}}(G(F))$ under the assumption of a suitable form of the local Langlands correspondence.

\subsection{A cohomology class attached to an enhanced $L$-parameter}
\label{subsection:A cohomology class attached to an enhanced $L$-parameter}

Let $G$ be a connected reductive group over $F$, and let $\varphi \colon W_{F} \times \SL_{2}(\mathbb{C}) \to {^LG}$ be a tempered $L$-parameter.
The basic form of the local Langlands conjecture predicts that the $\widehat{G}$-conjugacy class $[\varphi]$ of $\varphi$ determines a finite set $\Pi_{[\varphi]}(G)$ of irreducible representations of $G(F)$, called an \emph{$L$-packet}.
Refinements of the local Langlands conjecture give parametrizations of $\Pi_{[\varphi]}(G)$ in terms of irreducible representations of certain variants of the centralizer $\Cent_{\widehat{G}}(\varphi)$ of the image of $\varphi$ in $\widehat{G}$.
There are several such formulations in the literature, each involving a slightly different modification of $\Cent_{\widehat{G}}(\varphi)$.
In this paper, we adopt the refinement based on rigid inner forms as developed by Kaletha and Dillery (\cite{MR3548533, MR4617942}). However, the arguments here are not specific to this choice and can be adapted to other formulations with minor changes.

We recall the refinement of the local Langlands conjecture in \cite{MR3548533, MR4617942} briefly.
Let $Z$ be a (sufficiently large) finite central subgroup of $G$.
We write $\overline{G} = G/Z$.
Then the isogeny $G \to \overline{G}$ dualizes to an isogeny $\widehat{\overline{G}} \to \widehat{G}$.
We write $S_{\varphi} = \Cent_{\widehat{G}}(\varphi)$ and denote by $Z(\widehat{\overline{G}})^{+}$ and $S_{\varphi}^{+}$ the preimages of $Z(\widehat{G})^{\Gamma_{F}}$ and $S_{\varphi}$ in $\widehat{\overline{G}}$, respectively.
We fix a realization of $G$ as a rigid inner twist of its quasi-split inner form $G_{\qs}$.
More precisely, we fix an inner twist $\psi \colon G_\qs \to G$ as well as a lift of the cohomology class of the 1-cocycle 
\[
\Gamma_{F} \ni \sigma \mapsto \psi^{-1} \sigma(\psi) \in G_\qs/Z(G_\qs)
\]
to the class $[x]$ in the set $H^1(\mathcal{E}, Z\to G_{\qs})$ defined in \cite[Section~3]{MR4617942}.
Such a lift always exists if $Z$ contains the center of the derived group of $G$ by \cite[Corollary~3.8]{MR3548533} if $F$ has characteristic zero and by \cite[Proposition~5.12]{MR4617942} if $F$ has positive characteristic.
Note that if $F$ has characteristic zero, the cohomology set
$H^1(\mathcal{E} ,Z\to G_{\qs})$ identifies with the set $H^1(u\to W,Z\to G_{\qs})$ defined in \cite[Section~3.2]{MR3548533} (see \cite[Proposition~7.1]{MR5007737}).
The cohomology class $[x] \in H^1(\mathcal{E}, Z\to G_{\qs})$ determines a character of $\pi_0(Z(\widehat{\overline{G}})^{+})$ via the pairing of \cite[Corollary~5.4]{MR3548533} in characteristic zero and \cite[Corollary~7.11]{MR4617942} in positive characteristic.
We use the same notation $[x]$ for this character.
The natural inclusion $Z(\widehat{\overline{G}})^{+} \subset S_{\varphi}^{+}$ induces a map $\pi_0(Z(\widehat{\overline{G}})^{+}) \to \pi_0(S_{\varphi}^{+})$ whose image is a central subgroup.
Let $\Irr\left(
\pi_0(S_{\varphi}^{+}), [x]
\right)$ denote the set of irreducible representations of $\pi_0(S_{\varphi}^{+})$ whose precompositions with the map $\pi_0(Z(\widehat{\overline{G}})^{+}) \to \pi_0(S_{\varphi}^{+})$ are $[x]$-isotypic.
The conjectural local Langlands correspondence predicts a bijection 
\[
\iota_{\mathfrak{w}} \colon \Pi_{[\varphi]}(G) \isoarrow \Irr\left(
\pi_0(S_{\varphi}^{+}), [x]
\right)
\]
depending on a Whittaker datum $\mathfrak{w}$ for $G_{\qs}$.
For further details, see \cite[Section~5.4]{MR3548533} and \cite[Section~7.3]{MR4617942}.

We call a pair $(\varphi, \rho)$ consisting of an $L$-parameter $\varphi$ and $\rho \in \Irr\left(
\pi_0(S_{\varphi}^{+})
\right)$ an \emph{enhanced $L$-parameter} for $^LG$.
If $\rho \in \Irr\left(
\pi_0(S_{\varphi}^{+}), [x]
\right)$, we say that the enhanced $L$-parameter $(\varphi, \rho)$ is \emph{$G$-relevant}.
We identify $\rho$ with an irreducible representation of $S_{\varphi}^{+}$ via the inflation.
The group $\widehat{G}$ acts on the set of enhanced $L$-parameters as follows.
For $g \in \widehat{G}$ and an $L$-parameter $\varphi$, we write ${^g\varphi} = \Ad(g) \circ \varphi$.
Then we have $\Cent_{\widehat{G}}(^g\varphi) = {^g\Cent_{\widehat{G}}(\varphi)}$.
Since $\widehat{\overline{G}} \to \widehat{G}$ is a central isogeny, the conjugation action of $\widehat{\overline{G}}$ on itself descends to an action of $\widehat{G}$ on $\widehat{\overline{G}}$.
With respect to this $\widehat{G}$-action on $\widehat{\overline{G}}$, we have $S_{{^g\varphi}}^{+} = {^gS_{\varphi}^{+}}$.
Now, we define the action of $\widehat{G}$ on the set of enhanced $L$-parameters as
\[
g \cdot (\varphi, \rho) = ({^g\!\varphi}, {^g\rho}),
\]
where $^g\rho$ is the representation of $\pi_0(S_{^g\varphi}^{+}) = {^g\pi_0(S_{\varphi}^{+})}$ given by $x \mapsto \rho({^{g^{-1}}x})$.
For an $L$-parameter $\varphi$ (resp.\ an enhanced $L$-parameter $(\varphi, \rho)$), we denote by $[\varphi]$ (resp.\ $[(\varphi, \rho)]$) its $\widehat{G}$-orbit.
Following \cite[Section~3]{2017arXiv170103593A}, let $\Phi_{e}(^LG)$ denote the set of $\widehat{G}$-orbits of enhanced $L$-parameters and $\Phi_{e}(G)$ denote the subset of $\Phi_{e}(^LG)$ consisting of $\widehat{G}$-orbits of $G$-relevant enhanced $L$-parameters.

In analogy with the automorphic side, we attach a cohomology class of a 2-cocycle to an enhanced parameter.
First, we recall the definition of the unramified twists for enhanced $L$-parameters from \cite[Section~3]{2017arXiv170103593A}.
In \cite[Section~3.3.1]{MR3444233}, Haines provided a natural isomorphism
\[
X_{\unr}(G) \isoarrow \left(
Z(\widehat{G})^{I_{F}}
\right)_{\Frob_{F}}^{\circ}.
\]
Let us recall the construction briefly.
In \cite[Section~7]{MR1485921}, Kottwitz constructed a functorial surjective homomorphism
\[
\kappa_{G_{F^{\unr}}} \colon G(F^{\unr}) \twoheadrightarrow X^{*}(Z(\widehat{G}))_{I_{F}}
\]
called the \emph{Kottwitz homomorphism}, whose restriction to $\Frob_{F}$-fixed points
\[
\kappa_{G} \colon G(F) \twoheadrightarrow \left(
X^{*}(Z(\widehat{G}))_{I_{F}}
\right)^{\Frob_{F}}
\]
is also surjective.
According to \cite[Lemma~11.5.2]{KalethaPrasad}, we have 
\[
\kappa_{G_{F^{\unr}}}^{-1}\left(
\left(X^{*}(Z(\widehat{G}))_{I_{F}}\right)_\mathrm{tor}
\right) = G(F^{\unr})^{1}
\qquad
\text{and}
\qquad
\kappa_{G}^{-1}\left(
\left(\left(X^{*}(Z(\widehat{G}))_{I_{F}}\right)^{\Frob_{F}}\right)_\mathrm{tor}
\right) = G(F)^{1}.
\]
Hence, the Kottwitz homomorphism induces isomorphisms
\[
\kappa_{G_{F^{\unr}}} \colon G(F^{\unr})/G(F^{\unr})^{1} \isoarrow X^{*}(Z(\widehat{G}))_{I_{F}}/\left(X^{*}(Z(\widehat{G}))_{I_{F}}\right)_\mathrm{tor}
\]
and
\[
\kappa_{G} \colon G(F)/G(F)^{1} \isoarrow \left(
X^{*}(Z(\widehat{G}))_{I_{F}}
\right)^{\Frob_{F}} /\left(\left(X^{*}(Z(\widehat{G}))_{I_{F}}\right)^{\Frob_{F}}\right)_\mathrm{tor}.
\]
Thus, we have natural isomorphisms 
\begin{align*}
X_{\unr}(G_{F^{\unr}}) &= \Hom_{\mathrm{grp}}\left(
G(F^{\unr})/G(F^{\unr})^{1}, \mathbb{C}^{\times}
\right) \\
&\simeq \Hom_{\mathrm{grp}}\left(
X^{*}(Z(\widehat{G}))_{I_{F}}/\left(X^{*}(Z(\widehat{G}))_{I_{F}}\right)_\mathrm{tor}, \mathbb{C}^{\times}
\right) \\
&\simeq
\Hom_{\mathrm{grp}}\left(
X^{*}(Z(\widehat{G}))_{I_{F}}, \mathbb{C}^{\times}
\right)^{\circ} \\
&\simeq
\Hom_{\mathrm{grp}}\left(
X^{*}(Z(\widehat{G})^{I_{F}}), \mathbb{C}^{\times}
\right)^{\circ} \\
&\simeq
\Hom_{\mathbb{C}-\mathrm{Alg}}\left(
\mathbb{C}[X^{*}(Z(\widehat{G})^{I_{F}})], \mathbb{C}
\right)^{\circ} \\
&\simeq Z(\widehat{G})^{I_{F}, \circ},
\end{align*}
where $\Hom_{\mathrm{grp}}$ (resp.\ $\Hom_{\mathbb{C}-\mathrm{Alg}}$) denotes the set of group homomorphisms (resp.\ $\mathbb{C}$-algebra homomorphisms).
Similarly, we have
\begin{align*}
X_{\unr}(G) &= \Hom_{\mathrm{grp}}\left(
G(F)/G(F)^{1}, \mathbb{C}^{\times}
\right) \\
&\simeq
\Hom_{\mathrm{grp}}\left(
\left(
X^{*}(Z(\widehat{G}))_{I_{F}}
\right)^{\Frob_{F}} /\left(\left(X^{*}(Z(\widehat{G}))_{I_{F}}\right)^{\Frob_{F}}\right)_\mathrm{tor}, \mathbb{C}^{\times}
\right) \\
&\simeq 
\Hom_{\mathrm{grp}}\left(
\left(
X^{*}(Z(\widehat{G}))_{I_{F}}
\right)^{\Frob_{F}} , \mathbb{C}^{\times}
\right)^{\circ} \\
&\simeq 
\Hom_{\mathrm{grp}}\left(
X^{*}\left(
\left(
Z(\widehat{G})^{I_{F}}
\right)_{\Frob_{F}}
\right) , \mathbb{C}^{\times}
\right)^{\circ} \\
&\simeq \left(
Z(\widehat{G})^{I_{F}}
\right)_{\Frob_{F}}^{\circ}.
\end{align*}
We write $X_{\unr}(^LG) = \left(
Z(\widehat{G})^{I_{F}}
\right)_{\Frob_{F}}^{\circ}$.
Then the above isomorphism is written as
\begin{equation}
\label{comparisonofunramifiedcharacters}
X_{\unr}(G) \simeq X_{\unr}(^LG).
\end{equation}
For $\chi \in X_{\unr}(G)$ (resp.\ $\chi \in X_{\unr}(G_{F^{\unr}})$), we denote by $z_{\chi}$ the corresponding element of $X_{\unr}(^LG)$ (resp.\ $Z(\widehat{G})^{I_{F}, \circ}$) via the isomorphisms above.
Similarly, for $z \in X_{\unr}(^LG)$ (resp.\ $z \in Z(\widehat{G})^{I_{F}, \circ}$), we denote by $\chi_{z}$ the corresponding element of $X_{\unr}(G)$ (resp.\ $X_{\unr}(G_{F^{\unr}})$).

We define the action of $X_{\unr}(^LG)$ on $\Phi_{e}(G)$ as follows.
For $z \in Z(\widehat{G})^{I_{F}}$, let $z \cdot \varphi$ be the unique $L$-parameter such that 
\[
(z \cdot \varphi)\restriction_{I_{F} \times \SL_2(\mathbb{C})} = \varphi\restriction_{I_{F} \times \SL_2(\mathbb{C})}
\qquad
\text{and}
\qquad
(z \cdot \varphi)(\Frob_{F}) = z \cdot \varphi(\Frob_{F}).
\]
Noting that $S^{+}_{\varphi} = S^{+}_{z \cdot \varphi}$, we define an action of $Z(\widehat{G})^{I_{F}}$ on the set of enhanced $L$-parameters as $z \cdot (\varphi, \rho) = (z \cdot \varphi, \rho)$.
This action commutes with the conjugation action of $\widehat{G}$ on the set of enhanced $L$-parameters.
Hence, we can define an action of $Z(\widehat{G})^{I_{F}}$ on $\Phi_{e}(G)$ by $z \cdot [(\varphi, \rho)] := [(z \cdot \varphi, \rho)]$.
We denote by $Z(\widehat{G})^{I_{F}}_{[\varphi]}$ (resp.\ $Z(\widehat{G})^{I_{F}}_{[(\varphi, \rho)]}$) the stabilizer of $[\varphi]$ (resp.\ $[(\varphi, \rho)]$) in $Z(\widehat{G})^{I_{F}}$.
\begin{lemma}
\label{lemma:unramifiedrwistgaloisside}
The action of $Z(\widehat{G})^{I_{F}}$ on $\Phi_{e}(G)$ factors through the projection $Z(\widehat{G})^{I_{F}} \twoheadrightarrow (Z(\widehat{G})^{I_{F}})_{\Frob_{F}}$.
\end{lemma}
\begin{proof}
Let $z \in Z(\widehat{G})^{I_{F}}$ and $[(\varphi, \rho)] \in \Phi_{e}(G)$.
Then we have 
\begin{align*}
(\Frob_{F}(z) \cdot \varphi)(\Frob_{F}) &= \Frob_{F}(z) \cdot \varphi(\Frob_{F}) \\
&= (z^{-1} \rtimes 1) \cdot (z \cdot \varphi(\Frob_{F})) \cdot (z \rtimes 1) \\
&= (z^{-1} \rtimes 1) \cdot (z \cdot \varphi)(\Frob_{F}) \cdot (z \rtimes 1).
\end{align*}
Since $z \in Z(\widehat{G})^{I_{F}}$, we also have 
\[
\Frob_{F}(z) \cdot \varphi \restriction_{I_{F} \times \SL_2(\mathbb{C})} = \varphi\restriction_{I_{F} \times \SL_2(\mathbb{C})} = \Ad(z^{-1}) \circ \varphi\restriction_{I_{F} \times \SL_2(\mathbb{C})} = \Ad(z^{-1}) \circ (z \cdot \varphi)\restriction_{I_{F} \times \SL_2(\mathbb{C})}.
\]
Thus, we obtain that $\Frob_{F}(z) \cdot \varphi = \Ad(z^{-1}) \circ (z \cdot \varphi)$.
Since the conjugation action of $Z(\widehat{G})$ on $\pi_0(S_{\varphi}^{+})$ is trivial, we conclude that
\[
\Frob_{F}(z) \cdot [(\varphi, \rho)] = [(\Frob_{F}(z) \cdot \varphi, \rho)] = [(\Ad(z^{-1}) \circ (z \cdot \varphi), \rho)] = [(z \cdot \varphi, \rho)] = z \cdot [(\varphi, \rho)],
\]
which implies the claim.
\end{proof}

By Lemma~\ref{lemma:unramifiedrwistgaloisside}, the group $X_{\unr}(^LG) \subset (Z(\widehat{G})^{I_{F}})_{\Frob_{F}}$ acts on $\Phi_{e}(G)$.
For $[(\varphi, \rho)] \in \Phi_e(G)$ and $z \in X_{\unr}(^LG)$, we call $z \cdot [(\varphi, \rho)]$ the \emph{unramified twist} of $[(\varphi, \rho)]$ by $z$.
We denote by $X_{\unr}(^LG)_{[(\varphi, \rho)]}$ the stabilizer of $[(\varphi, \rho)]$ in $X_{\unr}(^LG)$.
We also define the subset $X_{\unr}(^LG) \cdot [(\varphi, \rho)]$ of $\Phi_e(G)$ as 
\[
X_{\unr}(^LG) \cdot [(\varphi, \rho)] = \left\{
z \cdot [(\varphi, \rho)] \mid z \in X_{\unr}(^LG)
\right\}.
\]

We now attach to an enhanced $L$-parameter $(\varphi, \rho)$ a cohomology class 
\[
[\mu_{(\varphi, \rho)}] \in H^2\left(
X_{\unr}(^LG)_{[(\varphi, \rho)]}, \mathbb{C}^{\times}
\right).
\]
Let $(\varphi, \rho)$ be an enhanced $L$-parameter.
We define a subgroup $N_{\varphi}$ of $\widehat{G}$ as 
\[
N_{\varphi} = \left\{
g \in \widehat{G} \mid {^g\varphi} = z \cdot \varphi \quad \text{for some $z \in Z(\widehat{G})^{I_{F}, \circ}$}
\right\}.
\]
For $n \in N_{\varphi}$, let ${^n\varphi}/\varphi \in Z(\widehat{G})^{I_{F}, \circ}$ be the element such that ${^n\varphi} = ({^n\varphi}/\varphi) \cdot \varphi$.
Then the map $n \mapsto {^n\varphi}/\varphi$ induces an exact sequence
\begin{equation}
\label{exactseqingGhat}
1 \to S_{\varphi} \to N_{\varphi} \to Z(\widehat{G})^{I_{F}, \circ}_{[\varphi]} \to 1.
\end{equation}
We denote by $N_{\varphi}^{+}$ the preimage of $N_{\varphi}$ in $\widehat{\overline{G}}$.
Pulling back \eqref{exactseqingGhat} along $\widehat{\overline{G}} \to \widehat{G}$, we obtain an exact sequence
\begin{equation}
\label{exactseqingGbarhat}
1 \to S_{\varphi}^{+} \to N_{\varphi}^{+} \to Z(\widehat{G})^{I_{F}, \circ}_{[\varphi]} \to 1.
\end{equation}
Let $N^{+}_{\varphi, \rho}$ denote the normalizer of $\rho$ in $N_{\varphi}^{+}$.
Then we have an exact sequence
\begin{equation}
\label{exactseqingGbarhatrhonormalize}
1 \to S_{\varphi}^{+} \to N_{\varphi, \rho}^{+} \to Z(\widehat{G})^{I_{F}, \circ}_{[(\varphi, \rho)]} \to 1.
\end{equation}
We define a cohomology class $[\mu_{(\varphi, \rho)}] \in H^2\left(
Z(\widehat{G})^{I_{F}, \circ}_{[(\varphi, \rho)]}, \mathbb{C}^{\times}
\right)$ as $[\mu_{(\varphi, \rho)}] = [\mu_{(N_{\varphi, \rho}^{+}, S_{\varphi}^{+}, \rho)}]$, where $[\mu_{(N_{\varphi, \rho}^{+}, S_{\varphi}^{+}, \rho)}]$ is the cohomology class in
\[
H^2\left(
N_{\varphi, \rho}^{+}/S_{\varphi}^{+}, \mathbb{C}^{\times}
\right) = H^2\left(
Z(\widehat{G})^{I_{F}, \circ}_{[(\varphi, \rho)]}, \mathbb{C}^{\times}
\right)
\]
defined in Appendix~\ref{sectionA cohomology class attached to a projective representation}.
\begin{proposition}
\label{prop:muforspecsideFunrvsF}
We have $[\mu_{(\varphi, \rho)}] \in H^2\left(
X_{\unr}(^LG)_{[(\varphi, \rho)]}, \mathbb{C}^{\times}
\right)$, where we regard $H^2\left(
X_{\unr}(^LG)_{[(\varphi, \rho)]}, \mathbb{C}^{\times}
\right)$ as a subset of $H^2\left(
Z(\widehat{G})^{I_{F}, \circ}_{[(\varphi, \rho)]}, \mathbb{C}^{\times}
\right)$ via the inflation map
\[
H^2\left(
X_{\unr}(^LG)_{[(\varphi, \rho)]}, \mathbb{C}^{\times}
\right)
\hookrightarrow H^2\left(
Z(\widehat{G})^{I_{F}, \circ}_{[(\varphi, \rho)]}, \mathbb{C}^{\times}
\right).
\]
\end{proposition}
\begin{proof}
Let $Z'$ denote the preimage of $Z(\widehat{G})^{I_{F}}$ via the isogeny $\widehat{\overline{G}} \to \widehat{G}$.
Since $S_{\varphi} \cap Z(\widehat{G})^{I_{F}} = Z(\widehat{G})^{\Gamma_{F}}$, we have $S_{\varphi}^{+} \cap Z' = Z(\widehat{\overline{G}})^{+}$.
Moreover, the image of $N_{\varphi}^{+} \cap Z'$ via the map $N_{\varphi}^{+} \to Z(\widehat{G})^{I_{F}, \circ}_{[\varphi]}$ in \eqref{exactseqingGbarhat} agrees with the kernel of the projection 
\[
Z(\widehat{G})^{I_{F}, \circ} \twoheadrightarrow \left(
Z(\widehat{G})^{I_{F}}
\right)_{\Frob_{F}}^{\circ} = X_{\unr}(^LG).
\]
Hence, the map $N_{\varphi}^{+} \to Z(\widehat{G})^{I_{F}, \circ}_{[\varphi]}$ induces an isomorphism
\[
N_{\varphi, \rho}^{+}/(S_{\varphi}^{+} \cdot (N_{\varphi}^{+} \cap Z')) \isoarrow X_{\unr}(^LG)_{[(\varphi, \rho)]}.
\]
We prove that the representation $\rho$ extends to a representation $\widetilde{\rho}$ of $S_{\varphi}^{+} \cdot (N_{\varphi}^{+} \cap Z')$.
Since the conjugation action of $Z' \subset Z(\widehat{\overline{G}})$ on $S_{\varphi}^{+}$ is trivial, we have 
\[
S_{\varphi}^{+} \cdot Z' \simeq (S_{\varphi}^{+} \times Z')/\left\{
(z, z^{-1}) \mid z \in Z(\widehat{\overline{G}})^{+}
\right\}.
\]
Since $Z'$ is abelian, the character $[x] = \rho|_{Z(\widehat{\overline{G}})^{+}}$ of $Z(\widehat{\overline{G}})^{+}$ extends to a character $\widetilde{[x]}$ of $Z'$.
Then the irreducible representation $\rho \boxtimes \widetilde{[x]}$ of $S_{\varphi}^{+} \times Z'$ descends to an irreducible representation $\widetilde{\rho}'$ of $S_{\varphi}^{+} \cdot Z'  = (S_{\varphi}^{+} \times Z')/Z(\widehat{\overline{G}})^{+}$.
We define the representation $\widetilde{\rho}$ as the restriction of $\widetilde{\rho}'$ to $S_{\varphi}^{+} \cdot (N_{\varphi}^{+} \cap Z')$.
By the construction, $\widetilde{\rho}$ is an extension of $\rho$.
Since the conjugation action of $N_{\varphi, \rho}^{+} \subset \widehat{\overline{G}}$ on $N_{\varphi}^{+} \cap Z' \subset Z(\widehat{\overline{G}})$ is trivial, the group $N_{\varphi, \rho}^{+}$ normalizes the representation $\widetilde{\rho}$.
Thus, we obtain that the cohomology class $[\mu_{(\varphi, \rho)}]$ is the image of the cohomology class 
\[
[\mu_{(N_{\varphi, \rho}^{+}, S_{\varphi}^{+} \cdot (N_{\varphi}^{+} \cap Z'), \widetilde{\rho})}] \in H^2\left(
N_{\varphi, \rho}^{+}/(S_{\varphi}^{+} \cdot (N_{\varphi}^{+} \cap Z')), \mathbb{C}^{\times}
\right) = H^2\left(
X_{\unr}(^LG)_{[(\varphi, \rho)]}, \mathbb{C}^{\times}
\right)
\]
defined in Appendix~\ref{sectionA cohomology class attached to a projective representation}.
\end{proof}

At the end of this subsection, we prove that $[\mu_{(\varphi, \rho)}]$ only depends on the $\widehat{G} \times Z(\widehat{G})^{I_{F}, \circ}$-orbit of $(\varphi, \rho)$.

\begin{proposition}
\label{prop:cohomologyclassingaloissideonlydependsoninertialequivalenceclass}
The cohomology class $[\mu_{(\varphi, \rho)}]$ only depends on $X_{\unr}(^LG) \cdot [(\varphi, \rho)]$.
\end{proposition}
\begin{proof}
Let $g \in \widehat{G}$ and $z \in Z(\widehat{G})^{I_{F}, \circ}$.
It suffices to show that $[\mu_{(\varphi, \rho)}] = [\mu_{(z \cdot (^g\varphi), {^g\rho})}]$.
By the definitions, we have 
\[
N_{z \cdot (^g\varphi)}^{+} = {^gN_{\varphi}^{+}}
\qquad
\text{and}
\qquad
Z(\widehat{G})^{I_{F}, \circ}_{[(z \cdot (^g\varphi), {^g\rho})]} = Z(\widehat{G})^{I_{F}, \circ}_{[(\varphi, \rho)]}.
\]
Hence, we have an isomorphism 
\[
N_{\varphi}^{+} \isoarrow N_{z \cdot (^g\varphi)}^{+}, \quad n \mapsto gng^{-1}.
\]
Since 
\[
{^{gng^{-1}}(z \cdot (^g\varphi))} = z \cdot {^{gn}\varphi} = z \cdot {^g(({^n\varphi}/\varphi) \cdot \varphi)} = z ({^n\varphi}/\varphi) \cdot {^g\varphi} = ({^n\varphi}/\varphi) \cdot (z \cdot {^g\varphi}),
\]
we have ${^{gng^{-1}}(z \cdot (^g\varphi))}/ (z \cdot (^g\varphi)) = {^n\varphi}/\varphi$.
Thus, the exact sequences \eqref{exactseqingGbarhatrhonormalize} for $(\varphi, \rho)$ and $(z \cdot (^g\varphi), {^g\rho})$ are isomorphic.
Since the representation $\rho$ corresponds to the representation $^g\rho$ of ${^gS_{\varphi}^{+}} = S_{z \cdot (^g\varphi)}^{+}$ under the isomorphism
\[
S_{\varphi}^{+} \isoarrow S_{z \cdot (^g\varphi)}^{+}, \quad n \mapsto gng^{-1},
\]
we obtain that $[\mu_{(\varphi, \rho)}] = [\mu_{(z \cdot (^g\varphi), {^g\rho})}]$.
\end{proof}

\subsection{A conjecture on the cohomology classes}

Let $\Rep^{\mathfrak{s}}(G(F))$ be a supercuspidal block.
Note that the inertial equivalence class $\mathfrak{s}$ is the set of irreducible representations in $\Rep^{\mathfrak{s}}(G(F))$.
In this subsection, we state a ``meta-conjecture'' relating the cohomology classes $[\mu_{\mathfrak{s}}]$ in Section~\ref{subsection:Azumaya algebra attached to a supercuspidal block} and $[\mu_{(\varphi, \rho)}]$ in Section~\ref{subsection:A cohomology class attached to an enhanced $L$-parameter} assuming a part of the local Langlands correspondence.
As the working assumption, we assume the following.

\begin{assumption}
\label{workingassumptionLLC}
We have an injective map
\[
\LLC^{\mathfrak{s}} \colon \mathfrak{s} \to \Phi_{e}(G)
\]
such that 
$
\LLC^{\mathfrak{s}}(\pi \otimes \chi) = z_{\chi} \cdot \LLC^{\mathfrak{s}}(\pi)
$ for all $\chi \in X_{\unr}(G)$.
\end{assumption}

Let $\mathfrak{s}^{\vee}$ denote the image of $\LLC^{\mathfrak{s}}$.
We define $[\mu_{\mathfrak{s}^{\vee}}] \in H^2\left(
X_{\unr}(^LG)_{[(\varphi, \rho)]}, \mathbb{C}^{\times}
\right)$ as $[\mu_{\mathfrak{s}^{\vee}}] = [\mu_{(\varphi, \rho)}]$, where $[(\varphi, \rho)]$ is any element in $\mathfrak{s}^{\vee}$ (see Proposition~\ref{prop:cohomologyclassingaloissideonlydependsoninertialequivalenceclass}).

\begin{lemma}
For any $\pi \in \mathfrak{s}$ and $[(\varphi, \rho)] \in \mathfrak{s}^{\vee}$, we have
\[
X_{\unr}(G)_{\pi} = X_{\unr}(^LG)_{[(\varphi, \rho)]}
\]
as subgroups of $X_{\unr}(G) \simeq X_{\unr}(^LG)$.
\end{lemma}
\begin{proof}
Since the group $X_{\unr}(G)_{\pi}$ does not depend on the choice of $\pi \in \mathfrak{s}$, and $X_{\unr}(^LG)_{[(\varphi, \rho)]}$ does not depend on the choice of $[(\varphi, \rho)]$, we may assume that $\LLC^{\mathfrak{s}}(\pi) = [(\varphi, \rho)]$.
Then Assumption~\ref{workingassumptionLLC} implies that for $\chi \in X_{\unr}(G)$, we have $\pi \otimes \chi \simeq \pi$ if and only if $z_{\chi} \cdot [(\varphi, \rho)] = [(\varphi, \rho)]$.
Thus, we conclude that $X_{\unr}(G)_{\pi} = X_{\unr}(^LG)_{[(\varphi, \rho)]}$.
\end{proof}

We now state the conjecture relating $[\mu_{\mathfrak{s}}]$ and $[\mu_{\mathfrak{s}^{\vee}}]$.

\begin{conjecture}
\label{conjecture:comparisonoftwococycles}
As elements of $H^2\left(
X_{\unr}(G)_{\pi}, \mathbb{C}^{\times}
\right) = H^2\left(
X_{\unr}(^LG)_{[(\varphi, \rho)]}, \mathbb{C}^{\times}
\right)$, we have
$[\mu_{\mathfrak{s}}] = [\mu_{\mathfrak{s}^{\vee}}]$.
\end{conjecture}

\section{The case of non-singular supercuspidal representations}
\label{sec:The case of non-singular supercuspidal representations}

In this section, we prove that the local Langlands correspondence for non-singular supercuspidal representations given by Kaletha (\cite{2019arXiv191203274K}) satisfies Conjecture~\ref{conjecture:comparisonoftwococycles}.
We explain the main result of this section more precisely.
In \cite{2019arXiv191203274K}, Kaletha defined and constructed a class of supercuspidal representations called non-singular supercuspidal representations.
He also gave a candidate for the local Langlands correspondence for non-singular supercuspidal representations, which depends on several choices (see \cite[Remark~4.5.9]{2019arXiv191203274K}).
In \cite{MR5056155}, Solleveld and Xu proved that these choices can be made so that the resulting local Langlands correspondence satisfies Assumption~\ref{workingassumptionLLC} for depth-zero non-singular supercuspidal representations.
In this section, we further prove that the choices can be made so that the resulting local Langlands correspondence satisfies Conjecture~\ref{conjecture:comparisonoftwococycles}.

In this section, we assume that $G$ splits over a tamely ramified extension of $F$.
We also assume that $p$ is an odd prime, is not a bad prime for any irreducible factor of the absolute root system of $G$ in the sense of \cite[4.1]{springer1970conjugacy}, does not divide the number of connected components of $Z(G)$, and does not divide the order of the fundamental group of $G_{\mathrm{der}}$.
These assumptions are required for the construction of the local Langlands correspondence in \cite{2019arXiv191203274K}. 
On the other hand, these assumptions are unnecessary when considering only depth-zero non-singular supercuspidal representations. 
For the definition of depth-zero non-singular supercuspidal representations, see Section~\ref{subsec:Computation of the cohomology class attached to a non-singular supercuspidal representation}.

\begin{remark}
\label{remarkaboutthecharacteristic}
In \cite{2019arXiv191203274K}, Kaletha considers non-archimedean local fields $F$ of both zero and positive characteristic.
However, his construction of the local Langlands correspondence relies on the theory of rigid inner forms developed in \cite{MR3548533} and the generalized Tate--Nakayama duality of \cite[Section~5]{2022arXiv221002519K}, both of which are developed under the assumption that $F$ has characteristic zero.
Dillery \cite{MR4617942} has extended the results of \cite{MR3548533} to arbitrary non-archimedean local fields by replacing Kaletha's cohomology set $H^1(u\to W, Z\to G_{\qs})$ with his cohomology set $H^1(\mathcal{E}, Z\to G_{\qs})$.
He has also announced an extension of the generalized Tate--Nakayama duality in \cite[Section~5]{2022arXiv221002519K} to this cohomology set.
With these extensions in place of the corresponding results of \cite{MR3548533,2022arXiv221002519K}, Kaletha's construction in \cite{2019arXiv191203274K} can therefore be carried out without assuming that $F$ has characteristic zero.
\end{remark}

\subsection{Computation of the cohomology class attached to a non-singular supercuspidal representation}
\label{subsec:Computation of the cohomology class attached to a non-singular supercuspidal representation}

In this subsection, we compute the cohomology class $[\mu_{\mathfrak{s}}]$ for a Bernstein block $\Rep^{\mathfrak{s}}(G(F))$ consisting of non-singular supercuspidal representations.
First, we recall the construction of a non-singular supercuspidal representation in \cite[Section~3.4]{2019arXiv191203274K} briefly.
Let $(S, \theta)$ be a tame $\resi$-non-singular elliptic pair, i.e., $S$ is a tame elliptic maximal torus of $G$, and $\theta$ is a character of $S(F)$ such that $(S,\theta)$ satisfies the conditions in \cite[Definition~3.4.1]{2019arXiv191203274K}. 
In the following subsection, we will assume that $(S, \theta)$ is a tame $F$-non-singular elliptic pair in the sense of \cite[Definition~3.4.1]{2019arXiv191203274K}, which is slightly stronger than the condition that $(S, \theta)$ is a tame $\resi$-non-singular elliptic pair and necessary for the construction of the local Langlands correspondence in \cite{2019arXiv191203274K}.
See \cite[Fact~3.1.5]{2019arXiv191203274K} for the difference between $\resi$-non-singular and $F$-non-singular.
From the $k_{F}$-non-singular pair $(S, \theta)$, Kaletha constructed a tuple
\begin{equation}
\label{reducibletuple}
\left((G^i)_{i=0}^{d}, y, (r_{i})_{i=0}^{d}, \kappa_{(S, \phi_{-1})}, (\phi_i)_{i=0}^{d}\right)
\end{equation}
that satisfies the conditions {\bf{D1}}, {\bf{D2}}, {\bf{D3}}, and {\bf{D5}} in \cite[Section~3]{Yu}, where $\kappa_{(S, \phi_{-1})}$ is a possibly reducible representation of $G^0(F)_{[y]}$ such that $\rho_{-1}\restriction_{G^0(F)_{y, 0}}$ is the inflation of a cuspidal representation of $\mathsf{G}^{0, \circ}_{[y]}(\resi) = G^{0}(F)_{[y],0}/ G^{0}(F)_{[y], 0+}$ for each irreducible constituent $\rho_{-1}$ of $\kappa_{(S, \phi_{-1})}$.
Hence, for such $\rho_{-1}$, the tuple
\[
\left((G^i)_{i=0}^{d}, y, (r_{i})_{i=0}^{d}, \rho_{-1}, (\phi_i)_{i=0}^{d}\right)
\]
satisfies the conditions {\bf{D1}} to {\bf{D5}} in \cite[Section~3]{Yu}.
Then we obtain an irreducible supercuspidal representation $\pi_{(S, \theta, \rho_{-1})}$ of $G(F)$ from this datum via the  twisted Yu construction introduced in \cite[Theorem~4.1.13]{FKS}, which is a modification of Yu's original construction~\cite{Yu}.
An irreducible supercuspidal representation of $G(F)$ obtained in this way is called non-singular.

We now explain in more detail the ``depth-zero part'' $(G^0, y, \kappa_{(S, \phi_{-1})})$ of \eqref{reducibletuple}.
The group $G^0$ is the connected reductive subgroup of $G$ with maximal torus $S$ and root system
\[
\Phi(G^0, S) = \left\{
\alpha \in \Phi(G, S) \mid (\theta \circ N_{E/F} \circ \alpha^{\vee})(E_{0+}^{\times})=1
\right\},
\]
where $E$ is the splitting field of $S$, and $E_{0+}^{\times}$ is the subgroup of $E^{\times}$ defined as
\[
E_{0+}^{\times} = \left\{
x \in E^{\times} \mid \ord_{E}(x-1) > 0
\right\}.
\]
Since $(S, \theta)$ is a tame $\resi$-non-singular elliptic pair, $S$ is a maximally unramified elliptic maximal torus of $G^0$, and $\theta$ is $\resi$-non-singular with respect to $G^0$ in the sense of \cite[Definition~3.1.1]{2019arXiv191203274K}.
The point $y \in \mathcal{B}(G^0, F)$ is a lift of the unique Frobenius-fixed point in $\mathcal{A}^{\red}(G^0, S, F^{\unr})$.
Since $N_{G^{0}}(S)(F)$ acts on $\mathcal{A}^{\red}(G^0, S, F^{\unr})$ Frobenius-equivariantly, we have $N_{G^0}(S)(F) \subset G^{0}(F)_{[y]}$.
The character $\phi_{-1}$ is a depth-zero character of $S(F)$ such that $\phi_{-1} \restriction_{S^{0}_{\mathrm{sc}}(F)} = \theta \restriction_{S^0_{\mathrm{sc}}(F)}$, where $S^0_{\mathrm{sc}}$ denotes the preimage of $S$ in $G^0_{\mathrm{sc}}$ and $\restriction_{S^0_{\mathrm{sc}}(F)}$ means the composition of $\phi_{-1}$ with $S^{0}_{\mathrm{sc}}(F) \to S(F)$.
Hence, $\phi_{-1}$ is a $\resi$-non-singular character of $S(F)$ with respect to $G^0$.
The representation $\kappa_{(S, \phi_{-1})}$ of $G^0(F)_{[y]}$ is constructed by using Deligne--Lusztig theory \cite{MR393266} and its extension to disconnected groups by Kaletha \cite[Section~2]{2019arXiv191203274K}.
We recall the construction of $\kappa_{(S, \phi_{-1})}$.
Let $\mathsf{U} \subset \mathsf{G}_{[y]}^{0, \circ}$ be the unipotent radical of a Borel $\overline{\resi}$-subgroup of $\mathsf{G}_{[y]}^{0, \circ}$ containing $\mathsf{S}^\circ$, and let $\Frob$ denote the Frobenius endomorphism of $\mathsf{G}^0_{[y]}$. We define the corresponding Deligne--Lusztig variety $Y_{\mathsf{U}}^{\mathsf{G}^0_{[y]}}$ as 
\[
Y_{\mathsf{U}}^{\mathsf{G}^0_{[y]}} = \left\{
g\mathsf{U} \in \mathsf{G}^0_{[y]} / \mathsf{U} \mid g^{-1}\Frob(g) \in \mathsf{U} \cdot \Frob(\mathsf{U})
\right\}.
\]
The group $\mathsf{G}^0_{[y]}(\resi)$ acts on $Y_{\mathsf{U}}^{\mathsf{G}^0_{[y]}}$ by the left multiplication and $\mathsf{S}(\resi)$ acts on $Y_{\mathsf{U}}^{\mathsf{G}^0_{[y]}}$ by the right multiplication, so these groups act on the $\ell$-adic cohomology group $H^{i}_{c} (Y_{\mathsf{U}}^{\mathsf{G}^0_{[y]}}, \overline{\mathbb{Q}_{\ell}})$.
We define
\[
H^{i}_{c} (Y_{\mathsf{U}}^{\mathsf{G}^0_{[y]}}, \overline{\mathbb{Q}_{\ell}})_{\phi_{-1}} = H^{i}_{c} (Y_{\mathsf{U}}^{\mathsf{G}^0_{[y]}}, \overline{\mathbb{Q}_{\ell}})\otimes_{\mathsf{S}(k_{F})}\phi_{-1}.
\]
Then $\mathsf{G}^0_{[y]}(k_F)$ acts on $H^{i}_{c} (Y_{\mathsf{U}}^{\mathsf{G}^0_{[y]}}, \overline{\mathbb{Q}_{\ell}})_{\phi_{-1}}$. 
According to \cite[Corollary~2.6.3]{2019arXiv191203274K}, this component vanishes unless $i = d(\mathsf{U})$, where $d(\mathsf{U})$ denotes the number of root hyperplanes separating the Weyl chambers of $\mathsf{U}$ and $\Frob(\mathsf{U})$.
We define the representation $\kappa_{(\mathsf{S}, \phi_{-1})}^{\mathsf{G}^0_{[y]}}$ of $\mathsf{G}^0_{[y]}(\resi)$ by the action of $\mathsf{G}^0_{[y]}(\resi)$ on $H^{d(\mathsf{U})}_{c} (Y_{\mathsf{U}}^{\mathsf{G}^{0}_{[y]}}, \overline{\mathbb{Q}_{\ell}})_{\phi_{-1}}$.
Now, we define the representation $\kappa_{(S, \phi_{-1})}$ of $G^0(F)_{[y]}$ as the inflation of $\kappa_{(\mathsf{S}, \phi_{-1})}^{\mathsf{G}^0_{[y]}}$ via the natural surjection $G^0(F)_{[y]} \twoheadrightarrow \mathsf{G}^0_{[y]}(\resi)$.
We regard these representations as representations over $\mathbb{C}$ via the fixed isomorphism $\mathbb{C} \simeq \overline{\mathbb{Q}_{\ell}}$.

For an irreducible subrepresentation $\rho_{-1}$ of $\kappa_{(S, \phi_{-1})}$, the compactly induced representation $\pi^{0}_{(S, \theta, \rho_{-1})} := \ind_{G^0(F)_{[y]}}^{G^0(F)}(\rho_{-1})$ is an irreducible supercuspidal representation by \cite[Proposition~6.6]{MR1371680}.
If $G = G^0$ and $\theta = \phi_{-1}$, we have $\pi_{(S, \theta, \rho_{-1})} = \pi^{0}_{(S, \theta, \rho_{-1})}$.
In this case, we call $\pi_{(S, \theta, \rho_{-1})}$ a depth-zero non-singular supercuspidal representation.
While the construction of $\pi_{(S, \theta, \rho_{-1})}$ is in general much more complicated than that of $\pi^{0}_{(S, \theta, \rho_{-1})}$, the corresponding Bernstein blocks are equivalent.
Let $\Rep^{\mathfrak{s}^{0}}(G^0(F))$ be the Bernstein block of $G^0(F)$ containing $\pi^{0}_{(S, \theta, \rho_{-1})}$, and let $\Rep^{\mathfrak{s}}(G(F))$ be the Bernstein block of $G(F)$ containing $\pi_{(S, \theta, \rho_{-1})}$.
We also define
\[
\pi^{0, \mathrm{ess}}_{(S, \theta, \rho_{-1})} := \ind_{G^0(F)_{[y]}}^{G^0(F)}(\rho_{-1}) \otimes \prod_{i=0}^{d} \phi_{i} \restriction_{G^0(F)},
\]
and let $\Rep^{\mathfrak{s}^{0, \mathrm{ess}}}(G^0(F))$ be the Bernstein block of $G^0(F)$ containing $\pi^{0, \mathrm{ess}}_{(S, \theta, \rho_{-1})}$.

\begin{proposition}[{\cite[Corollary~4.5.1]{2024arXiv240807805A}}]
\label{prop:reductiontodepthzerounramifiedtwist}
We have equivalences of categories 
\[
\Red^{\mathfrak{s}} \colon \Rep^{\mathfrak{s}}(G(F)) \isoarrow \Rep^{\mathfrak{s}^{0}}(G^0(F)) 
\qquad
\text{and}
\qquad
\Red^{\mathfrak{s}, \mathrm{ess}} \colon \Rep^{\mathfrak{s}}(G(F)) \isoarrow \Rep^{\mathfrak{s}^{0, \mathrm{ess}}}(G^0(F))
\]
such that for $\pi \in \Rep^{\mathfrak{s}}(G(F))$ and $\chi \in X_{\unr}(G)$, we have 
\[
\Red^{\mathfrak{s}}(\pi \otimes \chi) = \Red^{\mathfrak{s}}(\pi) \otimes \chi\restriction_{G^0(F)}
\qquad
\text{and}
\qquad
\Red^{\mathfrak{s}, \mathrm{ess}}(\pi \otimes \chi) = \Red^{\mathfrak{s}, \mathrm{ess}}(\pi) \otimes \chi\restriction_{G^0(F)}.
\]
\end{proposition}
\begin{proof}
The existence of an equivalence $\Red^{\mathfrak{s}}$ is a special case of \cite[Corollary~4.5.1]{2024arXiv240807805A}.
We also obtain the equivalence $\Red^{\mathfrak{s}, \mathrm{ess}}$ by combining $\Red^{\mathfrak{s}}$ with the equivalence 
\[
\Rep^{\mathfrak{s}^{0}}(G^0(F)) \isoarrow \Rep^{\mathfrak{s}^{0, \mathrm{ess}}}(G^0(F)), \qquad \pi^0 \mapsto \pi^0 \otimes \prod_{i=0}^{d} \phi_{i} \restriction_{G^0(F)}.
\]
The compatibility with the unramified twists follows from the construction of the equivalences.
\end{proof}

Since $S \subset G^0$ and $S$ is an elliptic maximal torus of $G$, the split ranks of $Z(G)$ and $Z(G^0)$ coincide.
Hence, we have $G^0(F)^{1} = G^0(F) \cap G(F)^{1}$ and the restriction map $X_{\unr}(G) \to X_{\unr}(G^0)$ is surjective.
Let $\pi \in \Rep^{\mathfrak{s}}(G(F))$ be an irreducible representation, and put $\pi^0 = \Red^{\mathfrak{s}}(\pi)$ and $\pi^{0, \mathrm{ess}} = \Red^{\mathfrak{s}, \mathrm{ess}}(\pi)$.
According to Proposition~\ref{prop:reductiontodepthzerounramifiedtwist}, the restriction map $X_{\unr}(G) \twoheadrightarrow X_{\unr}(G^0)$ induces a surjection $X_{\unr}(G)_{\pi} \twoheadrightarrow X_{\unr}(G^0)_{\pi^0} = X_{\unr}(G^0)_{\pi^{0, \mathrm{ess}}}$.
We regard $H^2\left(
X_{\unr}(G^0)_{\pi^0}, \mathbb{C}^{\times}
\right)$ as a subset of $H^2\left(
X_{\unr}(G)_{\pi}, \mathbb{C}^{\times}
\right)$ via the inflation map
\[
H^2\left(
X_{\unr}(G^0)_{\pi^0}, \mathbb{C}^{\times}
\right) 
\hookrightarrow
H^2\left(
X_{\unr}(G)_{\pi}, \mathbb{C}^{\times}
\right).
\]

\begin{corollary}
\label{corollary:reductionofcomputationofmutodepthzeroautm}
We have 
\[
[\mu_{\mathfrak{s}^0}] = [\mu_{\mathfrak{s}^{0, \mathrm{ess}}}] = [\mu_{\mathfrak{s}}] \in 
H^2\left(
X_{\unr}(G^0)_{\pi^{0}}, \mathbb{C}^{\times}
\right)
\subset
H^2\left(
X_{\unr}(G)_{\pi}, \mathbb{C}^{\times}
\right).
\]
\end{corollary}
\begin{proof}
The claim follows from Proposition~\ref{prop:reductiontodepthzerounramifiedtwist} and the definition of $[\mu_{\mathfrak{s}}]$.
\end{proof}

Hence, in order to compute $[\mu_{\mathfrak{s}}]$, we may assume that $G = G^0$ and $\theta = \phi_{-1}$.
From now on, throughout this subsection, we assume that $\theta$ is a $\resi$-non-singular depth-zero character of $S(F)$ with respect to $G$.
In the general setting, we apply the following arguments to $(G^0, \phi_{-1})$.
We recall the parametrization of irreducible subrepresentations of $\kappa_{(S, \theta)}$ in \cite[Section~3.3]{2019arXiv191203274K}.
Let $[\kappa_{(S, \theta)}]$ denote the set of irreducible constituents of $\kappa_{(S, \theta)}$.
We note that $[\kappa_{(S, \theta)}]$ does not denote a multiset.
We set $\theta^{1} = \theta \restriction_{S(F)^1}$ and define subgroups $N_{G}(S)(F)_{\theta}$ and $N_{G}(S)(F)_{\theta^{1}}$ of $N_{G}(S)(F)$ as
\[
N_{G}(S)(F)_{\theta} = \left\{
n \in N_{G}(S)(F) \mid {^n\theta} = \theta
\right\}
\qquad
\text{and}
\qquad
N_{G}(S)(F)_{\theta^{1}} = \left\{
n \in N_{G}(S)(F) \mid {^n\theta^{1}} = \theta^{1}
\right\}.
\]
In \cite[Section~2.7]{2019arXiv191203274K}, Kaletha defined an action $R_{\mathsf{U}}^{\mathsf{G}_{[y]}, \epsilon}$ of $N_{G}(S)(F)_{\theta}$ on $H^{d(\mathsf{U})}_{c} (Y_{\mathsf{U}}^{\mathsf{G}_{[y]}}, \overline{\mathbb{Q}_{\ell}})_{\theta}$, which we now recall.
We define a set $\mathcal{U}$ of unipotent radicals of Borel $\overline{\resi}$-subgroups of $\mathsf{G}_{[y]}^\circ$ containing $\mathsf{S}^\circ$ as
\[
\mathcal{U} = \left\{
{^n\mathsf{U}} \mid n \in  N_{G}(S)(F)_{\theta^{1}}
\right\}.
\]
For each $\mathsf{U}_1, \mathsf{U}_2 \in \mathcal{U}$, let 
\[
\Psi^{\mathsf{G}_{[y]}}_{\mathsf{U}_2, \mathsf{U}_1} \colon H^{d(\mathsf{U}_1)}_{c} (Y_{\mathsf{U}_1}^{\mathsf{G}_{[y]}}, \overline{\mathbb{Q}_{\ell}})_{\theta} \isoarrow H^{d(\mathsf{U}_2)}_{c} (Y_{\mathsf{U}_2}^{\mathsf{G}_{[y]}}, \overline{\mathbb{Q}_{\ell}})_{\theta}
\]
denote the $\mathsf{G}_{[y]}(\resi)$-equivariant isomorphism defined in \cite[Section~2.6]{2019arXiv191203274K}.
According to \cite[Fact~2.7.2, Proposition~3.2.4]{2019arXiv191203274K}, there exists a homogeneous 1-cochain 
\[
\epsilon \in C^1\left(
N_{G}(S)(F)_{\theta^{1}}/S(F), \overline{\mathbb{Q}_{\ell}}^{\times}
\right)
\]
 such that 
the collection  
\[
\epsilon \cdot \Psi^{\mathsf{G}_{[y]}} =
\left\{
[\epsilon \cdot \Psi^{\mathsf{G}_{[y]}}]_{\mathsf{U}_2, \mathsf{U}_1} \colon H^{d(\mathsf{U}_1)}_{c} (Y_{\mathsf{U}_1}^{\mathsf{G}_{[y]}}, \overline{\mathbb{Q}_{\ell}})_{\theta} \isoarrow H^{d(\mathsf{U}_2)}_{c} (Y_{\mathsf{U}_2}^{\mathsf{G}_{[y]}}, \overline{\mathbb{Q}_{\ell}})_{\theta}
\right\}_{\mathsf{U}_1, \mathsf{U}_2 \in \mathcal{U}}
\]
is normalized and $N_{G}(S)(F)_{\theta}$-equivariant in the sense of \cite[Definition~2.4.4]{2019arXiv191203274K}.
Here, a homogeneous 1-cochain means the map
\[
\epsilon \colon N_{G}(S)(F)_{\theta^{1}}/S(F) \times N_{G}(S)(F)_{\theta^{1}}/S(F) \to \overline{\mathbb{Q}_{\ell}}^{\times}
\]
such that $\epsilon(w w_1, w w_2) = \epsilon(w_1, w_2)$ for all $w, w_1, w_2 \in N_{G}(S)(F)_{\theta^{1}}/S(F)$, and $[\epsilon \cdot \Psi^{\mathsf{G}_{[y]}}]$ is defined as 
\[
[\epsilon \cdot \Psi^{\mathsf{G}_{[y]}}]_{{^{w_2}\mathsf{U}}, {^{w_1}\mathsf{U}}} = \epsilon(w_2, w_1) \Psi^{\mathsf{G}_{[y]}}_{{^{w_2}\mathsf{U}}, {^{w_1}\mathsf{U}}}
\]
for $w_1, w_2 \in N_{G}(S)(F)_{\theta^{1}}/S(F)$.
We note that such a cochain $\epsilon$ is not canonical (see \cite[Remark~2.7.3]{2019arXiv191203274K}).
We denote by $[\theta]$ the set of unramified twists of $\theta$:
\[
[\theta] = \left\{
\theta \otimes \chi \mid \chi \in X_{\unr}(S)
\right\}.
\]
Replacing $\theta$ by $\theta' \in [\theta]$ in the construction of $\Psi^{\mathsf{G}_{[y]}}_{\mathsf{U}_2, \mathsf{U}_1}$, 
we also have a $\mathsf{G}_{[y]}(\resi)$-equivariant isomorphism 
\[
\Psi^{\mathsf{G}_{[y]}, \theta'}_{\mathsf{U}_2, \mathsf{U}_1} \colon H^{d(\mathsf{U}_1)}_{c} (Y_{\mathsf{U}_1}^{\mathsf{G}_{[y]}}, \overline{\mathbb{Q}_{\ell}})_{\theta'} \isoarrow H^{d(\mathsf{U}_2)}_{c} (Y_{\mathsf{U}_2}^{\mathsf{G}_{[y]}}, \overline{\mathbb{Q}_{\ell}})_{\theta'}
\]
for each $\mathsf{U}_1, \mathsf{U}_2 \in \mathcal{U}$.
According to \cite[Corollary~2.5.8]{2019arXiv191203274K}, the collection 
\[
\epsilon \cdot \Psi^{\mathsf{G}_{[y]}, \theta'} =
\left\{
[\epsilon \cdot \Psi^{\mathsf{G}_{[y]}, \theta'}]_{\mathsf{U}_2, \mathsf{U}_1} \colon H^{d(\mathsf{U}_1)}_{c} (Y_{\mathsf{U}_1}^{\mathsf{G}_{[y]}}, \overline{\mathbb{Q}_{\ell}})_{\theta'} \isoarrow H^{d(\mathsf{U}_2)}_{c} (Y_{\mathsf{U}_2}^{\mathsf{G}_{[y]}}, \overline{\mathbb{Q}_{\ell}})_{\theta'}
\right\}_{\mathsf{U}_1, \mathsf{U}_2 \in \mathcal{U}}
\]
is also normalized and $N_{G}(S)(F)_{\theta'}$-equivariant.
For $n \in N_{G}(S)(F)_{\theta^{1}}$ and $\theta' \in [\theta]$, we define a $\mathsf{G}_{[y]}(\resi)$-equivariant isomorphism 
\[
R_{\mathsf{U}}^{\mathsf{G}_{[y]}, \epsilon, \theta'}(n) \colon H^{d(\mathsf{U})}_{c} (Y_{\mathsf{U}}^{\mathsf{G}_{[y]}}, \overline{\mathbb{Q}_{\ell}})_{\theta'} \isoarrow H^{d(\mathsf{U})}_{c} (Y_{\mathsf{U}}^{\mathsf{G}_{[y]}}, \overline{\mathbb{Q}_{\ell}})_{^n\theta'}
\]
as 
\[
R_{\mathsf{U}}^{\mathsf{G}_{[y]}, \epsilon, \theta'}(n) = \Ad(n)_{*} \circ [\epsilon \cdot \Psi^{\mathsf{G}_{[y]}, \theta'}]_{{}^{n^{-1}}\mathsf{U}, \mathsf{U}} \circ \kappa_{(\mathsf{S}, \theta')}^{\mathsf{G}_{[y]}}(n)^{-1},
\]
where $\Ad(n)_{*}$ is the map induced by $\Ad(n)$.
For $\theta' = \theta$, we simply write $R_{\mathsf{U}}^{\mathsf{G}_{[y]}, \epsilon, \theta}(n) = R_{\mathsf{U}}^{\mathsf{G}_{[y]}, \epsilon}(n)$.

\begin{remark}
In \cite[Section~2.7]{2019arXiv191203274K}, the operator $R_{\mathsf{U}}^{\mathsf{G}_{[y]}, \epsilon}(n)$ is defined as 
\[
R_{\mathsf{U}}^{\mathsf{G}_{[y]}, \epsilon}(n) =  [\epsilon \cdot \Psi^{\mathsf{G}_{[y]}}]_{\mathsf{U}, {^n\mathsf{U}}}  \circ r_{n^{-1}} \colon H^{d(\mathsf{U})}_{c} (Y_{\mathsf{U}}^{\mathsf{G}_{[y]}}, \overline{\mathbb{Q}_{\ell}})_{\theta} \to H^{d(\mathsf{U})}_{c} (Y_{\mathsf{U}}^{\mathsf{G}_{[y]}}, \overline{\mathbb{Q}_{\ell}})_{\theta}
\]
for $n \in N_{G}(S)(F)_{\theta}$, where 
\[
r_{n^{-1}} \colon H^{d(\mathsf{U})}_{c} (Y_{\mathsf{U}}^{\mathsf{G}_{[y]}}, \overline{\mathbb{Q}_{\ell}})_{\theta} \isoarrow H^{d({}^{n}\mathsf{U})}_{c} (Y_{{}^{n}\mathsf{U}}^{\mathsf{G}_{[y]}}, \overline{\mathbb{Q}_{\ell}})_{^n\theta}
\]
is a $\mathsf{G}_{[y]}(\resi)$-equivariant isomorphism induced from the isomorphism of varieties 
\[
Y_{\mathsf{U}}^{\mathsf{G}_{[y]}} \isoarrow Y_{^n\mathsf{U}}^{\mathsf{G}_{[y]}}, \qquad g \mathsf{U} \mapsto g \mathsf{U} n^{-1} = gn^{-1} \cdot {^n\mathsf{U}}.
\]
Let us verify that our definition of $R_{\mathsf{U}}^{\mathsf{G}_{[y]}, \epsilon}(n)$ is compatible with the definition in \cite{2019arXiv191203274K}.
By the definitions of $r_{n^{-1}}$ and $\kappa_{(\mathsf{S}, \theta)}^{\mathsf{G}_{[y]}}$, we have $r_{n^{-1}} = \Ad(n)_{*} \circ \kappa_{(\mathsf{S}, \theta)}^{\mathsf{G}_{[y]}}(n)^{-1}$.
Moreover, since $\epsilon \cdot \Psi^{\mathsf{G}_{[y]}}$ is $N_{G}(S)(F)_{\theta}$-equivariant, we have 
\[
\Ad(n)_{*} \circ [\epsilon \cdot \Psi^{\mathsf{G}_{[y]}}]_{n^{-1}\mathsf{U}n, \mathsf{U}} \circ  \Ad(n)_{*}^{-1} = [\epsilon \cdot \Psi^{\mathsf{G}_{[y]}}]_{\mathsf{U}, {^n\mathsf{U}}}.
\]
Thus, we obtain that 
\begin{align*}
[\epsilon \cdot \Psi^{\mathsf{G}_{[y]}}]_{\mathsf{U}, {^n\mathsf{U}}}  \circ r_{n^{-1}}  &= [\epsilon \cdot \Psi^{\mathsf{G}_{[y]}}]_{\mathsf{U}, {^n\mathsf{U}}}  \circ \Ad(n)_{*} \circ \kappa_{(\mathsf{S}, \theta)}^{\mathsf{G}_{[y]}}(n)^{-1} \\
&= \Ad(n)_{*} \circ [\epsilon \cdot \Psi^{\mathsf{G}_{[y]}}]_{n^{-1}\mathsf{U}n, \mathsf{U}} \circ \kappa_{(\mathsf{S}, \theta)}^{\mathsf{G}_{[y]}}(n)^{-1}.
\end{align*}
\end{remark}

Since the collection $\epsilon \cdot \Psi^{\mathsf{G}_{[y]}}$ is normalized and $N_{G}(S)(F)_{\theta}$-equivariant, we obtain the following multiplicativity of $R_{\mathsf{U}}^{\mathsf{G}_{[y]}, \epsilon}$.

\begin{lemma}
\label{lemma:normalizedR}
For $n_1, n_2 \in N_{G}(S)(F)_{\theta^{1}}$ and $\theta' \in [\theta]$, we have 
\[
R_{\mathsf{U}}^{\mathsf{G}_{[y]}, \epsilon, {^{n_2}\theta'}}(n_{1}) \circ R_{\mathsf{U}}^{\mathsf{G}_{[y]}, \epsilon, \theta'}(n_{2}) = R_{\mathsf{U}}^{\mathsf{G}_{[y]}, \epsilon, \theta'}(n_{1} n_{2}).
\]
\end{lemma}
\begin{proof}
According to \cite[Lemma~2.5.3, Fact~2.6.12]{2019arXiv191203274K}, we have 
\[
\Ad(n_2)_{*}^{-1} \circ [\epsilon \cdot \Psi^{\mathsf{G}_{[y]}, {^{n_2}\theta'}}]_{n_1^{-1}\mathsf{U}n_1, \mathsf{U}}  \circ \Ad(n_2)_{*} = [\epsilon \cdot \Psi^{\mathsf{G}_{[y]}, \theta'}]_{(n_1 n_2)^{-1}\mathsf{U}n_1 n_2, n_{2}^{-1}\mathsf{U} n_2}. 
\]
Moreover, since $\epsilon \cdot \Psi^{\mathsf{G}_{[y]}, \theta'}$ is normalized, we have 
\[
 [\epsilon \cdot \Psi^{\mathsf{G}_{[y]}, \theta'}]_{(n_1 n_2)^{-1}\mathsf{U}n_1 n_2, n_{2}^{-1}\mathsf{U} n_2} \circ  [\epsilon \cdot \Psi^{\mathsf{G}_{[y]}, \theta'}]_{n_2^{-1}\mathsf{U} n_2, \mathsf{U}} =   [\epsilon \cdot \Psi^{\mathsf{G}_{[y]}, \theta'}]_{(n_1 n_2)^{-1}\mathsf{U}n_1 n_2, \mathsf{U}}.
\]
We also have, by the definitions, 
\[
\kappa_{(\mathsf{S}, {^{n_2}\theta'})}^{\mathsf{G}_{[y]}}(n_1)^{-1}
\circ 
\Ad(n_2)_{*}
= \Ad(n_2)_{*} \circ \kappa_{(\mathsf{S}, \theta')}^{\mathsf{G}_{[y]}}(n_{2}^{-1} n_1 n_2)^{-1}
\]
Thus, we obtain that 
\begin{align*}
&\qquad R_{\mathsf{U}}^{\mathsf{G}_{[y]}, \epsilon, {^{n_2}\theta'}}(n_{1}) \circ R_{\mathsf{U}}^{\mathsf{G}_{[y]}, \epsilon, \theta'}(n_{2}) \\
&=
\left(
\Ad(n_1)_{*} \circ [\epsilon \cdot \Psi^{\mathsf{G}_{[y]}, {^{n_2}\theta'}}]_{{}^{n_{1}^{-1}}\mathsf{U}, \mathsf{U}} \circ \kappa_{(\mathsf{S}, {^{n_2}\theta'})}^{\mathsf{G}_{[y]}}(n_1)^{-1}
\right) \circ \left(
\Ad(n_2)_{*} \circ [\epsilon \cdot \Psi^{\mathsf{G}_{[y]}, \theta'}]_{{}^{n_{2}^{-1}}\mathsf{U}, \mathsf{U}} \circ \kappa_{(\mathsf{S}, \theta')}^{\mathsf{G}_{[y]}}(n_2)^{-1}
\right) \\
&= \Ad(n_1 n_2)_{*} \circ  [\epsilon \cdot \Psi^{\mathsf{G}_{[y]}, \theta'}]_{{}^{(n_{1} n_{2})^{-1}}\mathsf{U}, \mathsf{U}} \circ  \kappa_{(\mathsf{S}, \theta')}^{\mathsf{G}_{[y]}}(n_{2}^{-1} n_1 n_2)^{-1} \circ \kappa_{(\mathsf{S}, \theta')}^{\mathsf{G}_{[y]}}(n_2)^{-1} \\
&= \Ad(n_1 n_2)_{*} \circ  [\epsilon \cdot \Psi^{\mathsf{G}_{[y]}, \theta'}]_{{}^{(n_1 n_2)^{-1}}\mathsf{U}, \mathsf{U}} \circ  \kappa_{(\mathsf{S}, \theta')}^{\mathsf{G}_{[y]}}(n_1 n_2)^{-1} \\
&= R_{\mathsf{U}}^{\mathsf{G}_{[y]}, \epsilon, \theta'}(n_1 n_2). 
\qedhere
\end{align*}

\end{proof}

In particular, the map $n \mapsto R_{\mathsf{U}}^{\mathsf{G}_{[y]}, \epsilon}(n)$ defines an $N_{G}(S)(F)_{\theta}$-representation on $H^{d(\mathsf{U})}_{c} (Y_{\mathsf{U}}^{\mathsf{G}_{[y]}}, \overline{\mathbb{Q}_{\ell}})_{\theta}$.
Moreover, since $R_{\mathsf{U}}^{\mathsf{G}_{[y]}, \epsilon}$ commutes with the action $\kappa^{\mathsf{G}_{[y]}}_{(\mathsf{S}, \theta)}$ of $\mathsf{G}_{[y]}(\resi)$, we obtain an action of $\mathsf{G}_{[y]}(\resi) \times N_{G}(S)(F)_{\theta}$ on $H^{d(\mathsf{U})}_{c} (Y_{\mathsf{U}}^{\mathsf{G}_{[y]}}, \overline{\mathbb{Q}_{\ell}})_{\theta}$.
Let $\kappa^{\mathsf{G}_{[y]}, \epsilon}_{(\mathsf{S}, \theta)}$ denote this $\mathsf{G}_{[y]}(\resi) \times N_{G}(S)(F)_{\theta}$-representation.
We record a consequence of Lemma~\ref{lemma:normalizedR} for later use.

\begin{lemma}
\label{lemma:intertwiningoperator}
For $n \in N_{G}(S)(F)_{\theta^{1}}$, we have 
\[
R_{\mathsf{U}}^{\mathsf{G}_{[y]}, \epsilon}(n)  \in \Hom_{\mathsf{G}_{[y]}(\resi) \times N_{G}(S)(F)_{\theta}}\left(
\kappa^{\mathsf{G}_{[y]}, \epsilon}_{(\mathsf{S}, \theta)}, {^{(1 \times n^{-1})}\kappa^{\mathsf{G}_{[y]}, \epsilon}_{(\mathsf{S}, {^n\theta})}}
\right).
\]
\end{lemma}
\begin{proof}
The claim follows from Lemma~\ref{lemma:normalizedR} and the fact that $R_{\mathsf{U}}^{\mathsf{G}_{[y]}, \epsilon}(n)$ commutes with the $\mathsf{G}_{[y]}(\resi)$-action.
\end{proof}

For $\rho \in \Irr\left(
N_{G}(S)(F)_{\theta}, \theta
\right)$, we define 
\[
\kappa_{(\mathsf{S}, \theta, \rho)}^{\mathsf{G}_{[y]}, \epsilon} = 
\Hom_{N_{G}(S)(F)_{\theta}}\left(
\rho^{\vee}, \kappa^{\mathsf{G}_{[y]}, \epsilon}_{(\mathsf{S}, \theta)}
\right),
\]
where $\rho^{\vee}$ is the contragredient representation of $\rho$.
It is a $\mathsf{G}_{[y]}(\resi)$-representation.
We denote by $\kappa_{(S, \theta, \rho)}^{\epsilon}$ the representation of $G(F)_{[y]}$ obtained by the inflation of $\kappa_{(\mathsf{S}, \theta, \rho)}^{\mathsf{G}_{[y]}, \epsilon}$ via $G(F)_{[y]} \twoheadrightarrow \mathsf{G}_{[y]}(\resi)$.

\begin{proposition}[{\cite[Theorem~2.7.7]{2019arXiv191203274K}}]
The map $\rho \mapsto \kappa_{(S, \theta, \rho)}^{\epsilon}$ is a bijection
\[
\Irr\left(
N_{G}(S)(F)_{\theta}, \theta
\right) \isoarrow [\kappa_{(S, \theta)}].
\]
\end{proposition}

For $\rho \in \Irr\left(
N_{G}(S)(F)_{\theta}, \theta
\right)$, we define a non-singular supercuspidal representation $\pi^{\epsilon}_{(S, \theta, \rho)}$ as 
\[
\pi^{\epsilon}_{(S, \theta, \rho)} = \ind_{G(F)_{[y]}}^{G(F)}(\kappa_{(S, \theta, \rho)}^{\epsilon}).
\]
We compute the cohomology class $[\mu_{\mathfrak{s}}]$, where $\Rep^{\mathfrak{s}}(G(F))$ is the Bernstein block containing $\pi^{\epsilon}_{(S, \theta, \rho)}$.
For this purpose, we construct an isomorphism $I_{\chi} \colon \pi^{\epsilon}_{(S, \theta, \rho)} \isoarrow \pi^{\epsilon}_{(S, \theta, \rho)} \otimes \chi$ for each $\chi \in X_{\unr}(G)_{\pi}$.
We first define an isomorphism $I^{\kappa_{(S, \theta, \rho)}^{\epsilon}}_{\chi} \colon \kappa_{(\mathsf{S}, \theta, \rho)}^{\mathsf{G}_{[y]}, \epsilon} \isoarrow \kappa_{(\mathsf{S}, \theta, \rho)}^{\mathsf{G}_{[y]}, \epsilon} \otimes \chi$.
Let $\chi \in X_{\unr}(G)_{\pi}$.
According to \cite[Proposition~3.4.6]{2019arXiv191203274K}, we have $\pi^{\epsilon}_{(S, \theta, \rho)} \otimes \chi \simeq \pi^{\epsilon}_{(S, \theta \otimes \chi, \rho \otimes \chi)}$.
Hence, we obtain from \cite[Corollary~3.3.3]{2019arXiv191203274K} that there exists $n_{\chi} \in N_{G}(S)(F)_{\theta^{1}}$ such that $^{n_{\chi}} \rho \simeq \rho \otimes \chi$.
We fix such an $n_{\chi}$ and an isomorphism
\[
I_{n_{\chi}} \colon ^{n_{\chi}} \rho \isoarrow \rho \otimes \chi.
\]
If $n'_{\chi} \in N_{G}(S)(F)_{\theta^{1}}$ is another element such that ${^{n'_{\chi}} \rho} \simeq \rho \otimes \chi$, we have ${^{n'_{\chi}}\theta} = {^{n_\chi}\theta} = \theta \otimes \chi$.
Hence, we have $n_{\chi}^{-1}n'_{\chi} \in N_{G}(S)(F)_{\theta}$.
We define an isomorphism
\[
I_{n'_{\chi}} \colon ^{n'_{\chi}} \rho \isoarrow \rho \otimes \chi
\]
as $I_{n'_{\chi}} = I_{n_{\chi}} \circ \rho(n_{\chi}^{-1}n'_{\chi})$.
For $\chi_1, \chi_2 \in X_{\unr}(G)_{\pi}$, the composition $I_{n_{\chi_1}} \circ I_{n_{\chi_2}}$ of $I_{n_{\chi_{1}}} \in \Hom_{N_{G}(S)(F)_{\theta}}\left(
{^{n_{\chi_1}}\rho}, \rho\otimes \chi_{1}
\right) = \Hom_{N_{G}(S)(F)_{\theta}}\left(
{^{n_{\chi_1}}\rho}\otimes\chi_{2}, \rho\otimes \chi_{1}\chi_{2}
\right)$ and 
\[
I_{n_{\chi_2}} \in \Hom_{N_{G}(S)(F)_{\theta}}\left(
{^{n_{\chi_2}}\rho}, \rho\otimes \chi_{2}
\right) = \Hom_{N_{G}(S)(F)_{\theta}}\left(
{^{n_{\chi_1}n_{\chi_2}}\rho}, {^{n_{\chi_1}}\rho\otimes \chi_{2}}
\right) \subset \End_{\mathbb{C}}(\rho)
\]
is a non-zero element of $\Hom_{N_{G}(S)(F)_{\theta}}\left(
{^{n_{\chi_1} n_{\chi_2}}\rho}, \rho\otimes \chi_{1}\chi_{2}
\right)$.
By Schur's lemma, there exists 
\[
\mu_{\rho}\left((n_{\chi_1}, \chi_1), (n_{\chi_2}, \chi_2)\right) \in \mathbb{C}^{\times}
\]
such that 
\[
I_{n_{\chi_1}} \circ I_{n_{\chi_2}} = \mu_{\rho}\left((n_{\chi_1}, \chi_1), (n_{\chi_2}, \chi_2)\right) \cdot I_{n_{\chi_1} n_{\chi_2}}.
\]

To construct $I^{\kappa_{(S, \theta, \rho)}^{\epsilon}}_{\chi}$, we recall another isomorphism.
According to \cite[Lemma~2.6.10]{2019arXiv191203274K}, for $\theta' \in [\theta]$ and $\chi \in X_{\unr}(G)$, we have a canonical isomorphism 
\[
T^{\mathsf{G}_{[y]}, \theta'}_{\chi} \colon \kappa^{\mathsf{G}_{[y]}}_{(\mathsf{S}, \theta' \otimes \chi)} \isoarrow \kappa^{\mathsf{G}_{[y]}}_{(\mathsf{S}, \theta')} \otimes \chi.
\]
For later use, we briefly recall the construction of $T^{\mathsf{G}_{[y]}, \theta'}_{\chi}$ and record one of its properties.
Replacing $\mathsf{G}_{[y]}$ by $\mathsf{G}'_{[y]} := \mathsf{S} \cdot \mathsf{G}^{\circ}_{[y]}$ and $\mathsf{G}^{\circ}_{[y]}$ in the construction of $\kappa_{(\mathsf{S}, \theta')}^{\mathsf{G}_{[y]}}$, we obtain representations $\kappa_{(\mathsf{S}, \theta')}^{\mathsf{G}'_{[y]}}$ of $\mathsf{G}'_{[y]}(\resi)$ and $\kappa_{(\mathsf{S}^{\circ}, \theta')}^{\mathsf{G}^{\circ}_{[y]}}$ of $\mathsf{G}^{\circ}_{[y]}(\resi)$, respectively.
Note that $\kappa_{(\mathsf{S}^{\circ}, \theta')}^{\mathsf{G}^{\circ}_{[y]}}$ is the usual Deligne--Lusztig representation of $\mathsf{G}^{\circ}_{[y]}(\resi)$ arising from the pair $(\mathsf{S}^{\circ}, \theta'\restriction_{\mathsf{S}^{\circ}(\resi)})$.
According to \cite[Proposition~1.25, Corollary~1.27]{MR393266}, the natural isomorphism 
\[
\ind_{A}^{B}(V \otimes \chi) \isoarrow \ind_{A}^{B}(V) \otimes \chi, \qquad f \mapsto \chi^{-1} \cdot f,
\]
where $A = \left[
\left(
\left(
\mathsf{G}_{[y]}^\circ
\right)_{\mathrm{sc}} \times \mathsf{S}^{\circ}_{\mathrm{sc}}
\right)/Z\left(
\left(\mathsf{G}_{[y]}^\circ\right)_{\mathrm{sc}}
\right)
\right](\resi)$, $B = \left[
\left(
\mathsf{G}_{[y]}^\circ \times \mathsf{S}^{\circ}
\right)/Z\left(\mathsf{G}_{[y]}^\circ
\right)\right](\resi)$, and $V$ is the $\ell$-adic cohomology group of a Deligne--Lusztig variety for $\left(\mathsf{G}_{[y]}^\circ\right)_{\mathrm{sc}}$, induces isomorphisms
\[
\kappa_{(\mathsf{S}^{\circ}, \theta' \otimes \chi)}^{\mathsf{G}^{\circ}_{[y]}} \isoarrow \kappa_{(\mathsf{S}^{\circ}, \theta')}^{\mathsf{G}^{\circ}_{[y]}} \otimes \chi 
\qquad
\text{and}
\qquad
\kappa_{(\mathsf{S}, \theta' \otimes \chi)}^{\mathsf{G}'_{[y]}} \isoarrow \kappa_{(\mathsf{S}, \theta')}^{\mathsf{G}'_{[y]}} \otimes \chi.
\]
For more details, see \cite[Section~2.5]{2019arXiv191203274K}.
Moreover, according to \cite[Corollary~2.6.4, Remark~2.6.5]{2019arXiv191203274K}, we have a natural isomorphism
\[
\kappa_{(\mathsf{S}, \theta')}^{\mathsf{G}_{[y]}} \simeq \ind_{\mathsf{G}'_{[y]}(\resi)}^{\mathsf{G}_{[y]}(\resi)}\left(
\kappa_{(\mathsf{S}, \theta')}^{\mathsf{G}'_{[y]}}
\right).
\]
Thus, we obtain a natural isomorphism
\[
T^{\mathsf{G}_{[y]}, \theta'}_{\chi} \colon \kappa^{\mathsf{G}_{[y]}}_{(\mathsf{S}, \theta' \otimes \chi)} = \ind_{\mathsf{G}'_{[y]}(\resi)}^{\mathsf{G}_{[y]}(\resi)}\left(
\kappa_{(\mathsf{S}, \theta' \otimes \chi)}^{\mathsf{G}'_{[y]}}
\right) \isoarrow \ind_{\mathsf{G}'_{[y]}(\resi)}^{\mathsf{G}_{[y]}(\resi)}\left(
\kappa_{(\mathsf{S}, \theta')}^{\mathsf{G}'_{[y]}} \otimes \chi
\right) = \kappa^{\mathsf{G}_{[y]}}_{(\mathsf{S}, \theta')} \otimes \chi.
\]

We record a multiplicative property of $T^{\mathsf{G}_{[y]}, \theta'}_{\chi}$.
Note that for $\theta' \in [\theta]$ and $\chi_{1}, \chi_{2} \in X_{\unr}(G)$, we have
\[
T^{\mathsf{G}_{[y]}, \theta'}_{\chi_{1}} \in 
\Hom_{\mathsf{G}_{[y]}(\resi)}\left(
\kappa^{\mathsf{G}_{[y]}}_{(\mathsf{S}, \theta' \otimes \chi_{1})}, \kappa^{\mathsf{G}_{[y]}}_{(\mathsf{S}, \theta')} \otimes \chi_{1}
\right) =
\Hom_{\mathsf{G}_{[y]}(\resi)}\left(
\kappa^{\mathsf{G}_{[y]}}_{(\mathsf{S}, \theta' \otimes \chi_{1})} \otimes \chi_{2}, \kappa^{\mathsf{G}_{[y]}}_{(\mathsf{S}, \theta')} \otimes \chi_{1}\chi_{2}
\right).
\]
\begin{lemma}
\label{lemma:multipricativityofT}
For $\theta' \in [\theta]$ and $\chi_{1}, \chi_{2} \in X_{\unr}(G)$, we have 
\[
T^{\mathsf{G}_{[y]}, \theta'}_{\chi_{1}} \circ T^{\mathsf{G}_{[y]}, \theta' \otimes \chi_{1}}_{\chi_{2}} = T^{\mathsf{G}_{[y]}, \theta'}_{\chi_{1} \chi_{2}}.
\]
\end{lemma}
\begin{proof}
The claim follows from the construction of $T^{\mathsf{G}_{[y]}, \theta'}_{\chi}$.
\end{proof}

We write $T^{\mathsf{G}_{[y]}, \theta}_{\chi} = T^{\mathsf{G}_{[y]}}_{\chi}$.
As explained in the proof of \cite[Theorem~2.7.7]{2019arXiv191203274K}, we obtain from \cite[Lemma~2.5.7]{2019arXiv191203274K} that
\begin{align*}
T^{\mathsf{G}_{[y]}}_{\chi} &\in \Hom_{\mathsf{G}_{[y]}(\resi) \times N_{G}(S)(F)_{\theta}}\left(
\kappa^{\mathsf{G}_{[y]}, \epsilon}_{(\mathsf{S}, \theta \otimes \chi)}, \kappa^{\mathsf{G}_{[y]}, \epsilon}_{(\mathsf{S}, \theta)} \otimes (\chi \boxtimes \chi^{-1})
\right) \\
&= \Hom_{\mathsf{G}_{[y]}(\resi) \times N_{G}(S)(F)_{\theta}}\left(
{^{(1 \times n_{\chi}^{-1})}\kappa^{\mathsf{G}_{[y]}, \epsilon}_{(\mathsf{S}, \theta \otimes \chi)}}, {^{(1 \times n_{\chi}^{-1})}\kappa^{\mathsf{G}_{[y]}, \epsilon}_{(\mathsf{S}, \theta)}} \otimes (\chi \boxtimes \chi^{-1})
\right).
\end{align*}
Combining $T^{\mathsf{G}_{[y]}}_{\chi}$ with Lemma~\ref{lemma:intertwiningoperator}, we obtain an isomorphism 
\[
T^{\mathsf{G}_{[y]}}_{\chi} \circ R_{\mathsf{U}}^{\mathsf{G}_{[y]}, \epsilon}(n_{\chi}) \in
\Hom_{\mathsf{G}_{[y]}(\resi) \times N_{G}(S)(F)_{\theta}}\left(
\kappa^{\mathsf{G}_{[y]}, \epsilon}_{(\mathsf{S}, \theta)}, {^{(1 \times n_{\chi}^{-1})}\kappa^{\mathsf{G}_{[y]}, \epsilon}_{(\mathsf{S}, \theta)}} \otimes (\chi \boxtimes \chi^{-1})
\right).
\]
We now define the isomorphism 
\begin{align*}
I^{\kappa_{(S, \theta, \rho)}^{\epsilon}}_{\chi} \colon \kappa_{(\mathsf{S}, \theta, \rho)}^{\mathsf{G}_{[y]}, \epsilon} 
&= 
\Hom_{N_{G}(S)(F)_{\theta}}\left(
\rho ^{\vee}, \kappa^{\mathsf{G}_{[y]}, \epsilon}_{(\mathsf{S}, \theta)}
\right) \\
&\isoarrow 
\Hom_{N_{G}(S)(F)_{\theta}}\left(
({^{n_{\chi}^{-1}}\rho} \otimes \chi)^{\vee}, {^{(1 \times n_{\chi}^{-1})}\kappa^{\mathsf{G}_{[y]}, \epsilon}_{(\mathsf{S}, \theta)}} \otimes (\chi \boxtimes \chi^{-1})
\right) \\
&= 
\Hom_{N_{G}(S)(F)_{\theta}}\left(
\rho^{\vee}, \kappa^{\mathsf{G}_{[y]}, \epsilon}_{(\mathsf{S}, \theta)} \otimes (\chi \boxtimes 1) 
\right)\\
& = \kappa_{(\mathsf{S}, \theta, \rho)}^{\mathsf{G}_{[y]}, \epsilon} \otimes \chi
\end{align*}
as 
\[
I^{\kappa_{(S, \theta, \rho)}^{\epsilon}}_{\chi}(\Phi) = \left(
T^{\mathsf{G}_{[y]}}_{\chi} \circ R_{\mathsf{U}}^{\mathsf{G}_{[y]}, \epsilon}(n_{\chi})
\right) \circ \Phi \circ I_{n_{\chi}}^{\vee}
\]
for $\Phi \in \Hom_{N_{G}(S)(F)_{\theta}}\left(
\rho ^{\vee}, \kappa^{\mathsf{G}_{[y]}, \epsilon}_{(\mathsf{S}, \theta)}
\right)$, where 
\[
I_{n_{\chi}}^{\vee} \in \Hom_{N_{G}(S)(F)_{\theta}}\left(
({^{n_{\chi}^{-1}}\rho} \otimes \chi)^{\vee}, \rho^{\vee}
\right)
\]
is the dual of 
\[
I_{n_{\chi}}  \in \Hom_{N_{G}(S)(F)_{\theta}}\left(
{^{n_{\chi}} \rho}, \rho \otimes \chi
\right) = 
 \Hom_{N_{G}(S)(F)_{\theta}}\left(
 \rho, {^{n_{\chi}^{-1}}\rho} \otimes \chi
\right).
\]
The isomorphism $I^{\kappa_{(S, \theta, \rho)}^{\epsilon}}_{\chi}$ does not depend on the choice of $n_{\chi}$ in the following sense.
\begin{lemma}
\label{lemma:independenceonlift}
Let $n'_{\chi} \in N_{G}(S)(F)_{\theta^{1}}$ such that ${^{n'_{\chi}} \rho} \simeq \rho \otimes \chi$.
Then we have
\[
\left(
T^{\mathsf{G}_{[y]}}_{\chi} \circ R_{\mathsf{U}}^{\mathsf{G}_{[y]}, \epsilon}(n_{\chi})
\right) \circ \Phi \circ I_{n_{\chi}}^{\vee} = \left(
T^{\mathsf{G}_{[y]}}_{\chi} \circ R_{\mathsf{U}}^{\mathsf{G}_{[y]}, \epsilon}(n'_{\chi})
\right) \circ \Phi \circ I_{n'_{\chi}}^{\vee}
\]
for all $\Phi \in \Hom_{N_{G}(S)(F)_{\theta}}\left(
\rho ^{\vee}, \kappa^{\mathsf{G}_{[y]}, \epsilon}_{(\mathsf{S}, \theta)}
\right)$.
\end{lemma}
\begin{proof}
We write $n = n^{-1}_{\chi} n'_{\chi}$.
According to Lemma~\ref{lemma:normalizedR}, we have 
\[
R_{\mathsf{U}}^{\mathsf{G}_{[y]}, \epsilon}(n'_{\chi}) = R_{\mathsf{U}}^{\mathsf{G}_{[y]}, \epsilon}(n_{\chi}) \circ R_{\mathsf{U}}^{\mathsf{G}_{[y]}, \epsilon}(n)
= R_{\mathsf{U}}^{\mathsf{G}_{[y]}, \epsilon}(n_{\chi}) \circ \kappa^{\mathsf{G}_{[y]}, \epsilon}_{(\mathsf{S}, \theta)}(1, n).
\]
We also obtain from the definition of $I_{n'_{\chi}}$ that 
\[
I_{n'_{\chi}}^{\vee} = \left(
I_{n_{\chi}} \circ \rho(n)
\right)^{\vee} = \rho^{\vee}(n)^{-1} \circ I_{n_{\chi}}^{\vee}.
\]
Thus, the claim follows from the fact $\Phi \in \Hom_{N_{G}(S)(F)_{\theta}}\left(
\rho ^{\vee}, \kappa^{\mathsf{G}_{[y]}, \epsilon}_{(\mathsf{S}, \theta)}
\right)$.
\end{proof}

We also regard $I^{\kappa_{(S, \theta, \rho)}^{\epsilon}}_{\chi}$ as an isomorphism $\kappa_{(S, \theta, \rho)}^{\epsilon} \isoarrow \kappa_{(S, \theta, \rho)}^{\epsilon} \otimes \chi$ via the inflation.
We study the relation between $I^{\kappa_{(S, \theta, \rho)}^{\epsilon}}_{\chi_1} \circ I^{\kappa_{(S, \theta, \rho)}^{\epsilon}}_{\chi_2}$ and $I^{\kappa_{(S, \theta, \rho)}^{\epsilon}}_{\chi_1 \chi_2}$ for $\chi_1, \chi_2 \in X_{\unr}(G)_{\pi}$.

\begin{lemma}
\label{lemma:RTvsTR}
For $\chi \in X_{\unr}(G)$ and $n \in N_{G}(S)(F)_{\theta^{1}}$, we have 
\[
R_{\mathsf{U}}^{\mathsf{G}_{[y]}, \epsilon}(n) \circ T^{\mathsf{G}_{[y]}}_{\chi} = \chi(n) \cdot T^{\mathsf{G}_{[y]}, {^n\theta}}_{\chi} \circ R_{\mathsf{U}}^{\mathsf{G}_{[y]}, \epsilon, \theta \otimes \chi}(n).
\]
\end{lemma}
\begin{proof}
Since $T^{\mathsf{G}_{[y]}}_{\chi} \in \Hom_{\mathsf{G}_{[y]}(\resi)}\left(
\kappa^{\mathsf{G}_{[y]}}_{(\mathsf{S}, \theta \otimes \chi)}, \kappa^{\mathsf{G}_{[y]}}_{(\mathsf{S}, \theta)} \otimes \chi
\right)$, we have 
\[
\kappa_{(\mathsf{S}, \theta)}^{\mathsf{G}_{[y]}}(n)^{-1} \circ T^{\mathsf{G}_{[y]}}_{\chi} = \chi(n) \cdot T^{\mathsf{G}_{[y]}}_{\chi} \circ \kappa_{(\mathsf{S}, \theta \otimes \chi)}^{\mathsf{G}_{[y]}}(n)^{-1}.
\]
Moreover, according to \cite[Lemma~2.5.7]{2019arXiv191203274K}, we have
\[
\Ad(n)_{*} \circ [\epsilon \cdot \Psi^{\mathsf{G}_{[y]}}]_{n^{-1}\mathsf{U}n, \mathsf{U}} \circ T^{\mathsf{G}_{[y]}}_{\chi} = T^{\mathsf{G}_{[y]}, {^n\theta}}_{\chi} \circ \Ad(n)_{*} \circ [\epsilon \cdot \Psi^{\mathsf{G}_{[y]}, \theta \otimes \chi}]_{n^{-1}\mathsf{U}n, \mathsf{U}}.
\]
Thus, we conclude that 
\begin{align*}
R_{\mathsf{U}}^{\mathsf{G}_{[y]}, \epsilon}(n) \circ T_{\chi}^{\mathsf{G}_{[y]}} &= \left(
\Ad(n)_{*} \circ [\epsilon \cdot \Psi^{\mathsf{G}_{[y]}}]_{n^{-1}\mathsf{U}n, \mathsf{U}} \circ \kappa_{(\mathsf{S}, \theta)}^{\mathsf{G}_{[y]}}(n)^{-1}
\right) \circ T_{\chi}^{\mathsf{G}_{[y]}} \\
&= \chi(n) \cdot T^{\mathsf{G}_{[y]}, {^n\theta}}_{\chi} \circ \Ad(n)_{*} \circ [\epsilon \cdot \Psi^{\mathsf{G}_{[y]}, \theta \otimes \chi}]_{n^{-1}\mathsf{U}n, \mathsf{U}} \circ \kappa_{(\mathsf{S}, \theta \otimes \chi)}^{\mathsf{G}_{[y]}}(n)^{-1} \\
&= \chi(n) \cdot T^{\mathsf{G}_{[y]}, {^n\theta}}_{\chi} \circ R_{\mathsf{U}}^{\mathsf{G}_{[y]}, \epsilon, \theta \otimes \chi}(n).
\qedhere
\end{align*}
\end{proof}

Using Lemma~\ref{lemma:RTvsTR}, we can compute $\left(
T^{\mathsf{G}_{[y]}}_{\chi_1} \circ R_{\mathsf{U}}^{\mathsf{G}_{[y]}, \epsilon}(n_{\chi_1})
\right) \circ \left(
T^{\mathsf{G}_{[y]}}_{\chi_2} \circ R_{\mathsf{U}}^{\mathsf{G}_{[y]}, \epsilon}(n_{\chi_2}) 
\right)$ for $\chi_1, \chi_2 \in X_{\unr}(G)_{\pi}$ as follows.

\begin{proposition}
\label{prop:compositionofTcircR}
For $\chi_1, \chi_2 \in X_{\unr}(G)_{\pi}$, we have 
\[
\left(
T^{\mathsf{G}_{[y]}}_{\chi_1} \circ R_{\mathsf{U}}^{\mathsf{G}_{[y]}, \epsilon}(n_{\chi_1})
\right) \circ \left(
T^{\mathsf{G}_{[y]}}_{\chi_2} \circ R_{\mathsf{U}}^{\mathsf{G}_{[y]}, \epsilon}(n_{\chi_2}) 
\right) = \chi_{2}(n_{\chi_1}) \cdot T^{\mathsf{G}_{[y]}}_{\chi_1 \chi_2} \circ R_{\mathsf{U}}^{\mathsf{G}_{[y]}, \epsilon}(n_{\chi_1} n_{\chi_2}).
\]
\end{proposition}
\begin{proof}
According to Lemma~\ref{lemma:RTvsTR}, we have 
\[
\left(
T^{\mathsf{G}_{[y]}}_{\chi_1} \circ R_{\mathsf{U}}^{\mathsf{G}_{[y]}, \epsilon}(n_{\chi_1})
\right) \circ \left(
T^{\mathsf{G}_{[y]}}_{\chi_2} \circ R_{\mathsf{U}}^{\mathsf{G}_{[y]}, \epsilon}(n_{\chi_2}) 
\right) = \chi_{2}(n_{\chi_1}) \cdot \left(
T^{\mathsf{G}_{[y]}}_{\chi_1} \circ T^{\mathsf{G}_{[y]}, \theta \otimes \chi_1}_{\chi_2} 
\right) \circ \left(
R_{\mathsf{U}}^{\mathsf{G}_{[y]}, \epsilon, \theta \otimes \chi_2}(n_{\chi_1}) \circ R_{\mathsf{U}}^{\mathsf{G}_{[y]}, \epsilon}(n_{\chi_2}) 
\right).
\]
Thus, the claim follows from Lemma~\ref{lemma:normalizedR} and Lemma~\ref{lemma:multipricativityofT}.
\end{proof}

We now obtain the following description of $I^{\kappa_{(S, \theta, \rho)}^{\epsilon}}_{\chi_1} \circ I^{\kappa_{(S, \theta, \rho)}^{\epsilon}}_{\chi_2}$ for $\chi_1, \chi_2 \in X_{\unr}(G)_{\pi}$.

\begin{corollary}
\label{corollary:2-cocycleforkappa}
For $\chi_1, \chi_2 \in X_{\unr}(G)_{\pi}$, we have 
\[
I^{\kappa_{(S, \theta, \rho)}^{\epsilon}}_{\chi_1} \circ I^{\kappa_{(S, \theta, \rho)}^{\epsilon}}_{\chi_2} = \chi_2(n_{\chi_1}) \mu_{\rho}\left((n_{\chi_1}, \chi_1), (n_{\chi_2}, \chi_2)\right) \cdot I^{\kappa_{(S, \theta, \rho)}^{\epsilon}}_{\chi_1 \chi_2}. 
\]
\end{corollary}
\begin{proof}
By Lemma~\ref{lemma:independenceonlift}, we may assume that $n_{\chi_1 \chi_2} = n_{\chi_1} n_{\chi_2}$.
Then the claim follows from Proposition~\ref{prop:compositionofTcircR} and the definitions of $I^{\kappa_{(S, \theta, \rho)}^{\epsilon}}_{\chi}$ and $\mu_{\rho}$.
\end{proof}

We define the isomorphism $I_{\chi} \colon \pi^{\epsilon}_{(S, \theta, \rho)}\xrightarrow{\sim}\pi^{\epsilon}_{(S, \theta, \rho)} \otimes \chi$ as the composition of the isomorphism 
\[
\ind_{G(F)_{[y]}}^{G(F)}\big(I^{\kappa_{(S, \theta, \rho)}^{\epsilon}}_{\chi}\big) \colon 
\pi^{\epsilon}_{(S, \theta, \rho)} = \ind_{G(F)_{[y]}}^{G(F)}(\kappa_{(S, \theta, \rho)}^{\epsilon}) \isoarrow \ind_{G(F)_{[y]}}^{G(F)}(\kappa_{(S, \theta, \rho)}^{\epsilon} \otimes \chi)
\]
and the natural isomorphism 
\[
\ind_{G(F)_{[y]}}^{G(F)}(\kappa_{(S, \theta, \rho)}^{\epsilon} \otimes \chi) \isoarrow \ind_{G(F)_{[y]}}^{G(F)}(\kappa_{(S, \theta, \rho)}^{\epsilon}) \otimes \chi, \qquad f \mapsto [g \mapsto \chi(g)^{-1} \cdot f(g)].
\]

\begin{corollary}
\label{corollary:2-cocycleforpiSthetarho}
For $\chi_1, \chi_2 \in X_{\unr}(G)_{\pi}$, we have 
\[
I_{\chi_1} \circ I_{\chi_2} = \chi_2(n_{\chi_1}) \mu_{\rho}\left((n_{\chi_1}, \chi_1), (n_{\chi_2}, \chi_2)\right) \cdot I_{\chi_1 \chi_2}. 
\]
\end{corollary}
\begin{proof}
According to Corollary~\ref{corollary:2-cocycleforkappa},
for $f \in \ind_{G(F)_{[y]}}^{G(F)}(\kappa_{(S, \theta, \rho)}^{\epsilon})$ and $g \in G(F)$, we have
\begin{align*}
\left(
(I_{\chi_1} \circ I_{\chi_2} )(f)
\right)(g) &= \chi_1(g)^{-1} \cdot I^{\kappa_{(S, \theta, \rho)}^{\epsilon}}_{\chi_1}\left(
(I_{\chi_2}(f))(g)
\right) \\
&= \chi_1(g)^{-1} \chi_2(g)^{-1} \cdot \left(
I^{\kappa_{(S, \theta, \rho)}^{\epsilon}}_{\chi_1} \circ I^{\kappa_{(S, \theta, \rho)}^{\epsilon}}_{\chi_2}
\right)(f(g)) \\
&= (\chi_{1} \chi_{2})(g)^{-1} \chi_2(n_{\chi_1}) \mu_{\rho}\left((n_{\chi_1}, \chi_1), (n_{\chi_2}, \chi_2)\right) \cdot I^{\kappa_{(S, \theta, \rho)}^{\epsilon}}_{\chi_1 \chi_2}(f(g)) \\
&=  \chi_2(n_{\chi_1}) \mu_{\rho}\left((n_{\chi_1}, \chi_1), (n_{\chi_2}, \chi_2)\right) \cdot I_{\chi_1 \chi_2}(f)(g).
\end{align*}
Thus, we obtain the claim.
\end{proof}

At the end of this subsection, we describe the cocycle $(\chi_1, \chi_2) \mapsto \chi_2(n_{\chi_1}) \mu_{\rho}\left((n_{\chi_1}, \chi_1), (n_{\chi_2}, \chi_2)\right)$ using the cocycle attached to a projective representation.
We define a group structure on the set 
\[
\left\{
(n, \chi, c) \in N_{G}(S)(F) \times X_{\unr}(G) \times \mathbb{C}^{\times} \mid {^n\theta} = \theta \otimes \chi
\right\}
\]
as
\[
(n_1, \chi_1, c_1) \cdot (n_2, \chi_2, c_2) = (n_1 n_2, \chi_1 \chi_2, c_1 c_2 \chi_1(n_2)).
\]
We can check that $\left\{
(s, 1, \theta(s)^{-1}) \mid s \in S(F)
\right\}$ is a normal subgroup.
We define 
\[
\widetilde{\Box}_{\mathrm{autm}} = \left\{
(n, \chi, c) \in N_{G}(S)(F) \times X_{\unr}(G) \times \mathbb{C}^{\times} \mid {^n\theta} = \theta \otimes \chi
\right\}/ \left\{
(s, 1, \theta(s)^{-1}) \mid s \in S(F)
\right\}.
\]
Note that the subgroup 
\[
\Box_{\mathrm{autm}} = \left\{
(n, 1, c) \mid n \in N_{G}(S)(F)_{\theta}
\right\}
/
\left\{
(s, 1, \theta(s)^{-1}) \mid s \in S(F)
\right\}
\]
of $\widetilde{\Box}_{\mathrm{autm}}$ is written as $\Box_{2}$ in \cite[Section~4.5]{2019arXiv191203274K}.
The second projection $\widetilde{\Box}_{\mathrm{autm}} \to X_{\unr}(G)$ induces an exact sequence
\[
1 \to \Box_{\mathrm{autm}} \to \widetilde{\Box}_{\mathrm{autm}} \to X_{\unr}(G).
\]
We define a representation $\rho_{\autm}$ of $\Box_{\autm}$ on the representation space of $\rho$ as 
\[
\rho_{\autm}(n, 1, c) = c \cdot \rho(n),
\]
which is well-defined since the restriction of $\rho$ to $S(F)$ is $\theta$-isotypic.
According to \cite[Lemma~C.5]{2019arXiv191203274K}, the correspondence $\rho \mapsto \rho_{\autm}$ is a bijection 
\[
\Irr\left(
N_{G}(S)(F)_{\theta}, \theta
\right) \isoarrow \Irr\left(
\Box_{\autm}, \id_{\C^{\times}}
\right).
\]
Let $\widetilde{\Box}_{\mathrm{autm}, \rho_{\autm}}$ denote the normalizer of $\rho_{\autm}$ in $\widetilde{\Box}_{\mathrm{autm}}$.
\begin{lemma}
\label{lemma:intertwiningrhovsrhoautm}
For $\chi \in X_{\unr}(G)_{\pi}$, we have
\[
I_{n_{\chi}} \in \Hom_{\Box_{\autm}}\left(
{^{[(n_{\chi}, \chi, 1)]}\rho_{\autm}}, \rho_{\autm}
\right) \subset \End_{\mathbb{C}}(\rho_{\autm}) = \End_{\mathbb{C}}(\rho).
\]
Moreover, the image of $\widetilde{\Box}_{\mathrm{autm}, \rho_{\autm}}$ in $X_{\unr}(G)$ via the second projection agrees with $X_{\unr}(G)_{\pi}$.
\end{lemma}
\begin{proof}
For $\chi \in X_{\unr}(G)_{\pi}$, we have
\begin{align*}
I_{n_{\chi}} \circ {^{[(n_{\chi}, \chi, 1)]}\rho_{\autm}}(n, 1, c) &= 
I_{n_{\chi}} \circ \rho_{\autm}(n_{\chi}^{-1} n n_{\chi}, 1, c \chi(n)^{-1}) \\
&= I_{n_{\chi}} \circ (c \chi(n)^{-1} \cdot {^{n_{\chi}}\rho}(n)) \\
&= c \chi(n)^{-1} \cdot (\rho \otimes \chi)(n) \circ I_{n_{\chi}} \\
&= c \cdot \rho(n) \circ I_{n_{\chi}} \\
&= \rho_{\autm}(n, 1, c) \circ I_{n_{\chi}}.
\end{align*}
Hence, we have $I_{n_{\chi}} \in \Hom_{\Box_{\autm}}\left(
{^{[(n_{\chi}, \chi, 1)]}\rho_{\autm}}, \rho_{\autm}
\right)$.
In particular, we have $[(n_{\chi}, \chi, 1)] \in \widetilde{\Box}_{\mathrm{autm}, \rho_{\autm}}$.
Thus, $\chi$ is contained in the image of $\widetilde{\Box}_{\mathrm{autm}, \rho_{\autm}}$.

Conversely, let $(n, \chi, c) \in \widetilde{\Box}_{\mathrm{autm}, \rho_{\autm}}$.
Then by reversing the above argument, we have ${^n\rho} \simeq \rho \otimes \chi$.
Hence, we obtain from \cite[Proposition~3.4.6]{2019arXiv191203274K} that 
\[
\pi^{\epsilon}_{(S, \theta, \rho)} \otimes \chi \simeq \pi^{\epsilon}_{(S, \theta \otimes \chi, \rho \otimes \chi)} \simeq \pi^{\epsilon}_{(S, \theta, \rho)}.
\]
Thus, we conclude that the image $\chi$ of $(n, \chi, c)$ in $X_{\unr}(G)$ is contained in $X_{\unr}(G)_{\pi}$.
\end{proof}

Now, we have a short exact sequence
\[
1 \to \Box_{\autm} \to \widetilde{\Box}_{\autm, \rho_{\autm}} \to X_{\unr}(G)_{\pi} \to 1.
\]
By the construction in Appendix~\ref{sectionA cohomology class attached to a projective representation}, we attach a cohomology class
\[
[\mu_{\rho_{\autm}}] := [\mu_{(\widetilde{\Box}_{\autm, \rho_{\autm}}, \Box_{\autm}, \rho_{\autm})}] \in H^2\left(
 X_{\unr}(G)_{\pi}, \mathbb{C}^{\times}
 \right)
\]
to $\rho_{\autm}$.

\begin{proposition}
\label{prop:describemusasmurhoautm}
We have $[\mu_{\mathfrak{s}}] = [\mu_{\rho_{\autm}}]$.
\end{proposition}
\begin{proof}
For each $\chi \in X_{\unr}(G)_{\pi}$, we define
\[
T_{[(n_{\chi}, \chi, 1)]} =  \chi(n_{\chi})^{-1} \cdot I_{n_{\chi}}.
\]
According to Lemma~\ref{lemma:intertwiningrhovsrhoautm}, we have 
\[
T_{[(n_{\chi}, \chi, 1)]} \in \Hom_{\Box_{\autm}}\left(
{^{[(n_{\chi}, \chi, 1)]}\rho_{\autm}}, \rho_{\autm}
\right).
\]
 For $[(n, 1, c)] \in \Box_{\autm}$, we also define 
\[
T_{[(n_{\chi}, \chi, 1)] \cdot [(n, 1, c)]} = \chi(n_{\chi})^{-1} \cdot I_{n_{\chi}} \circ \rho_{\autm}(n, 1, c).
\]
Then the family 
\[
\mathcal{T} = 
\left\{
T_{g} \in \Hom_{\Box_{\autm}}\left(
{^g\rho_{\autm}}, \rho_{\autm}
\right)
\right\}_{g \in \widetilde{\Box}_{\autm, \rho_{\autm}}}
\]
satisfies Condition~\ref{conditionforcTinappendix} in Appendix~\ref{sectionA cohomology class attached to a projective representation}.
Then we have
\begin{align*}
T_{[(n_{\chi_1}, \chi_1, 1)]} \circ T_{[(n_{\chi_2}, \chi_2, 1)]} 
&= \chi_1(n_{\chi_1})^{-1} \chi_2(n_{\chi_2})^{-1} \cdot I_{n_{\chi_1}} \circ I_{n_{\chi_2}} \\
&= \chi_1(n_{\chi_1})^{-1} \chi_2(n_{\chi_2})^{-1} \mu_{\rho}\left(
(n_{\chi_1}, \chi_1), (n_{\chi_2}, \chi_2)
\right) \cdot I_{n_{\chi_1} n_{\chi_2}}.
\end{align*}
On the other hand, we have 
\begin{align*}
T_{[(n_{\chi_1}, \chi_1, 1)] \cdot [(n_{\chi_2}, \chi_2, 1)]} &= T_{[(n_{\chi_1} n_{\chi_2}, \chi_1 \chi_2, \chi_1(n_{\chi_2}))]} \\
&= T_{[(n_{\chi_1 \chi_2}, \chi_1 \chi_2, 1)] \cdot [(n_{\chi_1 \chi_2}^{-1} n_{\chi_1} n_{\chi_2}, 1, \chi_1(n_{\chi_2}) (\chi_1 \chi_2)(n_{\chi_1 \chi_2}^{-1} n_{\chi_1} n_{\chi_2})^{-1})]} \\
&= (\chi_1 \chi_2)(n_{\chi_1 \chi_2})^{-1} \cdot I_{n_{\chi_1 \chi_2}} \circ \rho_{\autm}(n_{\chi_1 \chi_2}^{-1} n_{\chi_1} n_{\chi_2}, 1, \chi_1(n_{\chi_2}) (\chi_1 \chi_2)(n_{\chi_1 \chi_2}^{-1} n_{\chi_1} n_{\chi_2})^{-1}) \\
&= (\chi_1 \chi_2)(n_{\chi_1 \chi_2})^{-1} \chi_1(n_{\chi_2}) (\chi_1 \chi_2)(n_{\chi_1 \chi_2}^{-1} n_{\chi_1} n_{\chi_2})^{-1} \cdot I_{n_{\chi_1 \chi_2}} \circ \rho(n_{\chi_1 \chi_2}^{-1} n_{\chi_1} n_{\chi_2}) \\
&= \chi_2(n_{\chi_1})^{-1} \chi_1(n_{\chi_1})^{-1} \chi_2(n_{\chi_2})^{-1} \cdot  I_{n_{\chi_1} n_{\chi_2}}.
\end{align*}
Thus, we obtain that 
\[
T_{[(n_{\chi_1}, \chi_1, 1)]} \circ T_{[(n_{\chi_2}, \chi_2, 1)]} = \chi_2(n_{\chi_1}) \mu_{\rho}\left(
(n_{\chi_1}, \chi_1), (n_{\chi_2}, \chi_2)
\right) \cdot T_{[(n_{\chi_1}, \chi_1, 1)] \cdot [(n_{\chi_2}, \chi_2, 1)]}.
\]
Now, the claim follows from Corollary~\ref{corollary:2-cocycleforpiSthetarho} and the definition of $[\mu_{\rho_{\autm}}]$.
\end{proof}

In order to compare $\widetilde{\Box}_{\autm}$ with its counterpart on the spectral side, we also introduce the $F^{\unr}$-version $\widetilde{\Box}_{\autm, F^{\unr}}$.
We define 
\[
\widetilde{\Box}_{\autm, F^{\unr}} = \left\{
(n, \chi, c) \in N_{G}(S)(F) \times X_{\unr}(G_{F^{\unr}}) \times \mathbb{C}^{\times} \mid {^n\theta} = \theta \otimes \chi
\right\}/ \left\{
(s, 1, \theta(s)^{-1}) \mid s \in S(F)
\right\},
\]
endowed with a group structure defined in a manner similar to that of $\widetilde{\Box}_{\autm}$.
Then we have an exact sequence 
\[
1 \to \Box_{\autm} \to \widetilde{\Box}_{\autm, F^{\unr}, \rho_{\autm}} \to X_{\unr}(G_{F^{\unr}})_{\pi} \to 1,
\]
where $\widetilde{\Box}_{\autm, F^{\unr}, \rho_{\autm}}$ denotes the normalizer of $\rho_{\autm}$ in $\widetilde{\Box}_{\autm, F^{\unr}}$, and $X_{\unr}(G_{F^{\unr}})_{\pi}$ denotes the inverse image of $X_{\unr}(G)_{\pi}$ via the restriction map $X_{\unr}(G_{F^{\unr}}) \twoheadrightarrow X_{\unr}(G)$.
By the construction in Appendix~\ref{sectionA cohomology class attached to a projective representation}, we obtain
\[
[\mu_{\rho_{\autm}, F^{\unr}}] := [\mu_{(\widetilde{\Box}_{\autm, F^{\unr}, \rho_{\autm}}, \Box_{\autm}, \rho_{\autm})}] \in H^2\left(
X_{\unr}(G_{F^{\unr}})_{\pi}, \mathbb{C}^{\times}
\right).
\]
We regard $H^2\left(
X_{\unr}(G)_{\pi}, \mathbb{C}^{\times}
\right)$ as a subset of $H^2\left(
X_{\unr}(G_{F^{\unr}})_{\pi}, \mathbb{C}^{\times}
\right)$ via the inflation map
\[
H^2\left(
X_{\unr}(G)_{\pi}, \mathbb{C}^{\times}
\right)
\hookrightarrow H^2\left(
X_{\unr}(G_{F^{\unr}})_{\pi}, \mathbb{C}^{\times}
\right).
\]

\begin{corollary}
\label{cor:describemusasmurhoautmFunrver}
We have $[\mu_{\mathfrak{s}}] = [\mu_{\rho_{\autm}}] = [\mu_{\rho_{\autm}, F^{\unr}}]$.
\end{corollary}
\begin{proof}
By the construction, we obtain that $[\mu_{\rho_{\autm}, F^{\unr}}]$ is the image of $[\mu_{\rho_{\autm}}]$ via the inflation map.
Thus, the claim follows from Proposition~\ref{prop:describemusasmurhoautm}.
\end{proof}

\subsection{Computation of the cohomology class attached to an enhanced $L$-parameter}
\label{subsec:Computation of the cohomology class attached to non-singular Lparameter}

Let $(S, \theta)$ be a tame $F$-non-singular elliptic pair in the sense of \cite[Definition~3.4.1]{2019arXiv191203274K}.
We use the same notation as in the previous subsection.
For each irreducible constituent $\rho_{-1}$ of $\kappa_{(S, \phi_{-1})}$, we have a non-singular supercuspidal representation $\pi_{(S, \theta, \rho_{-1})}$.
We also recall that the irreducible constituents of $\kappa_{(S, \phi_{-1})}$ are indexed by $\Irr\left(
N_{G^0}(S)(F)_{\phi_{-1}}, \phi_{-1}
\right)$, which is in bijection with $\Irr\left(
\Box_{\autm}, \id_{\C^{\times}}
\right)$.
If $\rho_{-1}$ corresponds to an irreducible representation $\rho_{\autm} \in \Irr\left(
\Box_{\autm}, \id_{\C^{\times}}
\right)$, we write $\pi_{(S, \theta, \rho_{-1})} = \pi_{(S, \theta, \rho_{\autm})}$.
In this subsection, we recall the construction of the enhanced parameter $[(\varphi_{\pi}, \rho_{\pi})]$ corresponding to a non-singular supercuspidal representation $\pi = \pi_{(S, \theta, \rho_{\autm})}$ via the local Langlands correspondence given in \cite{2019arXiv191203274K}.
We then describe the associated cohomology class $[\mu_{(\varphi_{\pi}, \rho_{\pi})}]$ in a manner similar to Proposition~\ref{prop:describemusasmurhoautm} and Corollary~\ref{cor:describemusasmurhoautmFunrver}.

First, we explain the construction of the (conjugacy class of) $L$-parameter $\varphi_{\pi} \colon W_{F} \to {^LG}$.
Recall that we have fixed an inner twist $\psi \colon G_{\qs} \to G$.
Let $j \colon S \hookrightarrow G$ denote the inclusion map.
We fix an embedding $j_{0} \colon S \hookrightarrow G_{\qs}$ over $F$ such that $\psi \circ j_{0}$ and $j$ are $G(\overline{F})$-conjugate.
We also assume that $j_0$ satisfies the following property.
Let $G^0_{\qs}$ be the connected reductive subgroup of $G_{\qs}$ with maximal torus $j_0(S)$ and root system
\[
\Phi(G^0_{\qs}, j_0(S)) = \left\{
\alpha \in \Phi(G_{\qs}, j_0(S)) \mid (\theta \circ N_{E/F} \circ \alpha^{\vee})(E_{0+}^{\times})=1
\right\},
\]
where $E$ is the splitting field of $S$, $E_{0+}^{\times}$ is the subgroup of $E^{\times}$ defined as
\[
E_{0+}^{\times} = \left\{
x \in E^{\times} \mid \ord_{E}(x-1) > 0
\right\},
\]
and we identify $\theta$ with the character $j_{0}(s) \mapsto \theta(s)$ of $j_{0}(S)(F)$.
The property of $j_0$ that we require is that $G^0_{\qs}$ is quasi-split and the unique Frobenius-fixed point in $\mathcal{A}^{\red}(G^0_{\qs}, j_0(S), F^{\unr})$ is an absolutely special vertex.
Such an embedding exists by \cite[Lemma~3.2.2, Lemma~3.4.12]{kal19}.
As explained in \cite[Section~4.4]{2019arXiv191203274K}, the choice of $j_0$ is closely related to the choice of a Whittaker datum for $G_{\qs}$.
By the construction in \cite[Section~5.1]{kal19}, the $G_{\qs}(\overline{F})$-conjugacy class of $j_0$ determines a $\widehat{G}$-conjugacy class of embeddings $\widehat{j} \colon \widehat{S} \hookrightarrow \widehat{G}$.
Moreover, we have a set of $\chi$-data $(\chi'')$ for $\Phi(G, S)$ computed in terms of $\theta$ as in \cite[Notation~4.3.4]{FKS}.
In this paper, we do not need the precise definition of $(\chi'')$.
We only need the fact that by using $\chi$-data $(\chi'')$, one can extend an embedding $\widehat{j} \colon \widehat{S} \hookrightarrow \widehat{G}$ to an embedding of $L$-groups ${^Lj} \colon {^LS} \hookrightarrow {^LG}$ as explained in \cite[Section~2.6]{MR909227}.
We now define the $L$-parameter $\varphi_{\pi} \colon W_{F} \to {^LG}$ by
\[
\varphi_{\pi} = {^Lj} \circ \varphi_{\theta},
\] 
where $\varphi_{\theta} \colon W_{F} \to {^LS}$ is the Langlands parameter associated to the character $\theta$ via the local Langlands correspondence for tori.

Next, we explain the construction of the enhancement $\rho_{\pi}$.
By the construction in \cite[Section~5.1]{MR3548533} and its positive-characteristic analogue in \cite[Section~7.1]{MR4617942}, the stably conjugate embeddings $j_0$ and $j$ determine a class $\mathrm{inv}(j_0,j)\in H^1(\mathcal{E},Z\to S)$, which corresponds to a character $\eta \colon \pi_0([\widehat{\overline{S}}]^{+}) \to \mathbb{C}^{\times}$ via the isomorphism in \cite[Corollary~5.4]{MR3548533} and \cite[Corollary~7.11]{MR4617942}.
Here, $\overline{S}=S/Z$, and $[\widehat{\overline{S}}]^+$ denotes the inverse image of $\widehat{S}^{\Gamma_F}$ under the isogeny $\widehat{\overline{S}}\to\widehat{S}$ dual to the quotient map $S\to\overline{S}$.
We also regard $\eta$ as a character of $[\widehat{\overline{S}}]^{+}$ via the inflation.
We identify $\widehat{S}^{\Gamma_{F}}$ with its image in $\widehat{G}$ via $\widehat{j} \colon \widehat{S} \hookrightarrow \widehat{G}$.
Then we have $\widehat{S}^{\Gamma_{F}} \subset S_{\varphi_{\pi}}$ and $[\widehat{\overline{S}}]^{+} \subset S_{\varphi_{\pi}}^{+}$.
Let $S_{\varphi_{\pi}, \eta}^{+}$ denote the stabilizer of $\eta$ in $S_{\varphi_{\pi}}^{+}$.
We define
\begin{align*}
\Box_{\spec} &= 
\left(
S_{\varphi_{\pi}, \eta}^{+} \times \mathbb{C}^{\times}
\right)/\left\{
(s, \eta(s)^{-1}) \mid s \in [\widehat{\overline{S}}]^{+}
\right\} \\
&\simeq 
\left(
\pi_0(S_{\varphi_{\pi}, \eta}^{+}) \times \mathbb{C}^{\times}
\right)/\left\{
(s, \eta(s)^{-1}) \mid s \in \pi_0([\widehat{\overline{S}}]^{+})
\right\}.
\end{align*}
We note that $\Box_{\spec}$ is written as $\Box_{1}$ in \cite[Section~4.5]{2019arXiv191203274K}.
According to \cite[Lemma~C.5]{2019arXiv191203274K}, we have a bijection 
\begin{equation}
\label{eq:thetaSvarphivsidBoxspec}
\Irr\left(
S_{\varphi_{\pi}}^{+}, \eta
\right) \isoarrow \Irr\left(
\Box_{\spec}, \id_{\C^{\times}}
\right).
\end{equation}
Here, $\Irr\left(
\Box_{\spec}, \id_{\C^{\times}}
\right)$ is the set of isomorphism classes of irreducible representations of $\Box_{\spec}$ whose restriction to the image of $\{1\}\times \C^{\times}$ is the identity character.
The bijection is obtained as follows.
As explained in the proof of \cite[Lemma~C.5]{2019arXiv191203274K}, the correspondence
\[
\Irr\left(
S_{\varphi_{\pi}, \eta}^{+}, \eta
\right) \to \Irr\left(
S_{\varphi_{\pi}}^{+}, \eta
\right), \qquad \sigma \mapsto \ind_{S_{\varphi_{\pi}, \eta}^{+}}^{S_{\varphi_{\pi}}^{+}}(\sigma)
\]
is a bijection.
On the other hand, for each $\sigma \in \Irr\left(
S_{\varphi_{\pi}, \eta}^{+}, \eta
\right)$, we can define an action $\sigma'$ of $\Box_{\spec}$ on the representation space of $\sigma$ as 
\[
\sigma'(s, c) = c \cdot \sigma(s),
\]
which gives a bijection
\[
\Irr\left(
S_{\varphi_{\pi}, \eta}^{+}, \eta
\right) \isoarrow \Irr\left(
\Box_{\spec}, \id_{\C^{\times}}
\right).
\]
Combining these two bijections, we obtain \eqref{eq:thetaSvarphivsidBoxspec}.

According to \cite[Proposition~4.5.1]{2019arXiv191203274K}, we have an isomorphism 
\[
\iota \colon \Box_{\autm} \isoarrow \Box_{\spec}.
\] 
We define the enhancement $\rho_{\pi}$ as the representation of $S_{\varphi_{\pi}}^{+}$ corresponding to $\rho_{\spec} := \iota(\rho_{\autm})$ via \eqref{eq:thetaSvarphivsidBoxspec}.
Note that the correspondence $\pi \mapsto \LLC^{\mathfrak{s}}_{\iota}(\pi) := [(\varphi_{\pi}, \rho_{\pi})]$ depends on $\iota$.

We now compute the cohomology class $[\mu_{(\varphi_{\pi}, \rho_{\pi})}]$.
Recall that $[\mu_{(\varphi_{\pi}, \rho_{\pi})}]$ is defined as the cohomology class $[\mu_{(N_{\varphi_{\pi}, \rho_{\pi}}^{+}, S_{\varphi_{\pi}}^{+}, \rho_{\pi})}]$ attached to the exact sequence 
\begin{equation}
\label{eq:exactseqSNZforpi}
1 \to S_{\varphi_{\pi}}^{+} \to N_{\varphi_{\pi}, \rho_{\pi}}^{+} \to Z(\widehat{G})^{I_{F}, \circ}_{[(\varphi_{\pi}, \rho_{\pi})]} \to 1
\end{equation}
and $\rho_{\pi} \in \Irr\left(
S_{\varphi_{\pi}}^{+}, \eta
\right)$.
We rewrite $[\mu_{(\varphi_{\pi}, \rho_{\pi})}]$ in terms of $\Box_{\spec}$ and $\rho_{\spec}$.
We write $\widehat{C} = \Cent_{\widehat{G}}(\varphi_{\pi}(I_{F}))^{\circ}$ and $\widehat{T} = \Cent_{\widehat{G}^0}(\widehat{C})$.
Since $N_{\varphi_{\pi}} \subset \Cent_{\widehat{G}}(\varphi_{\pi}(I_{F}))$, $N_{\varphi_{\pi}}$ normalizes $\widehat{C}$.
Moreover, according to \cite[Lemma~4.1.9]{2019arXiv191203274K}, we have $\widehat{G}^0 = \Cent_{\widehat{G}}(\varphi_{\pi}(P_{F}))$, which is also normalized by $N_{\varphi_{\pi}} \subset \Cent_{\widehat{G}}(\varphi_{\pi}(I_{F}))$.
Thus, we have $N_{\varphi_{\pi}} \subset N_{\widehat{G}}(\widehat{T})$.
According to \cite[Lemma~5.2.2]{kal19} and \cite[Lemma~4.1.3, Proposition~4.1.8]{2019arXiv191203274K}, $\widehat{T}$ is a maximal torus of $\widehat{G}$.
Replacing $\varphi_{\pi}$ by its $\widehat{G}$-conjugate, we assume that $\widehat{j}(\widehat{S}) = \widehat{T}$.
Then $\widehat{j}$ identifies $\widehat{S}^{\Gamma_{F}}$ as a normal subgroup of $N_{\varphi_{\pi}}$.
We also regard $[\widehat{\overline{S}}]^{+}$ as a normal subgroup of $N_{\varphi_{\pi}}^{+}$, and let $N^{+}_{\varphi_{\pi}, \eta}$ denote the normalizer of $\eta$ in $N_{\varphi_{\pi}}^{+}$.
We define a group $\widetilde{\Box}_{\spec, F^{\unr}}$ as 
\[
\widetilde{\Box}_{\spec, F^{\unr}} = \left\{
(n, z, c) \in N_{\varphi_{\pi}, \eta}^{+} \times Z(\widehat{G})^{I_{F}, \circ} \times \mathbb{C}^{\times} \mid {^n\varphi_{\pi}} = z \cdot \varphi_{\pi}
\right\}/\left\{
(s, 1, \eta(s)^{-1}) \mid s \in [\widehat{\overline{S}}]^{+}
\right\}.
\]
We identify $\Box_{\spec}$ with the subgroup 
\[
 \left\{
(n, 1, c) \mid n \in S_{\varphi_{\pi}, \eta}^{+}
\right\}/\left\{
(s, 1, \eta(s)^{-1}) \mid s \in [\widehat{\overline{S}}]^{+}
\right\}
\]
of $\widetilde{\Box}_{\spec, F^{\unr}}$.
The second projection of $\widetilde{\Box}_{\spec, F^{\unr}}$ induces an exact sequence 
\[
1 \to \Box_{\spec} \to \widetilde{\Box}_{\spec, F^{\unr}} \to Z(\widehat{G})^{I_{F}, \circ}.
\]
Let $\widetilde{\Box}_{\spec, F^{\unr}, \rho_{\spec}}$ denote the normalizer of $\rho_{\spec}$ in $\widetilde{\Box}_{\spec, F^{\unr}}$.

\begin{lemma}
The image of the second projection $\widetilde{\Box}_{\spec, F^{\unr}, \rho_{\spec}} \to Z(\widehat{G})^{I_{F}, \circ}$ agrees with $Z(\widehat{G})^{I_{F}, \circ}_{[(\varphi_{\pi}, \rho_{\pi})]}$.
\end{lemma}
\begin{proof}
Let $[(n, z, c)] \in \widetilde{\Box}_{\spec, F^{\unr}, \rho_{\spec}}$.
Then we have ${^n\varphi_{\pi}} = z \cdot \varphi_{\pi}$.
Moreover, by the construction of \eqref{eq:thetaSvarphivsidBoxspec}, we have ${^n\rho_{\pi}} \simeq \rho_{\pi}$.
Hence, we obtain that 
\[
z \cdot [(\varphi_{\pi}, \rho_{\pi})] = [(z \cdot \varphi_{\pi}, \rho_{\pi})] = [({^n\varphi_{\pi}}, \rho_{\pi})] = [({^n\varphi_{\pi}}, {^n\rho_{\pi}})] = [(\varphi_{\pi}, \rho_{\pi})]. 
\]
Conversely, suppose that $z \in Z(\widehat{G})^{I_{F}, \circ}_{[(\varphi_{\pi}, \rho_{\pi})]}$.
Then there exists $n \in N^{+}_{\varphi_{\pi}}$ such that $({^n\varphi_{\pi}}, {^n\rho_{\pi}}) = (z \cdot \varphi_{\pi}, \rho_{\pi})$.
Since $\rho_{\pi} \in \Irr\left(
S_{\varphi_{\pi}}^{+}, \eta
\right)$ and $n$ normalizes $\rho_{\pi}$, there exists $s \in S_{\varphi_{\pi}}^{+}$ such that $ns$ normalizes $\eta$.
Then we obtain that $[(ns, z, 1)] \in \widetilde{\Box}_{\spec, F^{\unr}, \rho_{\spec}}$.
\end{proof}

Now, we have a short exact sequence
\[
1 \to \Box_{\spec} \to \widetilde{\Box}_{\spec, F^{\unr}, \rho_{\spec}} \to Z(\widehat{G})^{I_{F}, \circ}_{[(\varphi_{\pi}, \rho_{\pi})]} \to 1.
\]
By the construction in Appendix~\ref{sectionA cohomology class attached to a projective representation}, we attach a cohomology class
\[
[\mu_{\rho_{\spec}}] := [\mu_{(\widetilde{\Box}_{\spec, F^{\unr}, \rho_{\spec}}, \Box_{\spec}, \rho_{\spec})}] \in H^2\left(
 Z(\widehat{G})^{I_{F}, \circ}_{[(\varphi_{\pi}, \rho_{\pi})]}, \mathbb{C}^{\times}
 \right)
\]
to $\rho_{\spec}$.

\begin{proposition}
\label{prop:describemusasmurhospec}
We have $[\mu_{(\varphi_{\pi}, \rho_{\pi})}] = [\mu_{\rho_{\spec}}]$.
\end{proposition}
\begin{proof}
We write $\rho_{\pi} = \ind_{S_{\varphi_{\pi}, \eta}^{+}}^{S_{\varphi_{\pi}}^{+}}(\rho_{\pi, 0})$ for $\rho_{\pi, 0} \in \Irr\left(
S_{\varphi_{\pi}, \eta}^{+}, \eta
\right)$.
For each $[(n, z, 1)] \in \widetilde{\Box}_{\spec, F^{\unr}, \rho_{\spec}}$, we have 
\[
\Hom_{S_{\varphi_{\pi}, \eta}^{+}}\left(
{^n\rho_{\pi, 0}}, \rho_{\pi, 0}
\right) = \Hom_{\Box_{\spec}}\left(
{^{[(n, z, 1)]}\rho_{\spec}}, \rho_{\spec}
\right)
\]
as subsets of $\End_{\mathbb{C}}(\rho_{\pi, 0}) = \End_{\mathbb{C}}(\rho_{\spec})$.
Moreover, for each non-zero element $T_{n} \in \Hom_{S_{\varphi_{\pi}, \eta}^{+}}\left(
{^n\rho_{\pi, 0}}, \rho_{\pi, 0}
\right)$, we have an $S_{\varphi_{\pi}, \eta}^{+}$-isomorphism
\[
\ind_{S_{\varphi_{\pi}, \eta}^{+}}^{S_{\varphi_{\pi}}^{+}}(T_{n}) \colon {^n(\ind_{S_{\varphi_{\pi}, \eta}^{+}}^{S_{\varphi_{\pi}}^{+}}(\rho_{\pi, 0}))} \isoarrow \ind_{S_{\varphi_{\pi}, \eta}^{+}}^{S_{\varphi_{\pi}}^{+}}(\rho_{\pi, 0}), \qquad [f \mapsto [x \mapsto T_{n}(f(n^{-1}xn))]].
\]
For $[(n_1, z_1, 1)], [(n_2, z_2, 1)] \in \widetilde{\Box}_{\spec, F^{\unr}, \rho_{\spec}}$, we have 
\[
\ind_{S_{\varphi_{\pi}, \eta}^{+}}^{S_{\varphi_{\pi}}^{+}}(T_{n_1}) \circ \ind_{S_{\varphi_{\pi}, \eta}^{+}}^{S_{\varphi_{\pi}}^{+}}(T_{n_2}) = \ind_{S_{\varphi_{\pi}, \eta}^{+}}^{S_{\varphi_{\pi}}^{+}}(T_{n_1} \circ T_{n_2}).
\]
Thus, the claim follows from the definitions of $[\mu_{(\varphi_{\pi}, \rho_{\pi})}]$ and $ [\mu_{\rho_{\spec}}]$.
\end{proof}

At the end of this subsection, we reduce the computation of $[\mu_{(\varphi_{\pi}, \rho_{\pi})}] = [\mu_{\rho_{\spec}}]$ to the case where $G = G^0$ and $\theta = \phi_{-1}$.
As explained in \cite[Section~4.4]{2019arXiv191203274K}, the $L$-parameter $\varphi_{\pi}$ decomposes as $\varphi_{\pi} = {^L{j_{G^0, G}}} \circ \varphi_{\pi^{0, \mathrm{ess}}}$, where 
$
{^Lj_{G^0, G}} \colon {^LG^0} \hookrightarrow {^LG}
$ is the extension of the embedding 
\[
\widehat{G}^0 = \Cent_{\widehat{G}}(\varphi_{\pi}(P_{F})) \subset \widehat{G}
\]
associated to $(\chi'')$, and $\varphi_{\pi^{0, \mathrm{ess}}}$ is the $L$-parameter associated to the non-singular representation 
$\pi^{0, \mathrm{ess}} = \pi^{0, \mathrm{ess}}_{(S, \theta, \rho_{-1})}$
of $G^0(F)$ defined in Section~\ref{subsec:Computation of the cohomology class attached to a non-singular supercuspidal representation}.
Recall that $\pi^{0, \mathrm{ess}}$ is defined as 
\[
\pi^{0, \mathrm{ess}} = \pi^{0} \otimes \prod_{i=0}^{d} \phi_{i} \restriction_{G^0(F)}, 
\]
where $\pi^0 = \ind_{G^0(F)_{[y]}}^{G^0(F)}(\rho_{-1})$.
We denote by $(\varphi_{\pi^0}, \rho_{\pi^0})$ the enhanced $L$-parameter associated to $\pi^0$.
Then we have
\[
\varphi_{\pi^{0, \mathrm{ess}}}(w) = z(w) \cdot \varphi_{\pi^0}(w)
\]
for $w \in W_{F}$, where $z \in Z^1(W_{F}, Z(\widehat{G}^0))$ is a lift of the class $[z] \in H^1(W_{F}, Z(\widehat{G}^0))$ associated to the character $\prod_{i=0}^{d} \phi_{i} \restriction_{G^0(F)}$ via the local Langlands correspondence for characters (\cite[\S10.2]{MR546608}).
Since $Z(\widehat{G}) \subset Z(\widehat{G}^0)$ and $N_{\varphi_{\pi}} \subset \Cent_{\widehat{G}}(\varphi_{\pi}(I_{F})) \subset \widehat{G}^0$, we have $S_{\varphi_{\pi^0}} = S_{\varphi_{\pi^{0, \mathrm{ess}}}} = S_{\varphi_{\pi}}$ and 
\[
N_{\varphi_{\pi^0}} = \left\{
g \in \widehat{G}^0 \mid {^g\varphi_{\pi^0}} = z \cdot \varphi_{\pi^0} \quad \text{for some $z \in Z(\widehat{G}^0)^{I_{F}, \circ}$}
\right\} \supset N_{\varphi_{\pi}}.
\]
Under the identification $S_{\varphi_{\pi^0}} = S_{\varphi_{\pi}}$, we have $\rho_{\pi^0} = \rho_{\pi} \in \Irr\left(
S_{\varphi_{\pi^0}}^{+}, \eta
\right)$, and the exact sequence \eqref{eq:exactseqSNZforpi} is the pullback of the exact sequence
\[
1 \to S_{\varphi_{\pi^0}}^{+} \to N_{\varphi_{\pi^0}, \rho_{\pi^0}}^{+} \to Z(\widehat{G}^0)^{I_{F}, \circ}_{[(\varphi_{\pi^0}, \rho_{\pi^0})]} \to 1
\]
attached to $(\varphi_{\pi^0}, \rho_{\pi^0})$ along the inclusion map $Z(\widehat{G})^{I_{F}, \circ}_{[(\varphi_{\pi}, \rho_{\pi})]} \hookrightarrow Z(\widehat{G}^0)^{I_{F}, \circ}_{[(\varphi_{\pi^0}, \rho_{\pi^0})]}$.
Thus, $[\mu_{(\varphi_{\pi}, \rho_{\pi})}]$ is the image of $[\mu_{(\varphi_{\pi^0}, \rho_{\pi^0})}]$ via the restriction map
\[
H^2\left(
 Z(\widehat{G}^{0})^{I_{F}, \circ}_{[(\varphi_{\pi^0}, \rho_{\pi^0})]}, \mathbb{C}^{\times}
 \right) \to H^2\left(
 Z(\widehat{G})^{I_{F}, \circ}_{[(\varphi_{\pi}, \rho_{\pi})]}, \mathbb{C}^{\times}
 \right).
\]
On the other hand, according to Proposition~\ref{prop:muforspecsideFunrvsF}, we have $[\mu_{(\varphi_{\pi}, \rho_{\pi})}] \in H^2\left(
X_{\unr}(^LG)_{[(\varphi_{\pi}, \rho_{\pi})]}, \mathbb{C}^{\times}
\right)$ and $[\mu_{(\varphi_{\pi^0}, \rho_{\pi^0})}] \in H^2\left(
X_{\unr}(^LG^{0})_{[(\varphi_{\pi^0}, \rho_{\pi^0})]}, \mathbb{C}^{\times}
\right)$.
Thus, we obtain the following.
\begin{proposition}
\label{prop:reductionofcomputationofmutodepthzerospec}
The cohomology class $[\mu_{(\varphi_{\pi}, \rho_{\pi})}] \in H^2\left(
X_{\unr}(^LG)_{[(\varphi_{\pi}, \rho_{\pi})]}, \mathbb{C}^{\times}
\right)$ agrees with the image of $[\mu_{(\varphi_{\pi^0}, \rho_{\pi^0})}] \in H^2\left(
X_{\unr}(^LG^{0})_{[(\varphi_{\pi^0}, \rho_{\pi^0})]}, \mathbb{C}^{\times}
\right)$ via the inflation map
\[
H^2\left(
X_{\unr}(^LG^{0})_{[(\varphi_{\pi^0}, \rho_{\pi^0})]}, \mathbb{C}^{\times}
\right) \hookrightarrow H^2\left(
X_{\unr}(^LG)_{[(\varphi_{\pi}, \rho_{\pi})]}, \mathbb{C}^{\times}
\right).
\]
\end{proposition}

Thus, as on the automorphic side, the computation of $[\mu_{(\varphi_{\pi}, \rho_{\pi})}]$ is reduced to the case where $\theta$ is an $F$-non-singular depth-zero character of $S(F)$ with respect to $G$.
In particular, by Corollary~\ref{corollary:reductionofcomputationofmutodepthzeroautm} and Proposition~\ref{prop:reductionofcomputationofmutodepthzerospec}, Conjecture~\ref{conjecture:comparisonoftwococycles} is valid for $\pi$ if and only if it is valid for $\pi^{0}$.

\subsection{Deformation of $\widetilde{\Box}_{\autm, F^{\unr}}$}
\label{subsec:deformationogtildeBoxautm}

In the remainder of this section, we prove that the groups $\widetilde{\Box}_{\autm, F^{\unr}}$ and $\widetilde{\Box}_{\spec, F^{\unr}}$ are isomorphic following a strategy similar to that of the proof of $\Box_{\autm} \simeq \Box_{\spec}$ in \cite{2019arXiv191203274K}.
In this subsection, we will describe $\widetilde{\Box}_{\autm, F^{\unr}}$ as a deformation of the ``base case'' $\widetilde{\Box}_{\autm, F^{\unr}, \base}$.
In Section~\ref{subsec:deformationofwidetildespec}, we will give a similar deformation of $\widetilde{\Box}_{\spec, F^{\unr}}$ and we will compare these deformations of $\widetilde{\Box}_{\autm, F^{\unr}}$ and $\widetilde{\Box}_{\spec, F^{\unr}}$ in Section~\ref{subsec:comparisonofextensions}.
As in Section~\ref{subsec:Computation of the cohomology class attached to a non-singular supercuspidal representation}, we assume that $\theta$ is a $\resi$-non-singular depth-zero character of $S(F)$ with respect to $G$.
We define the group $\widetilde{\Box}_{\autm, F^{\unr}, \base}$ as
\[
\widetilde{\Box}_{\autm, F^{\unr}, \base} = \left\{
(n, \chi, c) \in N_{G_{\qs}}(S)(F) \times X_{\unr}(G_{\qs, F^{\unr}}) \times \mathbb{C}^{\times} \mid {^n\theta} = \theta \otimes \chi
\right\}/ \left\{
(s, 1, \theta(s)^{-1}) \mid s \in S(F)
\right\},
\]
where we identify $S$ with the maximal torus $j_0(S)$ of $G_{\qs}$.
We also define 
\begin{align*}
W_{\autm, \base} &= \widetilde{\Box}_{\autm, F^{\unr}, \base}/\mathbb{C}^{\times} \\
&\simeq \left\{
(w, \chi) \in N_{G_{\qs}}(S)(F)/S(F) \times X_{\unr}(G_{\qs, F^{\unr}}) \mid {^w\theta} = \theta \otimes \chi
\right\}.
\end{align*}
Then we have an exact sequence 
\begin{equation}
\label{splittingseqinautomside}
1 \to \mathbb{C}^{\times} \to \widetilde{\Box}_{\autm, F^{\unr}, \base} \to W_{\autm, \base} \to 1.
\end{equation}

The following proposition is an extension of \cite[Corollary~4.5.7]{2019arXiv191203274K}.

\begin{proposition}
\label{prop:splittingogBoxautmbase}
We have a section $\spl_{\autm} \colon \widetilde{\Box}_{\autm, F^{\unr}, \base} \to \mathbb{C}^{\times}$ of the natural inclusion $\mathbb{C}^{\times} \to \widetilde{\Box}_{\autm, F^{\unr}, \base}$.
Hence, we have an isomorphism
\[
\widetilde{\Box}_{\autm, F^{\unr}, \base} \simeq W_{\autm, \base} \times \mathbb{C}^{\times}.
\]
\end{proposition}
\begin{proof}
Let $y_{0}$ be the unique Frobenius-fixed point in $\mathcal{A}^{\red}(G_{\qs}, S, F^{\unr})$.
By the assumption on $j_0$, the point $y_{0}$ is an absolutely special vertex.
Hence, according to \cite[Lemma~3.4.10]{kal19}, we have 
\[
N_{G_{\qs}}(S)(F) = \left(
N_{G_{\qs}}(S)(F) \cap G_{\qs}(F)_{y_{0}, 0}
\right)
\cdot S(F).
\]
We denote by $\theta_{0}$ the restriction of $\theta$ to $S(F)_{0}$.
According to \cite[Proposition~2.3.3]{2019arXiv191203274K}, the character $\theta_0$ extends to a character $\widetilde{\theta}_0$ of $N_{G_{\qs}}(S)(F)_{\theta_0} \cap G_{\qs}(F)_{y_{0}, 0}$.
We define a map $\spl_{\autm} \colon \widetilde{\Box}_{\autm, F^{\unr}, \base} \to \mathbb{C}^{\times}$ as follows.
Let $[(n, \chi, c)] \in \widetilde{\Box}_{\autm, F^{\unr}, \base}$.
We write $n = s n_0$ for $s \in S(F)$ and $n_{0} \in N_{G_{\qs}}(S)(F)_{\theta^{1}} \cap G_{\qs}(F)_{y_{0}, 0}$.
We define 
\[
\spl_{\autm}([(n, \chi, c)]) = c \cdot \theta(s) \widetilde{\theta}_0(n_0).
\]
Since the restrictions of $\theta$ and $\widetilde{\theta}_0$ to $S(F)_{0}$ are $\theta_0$, the map $\spl_{\autm}$ does not depend on the expression $n = s n_0$.
We prove that $\spl_{\autm}$ is a character of $\widetilde{\Box}_{\autm, F^{\unr}, \base}$.
Let $[(n, \chi, c)], [(n', \chi', c')] \in \widetilde{\Box}_{\autm, F^{\unr}, \base}$.
We write $n = s n_0$ and $n' = s' n'_0$ as above.
Note that $nn' = s (n_0 s' n_0^{-1}) \cdot (n_{0} n'_0)$ and $\chi(n') = \chi(s')$ since $\chi$ is trivial on $G_{\qs}(F)_{y_{0}, 0}$.
We have 
\begin{align*}
\spl_{\autm}(
[(n, \chi, c)] \cdot [(n', \chi', c')]
) & = \spl_{\autm}(
[(nn', \chi \chi', c c' \chi(n'))]
) \\
&= c c' \chi(n') \theta(s (n_{0} s' n_{0}^{-1})) \widetilde{\theta}_{0}(n_{0} n'_{0}) \\
&= c c' \chi(n') \theta(s) ({}^{n_{0}^{-1}}\theta(s')) \widetilde{\theta}_{0}(n_{0}) \widetilde{\theta}_{0}(n'_{0}) \\
&= c c' \chi(s') \theta(s)(\theta \otimes \chi^{-1})(s') \widetilde{\theta}_{0}(n_{0}) \widetilde{\theta}_{0}(n'_{0}) \\
&= c c'\theta(s)\theta(s') \widetilde{\theta}_{0}(n_{0}) \widetilde{\theta}_{0}(n'_{0}) \\
&= \spl_{\autm}([(n, \chi, c)]) \spl_{\autm}([(n', \chi', c')]).
\end{align*}
Thus, $\spl_{\autm} \colon \widetilde{\Box}_{\autm, F^{\unr}, \base} \to \mathbb{C}^{\times}$ gives a splitting of \eqref{splittingseqinautomside}.
\end{proof}

We now compare $\widetilde{\Box}_{\autm, F^{\unr}, \base}$ and $\widetilde{\Box}_{\autm, F^{\unr}}$.
Replacing the inner twist $\psi \colon G_{\qs} \to G$ by $\psi \circ \Ad(g)$ for some $g \in G_{\qs}(\overline{F})$, we assume that $j = \psi \circ j_{0}$.
Then the 1-cocycle 
\[
\Gamma_{F} \ni \sigma \mapsto \psi^{-1} \sigma(\psi) \in G_\qs/Z(G_\qs)
\]
is contained in $Z^1\left(
\Gamma_{F}, S/Z(G_\qs)
\right)$.
Since $H^1\left(
I_{F}, S/Z(G_\qs)
\right) = \{0\}$ by \cite[Theorem~1.9]{MR180554}, we may further replace $\psi$ by its $S(\overline{F})$-conjugate so that 
\[
[\sigma \mapsto \psi^{-1} \sigma(\psi)] \in Z^1\left(
\Gamma_{F}/I_{F}, (S/Z(G_\qs))(F^{\unr})
\right).
\]
In particular, $\psi$ induces isomorphisms $G_{\qs}(F^{\unr}) \isoarrow G(F^{\unr})$, $N_{G_{\qs}}(S)(F^{\unr}) \isoarrow N_{G}(S)(F^{\unr})$, and $X_{\unr}(G_{\qs, F^{\unr}}) \isoarrow X_{\unr}(G_{F^{\unr}})$.
Since $j = \psi \circ j_{0}$, the induced isomorphism $W_{G_{\qs}}(S)(F^{\unr}) \isoarrow W_{G}(S)(F^{\unr})$ is the identity map, where we regard $W_{G_{\qs}}(S)(F^{\unr})$ and $W_{G}(S)(F^{\unr})$ as subgroups of the automorphism group of $S$.
Moreover, according to \cite[Lemma~3.4.10]{kal19} and \cite[Lemma~E.1]{2019arXiv191203274K}, we have 
\[
N_{G_{\qs}}(S)(F)/S(F) = W_{G_{\qs}}(S)(F)
\qquad
\text{and}
\qquad
N_{G}(S)(F)/S(F) = W_{G}(S)(F)_{\eta},
\]
where $W_{G}(S)(F)_{\eta}$ denotes the stabilizer of $\eta$ in $W_{G}(S)(F)$ with respect to the action of $W_{G}(S)(F)$ on $[\widehat{\overline{S}}]^{+}$ induced by the natural Galois-equivariant action of $W_{G}(S)(F)$ on $S$.
For each $w \in W_{G}(S)(F)_{\eta} \subset W_{G_{\qs}}(S)(F)$, we fix lifts $n_{w} \in N_{G}(S)(F)$ and $n_{\base, w} \in N_{G_{\qs}}(S)(F)$.
We assume that $n_{1} = 1$ and $n_{\base, 1} = 1$.
We define $s_{\autm}(w) := n_{w} \cdot \psi(n_{\base, w})^{-1}$.
For $w_1, w_2 \in W_{G}(S)(F)_{\eta}$, we define $\alpha_{\autm}(w_1, w_2) = s_{\autm}(w_1) (w_1s_{\autm}(w_2) w_{1}^{-1}) s_{\autm}(w_1 w_2)^{-1}$.
\begin{lemma}
\label{lemma:sautminSFunralphaautminSF}
We have $s_{\autm}(w) \in S(F^{\unr})$ and $\alpha_{\autm}(w_1, w_2) \in S(F)$.
Moreover, $\alpha_{\autm}$ is a 2-cocycle $W_{G}(S)(F)_{\eta} \times W_{G}(S)(F)_{\eta} \to S(F)$.
\end{lemma}
\begin{proof}
Since $n_{w}$ and $n_{\base, w}$ are lifts of $w$, we have $s_{\autm}(w) \in S(\overline{F})$.
Moreover, since $n_{w} \in N_{G}(S)(F)$, $n_{\base, w} \in N_{G_{\qs}}(S)(F)$, and $\psi$ is defined over $F^{\unr}$, we obtain that $s_{\autm}(w) \in S(F^{\unr})$.
We prove that $\alpha_{\autm}(w_1, w_2) \in S(F)$.
For $w_1, w_2 \in W_{G}(S)(F)_{\eta}$, we have
\begin{align*}
s_{\autm}(w_1) (w_1s_{\autm}(w_2) w_{1}^{-1}) &= n_{w_1} \cdot \psi(n_{\base, w_1})^{-1} (w_1(n_{w_2} \cdot \psi(n_{\base, w_2})^{-1}) w_{1}^{-1}) \\
&= n_{w_1} \cdot \psi(n_{\base, w_1})^{-1} (\psi(n_{\base, w_1}) (n_{w_2} \cdot \psi(n_{\base, w_2})^{-1}) \psi(n_{\base, w_1})^{-1}) \\
&= n_{w_1} n_{w_2} \psi(n_{\base, w_1} n_{\base, w_2})^{-1} \\
&= (n_{w_1} n_{w_2} n_{w_1 w_2}^{-1}) n_{w_1 w_2} \psi(n_{\base, w_1 w_2})^{-1} \psi(n_{\base, w_1} n_{\base, w_2} n_{\base, w_1 w_2}^{-1})^{-1} \\
&= (n_{w_1} n_{w_2} n_{w_1 w_2}^{-1}) s_{\autm}(w_1 w_2) \psi(n_{\base, w_1} n_{\base, w_2} n_{\base, w_1 w_2}^{-1})^{-1}.
\end{align*}
Since $n_{w_1} n_{w_2} n_{w_1 w_2}^{-1}, n_{\base, w_1} n_{\base, w_2} n_{\base, w_1 w_2}^{-1} \in S(F)$, we obtain that 
\begin{align*}
\alpha_{\autm}(w_1, w_2) &= s_{\autm}(w_1) (w_1s_{\autm}(w_2) w_{1}^{-1}) s_{\autm}(w_1 w_2)^{-1} \\
&= (n_{w_1} n_{w_2} n_{w_1 w_2}^{-1}) \psi(n_{\base, w_1} n_{\base, w_2} n_{\base, w_1 w_2}^{-1})^{-1} \in S(F).
\end{align*}
The claim that $\alpha_{\autm}$ is a 2-cocycle follows from the definition of it.
\end{proof}
We write $\overline{\alpha}_{\autm}(w_1, w_2) = \theta(\alpha_{\autm}(w_1, w_2))$.
We define
\begin{align*}
W_{\autm} &= \widetilde{\Box}_{\autm, F^{\unr}}/\mathbb{C}^{\times} \\
&\simeq \left\{
(w, \chi) \in N_{G}(S)(F)/S(F) \times X_{\unr}(G_{F^{\unr}}) \mid {^w\theta} = \theta \otimes \chi
\right\} \\
&= 
\left\{
(w, \chi) \in W_{G}(S)(F)_{\eta} \times X_{\unr}(G_{F^{\unr}}) \mid {^w\theta} = \theta \otimes \chi
\right\}.
\end{align*}
Then we have an exact sequence 
\begin{equation}
\label{exactseqinautomsidedeformedver}
1 \to \mathbb{C}^{\times} \to \widetilde{\Box}_{\autm, F^{\unr}} \to W_{\autm} \to 1.
\end{equation}
Let $\left(
W_{\autm} \times \mathbb{C}^{\times}
\right)^{s_{\autm}}$ be the group with underlying set $W_{\autm} \times \mathbb{C}^{\times}$ and group law given by
\[
(w_1, \chi_1, c_1) \cdot (w_2, \chi_2, c_2) = (w_1 w_2, \chi_1 \chi_2, c_1 c_2 \overline{\alpha}_{\autm}(w_1, w_2) \chi_1(s_{\autm}(w_2))).
\]
Then we also have an exact sequence 
\begin{equation}
\label{twistedexactseqwautm}
1 \to \mathbb{C}^{\times} \to \left(
W_{\autm} \times \mathbb{C}^{\times}
\right)^{s_{\autm}} \to W_{\autm} \to 1.
\end{equation}

\begin{proposition}
\label{prop:describetildeboxasadeformationautm}
The extensions \eqref{exactseqinautomsidedeformedver} and \eqref{twistedexactseqwautm} are isomorphic.
\end{proposition}
\begin{proof}
Let $\widetilde{\Box}_{\autm, F^{\unr}, \base, \eta}$ denote the inverse image of $W_{G}(S)(F)_{\eta}$ in $\widetilde{\Box}_{\autm, F^{\unr}, \base}$.
Then we have a bijection of sets 
\[
\widetilde{\Box}_{\autm, F^{\unr}, \base, \eta} \isoarrow \widetilde{\Box}_{\autm, F^{\unr}}, \qquad [(n, \chi, c)] \mapsto [(s_{\autm}(\overline{n})\psi(n), \chi \circ \psi^{-1}, c)],
\]
where $\overline{n}$ denotes the image of $n$ in $W_{G}(S)(F)_{\eta}$.
Moreover, according to Proposition~\ref{prop:splittingogBoxautmbase}, we have an isomorphism of groups
\[
\widetilde{\Box}_{\autm, F^{\unr}, \base, \eta} \isoarrow W_{\autm} \times \mathbb{C}^{\times}, \qquad [(n, \chi, c)] \mapsto (\overline{n}, \chi \circ \psi^{-1}, \spl_{\autm}([(n, \chi, c)])).
\]
Combining these, we have a bijection of sets 
\[
I \colon \widetilde{\Box}_{\autm, F^{\unr}} \isoarrow W_{\autm} \times \mathbb{C}^{\times}, \qquad [(n, \chi, c)] \mapsto (\overline{n}, \chi, \spl_{\autm}([(\psi^{-1}(s_{\autm}(\overline{n})^{-1}n), \chi\circ\psi, c)])).
\]
Let $[(n_{1}, \chi_{1}, c_{1})], [(n_{2}, \chi_{2}, c_{2})] \in \widetilde{\Box}_{\autm, F^{\unr}}$.
We write $w_{1} = \overline{n_{1}}$ and $w_{2} = \overline{n_{2}}$.
Then we have
\begin{align*}
I([(n_{1}, \chi_{1}, c_{1})] \cdot [(n_{2}, \chi_{2}, c_{2})]) &= I([(n_{1}n_{2}, \chi_{1}\chi_{2}, c_{1} c_{2} \chi_{1}(n_{2}))]) \\
&= (w_{1} w_{2}, \chi_{1} \chi_{2}, \spl_{\autm}([(\psi^{-1}(s_{\autm}(w_{1} w_{2})^{-1}n_{1} n_{2}), (\chi_{1}\chi_{2})\circ\psi, c_{1} c_{2} \chi_{1}(n_{2}))])).
\end{align*}
Since
\begin{align*}
s_{\autm}(w_{1} w_{2})^{-1} &= (w_{1} s_{\autm}(w_{2})^{-1} w_{1}^{-1})s_{\autm}(w_{1})^{-1}\alpha_{\autm}(w_{1}, w_{2}) \\
&= \alpha_{\autm}(w_{1}, w_{2}) s_{\autm}(w_{1})^{-1} (w_{1} s_{\autm}(w_{2})^{-1} w_{1}^{-1}) \\
&= \alpha_{\autm}(w_{1}, w_{2}) s_{\autm}(w_{1})^{-1} (n_{1} s_{\autm}(w_{2})^{-1} n_{1}^{-1}),
\end{align*}
we have 
\begin{align*}
&\quad [(\psi^{-1}(s_{\autm}(w_{1} w_{2})^{-1}n_{1} n_{2}), (\chi_{1}\chi_{2})\circ\psi, c_{1} c_{2} \chi_{1}(n_{2}))] \\
& = [(\psi^{-1}(s_{\autm}(w_{1})^{-1} (n_{1} s_{\autm}(w_{2})^{-1} n_{1}^{-1})n_{1} n_{2}), (\chi_{1}\chi_{2})\circ\psi, c_{1} c_{2} \chi_{1}(n_{2}) \theta(\alpha_{\autm}(w_{1}, w_{2})))] \\
& = [(\psi^{-1}(s_{\autm}(w_{1})^{-1} n_{1}), (\chi_{1})\circ\psi, c_{1})] \cdot [(\psi^{-1}(s_{\autm}(w_{2})^{-1} n_{2}), (\chi_{2}) \circ\psi, c_{2})] \\
& \qquad\qquad\cdot [(1, 1, \overline{\alpha}_{\autm}(w_{1}, w_{2}) \chi_{1}(s_{\autm}(w_{2})))].
\end{align*}
Thus, we conclude that $I$ is an isomorphism of groups $\widetilde{\Box}_{\autm, F^{\unr}} \isoarrow \left(
W_{\autm} \times \mathbb{C}^{\times}
\right)^{s_{\autm}}$.
\end{proof}

\subsection{Deformation of $\widetilde{\Box}_{\spec, F^{\unr}}$}
\label{subsec:deformationofwidetildespec}

In this subsection, we will describe $\widetilde{\Box}_{\spec, F^{\unr}}$ as a deformation of the ``base case'' $\widetilde{\Box}_{\spec, F^{\unr}, \base}$.
The argument is parallel to that in the previous subsection.
We assume that $\theta$ is an $F$-non-singular depth-zero character of $S(F)$ with respect to $G$.
We define a group homomorphism $\varphi_{\base} \colon W_{F} \to {^LG}$ as $\varphi_{\base}(w) = {^Lj}(1 \rtimes w)$ for $w \in W_{F}$.
We define a subgroup $N_{\widehat{G}}(\widehat{T})^{\varphi_{\base}(W_F)}$ of $N_{\widehat{G}}(\widehat{T})$ as
\[
N_{\widehat{G}}(\widehat{T})^{\varphi_{\base}(W_F)} = \left\{
n \in N_{\widehat{G}}(\widehat{T}) \mid {^n\varphi_{\base}} = \varphi_{\base} 
\right\}.
\]
According to \cite[Lemma~4.5.3]{2019arXiv191203274K}, the natural short exact sequence
\[
1 \to \widehat{S} \xrightarrow{\widehat{j}} N_{\widehat{G}}(\widehat{T}) \to W_{\widehat{G}}(\widehat{T}) \simeq W_{G}(S) \to 1
\]
restricts to a short exact sequence 
\[
1 \to \widehat{S}^{\Gamma_{F}} \to N_{\widehat{G}}(\widehat{T})^{\varphi_{\base}(W_F)}  \to W_{G}(S)(F) \to 1.
\]
We denote by $N_{\widehat{\overline{G}}}(\widehat{\overline{T}})^{+}_{\varphi_{\base}}$ the preimage of $N_{\widehat{G}}(\widehat{T})^{\varphi_{\base}(W_F)}$ in $\widehat{\overline{G}}$.
Then we have a short exact sequence 
\[
1 \to [\widehat{\overline{S}}]^{+} \to N_{\widehat{\overline{G}}}(\widehat{\overline{T}})^{+}_{\varphi_{\base}} \to W_{G}(S)(F) \to 1.
\]
Let $N_{\widehat{\overline{G}}}(\widehat{\overline{T}})^{+}_{\varphi_{\base}, \eta}$ denote the stabilizer of $\eta$ in $N_{\widehat{\overline{G}}}(\widehat{\overline{T}})^{+}_{\varphi_{\base}}$.
Then we have a short exact sequence 
\begin{equation}
\label{eq:splittingseqspecsidebeforepushoutbyeta}
1 \to [\widehat{\overline{S}}]^{+} \to N_{\widehat{\overline{G}}}(\widehat{\overline{T}})^{+}_{\varphi_{\base}, \eta} \to W_{G}(S)(F)_{\eta} \to 1.
\end{equation}
We define the group $\widetilde{\Box}_{\spec, F^{\unr}, \base}$ as
\[
\widetilde{\Box}_{\spec, F^{\unr}, \base} = \left\{
(n, z, c) \in N_{\widehat{\overline{G}}}(\widehat{\overline{T}})^{+}_{\varphi_{\base}, \eta} \times Z(\widehat{G})^{I_{F}, \circ} \times \mathbb{C}^{\times} \mid {^{\overline{n}}\theta} = \theta \otimes \chi_{z}
\right\}/\left\{
(s, 1, \eta(s)^{-1}) \mid s \in [\widehat{\overline{S}}]^{+}
\right\},
\]
where $\overline{n}$ denotes the image of $n$ in $W_{G}(S)(F)_{\eta}$.
We also define 
\begin{align*}
W_{\spec, \base} &= \widetilde{\Box}_{\spec, F^{\unr}, \base}/\mathbb{C}^{\times} \\
&\simeq \left\{
(w, z) \in N_{\widehat{\overline{G}}}(\widehat{\overline{T}})^{+}_{\varphi_{\base}, \eta}/[\widehat{\overline{S}}]^{+} \times Z(\widehat{G})^{I_{F}, \circ} \mid {^w\theta} = \theta \otimes \chi_{z}
\right\} \\
&\simeq \left\{
(w, z) \in W_{G}(S)(F)_{\eta} \times Z(\widehat{G})^{I_{F}, \circ} \mid {^w\theta} = \theta \otimes \chi_{z}
\right\}.
\end{align*}
Then we have an exact sequence 
\begin{equation}
\label{splittingseqinspecside}
1 \to \mathbb{C}^{\times} \to \widetilde{\Box}_{\spec, F^{\unr}, \base} \to W_{\spec, \base} \to 1.
\end{equation}

\begin{proposition}
\label{prop:splittingogBoxspecbase}
We have a section $\spl_{\spec} \colon \widetilde{\Box}_{\spec, F^{\unr}, \base} \to \mathbb{C}^{\times}$ of the natural inclusion $\mathbb{C}^{\times} \to \widetilde{\Box}_{\spec, F^{\unr}, \base}$.
Hence, we have an isomorphism
\[
\widetilde{\Box}_{\spec, F^{\unr}, \base} \simeq W_{\spec, \base} \times \mathbb{C}^{\times}.
\]
\end{proposition}
\begin{proof}
According to \cite[Corollary~4.5.5]{2019arXiv191203274K}, the push-out of \eqref{eq:splittingseqspecsidebeforepushoutbyeta} along $\eta \colon [\widehat{\overline{S}}]^{+} \to \mathbb{C}^{\times}$ is split.
Hence, we have an extension $\widetilde{\eta}$ of $\eta$ to $N_{\widehat{\overline{G}}}(\widehat{\overline{T}})^{+}_{\varphi_{\base}, \eta}$.
Using $\widetilde{\eta}$, we can define the splitting $\spl_{\spec} \colon \widetilde{\Box}_{\spec, F^{\unr}, \base} \to \mathbb{C}^{\times}$ as
\[
\spl_{\spec}([(n, z, c)]) = c \cdot \widetilde{\eta}(n).
\qedhere
\]
\end{proof}

We now compare $\widetilde{\Box}_{\spec, F^{\unr}, \base}$ and $\widetilde{\Box}_{\spec, F^{\unr}}$.
We recall that $N_{\varphi_{\pi}} \subset N_{\widehat{G}}(\widehat{T})$.
The following lemma is a refinement of \cite[Lemma~4.3.3]{2019arXiv191203274K}.
\begin{lemma}
Let $n \in N_{\varphi_{\pi}}$.
Then the image $w$ of $n$ in $W_{\widehat{G}}(\widehat{T}) \simeq W_{G}(S)$ is contained in $W_{G}(S)(F)$.
Moreover, for $w \in W_{G}(S)(F)$ and $z \in Z(\widehat{G})^{I_{F}, \circ}$, the following conditions are equivalent:
\begin{enumerate}
\item
We have ${^w\theta} = \theta \otimes \chi_{z}$.
\item
There exists a lift $n \in N_{\varphi_{\pi}}$ of $w$ such that ${^n\varphi_{\pi}} = z \cdot \varphi_{\pi}$.
\end{enumerate}
\end{lemma}
\begin{proof}
We prove the first claim.
Since the embedding $\widehat{j} \colon \widehat{S} \isoarrow \widehat{T} \subset \widehat{G}$ is $W_{F}$-equivariant with respect to the natural action of $W_F$ on $\widehat{S}$ and the action $w \mapsto \Ad(\varphi_{\pi}(w))$ of $W_{F}$ on $\widehat{G}$, it suffices to show that 
\[
\Ad(\varphi_{\pi}(w) \cdot (n \rtimes 1)) = \Ad((n \rtimes 1) \cdot \varphi_{\pi}(w))
\]
on $\widehat{T}$ for all $w \in W_{F}$.
Since $n \in N_{\varphi_{\pi}}$, we have 
\[
(\varphi_{\pi}(w) \cdot (n \rtimes 1)) \cdot ((n \rtimes 1) \cdot \varphi_{\pi}(w))^{-1} \in Z(\widehat{G})^{I_{F}, \circ},
\]
which acts trivially on $\widehat{T}$.
Thus, we obtain the first claim.

We prove the second claim.
Let $w \in W_{G}(S)(F)$.
According to \cite[(2.6.2)]{MR909227}, there exists a lift $n_0 \in N_{\widehat{G}}(\widehat{T})$ of $w$ such that $\Ad(n_0) \circ {^Lj} = {^Lj} \circ \Ad(w)$ on $^LS$.
Then we have 
\[
{^{n_{0}}\varphi_{\pi}} = {^{n_{0}}({^Lj} \circ \varphi_{\theta})} = {^Lj} \circ (\Ad(w) \circ \varphi_{\theta}).
\]
Suppose that ${^w\theta} = \theta \otimes \chi_{z}$.
Then we have 
$
[\Ad(w) \circ \varphi_{\theta}] = [\varphi_{{\theta \otimes \chi_{z}}}] = [z \cdot \varphi_{\theta}]
$.
Hence, the $L$-parameters $\Ad(w) \circ \varphi_{\theta}$ and $z \cdot \varphi_{\theta}$ are $\widehat{S}$-conjugate.
Thus, the $L$-parameters ${^{n_0}\varphi_{\pi}} = {^Lj} \circ (\Ad(w) \circ \varphi_{\theta})$ and $z \cdot \varphi_{\pi}$ are $\widehat{T}$-conjugate.
Replacing $n_0$ by $n = t n_0$ for some $t \in \widehat{T}$, we obtain a lift $n \in N_{\varphi_{\pi}}$ of $w$ such that ${^n\varphi_{\pi}} = z \cdot \varphi_{\pi}$.
Conversely, suppose that we have a lift $n$ of $w$ such that ${^n\varphi_{\pi}} = z \cdot \varphi_{\pi}$.
Then there exists $t \in \widehat{T}$ such that $n_{0} = tn$.
We write $s = \widehat{j}^{-1}(t) \in \widehat{S}$.
We have
\begin{align*}
{^Lj} \circ (\Ad(w) \circ \varphi_{\theta}) &= {^{n_0}\varphi_{\pi}} \\
&= \Ad(t) \circ ({}^n\varphi_{\pi}) \\
&= \Ad(t) \circ (z \cdot \varphi_{\pi}) \\
&= {^Lj} \circ (\Ad(s) \circ (z \cdot \varphi_{\theta})).
\end{align*}
Hence, the $L$-parameters $\Ad(w) \circ \varphi_{\theta}$ and $z \cdot \varphi_{\theta}$ are $\widehat{S}$-conjugate.
Thus, we conclude that ${^w}\theta = \theta \otimes \chi_{z}$.
\end{proof}

We define 
\[
W_{\spec} = \widetilde{\Box}_{\spec, F^{\unr}}/\mathbb{C}^{\times}
\]
which is isomorphic to 
\begin{align*}
&\left\{
(w, z) \in N_{\varphi_{\pi}, \eta}^{+}/[\widehat{\overline{S}}]^{+} \times Z(\widehat{G})^{I_{F}, \circ} \mid {^w\varphi_{\pi}} = z \cdot \varphi_{\pi}
\right\} \\
&\simeq \left\{
(w, z) \in W_{G}(S)(F)_{\eta} \times Z(\widehat{G})^{I_{F}, \circ} \mid {^w\theta} = \theta \otimes \chi_{z}
\right\} \\
&= W_{\spec, \base}.
\end{align*}

For each $(w, \chi) \in W_{\spec} = W_{\spec, \base}$, we fix lifts $n_{(w, \chi)} \in N_{\varphi_{\pi}, \eta}^{+}$ and $n_{\base, (w, \chi)} \in N_{\widehat{\overline{G}}}(\widehat{\overline{T}})^{+}_{\varphi_{\base}, \eta}$ such that $[(n_{(w, \chi)}, \chi, 1)] \in \widetilde{\Box}_{\spec, F^{\unr}}$ and $[(n_{\base, (w, \chi)}, \chi, 1)] \in \widetilde{\Box}_{\spec, F^{\unr}, \base}$, respectively.
We define $s_{\spec}(w, \chi) := n_{(w, \chi)} \cdot n_{\base, (w, \chi)}^{-1} \in \widehat{\overline{S}}$, where we regard $\widehat{\overline{S}}$ as a subgroup of $\widehat{\overline{G}}$ using the isomorphism $\widehat{j} \colon \widehat{S} \isoarrow \widehat{T}$.
For $(w_1, \chi_1), (w_2, \chi_2) \in W_{\spec}$, we define 
\[
\alpha_{\spec}((w_1, \chi_1), (w_2, \chi_2)) = s_{\spec}(w_1, \chi_1) (w_1s_{\spec}(w_2, \chi_2) w_{1}^{-1}) s_{\spec}(w_1 w_2, \chi_1 \chi_2)^{-1}.
\]

\begin{lemma}
We have $\alpha_{\spec}((w_{1}, \chi_{1}), (w_{2}, \chi_{2})) \in [\widehat{\overline{S}}]^{+}$.
\end{lemma}
\begin{proof}
For $(w_{1}, \chi_{1}), (w_{2}, \chi_{2}) \in W_{\spec}$, we have
\begin{align*}
\quad&s_{\spec}(w_{1}, \chi_{1}) (w_{1}s_{\spec}(w_{2}, \chi_{2}) w_{1}^{-1}) \\
&= n_{(w_{1}, \chi_{1})} \cdot n_{\base, (w_{1}, \chi_{1})}^{-1} (w_{1} ( n_{(w_{2}, \chi_{2})} \cdot n_{\base, (w_{2}, \chi_{2})}^{-1}) w_{1}^{-1}) \\
&= n_{(w_{1}, \chi_{1})} \cdot n_{\base, (w_{1}, \chi_{1})}^{-1} (n_{\base, (w_{1}, \chi_{1})} ( n_{(w_{2}, \chi_{2})} \cdot n_{\base, (w_{2}, \chi_{2})}^{-1}) n_{\base, (w_{1}, \chi_{1})}^{-1}) \\
&= n_{(w_{1}, \chi_{1})} n_{(w_{2}, \chi_{2})} (n_{\base, (w_{1}, \chi_{1})} n_{\base, (w_{2}, \chi_{2})})^{-1} \\
&= (n_{(w_{1}, \chi_{1})} n_{(w_{2}, \chi_{2})} n_{(w_{1} w_{2}, \chi_{1} \chi_{2})}^{-1}) n_{(w_{1} w_{2}, \chi_{1} \chi_{2})} n_{\base, (w_{1} w_{2}, \chi_{1} \chi_{2})}^{-1} (n_{\base, (w_{1}, \chi_{1})} n_{\base, (w_{2}, \chi_{2})} n_{\base, (w_{1} w_{2}, \chi_{1} \chi_{2})}^{-1})^{-1} \\
&= (n_{(w_{1}, \chi_{1})} n_{(w_{2}, \chi_{2})} n_{(w_{1} w_{2}, \chi_{1} \chi_{2})}^{-1}) s_{\spec}(w_{1} w_{2}, \chi_{1} \chi_{2}) (n_{\base, (w_{1}, \chi_{1})} n_{\base, (w_{2}, \chi_{2})} n_{\base, (w_{1} w_{2}, \chi_{1} \chi_{2})}^{-1})^{-1}.
\end{align*}
Since $n_{(w_{1}, \chi_{1})} n_{(w_{2}, \chi_{2})} n_{(w_{1} w_{2}, \chi_{1} \chi_{2})}^{-1}$ commutes with $\varphi_{\pi}(W_{F})$, we have $n_{(w_{1}, \chi_{1})} n_{(w_{2}, \chi_{2})} n_{(w_{1} w_{2}, \chi_{1} \chi_{2})}^{-1} \in [\widehat{\overline{S}}]^{+}$.
Similarly, we have $n_{\base, (w_{1}, \chi_{1})} n_{\base, (w_{2}, \chi_{2})} n_{\base, (w_{1} w_{2}, \chi_{1} \chi_{2})}^{-1} \in [\widehat{\overline{S}}]^{+}$.
Thus, we obtain that 
\begin{align*}
\alpha_{\spec}((w_{1}, \chi_{1}), (w_{2}, \chi_{2})) &= s_{\spec}(w_{1}, \chi_{1}) (w_{1}s_{\spec}(w_{2}, \chi_{2}) w_{1}^{-1}) s_{\spec}(w_{1} w_{2}, \chi_{1} \chi_{2})^{-1} \\
&= (n_{(w_{1}, \chi_{1})} n_{(w_{2}, \chi_{2})} n_{(w_{1} w_{2}, \chi_{1} \chi_{2})}^{-1}) (n_{\base, (w_{1}, \chi_{1})} n_{\base, (w_{2}, \chi_{2})} n_{\base, (w_{1} w_{2}, \chi_{1} \chi_{2})}^{-1})^{-1} \\
&\in [\widehat{\overline{S}}]^{+}.
\qedhere
\end{align*}
\end{proof}
We write $\overline{\alpha}_{\spec}((w_{1}, \chi_{1}), (w_{2}, \chi_{2})) = \eta(\alpha_{\spec}((w_{1}, \chi_{1}), (w_{2}, \chi_{2})))$.
We have an exact sequence 
\begin{equation}
\label{exactseqinspecsidedeformedver}
1 \to \mathbb{C}^{\times} \to \widetilde{\Box}_{\spec, F^{\unr}} \to W_{\spec} \to 1.
\end{equation}
Let $\left(
W_{\spec} \times \mathbb{C}^{\times}
\right)^{s_{\spec}}$ be the group with underlying set $W_{\spec} \times \mathbb{C}^{\times}$ and group law given by
\[
(w_{1}, \chi_{1}, c_{1}) \cdot (w_{2}, \chi_{2}, c_{2}) = (
w_{1} w_{2}, \chi_{1} \chi_{2}, c_{1} c_{2} \overline{\alpha}_{\spec}(
(w_{1}, \chi_{1}), (w_{2}, \chi_{2})
)
).
\]
Then we also have an exact sequence 
\begin{equation}
\label{twistedexactseqwspec}
1 \to \mathbb{C}^{\times} \to \left(
W_{\spec} \times \mathbb{C}^{\times}
\right)^{s_{\spec}} \to W_{\spec} \to 1.
\end{equation}

\begin{proposition}
The extensions \eqref{exactseqinspecsidedeformedver} and \eqref{twistedexactseqwspec} are isomorphic.
\end{proposition}
\begin{proof}
The proposition follows from the same argument as in the proof of Proposition~\ref{prop:describetildeboxasadeformationautm}.
\end{proof}

\subsection{Comparison of extensions}
\label{subsec:comparisonofextensions}

We continue to assume that $\theta$ is an $F$-non-singular depth-zero character of $S(F)$ with respect to $G$.
In this subsection, we will prove that 
\[
\left(
W_{\autm} \times \mathbb{C}^{\times}
\right)^{s_{\autm}} \simeq \left(
W_{\spec} \times \mathbb{C}^{\times}
\right)^{s_{\spec}}.
\]
We note that the isomorphism 
\[
X_{\unr}(G_{F^{\unr}}) \isoarrow Z(\widehat{G})^{I_{F}, \circ}, \qquad \chi \mapsto z_{\chi}
\]
induces an isomorphism $W_{\autm} \isoarrow W_{\spec}$.
Hence, one can identify the underlying sets of $\left(
W_{\autm} \times \mathbb{C}^{\times}
\right)^{s_{\autm}} $ and $\left(
W_{\spec} \times \mathbb{C}^{\times}
\right)^{s_{\spec}}$.
We prove that there exists a bijection between $\left(
W_{\autm} \times \mathbb{C}^{\times}
\right)^{s_{\autm}}$ and $\left(
W_{\spec} \times \mathbb{C}^{\times}
\right)^{s_{\spec}}$ which is compatible with the multiplication laws.
We follow the argument in the proof of \cite[Proposition~8.2.1]{2022arXiv221002519K}.
\begin{remark}
The proof of \cite[Proposition~8.2.1]{2022arXiv221002519K} relies on the generalized Tate--Nakayama duality developed in \cite[Section~5]{2022arXiv221002519K}, where $F$ is assumed to have characteristic zero.
As noted in Remark~\ref{remarkaboutthecharacteristic}, Dillery has announced an extension of this duality to arbitrary non-archimedean local fields.
Assuming this extension, the argument below carries over to positive characteristic after replacing Kaletha's cohomology set $H^1(u\to W, Z\to S)$ with Dillery's cohomology set $H^1(\mathcal{E}, Z\to S)$ and making the corresponding modifications.
Since this extension has not yet appeared in published form, we follow \cite{2022arXiv221002519K} and present the argument in Kaletha's framework.
Thus, the proof as written establishes the result in characteristic zero, while the positive-characteristic case is conditional on Dillery's announced extension.
\end{remark}
Let $\varphi_{\theta, 0} \in Z^1(W_{F}, \widehat{S})$ be the map such that $\varphi_{\theta}(\sigma) = \varphi_{\theta, 0}(\sigma) \rtimes \sigma$ for $\sigma \in W_{F}$.
We fix a lift $x \in Z^1(u \to W, Z \to S)$ of $[x]$.
As explained in \cite[Section~8.2]{2022arXiv221002519K}, for each $(w, \chi) \in W_{\autm}$, the pair $(x^{-1}, s_{\autm}(w))$ is an element of the set $Z^1(u \to W,\; Z \to S \xrightarrow{1-w} S)$ defined in \cite[Section~5.1]{2022arXiv221002519K}.
As in the proof of \cite[Proposition~8.2.1]{2022arXiv221002519K}, for each $(w, \chi) \in W_{\autm}$, we fix an element $(\overline{\lambda}_{w}, \mu_{w})$ in the set $Z_{0}(W_{K/F}, X_{*}(\overline{S}) \xrightarrow{1-w} X_{*}(S))_{0}$ defined in \cite[Section~5.1]{2022arXiv221002519K} such that the image of $(\overline{\lambda}_{w}, \mu_{w})$ in $Z^1(u \to W,\; Z \to S \xrightarrow{1-w} S)$ under the isomorphism \cite[(5.2.3)]{2022arXiv221002519K} agrees with $(x^{-1}, s_{\autm}(w))$.
We give a more precise description of $(\overline{\lambda}_{w}, \mu_{w}) \in Z_{0}(W_{K/F}, X_{*}(\overline{S}) \xrightarrow{1-w} X_{*}(S))_{0}$.
The field $K$ is a suitably large Galois extension of $F$.
In particular, $S$ splits over $K$, and $z_{\chi} \in Z(\widehat{G})^{W_{K}}$ for all $(w, \chi) \in W_{\autm}$.
We define $W_{K/F} := W_{F}/[W_K, W_K]$.
 Then we have an exact sequence 
 \[
 1 \to K^{\times} \simeq W_{K}/[W_{K}, W_{K}] \to W_{K/F} \to W_{F}/W_{K} = \Gamma_{K/F} \to 1,
 \]
 where the isomorphism $K^{\times} \simeq W_{K}/[W_{K}, W_{K}]$ comes from the class field theory normalized so that $\varpi_{K}$ corresponds to $\Frob_{K}$.

The first component $\overline{\lambda}_{w}$ of $(\overline{\lambda}_{w}, \mu_{w})$ is an element of $X_{*}(\overline{S})$ such that $N_{K/F}(\overline{\lambda}_{w}) = 0$ and the second component $\mu_{w}$ is a function $W_{K/F} \to X_{*}(S)$ such that 
\[
\prod_{(\sigma_1, \sigma_2, x) \in \Gamma_{K/F} \times \Gamma_{K/F} \times K^{\times}} \sigma_{1}(\mu_{w}(x s(\sigma_2)))(c(\sigma_1, \sigma_2)^{-1} \sigma_{1}(x)^{-1}) = s_{\autm}(w)
\]
up to coboundary.
Here, $s \colon \Gamma_{K/F} \to W_{K/F}$ is a section chosen as in \cite[Section~4.4]{MR3548533}, and $c(\sigma_{1},\sigma_{2})$ is defined as
\[
c(\sigma_1, \sigma_2) = s(\sigma_1) s(\sigma_2) s(\sigma_1 \sigma_2)^{-1} \in W_{K}/[W_{K}, W_{K}] \simeq K^{\times}.
\]
Moreover, the pair $(\overline{\lambda}_{w}, \mu_{w})$ satisfies
\[
\overline{\lambda}_{w} - w(\overline{\lambda}_{w}) = \sum_{\sigma \in W_{K/F}} (\sigma^{-1}(\mu_{w}(\sigma)) - \mu_{w}(\sigma)).
\]
As explained in the proof of \cite[Proposition~8.2.1]{2022arXiv221002519K}, we may assume that $\overline{\lambda}_{w}$ is independent of $w$.
Let $\overline{\lambda}$ denote this common element of $X_{*}(\overline{S})$.
For $(w, \chi) \in W_{\autm}$, we put
\[
h'(w, \chi) = \langle \overline{\lambda}, s_{\spec}(w, z_{\chi}) \rangle \cdot \prod_{\sigma \in W_{K/F}} \langle \mu_{w^{-1}}(\sigma), \varphi_{\theta, 0}(\sigma)\rangle
\]
and $h(w, \chi) = \overline{\alpha}_{\autm}(w^{-1}, w) \chi(s_{\autm}(w))^{-1} \cdot h'(w, \chi)$,
where the angle brackets denote the canonical pairings $X_{*}(\overline{S}) \times \widehat{\overline{S}} \to \mathbb{C}^{\times}$ and $X_{*}(S) \times \widehat{S} \to \mathbb{C}^{\times}$.
We prove that the differential of the 1-cochain $(w, \chi) \mapsto h(w, \chi)$ coincides with the 2-cocycle 
\[
\left(
(w_1, \chi_1), (w_2, \chi_2)
\right) \mapsto 
\overline{\alpha}_{\autm}(w_1, w_2) \overline{\alpha}_{\spec}((w_1, z_{\chi_1}), (w_2, z_{\chi_2}))^{-1} \chi_1(s_{\autm}(w_2)).
\]
We begin with a technical lemma.
For $\chi \in X_{\unr}(G_{F^{\unr}})$, let $\varphi_{\chi, 0} \in Z^1(W_F, \widehat{S})$ be the 1-cocycle such that $\varphi_{\chi, 0}(\sigma) = 1$ for all $\sigma \in I_{F}$ and $\varphi_{\chi, 0}(\Frob_{F}) = z_{\chi} \in Z(\widehat{G})^{I_F, \circ} \subset \widehat{S}$.
\begin{lemma}
\label{lemma:replacingbyacoboundary}
Let $(w_1, \chi_1), (w_2, \chi_2) \in W_{\autm}$.
Then we have
\begin{align*}
&\quad  \langle \overline{\lambda}, s_{\spec}(w_1, z_{\chi_1}) ({^{w_1}s_{\spec}(w_2, z_{\chi_2})}) s_{\spec}(w_2, z_{\chi_2})^{-1} \rangle \cdot \prod_{\sigma \in W_{K/F}} \langle \mu_{w_1^{-1}}(\sigma), {^{w_2}\varphi_{\theta, 0}}(\sigma)\rangle \\
&= \left(
\langle \overline{\lambda}, s_{\spec}(w_1, z_{\chi_1}) \rangle \cdot \prod_{\sigma \in W_{K/F}} \langle \mu_{w_1^{-1}}(\sigma), \varphi_{\theta, 0}(\sigma)\rangle
\right)
 \cdot \left(
\prod_{\sigma \in W_{K/F}} \langle \mu_{w_1^{-1}}(\sigma), \varphi_{\chi_2, 0}(\sigma)\rangle
\right).
\end{align*}
\end{lemma}
\begin{proof}
We write $s_1 = s_{\spec}(w_1, z_{\chi_1})$ and $s_2 = s_{\spec}(w_2, z_{\chi_2})$.
By the definitions of $\varphi_{\pi}$ and $\varphi_{\base}$, we have $\varphi_{\pi}(\sigma) = \widehat{j}(\varphi_{\theta, 0}(\sigma)) \varphi_{\base}(\sigma)$ for $\sigma \in W_{F}$.
Since ${^{n_{(w_2, z_{\chi_2})}}\varphi_{\pi}} = z_{\chi_2} \cdot \varphi_{\pi}$ and ${^{n_{\base, (w_2, z_{\chi_2})}}\varphi_{\base}} = \varphi_{\base}$, for $\sigma \in W_{F}$, we have
\begin{align*}
\widehat{j}(\varphi_{\theta, 0}(\sigma)) \varphi_{\base}(\sigma) &= \varphi_{\pi}(\sigma) \\
&= (z_{\chi_{2}}^{-1} \cdot {^{n_{(w_2, z_{\chi_2})}}\varphi_{\pi}})(\sigma) \\
&= \varphi_{\chi_2, 0}(\sigma)^{-1} \cdot (\widehat{j}({^{w_2}\varphi_{\theta, 0}}(\sigma)) {{}^{n_{(w_2, z_{\chi_2})}}\varphi_{\base}(\sigma)}) \\
&= \varphi_{\chi_2, 0}(\sigma)^{-1} \cdot (\widehat{j}({^{w_2}\varphi_{\theta, 0}}(\sigma)) {{}^{s_2 n_{\base, (w_2, z_{\chi_2})}}\varphi_{\base}(\sigma)}) \\
&= \varphi_{\chi_2, 0}(\sigma)^{-1} \cdot (\widehat{j}({^{w_2}\varphi_{\theta, 0}}(\sigma)) {{}^{s_2}\varphi_{\base}(\sigma)}).
\end{align*}
Thus, we obtain that 
\[
({^{w_2}\varphi_{\theta, 0}}(\sigma)) \varphi_{\theta, 0}(\sigma)^{-1} = \varphi_{\chi_2, 0}(\sigma) \cdot \sigma(s_2) s_2^{-1}.
\]
Using this, we have
\begin{align*}
& \quad \langle \overline{\lambda}, s_{\spec}(w_1, z_{\chi_1}) ({^{w_1}s_{\spec}(w_2, z_{\chi_2})}) s_{\spec}(w_2, z_{\chi_2})^{-1} \rangle \cdot \prod_{\sigma \in W_{K/F}} \langle \mu_{w_1^{-1}}(\sigma), {^{w_2}\varphi_{\theta, 0}}(\sigma)\rangle \\
&= \langle \overline{\lambda}, s_{1} ({^{w_1}s_2}) s_{2}^{-1} \rangle \cdot \prod_{\sigma \in W_{K/F}} \langle \mu_{w_1^{-1}}(\sigma), \varphi_{\theta, 0}(\sigma) \varphi_{\chi_2, 0}(\sigma) \sigma(s_2) s_2^{-1}\rangle.
\end{align*}
We compute each of $\langle \overline{\lambda}, s_{1}({^{w_1}s_2}) s_{2}^{-1} \rangle$ and $\prod_{\sigma \in W_{K/F}} \langle \mu_{w_1^{-1}}(\sigma), \varphi_{\theta, 0}(\sigma) \varphi_{\chi_2, 0}(\sigma) \sigma(s_2) s_2^{-1}\rangle$.
Since 
\[
\overline{\lambda} - w_{1}^{-1}(\overline{\lambda}) = \sum_{\sigma \in W_{K/F}} (\sigma^{-1}(\mu_{w_{1}^{-1}}(\sigma)) - \mu_{w_{1}^{-1}}(\sigma)),
\]
we have 
\begin{align*}
\langle \overline{\lambda}, s_{1} ({^{w_1}s_2}) s_{2}^{-1} \rangle &= \langle \overline{\lambda}, s_{\spec}(w_{1}, z_{\chi_{1}}) \rangle \cdot \langle w_{1}^{-1}(\overline{\lambda}) - \overline{\lambda} , s_{2} \rangle \\
&= \langle \overline{\lambda}, s_{1} \rangle \cdot  \prod_{\sigma \in W_{K/F}} \langle (\sigma^{-1}(\mu_{w_{1}^{-1}}(\sigma)) - \mu_{w_{1}^{-1}}(\sigma)), s_{2} \rangle^{-1} \\
&= \langle \overline{\lambda}, s_{1} \rangle \cdot  \prod_{\sigma \in W_{K/F}} \langle \mu_{w_{1}^{-1}}(\sigma), \sigma(s_{2}) s_{2}^{-1} \rangle^{-1}.
\end{align*}
On the other hand, we have 
\begin{align*}
&\quad \prod_{\sigma \in W_{K/F}} \langle \mu_{w_1^{-1}}(\sigma), \varphi_{\theta, 0}(\sigma) \varphi_{\chi_2, 0}(\sigma) \sigma(s_2) s_2^{-1}\rangle \\
&= \left(
\prod_{\sigma \in W_{K/F}} \langle \mu_{w_1^{-1}}(\sigma), \varphi_{\theta, 0}(\sigma) \rangle
\right) \cdot 
\left(
\prod_{\sigma \in W_{K/F}} \langle \mu_{w_1^{-1}}(\sigma), \varphi_{\chi_2, 0}(\sigma)\rangle
\right)
\cdot 
\left(
\prod_{\sigma \in W_{K/F}} \langle \mu_{w_{1}^{-1}}(\sigma), \sigma(s_2) s_{2}^{-1} \rangle
\right).
\end{align*}
Combining these, we obtain the claim.
\end{proof}

We now compute the differential of $h(w, \chi)$.

\begin{proposition}
\label{prop:comparisonofcocyckeuptovarphichi0factor}
For $(w_1, \chi_1), (w_2, \chi_2) \in W_{\autm}$, we have 
\begin{align*}
&\quad h(w_1, \chi_1) h(w_2, \chi_2) h(w_1 w_2, \chi_1 \chi_2)^{-1} \cdot \left(
\prod_{\sigma \in W_{K/F}} \langle \mu_{w_1^{-1}}(\sigma), \varphi_{\chi_2, 0}(\sigma)\rangle
\right)\\
&= \overline{\alpha}_{\autm}(w_1, w_2) \overline{\alpha}_{\spec}((w_1, z_{\chi_1}), (w_2, z_{\chi_2}))^{-1} \chi_1(s_{\autm}(w_2)) \chi_2(s_{\autm}(w_{1}^{-1}))^{-1}.
\end{align*}
\end{proposition}
\begin{proof}
We write $z_{1} = z_{\chi_1}$ and $z_{2} = z_{\chi_2}$.
According to Lemma~\ref{lemma:replacingbyacoboundary}, we have
\begin{align*}
&\quad h'(w_1, \chi_1) h'(w_2, \chi_2) h'(w_1 w_2, \chi_1 \chi_2)^{-1} \cdot \left(
\prod_{\sigma \in W_{K/F}} \langle \mu_{w_1^{-1}}(\sigma), \varphi_{\chi_2, 0}(\sigma)\rangle
\right)\\
&=  \langle \overline{\lambda}, s_{\spec}(w_1, z_1) ({^{w_1}s_{\spec}(w_2, z_2)}) s_{\spec}(w_2, z_2)^{-1} \rangle \cdot \prod_{\sigma \in W_{K/F}} \langle \mu_{w_1^{-1}}(\sigma), {^{w_2}\varphi_{\theta, 0}}(\sigma)\rangle \cdot h'(w_2, \chi_2) h'(w_1 w_2, \chi_1 \chi_2)^{-1} \\
&= 
\langle \overline{\lambda}, s_{\spec}(w_1, z_1) ({^{w_1}s_{\spec}(w_2, z_2)}) s_{\spec}(w_1 w_2, z_1 z_2)^{-1} \rangle \cdot \prod_{\sigma \in W_{K/F}} \langle ({^{{w_2}^{-1}}}\mu_{w_1^{-1}} + \mu_{w_{2}^{-1}} - \mu_{(w_1 w_2)^{-1}})(\sigma), \varphi_{\theta, 0}(\sigma)\rangle.
\end{align*}
As explained in the proof of \cite[Proposition~8.2.1]{2022arXiv221002519K}, we obtain from the functoriality of the isomorphism \cite[(5.2.3)]{2022arXiv221002519K} and \cite[Corollary~5.2.6]{2022arXiv221002519K} that 
\begin{align*}
&\quad \langle \overline{\lambda}, s_{\spec}(w_1, z_1) ({^{w_1}s_{\spec}(w_2, z_2)}) s_{\spec}(w_1 w_2, z_1 z_2)^{-1} \rangle \cdot \prod_{\sigma \in W_{K/F}} \langle ({^{{w_2}^{-1}}}\mu_{w_1^{-1}} + \mu_{w_{2}^{-1}} - \mu_{(w_1 w_2)^{-1}})(\sigma), \varphi_{\theta, 0}(\sigma)\rangle \\
&= \overline{\alpha}_{\autm}(w_{2}^{-1}, w_1^{-1})^{-1} \overline{\alpha}_{\spec}((w_1, z_1), (w_2, z_2))^{-1}.
\end{align*}
Thus, we obtain that 
\begin{align*}
&\quad h(w_1, \chi_1) h(w_2, \chi_2) h(w_1 w_2, \chi_1 \chi_2)^{-1} \cdot \left(
\prod_{\sigma \in W_{K/F}} \langle \mu_{w_1^{-1}}(\sigma), \varphi_{\chi_2, 0}(\sigma)\rangle
\right) \\
&= \overline{\alpha}_{\autm}(w_1^{-1}, w_1) \overline{\alpha}_{\autm}(w_2^{-1}, w_2) \overline{\alpha}_{\autm}((w_1 w_2)^{-1}, w_1 w_2)^{-1} \overline{\alpha}_{\autm}(w_{2}^{-1}, w_1^{-1})^{-1} \overline{\alpha}_{\spec}((w_1, z_1), (w_2, z_2))^{-1} \\
&\times \chi_1(s_{\autm}(w_1))^{-1} \chi_2(s_{\autm}(w_2))^{-1} (\chi_1 \chi_2)(s_{\autm}(w_1 w_2)).
\end{align*}
Since $\alpha_{\autm}$ is a 2-cocycle, we have 
\begin{align*}
&\quad \alpha_{\autm}(w_{2}^{-1}, w_1^{-1}) ({^{(w_1 w_2)^{-1}}\alpha_{\autm}(w_1, w_2)}) \alpha_{\autm}((w_1 w_2)^{-1}, w_1 w_2) \\
&= \alpha_{\autm}(w_{2}^{-1}, w_1^{-1}) \alpha_{\autm}((w_1 w_2)^{-1}, w_1) \alpha_{\autm}(w_2^{-1}, w_2) \\
&= ({^{w_2^{-1}}\alpha_{\autm}(w_{1}^{-1}, w_1)}) \alpha_{\autm}(w_2^{-1}, 1) \alpha_{\autm}(w_2^{-1}, w_2) \\
&= ({^{w_2^{-1}}\alpha_{\autm}(w_{1}^{-1}, w_1)})  \alpha_{\autm}(w_2^{-1}, w_2).
\end{align*}
Applying $\theta$ to these, we have
\begin{align*}
&\quad \overline{\alpha}_{\autm}(w_{2}^{-1}, w_1^{-1}) \overline{\alpha}_{\autm}(w_1, w_2) \overline{\alpha}_{\autm}((w_1 w_2)^{-1}, w_1 w_2) (\chi_1 \chi_2)(\alpha_{\autm}(w_1, w_2)) \\
&= \overline{\alpha}_{\autm}(w_1^{-1}, w_1)  \overline{\alpha}_{\autm}(w_2^{-1}, w_2) \chi_2(\alpha_{\autm}(w_{1}^{-1}, w_1)).
\end{align*}
Hence, we obtain that
\begin{align*}
&\quad\overline{\alpha}_{\autm}(w_1^{-1}, w_1) \overline{\alpha}_{\autm}(w_2^{-1}, w_2) \overline{\alpha}_{\autm}((w_1 w_2)^{-1}, w_1 w_2)^{-1} \overline{\alpha}_{\autm}(w_{2}^{-1}, w_1^{-1})^{-1} \\
&= \overline{\alpha}_{\autm}(w_1, w_2) (\chi_1 \chi_2)(\alpha_{\autm}(w_1, w_2))\chi_2(\alpha_{\autm}(w_{1}^{-1}, w_1))^{-1} \\
&= \overline{\alpha}_{\autm}(w_1, w_2) \chi_1(s_{\autm}(w_1) s_{\autm}(w_2) s_{\autm}(w_1 w_2)^{-1}) \chi_2(s_{\autm}(w_2)s_{\autm}(w_1^{-1})^{-1} s_{\autm}(w_1 w_2)^{-1}).
\end{align*}
Substituting this into the above equation, we obtain the claim.
\end{proof}

We now obtain the following description of $\prod_{\sigma \in W_{K/F}} \langle \mu_{w_1^{-1}}(\sigma), \varphi_{\chi_2, 0}(\sigma)\rangle$.

\begin{proposition}
\label{prop:calculationofvarphichi0factor}
For $(w_{1},\chi_{1}),(w_{2},\chi_{2}) \in W_{\autm}$, we have 
\[
\prod_{\sigma \in W_{K/F}} \langle \mu_{w_1^{-1}}(\sigma), \varphi_{\chi_2, 0}(\sigma)\rangle = \chi_2(s_{\autm}(w_{1}^{-1}))^{-1}.
\]
\end{proposition}
\begin{proof}
We write $w = w_{1}^{-1}$ and $\chi = \chi_2$.
In the notation of \cite[Section~5.1]{2022arXiv221002519K},
the cocycle $(\varphi_{\chi, 0},1) \in Z^{1}_{\cts}(W_{F}, \widehat{Z}\to \widehat{\overline{S}}\leftarrow \widehat{S})$ defines a cohomology class
$
[(\varphi_{\chi,0},1)]
\in H^{1}_{\cts}
(W_F,\widehat{Z}\to\widehat{\overline{S}}\leftarrow\widehat{S})
$.
Kaletha \cite[(5.2.1)]{2022arXiv221002519K} defines a pairing
\[
\langle \bullet,\bullet\rangle\colon H^{1}(u\to W,Z\to S\xrightarrow{1 - w}S) \times H^{1}_{\cts}(W_{F},\widehat{Z}\to \widehat{\overline{S}}\leftarrow \widehat{S})\to \C^{\times},
\]
and by the definition of this pairing in \cite[Section~5.2]{2022arXiv221002519K}, the left-hand side of the identity in the proposition is equal to
\[
\langle [(x^{-1},s_{\autm}(w))], [(\varphi_{\chi, 0},1)]\rangle^{-1}.
\]

Consider the following diagram of complexes:
\[
\begin{tikzcd}
Z \arrow[r]\arrow[d] & S\arrow[r,"1 - w"]\arrow[d] & S\arrow[d]\\
1 \arrow[r] & 1\arrow[r] & G/G_{\der}.
\end{tikzcd}
\]
Since $(1 - w)s = sws^{-1}w^{-1}\in G_{\der}$ for every $s\in S$, the diagram is commutative.
Hence, it induces a map $H^{1}(u\to W,Z\to S\to S)\to H^{1}(u\to W, 1\to 1\to G/G_{\der})$.
Let $\overline{s_{\autm}}(w)$ denote the image of $s_{\autm}(w)$ in $G/G_{\der}$.
Then the image of $[(x^{-1},s_{\autm}(w))]$ in $H^{1}(u\to W, 1\to 1\to G/G_{\der})$ is $[(1,\overline{s_{\autm}}(w))]$.
We also note that $[(1,\overline{s_{\autm}}(w))]$ is the image of $\overline{s_{\autm}}(w)\in H^{0}(\Gamma_{F},G/G_{\der})$ under the natural map $H^{0}(\Gamma_F,G/G_{\der}) \to H^{1}(u\to W,1\to 1\to G/G_{\der})$
appearing in \cite[Corollary~5.2.6(1)]{2022arXiv221002519K}.

Dualizing the diagram above gives
\[
\begin{tikzcd}
\widehat{Z} \arrow[r]
  & \widehat{\overline{S}}
  & \widehat{S} \arrow[l] \\
1 \arrow[r] \arrow[u]
  & 1 \arrow[u]
  & Z(\widehat{G})^{\circ} \arrow[l] \arrow[u],
\end{tikzcd}
\]
where $\widehat{Z}$ is the kernel of the isogeny
$\widehat{\overline{S}}\to\widehat{S}$
dual to the natural quotient map $S\to\overline{S}$.
This induces a map $H^{1}_{\cts}(W_{F},1\to 1\leftarrow Z(\widehat{G})^{\circ})\to H^{1}_{\cts}(W_{F},\widehat{Z}\to \widehat{\overline{S}}\leftarrow \widehat{S})$.
Since $\varphi_{\chi,0}$ takes values in $Z(\widehat{G})^{\circ}$,
the cocycle $(\varphi_{\chi,0},1)$ also defines a class in $H^{1}_{\cts}
(W_F,1\to 1\leftarrow Z(\widehat{G})^{\circ})$,
whose image is the class $[(\varphi_{\chi,0},1)]$ considered above.
We denote these two classes by the same symbol.

Functoriality of the pairing now gives
\[
\left\langle
[(x^{-1},s_{\autm}(w))],
[(\varphi_{\chi,0},1)]
\right\rangle^{-1}
=
\left\langle
[(1,\overline{s_{\autm}}(w))],
[(\varphi_{\chi,0},1)]
\right\rangle^{-1}.
\]
Finally, by the compatibility with the Langlands pairing established in
\cite[Corollary~5.2.6(1)]{2022arXiv221002519K}, we obtain
\[
\left\langle
[(1,\overline{s_{\autm}}(w))],
[(\varphi_{\chi,0},1)]
\right\rangle^{-1}
=
\chi\bigl(\overline{s_{\autm}}(w)\bigr)^{-1}
=
\chi\bigl(s_{\autm}(w)\bigr)^{-1}.
\qedhere
\]
\end{proof}

\begin{corollary}
\label{corollary:comparisonofcocyckecompleted}
For $(w_1, \chi_1), (w_2, \chi_2) \in W_{\autm} = W_{\spec}$, we have 
\[
h(w_1, \chi_1) h(w_2, \chi_2) h(w_1 w_2, \chi_1 \chi_2)^{-1} 
= \overline{\alpha}_{\autm}(w_1, w_2) \overline{\alpha}_{\spec}((w_1, z_{\chi_1}), (w_2, z_{\chi_2}))^{-1} \chi_1(s_{\autm}(w_2)).
\]
\end{corollary}
\begin{proof}
The claim follows from Propositions~\ref{prop:comparisonofcocyckeuptovarphichi0factor} and \ref{prop:calculationofvarphichi0factor}.
\end{proof}

\begin{theorem}
\label{mainthm:comparisonoftildeautmandtildespec}
The exact sequences \eqref{twistedexactseqwautm} and \eqref{twistedexactseqwspec} are isomorphic.
Thus, we have isomorphisms
\[
\widetilde{\iota} \colon \widetilde{\Box}_{\autm, F^{\unr}} \isoarrow \left(
W_{\autm} \times \mathbb{C}^{\times}
\right)^{s_{\autm}} \isoarrow \left(
W_{\spec} \times \mathbb{C}^{\times}
\right)^{s_{\spec}} \isoarrow \widetilde{\Box}_{\spec, F^{\unr}}.
\]
\end{theorem}
\begin{proof}
We define a bijection of sets $\widetilde{\iota} \colon \left(
W_{\autm} \times \mathbb{C}^{\times}
\right)^{s_{\autm}} \isoarrow \left(
W_{\spec} \times \mathbb{C}^{\times}
\right)^{s_{\spec}}$ as
\[
\widetilde{\iota}(w, \chi, c) = (w, z_{\chi}, c \cdot h(w, \chi)).
\]
According to Corollary~\ref{corollary:comparisonofcocyckecompleted}, for $(w_{1}, \chi_{1}, c_{1}), (w_{2}, \chi_{2}, c_{2}) \in \left(
W_{\autm} \times \mathbb{C}^{\times}
\right)^{s_{\autm}}$, we have
\begin{align*}
&\quad \widetilde{\iota}(w_{1}, \chi_{1}, c_{1}) \cdot \widetilde{\iota}(w_{2}, \chi_{2}, c_{2}) \\
&= (w_{1}, z_{\chi_{1}}, c_{1} \cdot h(w_{1}, \chi_{1})) \cdot (w_{2}, z_{\chi_{2}}, c_{2} \cdot h(w_{2}, \chi_{2})) \\
&= (w_{1} w_{2}, z_{\chi_{1}} z_{\chi_{2}}, c_{1} c_{2} h(w_{1}, \chi_{1}) h(w_{2}, \chi_{2}) \overline{\alpha}_{\spec}((w_{1}, z_{\chi_{1}}), (w_{2}, z_{\chi_{2}}))) \\
&= (w_{1} w_{2}, z_{\chi_{1}} z_{\chi_{2}}, c_{1} c_{2} h(w_{1} w_{2}, \chi_{1} \chi_{2}) \overline{\alpha}_{\autm}(w_{1}, w_{2})\chi_{1}(s_{\autm}(w_{2}))) \\
&=\widetilde{\iota}\left(
(w_{1}, \chi_{1}, c_{1}) \cdot (w_{2}, \chi_{2}, c_{2})
\right).
\end{align*}
Thus, the map $\widetilde{\iota}$ is an isomorphism of groups.
\end{proof}

\begin{corollary}
\label{corollary:Kaletha'sLLCsatisfiesmainconjecture}
Let $(S, \theta)$ be a (not necessarily depth-zero) tame $F$-non-singular elliptic pair.
Let $\iota \colon \Box_{\autm} \isoarrow \Box_{\spec}$ be the restriction of an isomorphism $\widetilde{\iota} \colon \widetilde{\Box}_{\autm, F^{\unr}} \isoarrow \widetilde{\Box}_{\spec, F^{\unr}}$ in Theorem~\ref{mainthm:comparisonoftildeautmandtildespec}.
Then the bijection
\[
\LLC^{\mathfrak{s}}_{\iota} \colon \mathfrak{s} \isoarrow \mathfrak{s}^{\vee}, \qquad \pi \mapsto [(\varphi_{\pi}, \rho_{\pi})]
\]
satisfies Conjecture~\ref{conjecture:comparisonoftwococycles}.
\end{corollary}
\begin{proof}
According to Corollary~\ref{corollary:reductionofcomputationofmutodepthzeroautm} and Proposition~\ref{prop:reductionofcomputationofmutodepthzerospec}, we may assume that $\theta$ is an $F$-non-singular depth-zero character of $S(F)$ with respect to $G$.
Then the claim follows from Corollary~\ref{cor:describemusasmurhoautmFunrver}, Proposition~\ref{prop:describemusasmurhospec}, and Theorem~\ref{mainthm:comparisonoftildeautmandtildespec}.
\end{proof}

\appendix

\section{A cohomology class attached to a projective representation}\label{sec:A cohomology class attached to a projective representation}
\label{sectionA cohomology class attached to a projective representation}
\counterwithin{equation}{section}

Let $G$ be a group and $H$ be a normal subgroup.
For $g \in G$, we denote by $\overline{g}$ the image of $g$ in $G/H$.
Let $\rho$ be an irreducible representation of $H$ that is normalized by $G$.
We assume that Schur's lemma holds for irreducible representations of $H$.
In particular, for each $g \in G$, we have
\[
\dim_{\mathbb{C}}\left(
\Hom_{H}\left(
{^g\rho}, \rho
\right)
\right) = 1.
\]
In this setting, the representation $\rho$ extends to a projective representation $\widetilde{\rho}$ of $G$, which defines a cohomology class of a 2-cocycle $[\mu_{(G, H, \rho)}] \in H^2(G/H, \mathbb{C}^{\times})$.
We briefly record the construction of $[\mu_{(G, H, \rho)}]$.
We fix a family of non-zero elements 
\[
\mathcal{T} = 
\left\{
T_{g} \in \Hom_{H}\left(
{^g\rho}, \rho
\right)
\right\}_{g \in G}
\]
that satisfies the following condition.
\begin{condition}
\label{conditionforcTinappendix}
For all $g \in G$ and $h \in H$, we have $T_{gh} = T_{g} \circ \rho(h)$.
\end{condition}

Let $g_1, g_2 \in G$.
Then as subspaces of $\End_{\mathbb{C}}(\rho)$, we have
\[
\Hom_{H}\left(
{^{g_2}\rho}, \rho
\right) = \Hom_{H}\left(
{^{g_1 g_2}\rho}, {^{g_1}\rho}
\right).
\]
Hence, we can compose the elements $T_{g_1}$ and $T_{g_2}$, and the composition $T_{g_1} \circ T_{g_2}$ is a non-zero element of $\Hom_{H}\left(
{^{g_1 g_2}\rho}, \rho
\right)$.
We define $\mu_{\mathcal{T}}(g_1, g_2) \in \mathbb{C}^{\times}$ by
\[
T_{g_1} \circ T_{g_2} = \mu_{\mathcal{T}}(g_1, g_2) \cdot T_{g_1 g_2}.
\]
By Condition~\ref{conditionforcTinappendix}, the number $\mu_{\mathcal{T}}(g_1, g_2)$ only depends on the images of $g_1$ and $g_2$ in $G/H$.
We can also check that $\mu_{\mathcal{T}}$ defines a 2-cocycle on $G/H$.
Let $\mathcal{T}' = \left\{T'_{g}\right\}_{g \in G}$ be another family that satisfies Condition~\ref{conditionforcTinappendix}.
Then there exists a family of non-zero complex numbers $\left\{c_{w}\right\}_{w \in G/H}$ such that $T'_{g} = c_{\overline{g}} \cdot T_{g}$ for all $g \in G$.
Then for $w_1, w_2 \in G/H$, we have
\[
\mu_{\mathcal{T}'}(w_1, w_2) = c_{w_1} c_{w_2} c_{w_1 w_2}^{-1} \cdot  \mu_{\mathcal{T}}(w_1, w_2).
\]
Thus, the cohomology class $[\mu_{\mathcal{T}}]$ of $\mu_{\mathcal{T}}$ does not depend on the choice of $\mathcal{T}$.
Conversely, for any 2-cocycle $\mu' \in Z^2\left(
G/H, \mathbb{C}^{\times}
\right)$ that is cohomologous to $\mu_{\mathcal{T}}$, there exists a family $\mathcal{T}'$ satisfying Condition~\ref{conditionforcTinappendix} such that $\mu_{\mathcal{T}'} = \mu'$. 
We define the cohomology class $[\mu_{(G, H, \rho)}] \in H^2(G/H, \mathbb{C}^\times)$ as $[\mu_{(G, H, \rho)}] = [\mu_{\mathcal{T}}]$.

\section{An example}\label{sec:An example}
\counterwithin{equation}{subsection}
In this appendix, we calculate the cohomology classes $[\mu_{\mathfrak{s}}]$ and $[\mu_{\mathfrak{s}^{\vee}}]$ for the examples considered in \cite[\S 4.4]{MR3331726}.

\subsection{The setting}\label{subsec: example, setting}

We denote by $q$ the cardinality of $\resi$ and assume that $4 \mid (q -3)$.
Then $q$ is odd and we have $-1\notin (\resi^{\times})^{2}$.
Let $E/F$ be a ramified quadratic extension, $n\ge 4$ be a power of $2$, $D$ be a central division algebra over $E$ of invariant $1/n$, $\mathcal{O}_{D}$ be the ring of integers of $D$, and $\resi[D]$ be the residue field of $D$.
As in the main body of this paper, we write $F_{n}$ (resp.\ $E_{n}$) for the unramified extension of $F$ (resp.\ $E$) of degree $n$.
Let $\varpi_{F}$ (resp.\ $\varpi_{E},\varpi_{F_{n}},\varpi_{E_{n}}$) be a uniformizer of $F$ (resp.\ $E,F_{n},E_{n}$).
The non-trivial Galois action on $E$ is denoted by $x \mapsto \overline{x}$.
We also use the same notation $x \mapsto \overline{x}$ for the non-trivial element of the Galois group $\Gamma_{E_{n}/F_{n}}$ of $E_{n}/F_{n}$.
For an integer $m$, we denote by $\ord_{2}(m)$ the $2$-adic valuation of $m$.
Set $r = \ord_{2}(q^{n} - 1)$.
We remark the following easy fact.
\begin{lemma}\label{lem:or_2(q^a-1)}
Let $a$ be a positive integer.
Then we have
\[
\ord_{2}(q^{a} - 1) =
\begin{cases}
\ord_{2}(a) + \ord_{2}(q + 1) & (2 \mid a), \\
1 & (2 \nmid a).
\end{cases}
\]
In particular, if $a < n$, we have $\ord_{2}(q^{a} - 1) < r$.
\end{lemma}
\begin{proof}
Take odd $b\in \Z$ and $m\in\Z_{\ge 0}$ such that $a = 2^{m}b$.
Then 
\[
q^{a} - 1 = (q^{2^{m}} - 1)(q^{2^{m}(b - 1)} + \cdots + q^{2^{m}} + 1).
\]
Since $(q^{2^{m}(b - 1)} + \cdots + q^{2^{m}} + 1)$ is odd, we have $\ord_{2}(q^{a} - 1) = \ord_{2}(q^{2^{m}} - 1)$.
Hence we may assume that $a = 2^{m}$.
If $m=0$, we have $\ord_{2}(q - 1) = 1$ by the assumption.
Suppose that $m \ge 1$.
Then we have $\ord_{2}(q^{2^{m}} + 1) = 1$.
Hence, we obtain that 
\[
\ord_{2}(q^{2^{m + 1}} - 1) = \ord_{2}(q^{2^{m}} - 1) + \ord_{2}(q^{2^{m}} + 1) = \ord_{2}(q^{2^{m}} - 1) + 1.
\]
Thus, we have 
\[
\ord_{2}(q^{2^{m}} - 1) = (m-1) + \ord_{2}(q^{2} - 1)  = (m-1) + \ord_{2}(q - 1) + \ord_{2}(q + 1) = m + \ord_{2}(q + 1)
\]
by induction on $m$.
In particular, we have
\[
r = \ord_{2}(q^{n} - 1) = \ord_{2}(n) + \ord_{2}(q + 1) > 1.
\]
Since $n$ is a power of $2$, for any $a < n$, we have $\ord_{2}(a) < \ord_{2}(n)$.
Thus, we conclude that $\ord_{2}(q^{a} - 1) < r$.
\end{proof}
Let $G$ be a connected reductive group over $F$ such that 
\[
G(F) = \{(d,x,y)\in D^{\times}\times E^{\times}\times E^{\times}\mid N_{D/E}(d)x^{2}(y/\overline{y}) = 1\},
\]
where $N_{D/E}$ is the reduced norm of $D$.
Let $S$ be the maximal torus of $G$ such that $S(F) = G(F)\cap (E_{n}^{\times}\times E^{\times}\times E^{\times})$.
The maximal compact subgroup $S(F)^{1}$ of $S(F)$ is given by $S(F)^{1} = S(F)\cap (\mathcal{O}_{E_{n}}^{\times}\times \mathcal{O}_{E}^{\times}\times \mathcal{O}_{E}^{\times})$.
Fix a character $\varepsilon \colon \resi[D]^{\times}\to \C^{\times}$ of order $2^r$.
Let $\theta\colon S(F)\to \C^{\times}$ be a character such that $\theta(d,x,y) = \varepsilon(d\pmod{\varpi_{E_{n}}})$ for $d\in \mathcal{O}_{E_{n}}^{\times},x,y\in\mathcal{O}_{E}^{\times}$ such that $(d,x,y)\in G(F)$.

The field $E_{n}$ is a splitting field of $D$, and we have $G_{E_{n}} \subset \GL_{n}^{2}\times \GL_{1}^{2}\times \GL_{1}^{2}$.
We fix this embedding as follows.
Let $\sigma_{n}\in \Gamma_{E_{n}/E}$ be the arithmetic Frobenius element, and we fix a uniformizer $\varpi_{D}\in D$ such that $D = \bigoplus_{i = 0}^{n - 1}\varpi_{D}^{i} E_{n}$, $\varpi_{D}^{n} = \varpi_{E}$, and $\varpi_{D} x\varpi_{D}^{-1} = \sigma_{n}(x)$ for $x\in E_{n}$.
We denote by $\End_{E_{n}}(D)$ the algebra of right $E_{n}$-module endomorphisms of $D$.
Then $E_{n}\otimes_{E}D\ni x\otimes d\mapsto (a\mapsto dax)\in \End_{E_{n}}(D)$ is an isomorphism $E_{n}\otimes_{E}D\simeq\End_{E_{n}}(D)$.
By using the basis $\{\varpi_{D}^{i}\}_{0\le i\le n - 1}$, we have $\End_{E_{n}}(D)\simeq M_{n}(E_{n})$.
Hence, we have an $E_n$-algebra isomorphism $E_{n}\otimes_{E}D \isoarrow M_{n}(E_{n})$.
For $y\in E_{n}$, we have $y\varpi_{D}^{i} = \varpi_{D}^{i}\sigma_{n}^{-i}(y)$.
Hence, the image of $E_{n}\otimes_{E} E_{n}\subset E_{n}\otimes_{E}D\simeq M_{n}(E_{n})$ is the subalgebra of diagonal matrices, and the embedding is given by $E_{n}\otimes_{E} E_{n} \ni x\otimes y\mapsto \diag(xy, x\sigma_{n}^{-1}(y), \ldots , x\sigma_{n}^{-(n-1)}(y)) \in M_{n}(E_{n})$, where $\diag(xy, x\sigma_{n}^{-1}(y), \ldots , x\sigma_{n}^{-(n-1)}(y))$ denotes the diagonal matrix in $M_{n}(E_{n})$ whose $(i, i)$-entry is $x \sigma_{n}^{-(i-1)}(y)$.

Let $\overline{D}$ (resp.\ $\overline{E_{n}}$) be the $E$-algebra obtained from $D$ (resp.\ $E_{n}$) by twisting its $E$-algebra structure by $x\mapsto \overline{x}$ and $\overline{M_{n}(E_{n})}$ be the $E_{n}$-algebra obtained from $M_{n}(E_{n})$ by twisting its $E_{n}$-algebra structure by $x\mapsto \overline{x}$.
Then $x\otimes d\mapsto (a\mapsto da\overline{x})$ induces an isomorphism $E_{n}\otimes_{E}\overline{D}\simeq \End_{E_{n}}(\overline{D})\simeq \overline{M_{n}(E_{n})}$.
Combining this with an $E_{n}$-algebra isomorphism
 \[
 \overline{M_{n}(E_{n})}\ni (a_{ij})\mapsto (\overline{a_{ij}})\in M_{n}(E_{n}),
 \]
 we obtain an isomorphism  $E_{n}\otimes_{E} \overline{D} \isoarrow M_{n}(E_{n})$.
We have an $E$-algebra embedding $\overline{E_{n}}\hookrightarrow \overline{D}$.
By the construction, the image of the subalgebra $E_{n}\otimes_{E}\overline{E_{n}}\subset E_{n}\otimes_{E} \overline{D}$ in $M_{n}(E_{n})$ is the subalgebra of diagonal matrices, and the embedding is given by $E_{n}\otimes_{E} \overline{E_{n}} \ni x\otimes y\mapsto \diag(x\overline{y}, x\sigma_{n}^{-1}(\overline{y}), \ldots , x\sigma_{n}^{-(n-1)}(\overline{y})) \in M_{n}(E_{n})$.

Combining the isomorphisms $E_{n}\otimes_{E}D \isoarrow M_{n}(E_{n})$ and $E_{n}\otimes_{E} \overline{D} \isoarrow M_{n}(E_{n})$ above, we obtain isomorphisms of $E_{n}$-algebras
\[
\begin{tikzcd}
E_{n}\otimes_{F}D\arrow[r,dash,"\sim"] & E_{n}\otimes_{E}E\otimes_{F}D \arrow[r,dash,"\sim"] & (E_{n}\otimes_{E}D) \times (E_{n} \otimes_{E} \overline{D}) \arrow[r,dash,"\sim"]& M_{n}(E_{n})^{2}\\
E_{n}\otimes_{F}E_{n}\arrow[r,dash,"\sim"]\arrow[u,phantom,"\subset" sloped] & E_{n}\otimes_{E}E\otimes_{F}E_{n} \arrow[r,dash,"\sim"]\arrow[u,phantom,"\subset" sloped] & (E_{n}\otimes_{E}E_{n}) \times (E_{n} \otimes_{E} \overline{E_{n}}) \arrow[r,dash,"\sim"]\arrow[u,phantom,"\subset" sloped] & (E_{n}^{n})^{2}.\arrow[u,phantom,"\subset" sloped]
\end{tikzcd}
\]
Here the isomorphism $E_{n}\otimes_{E}\overline{E_{n}}\simeq E_{n}^{n}$ is given by $E_{n}\otimes_{E}\overline{E_{n}} \xrightarrow{x\otimes y\mapsto x\otimes\overline{y}}E_{n}\otimes_{E}E_{n}\simeq E_{n}^{n}$.
Hence, we also have an isomorphism of $E_{n}$-algebras $E_{n}\otimes_{F}(D\times E\times E)\isoarrow (M_{n}(E_{n})\times E_{n}\times E_{n})^{2}$, which gives an embedding $G_{E_{n}} \subset (\GL_{n} \times \GL_{1} \times \GL_{1})^{2}$.
As a subgroup of $G_{E_{n}}$, we have $S_{E_{n}} \subset (\GL_{1}^{n}\times \GL_{1}\times \GL_{1})^{2}$, where we identify $\GL_{1}^{n}$ as the diagonal torus of $\GL_{n}$.
By the construction above, the first component of the embedding 
\[
S(E_{n}) \subset 
(E_{n} \otimes_{F} E_{n})^{\times} \times (E_{n} \otimes_{F} E)^{\times} \times (E_n \otimes_{F} E)^{\times}
 \xrightarrow{\sim} ((E_{n}^{\times})^{n} \times E_{n}^{\times} \times E_{n}^{\times})^{2}
\]
is given by
\[
(E_{n} \otimes_{F} E_{n})^{\times} \xrightarrow{\sim} ((E_{n}^{\times})^{n})^{2}, \qquad x \otimes y \mapsto ((x\sigma_{n}^{-(i-1)}(y))_{i},(x\sigma_{n}^{-(i-1)}(\overline{y}))_{i}).
\]
We also note that the first component of the norm map $N_{E_n/F} \colon S(E_{n}) \to S(F)$ is given by the restriction of
\[
(E_{n}\otimes_{F}E_{n})^{\times} \to E_{n}^{\times}, \qquad x\otimes y\mapsto N_{E_{n}/F}(x)y^{2n}.
\]
The composition of this norm map with the isomorphism $((E_{n}^{\times})^{n})^{2}\simeq (E_{n} \otimes_{F} E_{n})^{\times}$ is
\[
((E_{n}^{\times})^{n})^{2} \to E_{n}^{\times}, \qquad ((a_{i}),(b_{i}))\mapsto \prod_{i=1}^{n}\sigma_{n}^{i-1}(a_{i}\overline{b_{i}}).
\]

We prove that $\theta$ is $F$-non-singular.
For each $k \in \{1, 2\}$ and $1 \le i, j \le n$ with $i \neq j$, we define a homomorphism 
\[
\alpha_{i, j}^{(k)\vee} \colon \GL_{1} \to (\GL_{1}^{n}\times \GL_{1}\times \GL_{1})^{2}
\]
as
\[
\alpha_{i, j}^{(1)\vee}(x) = \left(
((x_{ij}), 1, 1), (1, 1,1)
\right)
\qquad
\text{and}
\qquad
\alpha_{i, j}^{(2)\vee}(x) = \left(
(1, 1,1), ((x_{ij}), 1, 1)
\right),
\]
where $(x_{ij}) \in \GL_{1}^{n}$ is the element with $i$-th component $x$, $j$-th component $x^{-1}$, and all other components equal to $1$.
For $\alpha \in \Phi(G, S)$, the composition of $\alpha^{\vee}$ with the embedding $S_{E_{n}} \subset (\GL_{1}^{n}\times \GL_{1}\times \GL_{1})^{2}$ agrees with $\alpha_{i, j}^{(k)\vee}$ for some $i, j, k$.
We prove that $\theta\circ N_{E_{n}/F}\circ\alpha_{i,j}^{(k)\vee}$ is non-trivial on $\mathcal{O}_{E_{n}}^{\times}$.
We only deal with the case $k = 1$ since the case $k = 2$ is similar.
By the description of the norm map above, for $t \in \mathcal{O}_{E_{n}}^{\times}$, we have
\[
\theta\circ N_{E_{n}/F}\circ\alpha_{i,j}^{(1)\vee}(t) = \varepsilon(\sigma_{n}^{i-1}(t) \sigma_{n}^{j-1}(t)^{-1}) = \varepsilon(t^{q^{i-1} - q^{j-1}}).
\]
According to Lemma~\ref{lem:or_2(q^a-1)}, we have
\[
\ord_{2}(q^{i-1} - q^{j-1}) = \ord_{2}(q^{\abs{i-j}} - 1) < r.
\]
Since the order of $\varepsilon$ is $2^{r}$, we obtain that $\theta\circ N_{E_{n}/F}\circ\alpha_{i,j}^{(k)\vee}$ is non-trivial on $\mathcal{O}_{E_{n}}^{\times}$.
Hence $\theta$ is $F$-non-singular.

\subsection{Automorphic side}

In this subsection, we compute the cohomology class $[\mu_{\mathfrak{s}}]$ for the Bernstein block $\Rep^{\mathfrak{s}}(G(F))$ containing the non-singular supercuspidal representation arising from $\theta$.
We use the same notation as in Sections~\ref{subsection:Azumaya algebra attached to a supercuspidal block} and \ref{subsec:Computation of the cohomology class attached to a non-singular supercuspidal representation}.
Since $\theta$ is a depth-zero character, we have $G^0 = G$.
We also note that $G(F)^{1} = G(F)\cap (\mathcal{O}_{D}^{\times}\times \mathcal{O}_{E}^{\times}\times \mathcal{O}_{E}^{\times})$, which is the unique maximal compact open subgroup $G(F)_{y,0}$ of $G(F)$.
The inclusion $S(F) \subset G(F)$ induces an isomorphism $S(F)^{1}/S(F)_{0+} \isoarrow G(F)^{1}/G(F)_{y, 0+}$.
Hence, the character $\theta^{1} := \theta \restriction_{S(F)^{1}}$ of $S(F)^{1}$ extends to a character of $G(F)^{1}$ that is trivial on $G(F)_{y, 0+}$. We use the same letter $\theta^{1}$ for this extension.
The character $\theta^{1}$ is also the representation $\pi^{1}$ of $G(F)^{1}$ in Section~\ref{subsection:Azumaya algebra attached to a supercuspidal block}.
First, we study the structure of $G(F)/G(F)^{1}$.

\begin{lemma}\label{lem:description of Henniart's group}
Let $a,b,c\in\Z,\ d\in \mathcal{O}_{D}^{\times},\ x,y\in\mathcal{O}_{E}^{\times}$ such that $(\varpi_{D}^{a}d,\varpi_{E}^{b}x,\varpi_{E}^{c}y)\in G(F)$.
Then we have $a + 2b = 0$, and one of the following conditions is satisfied:
\begin{itemize}
\item $c \in 2 \Z$ and $d\in (\resi[D]^{\times})^{2}$;
\item $c \in 2\Z + 1$ and $d\notin (\resi[D]^{\times})^{2}$,
\end{itemize}
where we identify $d$ with its image in $\resi[D]$ by abuse of notation.
Conversely, suppose that $a, b,c \in \Z$ and $d\in \mathcal{O}_{D}^{\times}$ satisfy the above conditions.
Then there exist $x,y\in \mathcal{O}_{E}^{\times}$ such that $(\varpi_{D}^{a}d,\varpi_{E}^{b}x,\varpi_{E}^{c}y)\in G(F)$.
\end{lemma}
\begin{proof}
Since $\varpi_{D}^{n} = \varpi_{E}$, we have $N_{D/E}(\varpi_{D}) = - \varpi_{E}$.
Hence, the element $(\varpi_{D}^{a}d,\varpi_{E}^{b}x,\varpi_{E}^{c}y)$ is contained in $G(F)$ if and only if 
\[
(-1)^{a} \varpi_{E}^{a + 2b}N_{D/E}(d)x^{2}(y/\overline{y})(\varpi_{E}/\overline{\varpi_{E}})^{c} = 1.
\]
Since $E/F$ is ramified and tamely ramified, we have $\overline{\varpi_{E}}/\varpi_{E} = -1$ and $y/\overline{y} = 1$ in $\resi[E]^{\times} = \resi^{\times}$.
Thus, we have $a + 2b = 0$ and $N_{D/E}(d)\in (-1)^{a + c}(\resi[E]^{\times})^{2}$ in $\resi[E]^{\times}$.
In this case, we have $a \in 2 \mathbb{Z}$, hence $(-1)^{a} = 1$.
Thus, we obtain from the assumption $-1 \notin (\resi^{\times})^{2}$ that $N_{D/E}(d)\in (\resi[E]^{\times})^{2}$ (resp.\ $N_{D/E}(d)\notin (\resi[E]^{\times})^{2}$) if $c\in 2\Z$ (resp.\ $c\in 2\Z + 1$).
We fix a generator $x$ of $\resi[D]^{\times}$.
Then $\resi[E]^{\times}$ is generated by $y := x^{(q^{n}-1)/(q-1)}$.
On the other hand, the norm map $N_{D/E}$ sends $x$ to
\[
x^{1 + q + \cdots + q^{n-1}} = x^{\frac{q^{n} - 1}{q-1}} = y.
\]
Hence, we have $N_{D/E}(d)\in (\resi[E]^{\times})^{2}$ if and only if $d\in (\resi[D]^{\times})^{2}$.
Thus, we obtain the first part of the lemma.

Conversely, suppose that $a, b,c \in \Z$ and $d\in \mathcal{O}_{D}^{\times}$ satisfy the conditions in the lemma.
By the above argument, there exists $x' \in \mathcal{O}_{E}^{\times}$ such that 
\[
N_{D/E}(\varpi_{D}^{a}d)(\varpi_{E}^{b}x')^{2}(\varpi_{E}/\overline{\varpi_{E}})^{c} = N_{D/E}(d) (x')^{2} (\varpi_{E}/\overline{\varpi_{E}})^{c} \equiv 1 \pmod{\varpi_{E}}.
\]
Since the map 
\[
1 + \varpi_{E} \mathcal{O}_{E} \to 1 + \varpi_{E} \mathcal{O}_{E}, \qquad z \mapsto z^{2}
\]
is an isomorphism, we can replace $x'$ by an element $x \in x'(1 + \varpi_{E} \mathcal{O}_{E} )$ so that 
\[
N_{D/E}(\varpi_{D}^{a}d)(\varpi_{E}^{b}x)^{2}(\varpi_{E}/\overline{\varpi_{E}})^{c} = 1.
\]
Then we have $(\varpi_{D}^{a}d,\varpi_{E}^{b}x,\varpi_{E}^{c})\in G(F)$.
\end{proof}

\begin{corollary}\label{cor:example, G/G(F)^1}
We have an isomorphism 
\[
G(F)/G(F)^{1} \isoarrow \left\{
(a, b, c) \in \mathbb{Z}^{3} \mid a + 2b = 0
\right\}, \qquad (\varpi_{D}^{a}d,\varpi_{E}^{b}x,\varpi_{E}^{c}y) \mapsto (a, b, c).
\]
\end{corollary}
\begin{proof}
The claim follows from Lemma~\ref{lem:description of Henniart's group}.
\end{proof}

We identify the image of $(\varpi_{D}^{a}d,\varpi_{E}^{b}x,\varpi_{E}^{c}y)$ in $G(F)/G(F)^{1}$ with the element $(a, b, c) \in \mathbb{Z}^{3}$ satisfying $a + 2b = 0$.
We write $N = N_{G(F)}(\theta^{1})= \{g\in G(F)\mid {}^{g}\theta^{1} = \theta^{1}\}$ as in Section~\ref{subsection:Azumaya algebra attached to a supercuspidal block}.
\begin{proposition}\label{prop:example. N/K}
We have $N/G(F)^{1} \simeq \{(a,b, c) \in G(F)/G(F)^{1}\mid a\in (n/2)\Z\}$.
\end{proposition}
\begin{proof}
Let $g = (\varpi_{D}^{a}d,\varpi_{E}^{b}x,\varpi_{E}^{c}y)\in G(F)$ and $k = (z,x',y')\in G(F)^{1}$.
By Lemma~\ref{lem:description of Henniart's group}, we have $z\in (\resi[D]^{\times})^{2}$, where we identify $z$ with its image in $\resi[D]^{\times}$.
On the other hand, for any $z\in (\resi[D]^{\times})^{2}$, we have an element $k = (z,x',y')\in G(F)^{1}$.
Since $\varpi_{D}$ acts on $\resi[D]$ via the Frobenius automorphism, we have $\theta^{1}(gkg^{-1}) = \varepsilon(z^{q^{a}})$.
Hence, $g\in N$ if and only if $\varepsilon(z^{q^{a}}) = \varepsilon(z)$ for all $z\in (\resi[D]^{\times})^{2}$.
Since $\varepsilon$ has order $2^{r}$, we obtain that $g \in N$ if and only if $2(q^{a} - 1)\equiv 0 \pmod{2^{r}}$, equivalently, $2a\in n\Z$ by Lemma~\ref{lem:or_2(q^a-1)}.
\end{proof}

We fix a family of non-zero elements 
\[
\mathcal{T} = 
\left\{
T_{n} \in \Hom_{G(F)^{1}}\left(
{^n\theta^{1}}, \theta^{1}
\right)
\right\}_{n \in N}
\]
that satisfies Condition~\ref{conditionforcT}.
According to Lemma~\ref{lemmacomputationofc}, for $n_{1}, n_{2} \in N$, we have 
\[
c_{\mu_{\mathcal{T}}}(\overline{n}_1, \overline{n}_2) = \theta^{1}(n_{1}^{-1} n_{2}^{-1} n_{1} n_{2}).
\]
Hence, the radical of $c_{\mu_{\mathcal{T}}}$ is the group $T/G(F)^{1}$, where
\[
T = \left\{
n \in N \mid \text{$\theta^{1}(n^{-1}h^{-1}nh) = 1$ for all $h \in N$}
\right\} = \left\{
n \in N \mid \text{$\theta^{1}(nhn^{-1} h^{-1}) = 1$ for all $h \in N$}
\right\}.
\]

\begin{proposition}\label{prop:example. T/K}
We have 
\[
T/G(F)^1\simeq \{(a,b, c)\in G(F)/G(F)^{1}\mid a\in n\Z,\ c\in 2\Z\}.
\]
\end{proposition}
\begin{proof}
Let $g = (\varpi_{D}^{a}d,\varpi_{E}^{b}x,\varpi_{E}^{c}y) \in T$.
Then for any $h = (\varpi_{D}^{a'}d',\varpi_{E}^{b'}x',\varpi_{E}^{c'}y')\in N$, we have $\theta^{1}(g hg^{-1} h^{-1}) = \varepsilon(d^{q^{a}(1 - q^{a'})}(d')^{q^{a'}(q^{a} - 1)}) = 1$.
By considering the case $a' = 0$, we have $\varepsilon((d')^{q^{a} - 1}) = 1$.
Since we can take any $d'$, we have $q^{a} - 1\equiv 0\pmod{2^{r}}$.
Hence, we obtain that $a\equiv 0\pmod{n}$ by Lemma~\ref{lem:or_2(q^a-1)}.
Then we have $1 = \varepsilon(d^{q^{a}(1 - q^{a'})}(d')^{q^{a'}(q^{a} - 1)}) = \varepsilon(d^{1 - q^{a'}})$.
Since $a'$ can be any element in $(n/2)\Z$, we obtain that $d\in (\resi[D]^{\times})^{2}$, equivalently $c\in 2\Z$ by Lemma~\ref{lem:description of Henniart's group}.
Conversely, any such element is contained in $T$.
\end{proof}

We take generators of $N/G(F)^1$ as follows.
Let $\delta\in D^{\times}$ such that $N_{D/E}(\delta) = \overline{\varpi_{E}}/\varpi_{E}$.
Since the image of $N_{D/E}(\delta)$ in $\resi[E]^{\times}$ is $-1$, from the proof of Lemma~\ref{lem:description of Henniart's group}, we have $\delta\notin (\resi[D]^{\times})^{2}$.
Put $g_{1} = (\varpi_{D}^{n/2},\varpi_{E}^{- n/4},1)$ and $g_{2} = (\delta,1,\varpi_{E})$.
According to Propositions~\ref{prop:example. N/K} and \ref{prop:example. T/K}, $N/G(F)^{1}$ is generated by the images of $g_1$ and $g_2$ in $N/G(F)^1$, and  $T/G(F)^{1}$ is generated by the images of $g_{1}^{2}, g_{2}^{2}$ in $T/G(F)^{1}$.
Replacing $\mathcal{T}$ if necessary, we assume that $T_{g_{1}^{a_{1}}g_{2}^{a_{2}}} = 1$ for $a_{1}, a_{2}\in \Z$.

\begin{proposition}
\label{prop:calculationofmuTforHV}
We have $\mu_{\mathcal{T}}(g_{1}^{a_{1}}g_{2}^{a_{2}}, g_{1}^{b_{1}}g_{2}^{b_{2}}) = (-1)^{a_{2}b_{1}}$ for $a_{1}, a_{2}, b_{1}, b_{2}\in \Z$.
\end{proposition}
\begin{proof}
By Condition~\ref{conditionforcT}, we have
\[
T_{g_{1}^{a_{1}}g_{2}^{a_{2}}g_{1}^{b_{1}}g_{2}^{b_{2}}} = \theta^{1}(g_{2}^{-a_2 - b_2} g_{1}^{-a_1 - b_1} g_{1}^{a_1} g_{2}^{a_2} g_{1}^{b_1} g_{2}^{b_2}) T_{g_{1}^{a_{1} + b_{1}}g_{2}^{a_{2} + b_{2}}} = \theta^{1}(g_{2}^{-a_2 - b_2} g_{1}^{-a_1 - b_1} g_{1}^{a_1} g_{2}^{a_2} g_{1}^{b_1} g_{2}^{b_2}).
\]
Since $g_1, g_2 \in N$, we obtain that 
\begin{align*}
\theta^{1}(g_{2}^{-a_2 - b_2} g_{1}^{-a_1 - b_1} g_{1}^{a_1} g_{2}^{a_2} g_{1}^{b_1} g_{2}^{b_2}) &= \theta^{1}(g_{2}^{-b_2} g_{2}^{- a_2} g_{1}^{-b_1} g_{2}^{a_2} g_{1}^{b_1} g_{2}^{b_2}) \\
&= \theta^{1}(g_{2}^{- a_2} g_{1}^{-b_1} g_{2}^{a_2} g_{1}^{b_1}) \\
& = \varepsilon(\delta^{a_{2}(q^{-(n/2)b_{1}} - 1)}).
\end{align*}
Since 
\[
1 = T_{g_{1}^{a_{1}}g_{2}^{a_{2}}} \circ T_{g_{1}^{b_{1}}g_{2}^{b_{2}}} = \mu_{\mathcal{T}}(g_{1}^{a_{1}}g_{2}^{a_{2}}, g_{1}^{b_{1}}g_{2}^{b_{2}}) \cdot T_{g_{1}^{a_{1}}g_{2}^{a_{2}}g_{1}^{b_{1}}g_{2}^{b_{2}}},
\]
we have $\mu_{\mathcal{T}}(g_{1}^{a_{1}}g_{2}^{a_{2}}, g_{1}^{b_{1}}g_{2}^{b_{2}}) = \varepsilon(\delta^{a_{2}(q^{-(n/2)b_{1}} - 1)})^{-1}$.
By Lemma~\ref{lem:or_2(q^a-1)}, we have 
\begin{align*}
\ord_{2}(a_{2}(q^{-(n/2)b_{1}} - 1))  &= \ord_{2}(a_2) +\ord_{2}((n/2) b_{1}) + \ord_{2}(q + 1) \\
&= \ord_{2}(a_2)  + \ord_{2}(b_{1}) + \ord_{2}(n) + \ord_{2}(q + 1) -1 \\
&= \ord_{2}(a_2)  + \ord_{2}(b_{1}) + r-1 \ge r-1.
\end{align*}
Since $\varepsilon$ has the order $2^{r}$ and $\delta\notin (\resi[D]^{\times})^{2}$, we have $\varepsilon(\delta^{a_{2}(q^{-(n/2)b_{1}} - 1)})\in \{\pm 1\}$, and $\varepsilon(\delta^{a_{2}(q^{-(n/2)b_{1}} - 1)}) = 1$ if and only if $\ord_{2}(a_2)  + \ord_{2}(b_{1}) \ge 1$.
Thus, we obtain the claim.
\end{proof}

Let $\pi$ be an irreducible supercuspidal representation of $G(F)$ containing $\theta^{1}$ as a $G(F)^{1}$-representation.
Let $\Rep^{\mathfrak{s}}(G(F))$ be the Bernstein block containing $\pi$.
According to Corollary~\ref{corollarymus=muvee}, the cohomology class $[\mu_{\mathfrak{s}}] \in H^{2}(X_{\unr}(G)_{\pi},\C^{\times})$ factors through the restriction map
\[
X_{\unr}(G)_{\pi} = (G(F)/T)^{*} \twoheadrightarrow (N/T)^{*}
\]
and is described as follows.
Let $\chi_1, \chi_2 \in (N/T)^{*}$ be the characters such that $\chi_{1}(g_{1}) = 1$, $\chi_{1}(g_{2}) = -1$, $\chi_{2}(g_{1}) = -1$, and $\chi_{2}(g_{2}) = 1$.
Then $(N/T)^{*}$ is generated by $\chi_{1},\chi_{2}$, and we obtain from Corollary~\ref{corollarymus=muvee} and Proposition~\ref{prop:calculationofmuTforHV} that the cohomology class $[\mu_{\mathfrak{s}}]$ is represented by
\begin{equation}
\label{eq:computationofmusexample}
\mu_{\mathfrak{s}}(\chi_{1}^{a_{1}}\chi_{2}^{a_{2}},\chi_{1}^{b_{1}}\chi_{2}^{b_{2}}) = (-1)^{a_{2}b_{1}}
\end{equation}
for $a_{1}, a_{2}, b_{1}, b_{2}\in \Z/2\Z$.

\subsection{Spectral side}

In this subsection, we compute the cohomology class $[\mu_{\mathfrak{s}^{\vee}}] = [\mu_{(\varphi_{\pi}, \rho_{\pi})}]$ for the enhanced $L$-parameter $(\varphi_{\pi}, \rho_{\pi})$ attached to $\pi$.
We use the same notation as in Sections~\ref{sec:aconjecture} and \ref{subsec:Computation of the cohomology class attached to non-singular Lparameter}.
Let $\sigma\in W_{F}$ and $\tau\in I_{F}$ be lifts of $\sigma_{n}\in\Gamma_{E_{n}/E}$ and $x\mapsto \overline{x}\in \Gamma_{E_{n}/F_{n}}$, respectively.
As we have seen in Section~\ref{subsec: example, setting}, the group $G$ splits over $E_{n}$.
The dual group $\widehat{G}$ of $G$ is given by
\[
\widehat{G} = (\GL_{n}(\C)\times \C^{\times}\times \C^{\times})^{2}/A
\]
where
\[
A = \langle ((t_{1}1_{n},t_{1}^{2},t_{1}),(1_{n},1,t_{1}^{-1})),((1_{n},1,t_{2}^{-1}),(t_{2}1_{n},t_{2}^{2},t_{2}))\mid t_{1},t_{2}\in \C^{\times}\rangle,
\]
and $1_{n}\in\GL_{n}(\C)$ is the identity matrix.
The action of $\sigma$ on $\widehat{G}$ is trivial (since $D^{\times}$ is an inner form of $\GL_{n}$) and $\tau$ swaps the two factors.

We also have
\[
\widehat{S} = (\GL_{1}(\C)^{n}\times \GL_{1}\times \GL_{1})^{2}/A.
\]
The element $\tau$ again swaps the two factors, and $\sigma$ acts on $\widehat{S}$ as 
\[
\sigma(((t_{1},\ldots,t_{n}),x,y),((t'_{1},\ldots,t'_{n}),x',y')) = (((t_{2},\ldots,t_{n},t_{1}),x,y),((t'_{2},\ldots,t'_{n},t'_{1}),x',y')).
\]

Note that $X_{\unr}(^LG) = Z(\widehat{G})^{I_{F},\circ}$ since the Frobenius element $\sigma$ acts trivially on $\widehat{G}$.
The subgroup $Z(\widehat{G})^{I_{F},\circ}$ of $\widehat{G}$ is given by
\begin{equation}\label{eq:example, Z(G^)^I_F}
X_{\unr}(^LG) = Z(\widehat{G})^{I_{F},\circ} = Z(\widehat{G})^{I_{F}} = \{((t1_{n},x,y),(t1_{n},x,y))\mid t,x,y\in \C^{\times}\}/A',
\end{equation}
where 
\[
A' = \{((t1_{n},t^{2},1),(t1_{n},t^{2},1))\mid t\in \C^{\times}\}.
\]
Indeed, let $g = ((g_{1},x_{1},y_{1}),(g_{2},x_{2},y_{2}))\in Z(\widehat{G})^{I_{F}}$.
Since $g$ is fixed by $\tau$, we have $g_{2} = g_{1}t_{1},\ g_{1} = g_{2}t_{2},\ x_{2} = x_{1}t_{1}^{2},\ x_{1} = x_{2}t_{2}^{2}$, and $y_{2} = t_{1}t_{2}^{-1}y_{1}$ for some $t_{1},t_{2}\in \C^{\times}$.
We obtain from $g_{2} = g_{1}t_{1}$ and $g_{1} = g_{2}t_{2}$ that $t_{1}t_{2} = 1$.
Thus, we obtain that 
\begin{align*}
g &= ((g_{1},x_{1},y_{1}),(g_{2},x_{2},y_{2})) \\
&= ((g_{1},x_{1},y_{1}),(g_{1} t_{1}, x_{1} t_{1}^{2} ,y_{1} t_{1}^{2})) \\
&= ((g_{1},x_{1},t_{1}y_{1}),(g_{1},x_{1},t_{1}y_{1})) \cdot ((1_{n}, 1, t_{1}^{-1}), (t_{1}1_{n}, t_{1}^{2}, t_{1})).
\end{align*}
Since $((1_{n}, 1, t_{1}^{-1}), (t_{1}1_{n}, t_{1}^{2}, t_{1})) \in A$, we may assume that $g_{2} = g_{1},\ x_{2} = x_{1}$, and $y_{2} = y_{1}$.
Since $g \in Z(\widehat{G})$, we have $g((h_{1},a_{1},b_{1}),(h_{2},a_{2},b_{2}))g^{-1} = ((h_{1},a_{1},b_{1}),(h_{2},a_{2},b_{2}))$ in $\widehat{G}$ for any $h_{1},h_{2}\in \GL_{n}(\C)$ and $a_{1},a_{2},b_{1},b_{2}\in \C^{\times}$.
Hence, there exist $t_{1},t_{2}\in \C^{\times}$ such that $g_{1}h_{1}g_{1}^{-1} = t_{1}h_{1},\ g_{1}h_{2}g_{1}^{-1} = t_{2}h_{2},\ a_{1} = t_{1}^{2}a_{1},\ a_{2} = t_{2}^{2}a_{2},\ b_{1} = t_{1}t_{2}^{-1}b_{1}$, and $b_{2} = t_{1}^{-1}t_{2}b_{2}$.
We obtain from $b_{1} = t_{1}t_{2}^{-1}b_{1}$ that $t_{1} = t_{2}$ and from $a_{1} = t_{1}^{2}a_{1}$ that $t_{1} = \pm 1$.
In particular, we have $t_{1}t_{2} = 1$.
Combining this with the conditions $g_{1}h_{1}g_{1}^{-1} = t_{1}h_{1}$ and $g_{1}h_{2}g_{1}^{-1} = t_{2}h_{2}$, we obtain that $g_{1}h_{1}h_{2}g_{1}^{-1} = h_{1}h_{2}$ for any $h_{1},h_{2}\in \GL_{n}$.
Thus, we conclude that $g_{1}\in Z(\GL_{n})(\C)$.
Therefore, we obtain the above description of $Z(\widehat{G})^{I_{F}}$.
Since this is connected, we have $Z(\widehat{G})^{I_{F}} = Z(\widehat{G})^{I_{F},\circ}$.

We also record an explicit description of the isomorphism $X_{\unr}(G) \simeq X_{\unr}(^LG)$.
Let $z = ((t1_{n},x,y),(t1_{n},x,y)) \in X_{\unr}(^LG)$ and $g \in G(F)$.
According to Lemma~\ref{lem:description of Henniart's group}, we have $g = (\varpi_{D}^{a}d,\varpi_{E}^{b}x',\varpi_{E}^{c}y')$ for some $a, b, c \in \mathbb{Z}$, $d\in \mathcal{O}_{D}^{\times}$, and $x',y'\in\mathcal{O}_{E}^{\times}$ such that $a + 2b = 0$.
Then the unramified character $\chi_{z}$ corresponding to $z$ is defined as
\begin{equation}
\label{eq:explicitKottwitzisom}
\chi_{z}(g) = t^{a} x^{b} y^{c},
\end{equation}
which is well-defined by the condition that $a + 2b = 0$.

We compute the $L$-parameter $\varphi_{\pi} = {^Lj} \circ \varphi_{\theta}$ attached to $\pi$.
First, we compute the $L$-parameter $\varphi_{\theta} \colon W_{F} \to {^LS}$ of $\theta$.
Since the computations of $[\mu_{\mathfrak{s}}]$ and $[\mu_{\mathfrak{s}^{\vee}}]$ only depend on the inertial equivalence class of $\pi$, we may freely replace $\theta$ with $\theta\chi$ for $\chi\in X_{\unr}(G)$.
Note that $X_{\unr}(G)\to X_{\unr}(S)$ is surjective since $S$ is elliptic.
Hence, we may assume that there exists $\theta^{D}\colon E_{n}^{\times}\to \C^{\times}$ such that $\theta^{D}(\varpi_{E_{n}}) = 1$ and $\theta(d,x,y) = \theta^{D}(d)$, where $\varpi_{E_n}$ is a uniformizer of $E_n$ that corresponds to $\sigma^{-n} \in W_{E_n} \subset W_F$ via the local class field theory.
We define $\varphi_{S}^{D}$ as the composition 
\[
W_{F}\xrightarrow{\text{$L$-parameter of $\theta^{D}$}}{}^{L}(\Res_{E_{n}/F}\GL_{1})\xrightarrow{\text{projection}}\widehat{\Res_{E_{n}/F}\GL_{1}}\simeq (\GL_{1}(\C)^{n})^{2},
\]
where $\Res_{E_{n}/F}\GL_{1}$ denotes the Weil restriction of $\GL_{1}$ to $F$.
Then the $L$-parameter $\varphi_{\theta}$ is given by
\[
\varphi_{\theta}(w) = ((\varphi_{S, 1}^{D}(w), 1, 1), (\varphi_{S, 2}^{D}(w), 1, 1)) \rtimes w
\]
for $w \in W_{F}$, where we write $\varphi_{S}^{D}(w) = (\varphi_{S, 1}^{D}(w), \varphi_{S, 2}^{D}(w))$.
The cohomology class $[\varphi_{S}^{D}]$ of the 1-cocycle $\varphi_{S}^{D}$ is given by Shapiro's lemma as follows.
The character $\theta^{D}$ is identified with a character $W_{E_{n}}\to \C^{\times}$ by the local class field theory, which we also regard as an element of $H^{1}(W_{E_{n}},\C^{\times})$.
By Shapiro's lemma, we have 
\[
H^{1}(W_{E_{n}},\C^{\times})\simeq H^{1}(W_{F},\Ind_{W_{E_{n}}}^{W_{F}}\C^{\times}) = H^1(W_{F}, \widehat{\Res_{E_{n}/F}\GL_{1}}).
\]
Through this isomorphism, the character $\theta^{D}$ corresponds to $[\varphi_{S}^{D}] \in H^1(W_{F}, \widehat{\Res_{E_{n}/F}\GL_{1}})$.
For later use, we record an explicit description of $\varphi_{S}^{D}$.
Define a map $s\colon W_{E_{n}}\backslash W_{F}\to W_{F}$ by $s(W_{E_{n}}\sigma^{k}) = \sigma^{k}\ (0\le k\le n - 1)$ and $s(W_{E_{n}}\sigma^{k}\tau) = \sigma^{k}\tau\ (0\le k\le n - 1)$.
For $x\in W_{E_{n}}\backslash W_{F}$, we define a map $s_{x}\colon W_{F}\to W_{E_{n}}$ by
\[
s_{x}(w)s(xw) = s(x)w.
\]
Then we have $\varphi_{S}^{D}(w)(x) = \theta^{D}(s_{x}(w))$.

We first consider the case $w\in I_{E_{n}}$.
For $x = W_{E_{n}} \sigma^{k}$, we have 
\[
s_{W_{E_{n}}\sigma^{k}}(w) = s(W_{E_n}\sigma^{k})ws(W_{E_n}\sigma^{k}w)^{-1} = \sigma^{k}w\sigma^{-k}.
\]
We note that the character $\theta^{D} \colon W_{E_{n}} \to \mathbb{C}^{\times}$ is trivial on the wild inertia subgroup $P_{E_{n}}$ of $W_{E_{n}}$.
Since $\sigma^{k}w\sigma^{-k} \equiv w^{q^{k}} \pmod{P_{E_n}}$, we have $\varphi_{S}^{D}(w)(W_{E_{n}}\sigma^{k}) = \theta^{D}(w^{q^{k}})$.
Similarly, since
\[
s_{W_{E_{n}}\sigma^{k}\tau}(w) = s(W_{E_n}\sigma^{k}\tau)w s(W_{E_n}\sigma^{k}\tau w)^{-1} = \sigma^{k}\tau w  \tau^{-1}\sigma^{-k} \equiv w^{q^{k}} \pmod{P_{E_n}},
\]
we obtain that $\varphi_{S}^{D}(w)(W_{E_n}\sigma^{k}\tau) = \theta^{D}(w^{q^{k}})$.

We next calculate $\varphi_{S}^{D}(\tau)$.
We have 
\[
s_{W_{E_n}\sigma^{k}}(\tau) = s(W_{E_n} \sigma^{k})\tau s( W_{E_n}\sigma^{k}\tau)^{-1} = 1
\]
and
\[
s_{W_{E_n}\sigma^{k}\tau}(\tau) = s(W_{E_n} \sigma^{k}\tau)\tau s(W_{E_n} \sigma^{k}\tau^{2})^{-1} = \sigma^{k}\tau^{2}\sigma^{-k}.
\]
Hence $\varphi_{S}^{D}(\tau)(W_{E_n} \sigma^{k}) = 1$ and $\varphi_{S}^{D}(\tau)(W_{E_n} \sigma^{k}\tau) = \theta^{D}(\tau^{2q^{k}})$.

Finally, we calculate $\varphi_{S}^{D}(\sigma)$.
We have 
\begin{align*}
s_{W_{E_{n}}\sigma^{k}}(\sigma) &= s(W_{E_{n}}\sigma^{k})\sigma s(W_{E_{n}}\sigma^{k}\sigma)^{-1} \\
&= \sigma^{k}\sigma s(W_{E_{n}}\sigma^{k + 1})^{-1} \\
&=
\begin{cases}
 1 & (k < n-1), \\
 \sigma^{n} & (k = n-1).
\end{cases}
\end{align*}
Since $\theta^{D}(\sigma^{n}) = \theta^{D}(\varpi_{E_{n}}^{-1}) = 1$, we obtain that
$
\varphi_{S}^{D}(\sigma)(W_{E_n} \sigma^{k})= 
1$.
Similarly, we have
\begin{align*}
s_{W_{E_n}\sigma^{k} \tau}(\sigma) &= s(W_{E_n}\sigma^{k}\tau)\sigma s(W_{E_n}\sigma^{k}\tau\sigma)^{-1} \\
&= 
\begin{cases}
\sigma^{k}\tau\sigma \tau^{-1}\sigma^{-k - 1} & (k < n-1), \\
\sigma^{n-1} \tau \sigma \tau^{-1} & (k = n-1). 
\end{cases}
\end{align*}
Hence, we obtain that 
\[
\varphi_{S}^{D}(\sigma)(W_{E_n} \sigma^{k}\tau)= 
\begin{cases}
\theta^{D}(\tau^{q^{k} - q^{k + 1}}) & (k < n-1), \\
\theta^{D}(\tau^{q^{n - 1} - q^{n}})\theta^{D}(\sigma^{n}) = \theta^{D}(\tau^{q^{n - 1} - q^{n}}) & (k = n-1).
\end{cases}
\]

Summarizing the calculations above, we have
\begin{align*}
\varphi_{S}^{D}(w) & = ((\theta^{D}(w),\theta^{D}(w)^{q},\ldots,\theta^{D}(w)^{q^{n - 1}}),(\theta^{D}(w),\theta^{D}(w)^{q},\ldots,\theta^{D}(w)^{q^{n - 1}}))\ (w\in I_{E_{n}}),\\
\varphi_{S}^{D}(\tau) & = ((1,1,\ldots,1),(\theta^{D}(\tau^{2}),\theta^{D}(\tau^{2q}),\ldots,\theta^{D}(\tau^{2q^{n - 1}}))), \\
\varphi_{S}^{D}(\sigma) & = ((1,1,\ldots,1),(\theta^{D}(\tau^{1 - q}),\ldots,\theta^{D}(\tau^{q^{n - 2} - q^{n - 1}}),\theta^{D}(\tau^{q^{n - 1} - q^{n}}))).
\end{align*}

We also record a lemma about $\theta^{D}$ for later use.
\begin{lemma}\label{lem:theta(tau^2)}
We have
$\theta^{D}(\tau^{2})^{2^{r - 1}} = \theta^{D}(\tau^{1 - q})^{2^{r - 1}} = -1$.
\end{lemma}
\begin{proof}
Let $W_{F_{n}}^{\mathrm{ab}}$ (resp.\ $W_{E_{n}}^{\mathrm{ab}}$) be the abelianization of $W_{F_{n}}$ (resp.\ $W_{E_{n}}$) and $\tau_{F_{n}}\in F_{n}^{\times}$ be the element corresponding to $\tau \in W_{F_n}$ via the local class field theory for $F_{n}$.
We also think of $\tau_{F_{n}}$ as an element of $E_{n}^{\times}$.
The transfer map $W_{F_{n}}^{\mathrm{ab}}\to W_{E_{n}}^{\mathrm{ab}}$ sends $\tau$ to $\tau^{2}$.
Hence $\tau_{F_{n}}\in E_{n}^{\times}$ corresponds to $\tau^{2}\in W_{E_{n}}^{\mathrm{ab}}$ via the local class field theory for $E_{n}$.
Therefore $\theta^{D}(\tau^{2}) = \theta^{D}(\tau_{F_{n}})$.
Since  $\tau\in W_{F_{n}}\setminus W_{E_{n}}$, we have $\tau_{F_{n}}\notin N_{E_{n}/F_{n}}(E_{n}^{\times})$.
Since the norm map induces $x\mapsto x^{2}$ on the residue field $\resi[F_n]^{\times} = \resi[E_n]^{\times}$, the image of $\tau_{F_{n}}$ in $\resi[F_{n}]^{\times} = \resi[E_{n}]^{\times}$ is not contained in $(\resi[E_{n}]^{\times})^{2}$.
By the assumptions, we have $r = \ord_{2}(\#\resi[E_{n}]^{\times})$ and $\theta^{D}$ has order $2^{r}$ on $\resi[E_{n}]^{\times}$.
Hence, we have $\theta^{D}(\tau^{2})^{2^{r - 1}} = \theta^{D}(\tau_{F_n})^{2^{r-1}} = -1$.
Since we assume that $4 \mid (q-3)$, we also obtain that $ \theta^{D}(\tau^{1 - q})^{2^{r - 1}} = -1$.
\end{proof}

Next, we calculate ${}^{L}j$ following \cite[2.3]{MR4768407}, which is based on \cite[Section~2.6]{MR909227}.
As in  Section~\ref{subsec: example, setting}, the set of roots is given by $\Phi(G, S) = \{\alpha_{i,j}^{(1)},\alpha_{i,j}^{(2)}\mid 1\le i\ne j\le n\}$.
We have $\sigma(\alpha_{i,j}^{(k)}) = \alpha_{i - 1,j - 1}^{(k)}$ (with indices taken modulo $n$) and $\tau(\alpha_{i,j}^{(k)}) = \alpha_{i,j}^{(3 - k)}$.
For $\alpha\in \Phi(G, S)$, let $W_{\alpha}$ (resp.\ $W_{\pm \alpha}$) be the stabilizer of $\alpha$ (resp.\ $\{\alpha,-\alpha\}$) in $W_{F}$.
Then we have 
\[
W_{\alpha_{i,j}^{(k)}} = W_{E_{n}},\quad 
W_{\pm\alpha_{i,j}^{(k)}} = \begin{cases}
W_{E_{n}}\{1,\sigma^{n/2}\} & \text{} \lvert i - j\rvert = n/2, \\
W_{E_{n}} & \text{otherwise.}
\end{cases}
\]
A root $\alpha$ is called symmetric if $W_{\alpha} \ne W_{\pm \alpha}$ and asymmetric otherwise. 
If $\alpha$ is symmetric, then $\overline{F}^{W_{\alpha}}\supset \overline{F}^{W_{\pm \alpha}}$ is a quadratic extension.
In this case, $\alpha$ is called unramified if this extension is unramified and is called ramified otherwise.
In our case, $\alpha_{i,j}^{(k)}$ is symmetric if and only if $\lvert i - j\rvert = n/2$ and in this case, it is unramified.

Let $\chi = \{\chi_{\alpha}\}_{\alpha\in \Phi(G, S)}$ be the unique minimally ramified $\chi$-data for $\Phi(G, S)$ in the sense of \cite[Definition~4.15]{MR4768407}.
In particular, $\chi_{\alpha_{i,j}^{(k)}}$ is trivial unless $\lvert i - j\rvert =  n/2$.
The set of $\chi$-data $\chi$ determines an $L$-embedding ${^LS} \to {^LG}$ as
\[
{}^{L}j(s \rtimes w) = \widehat{j}(sr_{\chi}(w))n(\omega_{S,G}(w)) \rtimes w
\]
for $s\in \widehat{S}$ and $w\in W_{F}$, where $r_{\chi}(w)\in\widehat{S}$ is a function defined in \cite[(4.10)]{MR4768407}, $\omega_{S,G}(w)$ is the element in the Weyl group $W_{\widehat{G}}(\widehat{j}(\widehat{S}))$ whose action on $\widehat{j}(\widehat{S})$ is given by $\widehat{j}(s)\mapsto \Ad(1 \rtimes w)^{-1}
(\widehat{j}(\Ad(1 \rtimes w)(s)))$, and $n$ is a Tits lift.
In our case, $n(\omega_{S,G}(w)) = 1$ for $w\in I_{F}$ and
\[
n(\omega_{S,G}(\sigma)) = ((\omega_{\sigma}, 1, 1), (\omega_{\sigma}, 1, 1)),
\]
where 
\[
\omega_{\sigma} = 
\begin{pmatrix}
0 & -1 & \cdots & 0 & 0\\
& & \ddots & & \\
0 & 0 & \cdots & 0 & -1\\
1 & 0 & \cdots & 0 & 0 
\end{pmatrix}.
\]

We calculate the function $r_{\chi}$.
To do this, we fix some data.
Let $\underline{\Phi} = \Phi(G, S)/W_{F}$ be the set of $W_{F}$-orbits in $\Phi(G, S)$.
Then $\underline{\Phi}\simeq \{\alpha_{1,k}^{(1)}\mid 2\le k\le n\}$ and $\underline{\Phi}/\{\pm 1\} \simeq \{\alpha_{1,k}^{(1)}\mid 2\le k\le n/2 + 1\}$.
For $\alpha = \alpha_{1,l}^{(1)}$, we have
\[
W_{\pm\alpha}\backslash W_{F}\simeq \{\sigma^{k},\sigma^{k}\tau\mid 0\le k\le N - 1\}
\]
where $N = n$ if $l\ne n/2 + 1$ and $N = n/2$ otherwise.
We define a section $u^{\alpha}\colon W_{\pm\alpha}\backslash W_{F}\to W_{F}$ by $u^{\alpha}(W_{\pm\alpha} \sigma^{k}\tau^{e}) = \sigma^{k}\tau^{e}$ for $0\le k\le N - 1$ and $e = 0,1$.
For $x \in W_{\pm\alpha}\backslash W_{F}$, we define $u^{\alpha}_{x} \colon W_{F} \to W_{\pm \alpha}$ by $u^{\alpha}_{x}(w)u^{\alpha}(xw) = u^{\alpha}(x)w$.
The quotient $W_{\alpha}\backslash W_{\pm \alpha}$ is not trivial only when $\alpha = \alpha_{i,j}^{(k)}$ with $\lvert i - j\rvert = n/2$.
In this case, we define $v^{\alpha}\colon W_{\alpha}\backslash W_{\pm \alpha}\to W_{\pm \alpha}$ by $v^{\alpha}(W_{\alpha}) = 1$ and $v^{\alpha}(W_{\alpha}\sigma^{n/2}) = \sigma^{n/2}$.
Otherwise, we set $v^{\alpha} = 1$.
For $x \in W_{\alpha}\backslash W_{\pm \alpha}$, we define $v^{\alpha}_{x} \colon W_{\pm \alpha} \to W_{\alpha}$ by $v^{\alpha}_{x}(w)v^{\alpha}(xw) = v^{\alpha}(x)w$.

Let $w\in W_{F}$.
The element $r_{\chi}(w)$ is defined as $r_{\chi}(w) = \prod_{\alpha,x}u^{\alpha}(x)^{-1}(\alpha)(\chi_{\alpha}(v^{\alpha}_{1}(u^{\alpha}_{x}(w))))$, where $\alpha$ runs through $\underline{\Phi}/\{\pm 1\}$ and $x$ runs through $W_{\pm\alpha}\backslash W_{F}$.
Since $\chi_{\alpha} = 1$ for asymmetric $\alpha$, we only have to consider symmetric $\alpha$.
In our case, we only consider $\alpha_{1,n/2 + 1}^{(1)}$.
For simplicity of the notation, we write $\alpha = \alpha_{1,n/2 + 1}^{(1)}$.
Then we have
\[
r_{\chi}(w) = \prod_{0\le k\le n/2 - 1,e = 0,1}u^{\alpha}(\sigma^{k}\tau^{e})^{-1}(\alpha)(\chi_{\alpha}(v_{1}^{\alpha}(u^{\alpha}_{\sigma^{k}\tau^{e}}(w)))).
\]
Note that $u^{\alpha}(\sigma^{k}\tau^{e})^{-1}(\alpha)\in \Phi(G, S) \subset X^{*}(S)$, and we regard it as an element in $X_{*}(\widehat{S})$.

\begin{lemma}
We have $r_{\chi}(w) = 1$ for all $w \in I_{F} \cup \{\sigma\}$.
\end{lemma}
\begin{proof}
First assume that $w\in I_{E_{n}}$.
Then we have 
\[
u^{\alpha}_{\sigma^{k}\tau^{e}}(w) = u^{\alpha}(\sigma^{k}\tau^{e})wu^{\alpha}(\sigma^{k}\tau^{e}w)^{-1} = \sigma^{k}\tau^{e}w(\sigma^{k}\tau^{e})^{-1}\in I_{E_{n}}.
\]
Hence, we have $v_{1}^{\alpha}(u^{\alpha}_{\sigma^{k}\tau^{e}}(w)) \in I_{E_{n}}$.
Since $\chi$ is minimally ramified, $\chi_{\alpha}$ is unramified.
Hence, we obtain that $\chi_{\alpha}(v_{1}^{\alpha}(u^{\alpha}_{\sigma^{k}\tau^{e}}(w))) = 1$.
Thus, we have $r_{\chi}(w) = 1$ in this case.
Similarly, since we have
\[
v_{1}^{\alpha}(u^{\alpha}_{\sigma^{k}\tau^{e}}(\tau)) =
\begin{cases}
1 & (e = 0), \\
\sigma^{k}\tau^{2}\sigma^{-k} & (e = 1),
\end{cases}
\]
we obtain that $r_{\chi}(\tau) = 1$.

We now calculate $r_{\chi}(\sigma)$.
In this case, we have
\[
u^{\alpha}_{\sigma^{k}}(\sigma) = 
\begin{cases}
1 & (k < n/2 - 1),\\
\sigma^{n/2} & (k = n/2 - 1),
\end{cases}
\quad
v^{\alpha}_{1}(1) = 
v^{\alpha}_{1}(\sigma^{n/2}) =1,
\]
and
\begin{align*}
& u^{\alpha}_{\sigma^{k}\tau}(\sigma) = 
(\sigma^{k}\tau\sigma^{-k})(\sigma^{k + 1}\tau^{-1}\sigma^{-k - 1})
\begin{cases}
1 & (k < n/2 - 1)\\
\sigma^{n/2} & (k = n/2 - 1)
\end{cases}\in 
I_{E_{n}}
\begin{cases}
1 & (k < n/2 - 1),\\
\sigma^{n/2} & (k = n/2 - 1),
\end{cases}
\\
& v^{\alpha}_{1}(I_{E_{n}}) =I_{E_{n}},\quad
v^{\alpha}_{1}(I_{E_{n}}\sigma^{n/2}) = I_{E_{n}}.
\end{align*}
Since $\chi$ is minimally ramified, $\chi_{\alpha}$ is unramified; hence $\chi_{\alpha}(I_{E_{n}}) = 1$.
Thus, we obtain the claim.
\end{proof}

We calculate the group $N_{\varphi_{\pi}}$ using the description of $\varphi_{\pi}$ above.
We begin with the calculation of the centralizer of $\varphi_{\pi}(I_{F})$ in $\widehat{G}$.
We define $\kappa \in \GL_{n}(\C)$ as $\kappa = \omega_{\sigma}^{n/2}$.
This element satisfies the following property.
\begin{lemma}\label{lem:example, permutation kappa and diagonal matrix}
Let $z\in \C^{\times}$ such that $z^{2^{r}} = 1$.
Then we have
\[
\kappa\diag(z,z^{q},\ldots,z^{q^{n - 1}})\kappa^{-1} = z^{2^{r-1}}\diag(z,z^{q},\ldots,z^{q^{n - 1}}).
\]
\end{lemma}
\begin{proof}
By a calculation, it suffices to show that $z^{q^{i}(q^{n/2} - 1)} = z^{q^{i}(1 - q^{n/2})} = z^{2^{r - 1}}$ for $0\le i\le n/2 - 1$.
We obtain from Lemma~\ref{lem:or_2(q^a-1)} and the assumption $z^{2^{r}} = 1$ that $z^{q^{n/2} - 1} \in \{1,-1\}$, and $z^{q^{n/2} - 1} = 1$ if and only if $z^{2^{r - 1}} = 1$.
Since $q$ is odd, the lemma follows.
\end{proof}

Using the element $\kappa$, we can describe the centralizer of $\varphi_{\pi}(I_{F})$ in $\widehat{G}$ as follows.

\begin{lemma}\label{lem:example, centralizer of varphi(I_F)}
The centralizer of $\varphi_{\pi}(I_{F})$ in $\widehat{G}$ is given by
\begin{align*}
\left(
\{((g,x,y),(g,x,y))\mid g\in \widehat{S}^{D}, x,y \in \C^{\times}\}
\cup 
\{((g,x,y),(g,x,- y))\mid g\in \widehat{S}^{D}\kappa, x,y \in \C^{\times}\}
\right)
/A',
\end{align*}
where $\widehat{S}^{D}$ is the subgroup of $\GL_{n}(\C)$ consisting of diagonal matrices.
\end{lemma}
\begin{proof}
Let $((g_{1},x_{1},y_{1}),(g_{2},x_{2},y_{2}))\in \widehat{G}$ be an element of the centralizer of $\varphi_{\pi}(I_{F})$.
Take $w\in I_{E_{n}} = I_{E}$ such that $\theta^{D}(w)$ has order $2^{r}$.
By Lemma~\ref{lem:or_2(q^a-1)}, we have $\theta^{D}(w)^{q^{i}}\ne \theta^{D}(w)^{q^{j}}$ for all $0\le i,j \le n - 1$ with $i \ne j$.
We also have the following: $\theta^{D}(w)^{q^{i} - q^{j}} = -1$ if and only if $i - j\equiv n/2\pmod{n}$.
Let $d = (\theta^{D}(w),\theta^{D}(w)^{q},\ldots,\theta^{D}(w)^{q^{n - 1}})$.
Then $\varphi_{\pi}(w) = ((d,1,1),(d,1,1)) \rtimes w$.
Since $((g_{1},x_{1},y_{1}),(g_{2},x_{2},y_{2}))$ commutes with $\varphi_{\pi}(w)$, there exist $t_{1},t_{2}\in \C^{\times}$ such that 
\[
((g_{1}dg_{1}^{-1},1,1),(g_{2}dg_{2}^{-1},1,1)) = ((dt_{1},t_{1}^{2},t_{1}t_{2}^{-1}),(dt_{2},t_{2}^{2},t_{1}^{-1}t_{2})).
\]
Hence, we have $t_{1}  = t_{2}\in \{1,-1\}$.
We also have $g_{1}dg_{1}^{-1} = dt_{1}$.
If $t_{1} = 1$, then it implies that $g_{1}$ is a diagonal matrix.
Assume that $t_{1} = -1$.
Then we have $d^{-1}g_{1}d = -g_{1}$.
The $(i,j)$-entry of $d^{-1}g_{1}d$ is $\theta^{D}(w)^{q^{j - 1} - q^{i - 1}}$ times the $(i,j)$-entry of $g_{1}$.
Hence, the $(i,j)$-entry of $g_{1}$ is non-zero only if $j - i\equiv n/2\pmod{n}$.
Thus, we obtain that $g_{1}\in \widehat{S}^{D}\kappa$.

By the computations of $\varphi_{\theta}(\tau)$ and $r_{\chi}$, we have $\varphi_{\pi}(\tau) = ((1,1,1),(d',1,1)) \rtimes \tau$, where $d' = (\theta^{D}(\tau^{2}),\theta^{D}(\tau^{2q}),\ldots,\theta^{D}(\tau^{2q^{n - 1}}))$.
According to Lemmas~\ref{lem:theta(tau^2)} and \ref{lem:example, permutation kappa and diagonal matrix}, we have 
\[
d' g_{1} (d')^{-1} =
\begin{cases}
g_{1} & (g_{1} \in \widehat{S}^{D}), \\
-g_{1} & (g_{1} \in \widehat{S}^{D} \kappa).
\end{cases}
\]
Since $((g_{1},x_{1},y_{1}),(g_{2},x_{2},y_{2}))$ commutes with $\varphi_{\pi}(\tau)$, we have 
\[
((g_{1}g_{2}^{-1},x_{1}x_{2}^{-1},y_{1}y_{2}^{-1}),(g_{2}d'g_{1}^{-1},x_{2}x_{1}^{-1},y_{2}y_{1}^{-1})) = ((s_{1}1_{n},s_{1}^{2},s_{1}s_{2}^{-1}),(d's_{2},s_{2}^{2},s_{1}^{-1}s_{2}))
\]
for some $s_{1},s_{2}\in \C^{\times}$.
We obtain from the conditions $g_{1}g_{2}^{-1} = s_{1}$ and $g_{2}d'g_{1}^{-1} = d's_{2}$ that 
\[
s_{1} = \begin{cases}
s_{2}^{-1} & (g_{1} \in \widehat{S}^{D}), \\
-s_{2}^{-1} & (g_{1} \in \widehat{S}^{D} \kappa).
\end{cases}
\]
Hence, we obtain that $x_{1} = x_{2}s_{1}^{2}$ and 
\[
y_{1} =  y_{2}s_{1}s_{2}^{-1} 
=
\begin{cases}
y_{2}s_{1}^{2} & (g_{1} \in \widehat{S}^{D}), \\
-y_{2}s_{1}^{2} & (g_{1} \in \widehat{S}^{D} \kappa).
\end{cases}
\]
Combining the above arguments, we obtain that 
\begin{align*}
((g_{1},x_{1},y_{1}),(g_{2},x_{2},y_{2})) &=
\begin{cases}
 ((g_1, x_1, y_1), (s_{1}^{-1} g_1, s_{1}^{-2} x_1, s_{1}^{-2} y_1)) &  (g_{1} \in \widehat{S}^{D}), \\
 ((g_1, x_1, y_1), (s_{1}^{-1} g_1, s_{1}^{-2} x_1, - s_{1}^{-2} y_1)) & (g_{1} \in \widehat{S}^{D} \kappa)
 \end{cases} \\
 &= 
 \begin{cases}
 ((g_1, x_1, s_1^{-1} y_1), (g_1, x_1, s_1^{-1} y_1)) \cdot ((1, 1, s_{1}), (s_{1}^{-1}, s_{1}^{-2}, s_1^{-1})) &  (g_{1} \in \widehat{S}^{D}), \\
 ((g_1, x_1, s_1^{-1} y_1), (g_1, x_1, - s_1^{-1} y_1)) \cdot ((1, 1, s_{1}), (s_{1}^{-1}, s_{1}^{-2}, s_1^{-1})) & (g_{1} \in \widehat{S}^{D} \kappa).
 \end{cases}
\end{align*}
Since $((1, 1, s_{1}), (s_{1}^{-1}, s_{1}^{-2}, s_1^{-1})) \in A$, these elements are contained in the set described in the lemma.
The reverse inclusion is easy to check.
\end{proof}

Now, we calculate the group $N_{\varphi_{\pi}}$.

\begin{proposition}
\label{prop:computationofNSZ_varphi}
Let $\mu_{n}(\C)$ denote the set of $n$-th roots of unity in $\C^{\times}$.
Then we have
\[
N_{\varphi_{\pi}} = (N_{\varphi_{\pi}, 1} \cup N_{\varphi_{\pi}, 2})/A',
\]
where 
\[
N_{\varphi_{\pi}, 1} = \{((g_{a, s},x,y),(g_{a, s},x,y))\mid a, x,y \in \C^{\times},s\in \mu_{n}(\C)\}
\]
and 
\[
N_{\varphi_{\pi}, 2} = \{
((g_{a, s}\kappa,x,y),(g_{a, s}\kappa,x,- y))\mid a, x,y \in \C^{\times},s\in\mu_{n}(\C)
\}
\]
with $g_{a,s}$ defined as
\[
g_{a, s} = a\diag(1,s^{-1},\ldots,s^{-(n - 1)}).
\]
The subgroup $S_{\varphi_{\pi}}$ is given by
\[
S_{\varphi_{\pi}} = \{((g_{a, s}, x,y), (g_{a, s}, x, y)) \mid a, x, y \in \C^{\times},\ s\in\{1,-1\}\}/
A'.
\]
Moreover, we have
\[
Z(\widehat{G})^{I_{F}, \circ}_{[\varphi_{\pi}]} = \left\{
((t1_{n},x,y),(t1_{n},x,y))\in Z(\widehat{G})^{I_{F}, \circ} \mid t^{n} = x^{n/2}, y \in \{1,-1\}
\right\}/A'.
\]
\end{proposition}
\begin{proof}
Let $((g,x,y),(g,x,\pm y))\in N_{\varphi_{\pi}} \subset \Cent_{\widehat{G}}(\varphi_{\pi}(I_{F}))$.
We write $\omega = \omega_{\sigma}$ and 
\[
\theta^{D}_{\tau} = (\theta^{D}(\tau^{1-q}),\theta^{D}(\tau^{1-q})^{q},\ldots,\theta^{D}(\tau^{1-q})^{q^{n-1}}).
\]
Then we have
\[
\varphi_{\pi}(\sigma) = ((\omega,1,1),(\theta^{D}_{\tau}\omega,1,1)) \rtimes \sigma.
\]
By \eqref{eq:example, Z(G^)^I_F}, there exist $t_{1},t_{2}\in \C^{\times}$ and $((t_{3}1_{n},x_{1},y_{1}),(t_{3}1_{n},x_{1},y_{1}))\in Z(\widehat{G})^{I_{F},\circ}$ such that
\[
((g\omega g^{-1},1,1),(g\theta^{D}_{\tau}\omega g^{-1},1,1)) = ((\omega t_{3}t_{1},x_{1}t_{1}^{2},y_{1}t_{1}t_{2}^{-1}),(\theta^{D}_{\tau}\omega t_{3}t_{2},x_{1}t_{2}^{2},y_{1}t_{1}^{-1}t_{2})).
\]
Set $s = t_{3}t_{1}$.
Then we have $g = s(\omega g\omega^{-1})$.

First, we assume that $g\in \widehat{S}^{D}$.
Since the conjugation action of $\omega$ sends $\diag(z_{1},\ldots,z_{n})$ to $\diag(z_{2},z_{3},\ldots,z_{n},z_{1})$, we obtain from the condition $g = s(\omega g\omega^{-1})$ that  $g = g_{a,s}$ for some $a\in \C^{\times}$ and that $s^{n} = 1$.
Since $g$ commutes with $\theta^{D}_{\tau}$, we obtain from $g\theta^{D}_{\tau}\omega g^{-1} = \theta^{D}_{\tau}\omega t_{3}t_{2}$ that $g = t_3 t_2 (\omega g \omega^{-1})$.
Thus, we have $t_{2} = s t_{3}^{-1} = t_{1}$, and hence $y_{1} = 1$.
Moreover, since $x_1 = t_1^{-2} = t_3^{2}s^{-2}$ and $s^{n} = 1$, we obtain that $x_1^{n/2} = t_{3}^{n}$.

Next, we consider the case $g\in \widehat{S}^{D}\kappa$.
Since $\kappa = \omega^{n/2}$ commutes with $\omega$, we obtain from $g = s(\omega g\omega^{-1})$ that $g = g_{a,s} \kappa$ for some $a\in \C^{\times}$ and that $s^{n} = 1$ by the same argument as above.
Since $g \theta^{D}_{\tau}g^{-1} = -\theta^{D}_{\tau}$ by Lemmas~\ref{lem:theta(tau^2)} and \ref{lem:example, permutation kappa and diagonal matrix}, we obtain from $g\theta^{D}_{\tau}\omega g^{-1} = \theta^{D}_{\tau}\omega t_{3}t_{2}$ that  $t_{1} = -t_{2}$ and $y_{1} = -1$.
Moreover, since $x_1 = t_1^{-2} = t_3^{2}s^{-2}$ and $s^{n} = 1$, we obtain that $x_1^{n/2} = t_{3}^{n}$.
From these arguments, we get the description of $N_{\varphi_{\pi}}$ and $Z(\widehat{G})^{I_{F},\circ}_{[\varphi_{\pi}]}$.

If $((g,x,y),(g,x,y))\in S_{\varphi_{\pi}}$, then $t_{3} = x_{1} = y_{1} = 1$.
Hence, $g$ should be in $\widehat{S}^{D}$ by the above argument.
In this case, we also have $t_{1}^{2} = 1$.
Hence $s = t_{1}\in \{1,-1\}$.
The reverse inclusions are immediate.
\end{proof}

Let $Z(D^1)$ denote the center of $D^1 := [D^{\times}, D^{\times}]$.
We define a finite central subgroup $Z$ of $G$ by
\[
Z(F) = Z(D^1) \times \{1\} \times \{1\} \subset \{(d,x,y)\in D^{\times}\times E^{\times}\times E^{\times}\mid N_{D/E}(d)x^{2}(y/\overline{y}) = 1\}.
\]
We write $\overline{S} = S/Z$ as in Section~\ref{subsec:Computation of the cohomology class attached to non-singular Lparameter}.
Then we have
\[
\widehat{\overline{S}} = \left\{
(((t_{1, i}), t_1, x_1, y_1), ((t_{2, i}), t_2, x_2, y_2)) \in (\GL_{1}(\C)^{n}\times \GL_{1} \times \GL_{1}\times \GL_{1})^{2} \mid \prod t_{1, i} = t_1^{n}, \prod t_{2, i} = t_2^{n}
\right\}/B,
\]
where 
\[
B  = \langle ((t_{1}1_{n},t_{1}, t_{1}^{2},t_{1}),(1_{n},1, 1,t_{1}^{-1})),((1_{n},1, 1,t_{2}^{-1}),(t_{2}1_{n},t_{2}, t_{2}^{2},t_{2}))\mid t_{1},t_{2}\in \C^{\times}\rangle.
\]
The dual isogeny of the natural surjection $S \to \overline{S}$ is given by
\[
\widehat{\overline{S}} \ni (((t_{1, i}), t_1, x_1, y_1), ((t_{2, i}), t_2, x_2, y_2)) \mapsto (((t_{1, i}), x_1, y_1), ((t_{2, i}), x_2, y_2)) \in \widehat{S}.
\]
According to Proposition~\ref{prop:computationofNSZ_varphi}, we have
\begin{align*}
[\widehat{\overline{S}}]^{+} &= S_{\varphi_{\pi}}^{+} \\
&= 
\left\{
((g_{a, s},t_{1}, x,y),(g_{a, s},t_{2}, x,y))\mid a, x, y,  t_1, t_2 \in \C^{\times},\ s\in\{1,-1\},\ a^{n} = t_1^{n} = t_{2}^{n}
\right\}/B',
\end{align*}
where 
\[
B' =  \langle ((t1_{n}, t, t^{2},1),(t1_{n}, t, t^{2}, 1))\mid t\in \C^{\times}\rangle.
\]
Thus, its group of connected components is described as
\begin{equation}
\label{pi0ofwidehatoverlineS+}
\pi_{0}([\widehat{\overline{S}}]^{+}) \isoarrow \mu_{n}(\C) \times \mu_{n}(\C) \times \{1,-1\}, \quad ((g_{a, s},t_{1}, x,y),(g_{a, s},t_{2}, x,y)) \mapsto (a/t_1, a/t_2, s).
\end{equation}
We identify the group $\pi_{0}\left(
Z\left(
\widehat{G/Z(G)}
\right)^{\Gamma_{F}}
\right)
$ with its image under the map 
\[
\pi_{0}\left(
Z\left(
\widehat{G/Z(G)}
\right)^{\Gamma_{F}}
\right)
 \to \pi_{0}([\widehat{\overline{S}}]^{+})
 \]
induced by the maps
\[
Z\left(
\widehat{G/Z(G)}
\right)^{\Gamma_{F}}\subset \widehat{(S/Z(G))}^{\Gamma_{F}} \to (\widehat{\overline{S}})^{\Gamma_{F}} \subset [\widehat{\overline{S}}]^{+}.
\]
Since 
\[
\pi_{0}\left(
Z\left(
\widehat{G/Z(G)}
\right)^{\Gamma_{F}}
\right) = \pi_0\left(
Z(\SL_{n}(\mathbb{C}) \times \SL_{n}(\mathbb{C}))^{\Gamma_{F}}
\right) \simeq \mu_{n}(\mathbb{C}),
\]
under the isomorphism in \eqref{pi0ofwidehatoverlineS+}, we have
\[
\pi_{0}\left(
Z\left(
\widehat{G/Z(G)}
\right)^{\Gamma_{F}}
\right) = \left\{
(z, z, 1) \mid z \in \mu_n(\mathbb{C})
\right\} \subset \mu_{n}(\C) \times \mu_{n}(\C) \times \{1,-1\}.
\]
We also obtain from Proposition~\ref{prop:computationofNSZ_varphi} that 
\[
N_{\varphi_{\pi}}^{+} = (N^{+}_{\varphi_{\pi}, 1} \cup N^{+}_{\varphi_{\pi}, 2})/B',
\]
where 
\[
N^{+}_{\varphi_{\pi}, 1} = \{((g_{a, s}, t_{1}, x,y),(g_{a, s}, t_{2}, x,y))\mid a, x,y, t_1, t_2 \in \C^{\times},\ s\in \mu_{n}(\C),\ a^{n}s^{\frac{n(n-1)}{2}} = t_1^{n} = t_2^{n}\}
\]
and 
\[
N^{+}_{\varphi_{\pi}, 2} = \{
((g_{a, s}\kappa, t_1, x,y),(g_{a, s}\kappa, t_2, x,- y))\mid a, x,y, t_1, t_2  \in \C^{\times},\ s\in\mu_{n}(\C),\ a^{n}s^{\frac{n(n-1)}{2}} = t_1^{n} = t_2^{n}
\}.
\]

Let $\eta \colon \pi_0([\widehat{\overline{S}}]^{+}) \to \mathbb{C}^{\times}$ be the character defined in Section~\ref{subsec:Computation of the cohomology class attached to non-singular Lparameter}.
Since $D$ is a division algebra, the precomposition of $\eta$ with the map
\[
\mu_{n}(\C)  \simeq \pi_{0}\left(
Z\left(
\widehat{G/Z(G)}
\right)^{\Gamma_{F}}
\right) \to \pi_{0}([\widehat{\overline{S}}]^{+})
\]
has order $n$.
Since $S_{\varphi_{\pi}}^{+} = [\widehat{\overline{S}}]^{+}$, we have $\rho_{\pi} = \eta$.

We now compute the cohomology class
\[
[\mu_{(\varphi_{\pi}, \rho_{\pi})}]
\in H^2\left(
N_{\varphi_{\pi}, \rho_{\pi}}^{+}/S_{\varphi_{\pi}}^{+}, \mathbb{C}^{\times}
\right) = H^2\left(
Z(\widehat{G})^{I_{F}, \circ}_{[(\varphi_{\pi}, \rho_{\pi})]}, \mathbb{C}^{\times}
\right).
\]
We fix a primitive $n$-th root of unity $\zeta \in \mu_{n}(\C)$ and a square root $\zeta^{\frac{n-1}{2}} \in \C^{\times}$ of $\zeta^{n-1}$.
We define an element $n_{\zeta} \in N^{+}_{\varphi_{\pi}, 1}$ as $n_{\zeta} = ((g_{1, \zeta}, \zeta^{\frac{n-1}{2}}, 1, 1), (g_{1, \zeta}, \zeta^{\frac{n-1}{2}}, 1, 1))$.
We also define an element $n_{\kappa} \in N^{+}_{\varphi_{\pi}, 2}$ as $n_{\kappa} = ((\kappa, 1, 1, 1), (\kappa, 1, 1, -1))$.
According to Proposition~\ref{prop:computationofNSZ_varphi}, the set 
\[
\left\{
n_{\kappa}^{a} n_{\zeta}^{b} \mid a \in \{0,1\},\ 0 \le b \le \frac{n}{2}-1
\right\}
\]
is a set of representatives for $N_{\varphi_{\pi}}^{+}/S_{\varphi_{\pi}}^{+}$.
We define a subgroup $N'$ of $N_{\varphi_{\pi}}^{+}$ as
\[
N' = \left \langle
S_{\varphi_{\pi}}^{+}, n_{\zeta}^{a} \mid a \in 2 \mathbb{Z}
\right \rangle.
\]
Then the set 
\[
\left\{
n_{\kappa}^{a} n_{\zeta}^{b} \mid a, b \in \{0, 1\}
\right\}
\]
is a set of representatives for $N_{\varphi_{\pi}}^{+} / N'$.

\begin{lemma}
\label{lemma:extensionofetaanditsnormalizer}
The character $\eta$ extends to a character $\widetilde{\eta}$ of $N'$.
Moreover, for any extension $\widetilde{\eta}$, the elements $n_{\zeta}$ and $n_{\kappa}$ normalize $\widetilde{\eta}$.
\end{lemma}
\begin{proof}
The first claim follows from the fact that $N' \subset \widehat{\overline{S}}$ is abelian.
Since $n_{\zeta} \in \widehat{\overline{S}}$, the element $n_{\zeta}$ normalizes $\widetilde{\eta}$.
We prove that $n_{\kappa}$ normalizes $\widetilde{\eta}$.
Let $((g_{a, s},t_{1}, x,y),(g_{a, s},t_{2}, x,y)) \in N'$.
Then we have
\[
n_{\kappa} ((g_{a, s},t_{1}, x,y),(g_{a, s},t_{2}, x,y)) n_{\kappa}^{-1} = ((g'_{a, s},t_{1}, x,y),(g'_{a, s},t_{2}, x,y)),
\]
where 
\[
g'_{a, s} = a\diag(s^{-\frac{n}{2}}, s^{-(\frac{n}{2} + 1)}, \ldots , s^{-(n-1)}, 1, s^{-1},\ldots, s^{-(\frac{n}{2} - 1)}).
\]
Since  $s^{n/2} = 1$ by the definition of $N'$, we obtain that $g_{a, s} = g'_{a, s}$.
Thus, we obtain that $n_{\kappa} ((g_{a, s},t_{1}, x,y),(g_{a, s},t_{2}, x,y)) n_{\kappa}^{-1} = ((g_{a, s},t_{1}, x,y),(g_{a, s},t_{2}, x,y))$.
In particular, $n_{\kappa}$ normalizes $\widetilde{\eta}$.
\end{proof}

According to Lemma~\ref{lemma:extensionofetaanditsnormalizer}, we have $N_{\varphi_{\pi}, \rho_{\pi}}^{+} = N_{\varphi_{\pi}}^{+}$, and the cohomology class $[\mu_{(\varphi_{\pi}, \rho_{\pi})}]$ agrees with the image of the cohomology class $[\mu_{N_{\varphi_{\pi}, \rho_{\pi}}^{+}, N', \widetilde{\eta}}] \in H^2\left(
N_{\varphi_{\pi}, \rho_{\pi}}^{+}/N', \C^{\times}
\right)$ via the inflation map
\[
H^2\left(
N_{\varphi_{\pi}, \rho_{\pi}}^{+}/N', \mathbb{C}^{\times}
\right)
\hookrightarrow H^2\left(
N_{\varphi_{\pi}, \rho_{\pi}}^{+}/S_{\varphi_{\pi}}^{+}, \mathbb{C}^{\times}
\right)
\]
induced by the natural surjection $N_{\varphi_{\pi}, \rho_{\pi}}^{+}/S_{\varphi_{\pi}}^{+} \twoheadrightarrow N_{\varphi_{\pi}, \rho_{\pi}}^{+}/N'$.

We compute $[\mu_{N_{\varphi_{\pi}, \rho_{\pi}}^{+}, N', \widetilde{\eta}}]$.
We fix a family of non-zero elements 
\[
\mathcal{T}' = 
\left\{
T'_{n} \in \Hom_{N'}\left(
{^n\widetilde{\eta}}, \widetilde{\eta}
\right)
\right\}_{n \in N_{\varphi_{\pi}, \rho_{\pi}}^{+}}
\]
satisfying Condition~\ref{conditionforcTinappendix} such that $T'_{1} = \id_{\mathbb{C}}$, $(T'_{n_{\zeta}})^{2} = \widetilde{\eta}(n_{\zeta}^{2})$, $(T'_{n_{\kappa}})^{2} = \widetilde{\eta}(n_{\kappa}^{2})$, and $T'_{n_{\kappa} n_{\zeta}} = T'_{n_{\kappa}} \circ T'_{n_{\zeta}}$.

\begin{proposition}
\label{prop:computationofmuT'}
Let $a_1, a_2, b_1, b_2 \in \{0, 1\}$.
Then we have
$\mu_{\mathcal{T}'}(n_{\kappa}^{a_{1}}n_{\zeta}^{a_{2}}, n_{\kappa}^{b_{1}}n_{\zeta}^{b_{2}}) = (-1)^{a_{2}b_{1}}$.
\end{proposition}
\begin{proof}
By Condition~\ref{conditionforcTinappendix}, we have
\[
T'_{n_{\kappa}^{a_{1}}n_{\zeta}^{a_{2}}n_{\kappa}^{b_{1}}n_{\zeta}^{b_{2}}} = \widetilde{\eta}(n_{\zeta}^{-a_2 - b_2} n_{\kappa}^{-a_1 - b_1} n_{\kappa}^{a_1} n_{\zeta}^{a_2} n_{\kappa}^{b_1} n_{\zeta}^{b_2})  \cdot T'_{n_{\kappa}^{a_{1} + b_{1}}n_{\zeta}^{a_{2} + b_{2}}}.
\]
Since $n_{\zeta}$ and $n_{\kappa}$ normalize $\widetilde{\eta}$, we have
\begin{align*}
\widetilde{\eta}(n_{\zeta}^{-a_2 - b_2} n_{\kappa}^{-a_1 - b_1} n_{\kappa}^{a_1} n_{\zeta}^{a_2} n_{\kappa}^{b_1} n_{\zeta}^{b_2})&= \widetilde{\eta}(n_{\zeta}^{-b_2} n_{\zeta}^{- a_2} n_{\kappa}^{-b_1} n_{\zeta}^{a_2} n_{\kappa}^{b_1} n_{\zeta}^{b_2}) \\
&= \widetilde{\eta}(n_{\zeta}^{- a_2} n_{\kappa}^{-b_1} n_{\zeta}^{a_2} n_{\kappa}^{b_1}) \\
&= 
\begin{cases}
\widetilde{\eta}((-1_{n}, 1, 1, 1), (- 1_{n}, 1, 1, 1)) & (a_2 = b_1 = 1), \\
1 & (\text{otherwise}).
\end{cases}
\end{align*}
Since the restriction of $\eta$ to 
$\mu_{n}(\C) \simeq \pi_{0}\left(
Z\left(
\widehat{G/Z(G)}
\right)^{\Gamma_{F}}
\right) \subset \pi_{0}([\widehat{\overline{S}}]^{+})$ has order $n$, we have 
\[
\eta((-1_{n}, 1, 1, 1), (- 1_{n}, 1, 1, 1)) = -1.
\]
Hence, we obtain that 
\[
\widetilde{\eta}(n_{\zeta}^{-a_2 - b_2} n_{\kappa}^{-a_1 - b_1} n_{\kappa}^{a_1} n_{\zeta}^{a_2} n_{\kappa}^{b_1} n_{\zeta}^{b_2}) = (-1)^{a_2 b_1}.
\]
We also obtain from the conditions of $\mathcal{T}'$ that $T'_{n_{\kappa}^{a_{1} + b_{1}}n_{\zeta}^{a_{2} + b_{2}}} = T'_{n_{\kappa}^{a_{1}}n_{\zeta}^{a_{2}}} \circ T'_{n_{\kappa}^{b_{1}}n_{\zeta}^{b_{2}}}$.
Thus, we conclude that 
\[
\mu_{\mathcal{T}'}(n_{\kappa}^{a_{1}}n_{\zeta}^{a_{2}}, n_{\kappa}^{b_{1}}n_{\zeta}^{b_{2}}) = \widetilde{\eta}(n_{\zeta}^{-a_2 - b_2} n_{\kappa}^{-a_1 - b_1} n_{\kappa}^{a_1} n_{\zeta}^{a_2} n_{\kappa}^{b_1} n_{\zeta}^{b_2})^{-1} = (-1)^{a_2 b_1}.
\qedhere
\]
\end{proof}

Now, we obtain a special case of Corollary~\ref{corollary:Kaletha'sLLCsatisfiesmainconjecture}.

\begin{corollary}
In the above setting, we have $[\mu_{\mathfrak{s}}] = [\mu_{\mathfrak{s}^{\vee}}]$.
\end{corollary}
\begin{proof}
Under the identification $X_{\unr}(G) \simeq X_{\unr}(^LG)$ in \eqref{eq:explicitKottwitzisom}, the subquotient $(N/T)^{*}$ of $X_{\unr}(G)$ corresponds to the subquotient $N_{\varphi_{\pi}, \rho_{\pi}}^{+}/N'$ of $X_{\unr}(^LG)$.
Moreover, the generators $\chi_1$ and $\chi_2$ in $(N/T)^{*}$ correspond to the generators $n_{\kappa}$ and $n_{\zeta}$ of $N_{\varphi_{\pi}, \rho_{\pi}}^{+}/N'$, respectively.
Thus, the claim follows from \eqref{eq:computationofmusexample} and Proposition~\ref{prop:computationofmuT'}.
\end{proof}

\end{document}